\documentclass{amsart}
\usepackage{graphics}
\usepackage[all]{xy}

\newtheorem{defn}[subsubsection]{Definition}
\newtheorem{thm}[subsubsection]{Theorem}
\newtheorem{lem}[subsubsection]{Lemma}

\newtheorem*{fact}{Fact}
\newtheorem{prop}[subsubsection]{Proposition}
\newtheorem*{fpbw}{The formal Poincare-Birkhoff-Witt theorem}
\newtheorem*{HW}{Homological Wheeling}

\newcommand{\ASreln}{\text{{AS}}}
\newcommand{\STUreln}{\text{{STU}}}
\newcommand{\IHXreln}{\text{{IHX}}}
\newcommand{\Aspace}{\mathcal{A}}
\newcommand{\Bspace}{\mathcal{B}}

\newcommand{\Wspace}{\mathcal{W}}

\newcommand{\Wbase}{\Wspace_{\text{basic}}}
\newcommand{\Wbasei}{\Wspace^i_{\text{basic}}}
\newcommand{\What}{\widehat{\mathcal{W}}}
\newcommand{\Whatiota}{\widehat{\mathcal{W}}_\iota}
\newcommand{\Whatwedge}{\widehat{\mathcal{W}}_\wedge}
\newcommand{\Whatwedgeiota}{\widehat{\mathcal{W}}_{\wedge\iota}}
\newcommand{\WhatF}{\widehat{\mathcal{W}}_{\mathrm{\bf F}}}
\newcommand{\WhatFiota}{\widehat{\mathcal{W}}_{\mathrm{\bf F}\iota}}
\newcommand{\baseFtobull}{\text{B}_{\mathrm{F}\rightarrow\bullet}}
\newcommand{\basebulltoF}{\text{B}_{\bullet\rightarrow\mathrm{F}}}
\newcommand{\ncw}{\widetilde{\mathcal{W}}}
\newcommand{\intoper}{\int_{\mathcal{T}_{\mathrm{dR}}}}
\newcommand{\di}{(d,\iota)}

\begin{document}
\title[The combinatorics of wheeling]{Non-commutative Chern-Weil theory and the combinatorics of Wheeling}
\author[A. Kricker]{Andrew Kricker}
\address{Division of Mathematical Sciences \\ School of
Physical and Mathematical Sciences \\ Nanyang Technological
University, Singapore, 637616} \email{ajkricker@ntu.edu.sg}
\begin{abstract}
This work applies the ideas of Alekseev and Meinrenken's
Non-commutative Chern-Weil Theory to describe a completely
combinatorial and constructive proof of the Wheeling Theorem. In
this theory, the crux of the proof is, essentially, the familiar
demonstration that a characteristic class does not depend on the
choice of connection made to construct it. To a large extent this
work may be viewed as an exposition of the details of some of
Alekseev and Meinrenken's theory written for Kontsevich integral
specialists. Our goal was a presentation with full combinatorial
detail in the setting of Jacobi diagrams -- to achieve this goal
certain key algebraic steps required replacement with
substantially different combinatorial arguments.

\end{abstract}
 \maketitle

\section{Introduction and Outline}
This paper is organized around a purely combinatorial proof of
what is called the Wheeling Isomorphism \cite{BGRT, BLT}. This is an algebra isomorphism between a pair of ``diagrammatic"
algebras, $\Aspace$ and $\Bspace$, $(\chi_{\Bspace}\circ\partial_\Omega)\, :\, \Bspace \to \Aspace$.

We will recall these spaces and maps shortly, beginning in Section
\ref{spacesintro}. The formal development of the theory begins in
Section \ref{commweil}. Our underlying aim is to use this proof to
introduce a combinatorial reconstruction of some elements of
Alekseev and Meinrenken's Non-commutative Chern-Weil theory
(see \cite{Mein} for an introduction).

The work falls naturally into two parts.
\begin{itemize}
\item{In the first part we'll describe and prove a statement we'll
call ``Homological Wheeling". This is a lifting of Wheeling
to the setting of a ``Non-commutative Weil complex for diagrams".
That complex will be introduced in Section \ref{noncommsect}.
After that, Homological Wheeling will be stated as Theorem
\ref{reexpress}. The homotopy equivalence at the heart of the
proof will be constructed in Section \ref{homotopyconstruct}.}
\item{The remainder of the paper will be occupied with describing
how the usual statement of Wheeling is recovered from Homological
Wheeling when certain relations are introduced.
This part consists mainly of the sort of delicate
``gluing legs to legs" combinatorics that Kontsevich integral
specialists will find very familiar (though here we have the extra
challenge of keeping track of permutations and signs). The details
of a key combinatorial identity (calculating how the ``Wheels"
power series appears in the theory) will appear in an accompanying
publication \cite{K}.}
\end{itemize}

Wheeling is a combinatorial strengthening of the Duflo isomorphism
of Lie theory. To recall Duflo: The classical
Poincare-Birkhoff-Witt theorem concerns $\mathfrak{g}$, an
arbitrary Lie algebra, $S(\mathfrak{g})^{\mathfrak{g}}$, the
algebra of $\mathfrak{g}$-invariants in the symmetric algebra of
$\mathfrak{g}$, and $U(\mathfrak{g})^{\mathfrak{g}}$, the algebra
of $\mathfrak{g}$-invariants in the universal enveloping algebra
of $\mathfrak{g}$. PBW says that the natural averaging map
actually gives a vector space isomorphism between the vector
spaces underlying these algebras. The averaging map, however, does
not respect the product structures. Duflo's fascinating discovery \cite{D}
was of an
explicit vector-space isomorphism of $S(\mathfrak{g})^{\mathfrak{g}}$ (in the case that $\mathfrak{g}$ is a finite-dimensional Lie algebra) whose composition with the averaging map gives an
isomorphism of {\bf algebras} $S(\mathfrak{g})^{\mathfrak{g}}\cong U(\mathfrak{g})^{\mathfrak{g}}$.

\subsection{The Kontsevich integral as the underlying motivation.}
Wheeling is a combinatorial analogue of the Duflo
isomorphism, but its origins and role are more interesting than
that fact suggests. When Bar-Natan, Garoufalidis, Rozansky and
Thurston conjectured that the Wheeling map was an algebra
isomorphism they were motivated, in part, by questions that arose
in the developing theory of the Kontsevich integral knot invariant (see
\cite{BGRT} and \cite{BLT} and discussions therein).
Wheeling is a fundamental topic in the study
of the Kontsevich integral because:
\begin{itemize}
\item{The Kontsevich integral of the unknot turns out to be {\it
precisely} the power series of graphs that is used in Wheeling.} \item{In addition to
the above fact, Bar-Natan, Le, and Thurston showed that the
Wheeling isomorphism was a consequence of the invariance of the
Kontsevich integral under a certain elementary isotopy (this is
their famous `1+1=2' argument). See \cite{BLT}.} \item{Wheeling
has proved to be an indispensable tool in almost all `hard'
theorems that can be proved about the Kontsevich integral (one
example is the rationality of the loop expansion, \cite{GK}).}
\end{itemize}

Here is the interesting thing: Wheeling has a {\it combinatorial}
statement, but the BLT proof of it has, at its heart, an analytic
machine -- the monodromy of the Knizhnik-Zamolodchikov
differential equations. The interest in a purely combinatorial
proof of Wheeling is whether we can reverse the above story and
use such a proof to learn something about the Kontsevich integral.
A purely combinatorial construction of the Kontsevich integral
seems to be an important goal. Perhaps the rich algebraic
structures and ideas that are associated with this new proof will
give a novel approach to this goal?


\subsection{A sketch of a few key constructions.}\label{sketch}
We'll begin by sketching out a few of the key players in the
proof. After that we'll take a step back and explain how these
structures are motivated by some recent developments in the
Chern-Weil theory of characteristic classes.

Recall that the space $\Aspace$ is generated by some combinatorial
gadgets called {\it Jacobi diagrams} (which we'll occasionally refer to as {\it ordered Jacobi diagrams}), and that $\Bspace$ is
generated by {\it symmetric Jacobi diagrams}. See Section
\ref{spacesintro} for the precise definitions we'll use. The
averaging map $\chi_\Bspace:\Bspace\rightarrow \Aspace$ is a
vector space isomorphism between the two spaces -- this is just a
formal version of the classical Poincare-Birkhoff-Witt theorem.
Here is an example of the averaging map:
\[
\chi_\Bspace\left(\raisebox{-0.8cm}{\scalebox{0.25}{\includegraphics{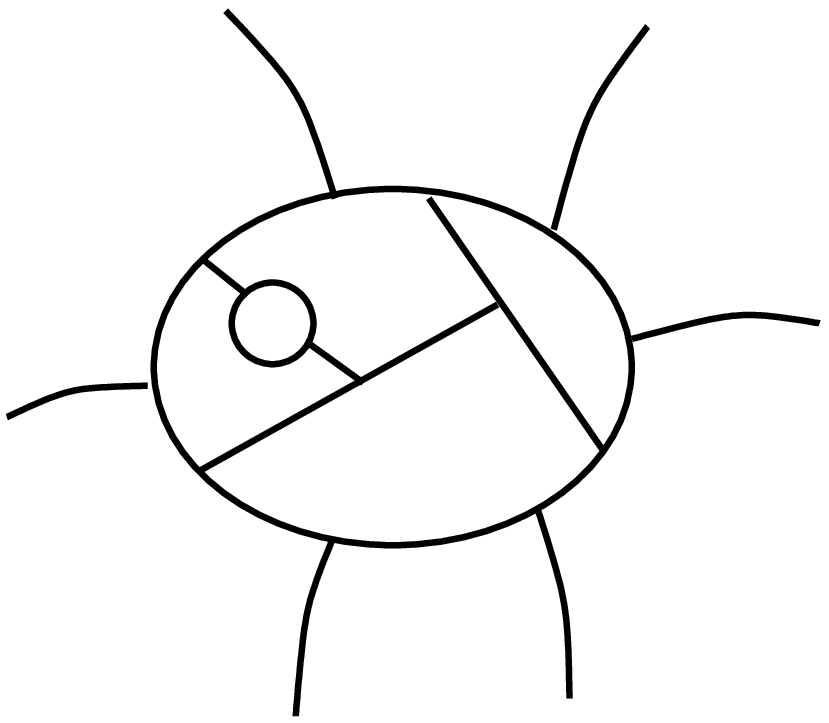}}}\right)=\frac{1}{6!}\sum_{\sigma\in\Sigma_6}\raisebox{-1cm}{\scalebox{0.22}{\includegraphics{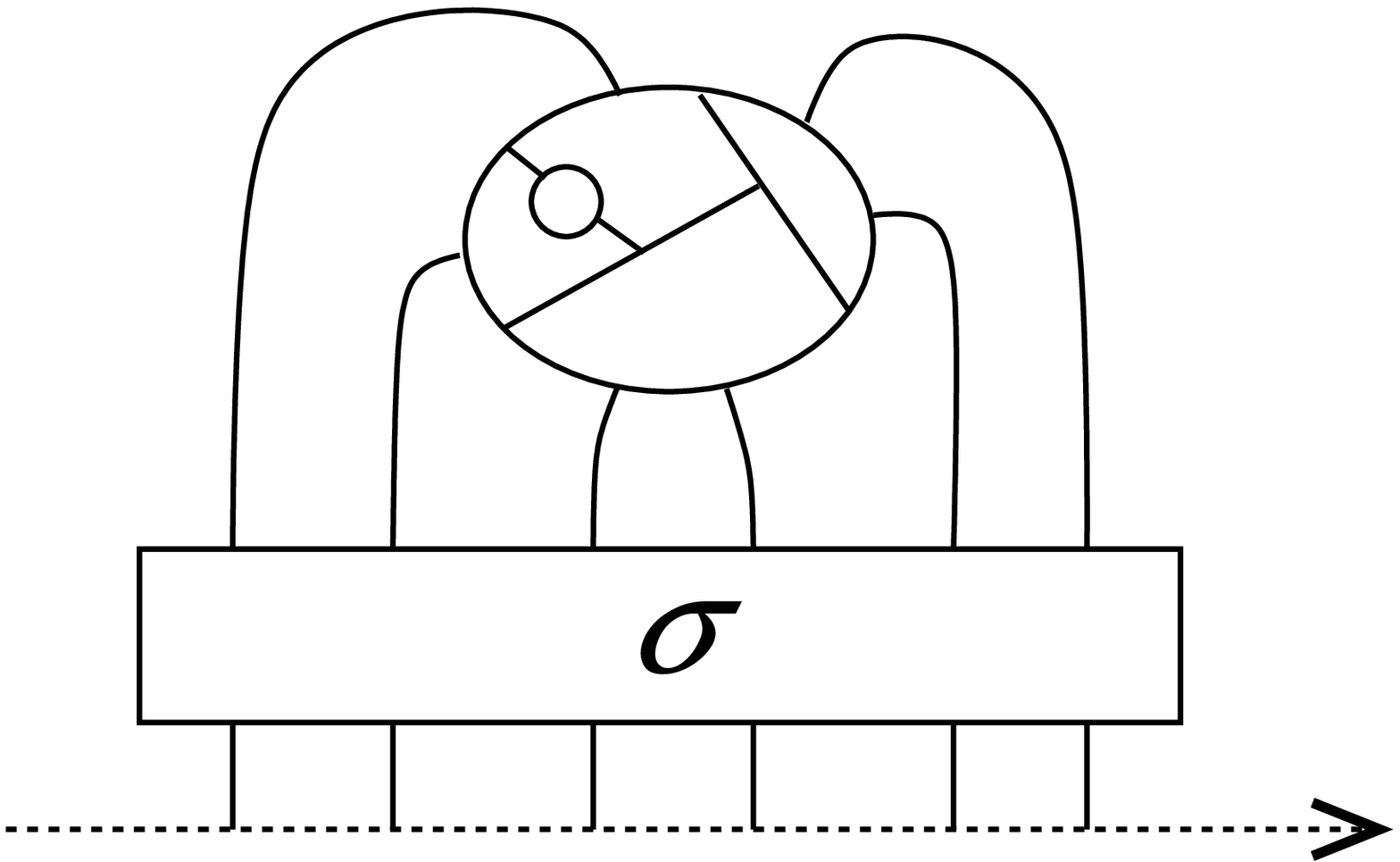}}}\ \in \Aspace.
\]
Of course, $\chi_\Bspace$ does not respect the natural product
structures on $\Aspace$ and $\Bspace$. The purpose of the
``Wheeling Isomorphism" $\partial_\Omega:\Bspace\rightarrow
\Bspace$ is to promote $\chi_\Bspace$ to an algebra isomorphism:
$(\chi_\Bspace\circ
\partial_\Omega): \Bspace \rightarrow \Aspace$.

Speaking generally, the thrust of this work is to study an
averaging map
\[
\chi_\Wspace : \Wspace \rightarrow \ncw
\]
between some algebraic systems $\Wspace$ and $\ncw$ that enlarge
the spaces $\Bspace$ and $\Aspace$ in certain ways.

First we'll focus on $\Wspace$. $\Wspace$ will be an example of
what we'll call an {\bf $\iota$-complex} (see Section \ref{ics}).
This is a diagram in the
category of vector spaces of the following form:
\[
\xymatrix{ 0 \ar[r] & \mathcal{W}^0 \ar[r]^{d^0}
\ar[dl]_{\iota^0}& \mathcal{W}^1
\ar[r]^{d^1} \ar[dl]_{\iota^1} & \mathcal{W}^2 \ar[dl]_{\iota^2} \ar[r]^{d^2} & \ldots \\
0 \ar[r] & \mathcal{W}^0_\iota \ar[r]^{d^0_\iota} &
\mathcal{W}^1_\iota \ar[r]^{d^1_\iota}
& \mathcal{W}^2_\iota \ar[r]^{d^2_\iota} & \ldots \\
}
\]
To be an $\iota$-complex: The rows of this diagram are co-chain
complexes, and the diagonal maps form a degree $-1$ map between
these complexes (in other words, $\iota^{i+1}\circ d^i =
-d^{i-1}_\iota\circ \iota^i$, for all $i\geq 0$).

The vector spaces $\Wspace^i$ are generated by some combinatorial
gadgets that we'll call {\bf Weil diagrams}. Weil diagrams are
very much like Jacobi diagrams, but Weil diagrams have two types
of legs: legs of leg-grade 1, and legs of leg-grade 2. The legs of
leg-grade 2 are distinguished by being drawn with a fat dot.
$\Wspace^i$ will be generated by Weil diagrams with total
leg-grade $i$. The legs of a Weil diagram are ordered, and when we
transpose the position of two adjacent legs, we pick up a sign
$(-1)^{xy}$, where $x$ and $y$ are the leg-grades of the involved
legs. For example:
\[
\raisebox{-3ex}{\scalebox{0.25}{\includegraphics{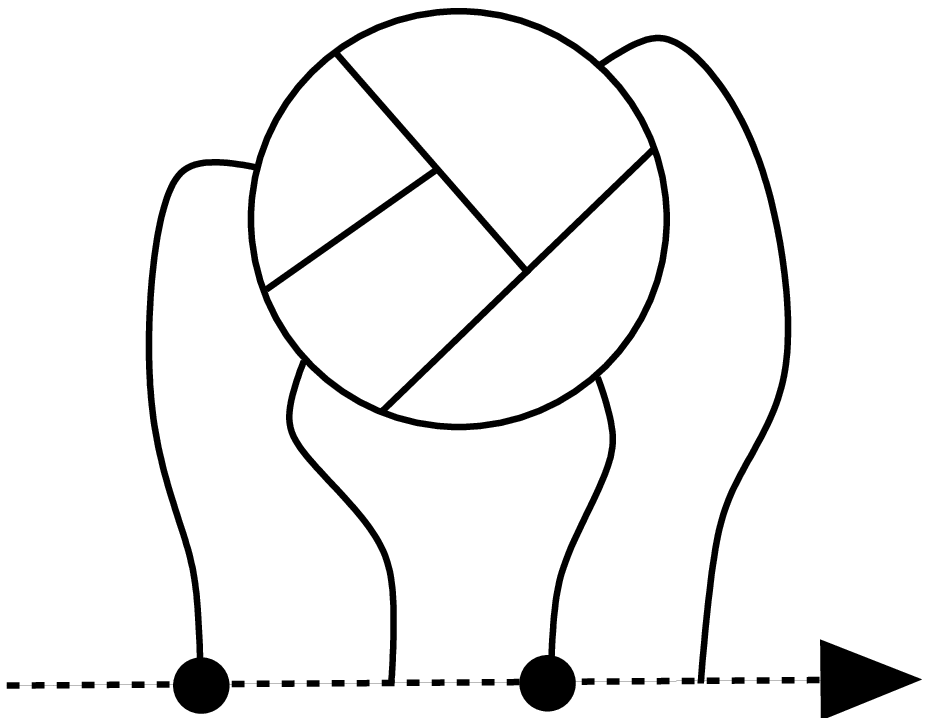}}} \ \
\ =\ \ \
\raisebox{-3ex}{\scalebox{0.25}{\includegraphics{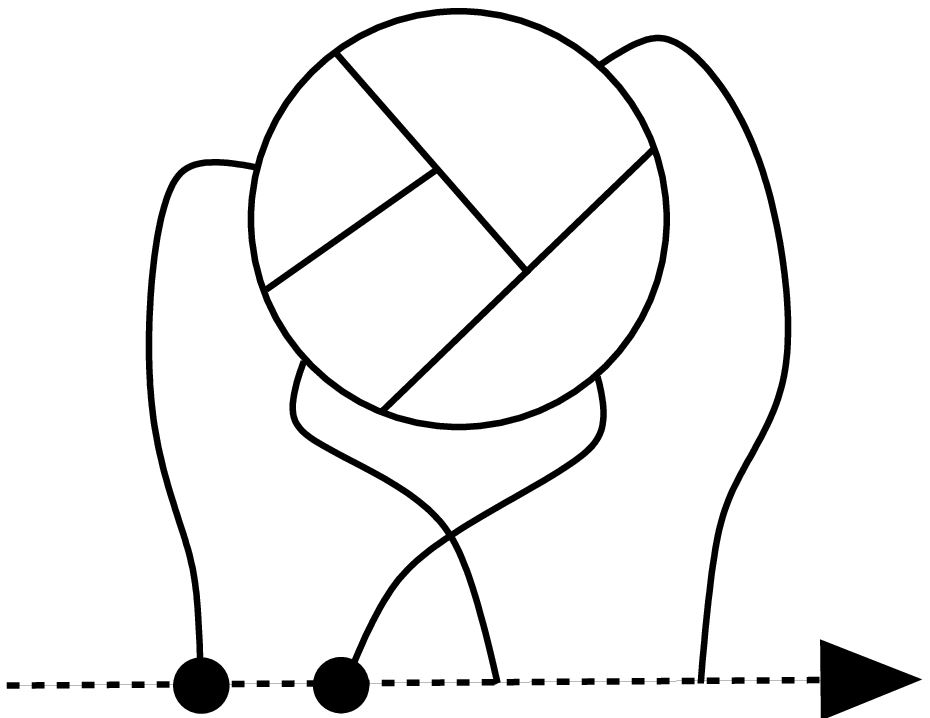}}} \ \
\ =\ \ \ -\
\raisebox{-3ex}{\scalebox{0.25}{\includegraphics{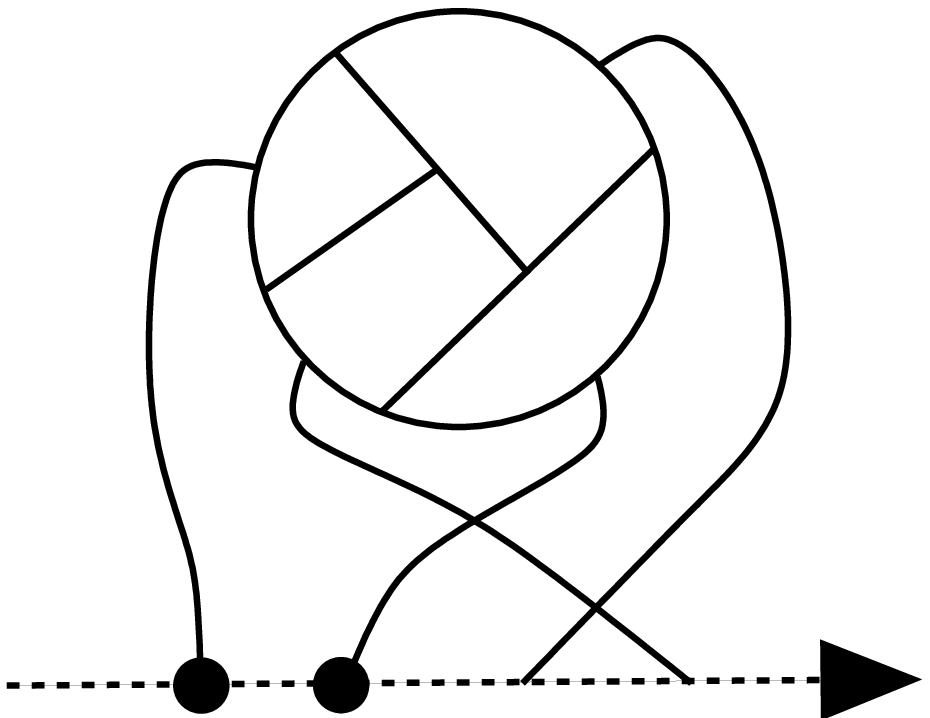}}}\ \
\ \ \ \mbox{in\ \ ${\mathcal{W}}^6$.}
\]
The spaces $\Wspace^i_\iota$ are exactly the same, but have
exactly one degree 1 vertex labelled $\iota$ which floats freely
(i.e. is not included in the ordering of the legs). See Section
\ref{commweil} for the full precise definition of this system
$\Wspace$: the {\bf commutative Weil complex for diagrams}.

Now it follows from the definition of an $\iota$-complex that if
we take the kernels of the maps $\iota$, then we get another
co-chain complex:
\[
\xymatrix{0 \ar[r] & \mathrm{Ker}\left(\iota^0\right) \ar[r]^{d^0}
& \mathrm{Ker}\left(\iota^1\right) \ar[r]^{d^1} &
\mathrm{Ker}\left(\iota^2\right) \ar[r]^{d^2} & \ldots}
\]
This is called $\Wbase$, the {\bf basic subcomplex} of the
original $\iota$-complex, and its cohomology is called the {\bf
basic cohomology}.

The relationship of this ``enlarged" structure, $\Wspace$, to the
original space, $\Bspace$, is the following:
\[
H^{i}(\Wbase) \cong \Bspace^i,
\]
where $\Bspace^i$ is the space\footnote{Our grading of $\Bspace$
is non-traditional. See Section \ref{spacesintro} for details.} of
symmetric Jacobi diagrams of leg-grade $i$.

Now consider $\ncw$, the target of the averaging map $\chi_\Wspace
: \Wspace \rightarrow \ncw$.  This will also be an example of an
$\iota$-complex, and $\chi_\Wspace$ will be a ``map of
$\iota$-complexes".

In Section \ref{noncommsect}, $\ncw$ will be defined in exactly
the same way as $\Wspace$, but without introducing the leg
transposition relations. So the relationship of $\ncw$ to
$\Wspace$ will be analogous to the relationship between a tensor
algebra $T(V)$, and a symmetric algebra $S(V)$ over some vector
space $V$.

The way that this structure enlarges $\Aspace$ is a bit different
to the way that $\Wspace$ enlarges $\Bspace$. In Section
\ref{gettinghw} we'll see that we can introduce some relations to
obtain:
\[
\left(\bigoplus_{i=0}^{\infty}\ncw^i\right)/(\text{relns.})\ \
\cong\ \Aspace\oplus\ldots
\]

The general flow of this work is:
\begin{itemize}
\item{To show that $\chi_\Wspace$ descends to a genuine
isomorphism of the basic cohomology rings of $\Wspace$ and
$\ncw$.} \item{To show that the usual statement of wheeling is
recovered from this fact when we pass to the ``smaller" structures
$\Aspace$ and $\Bspace$.}
\end{itemize}

\subsection{The origins of this work in Alekseev-Meinrenken's
Non-commutative Chern-Weil theory} This work originated in a
course given by E. Meinrenken at the University of Toronto in the
Spring of 2003 in which
Alekseev and Meinrenken's work on a Non-commutative Chern-Weil
theory was introduced \cite{ML}. See \cite{Mein} for a short review of the
theory, \cite{AleksMein} for the original paper, and
\cite{AleksMein2} for later developments.

The goal of this work is a self-contained reconstruction and exploration of
their technology in the combinatorial setting of the Jacobi
diagrams of Quantum Topology. Before beginning that discussion in the
next section, we'll give here a very brief sketch of the outlines of classical
Chern-Weil theory.
A more detailed
discussion of how our combinatorial definitions fit into the
classical picture will be given in Section \ref{lengthyaside},
once those definitions have been made.

So let $G$ be a compact Lie group with Lie algebra $\mathfrak{g}$.
Chern-Weil theory is a theory of {\bf characteristic classes} in
de Rham cohomology of smooth principal $G$-bundles. Recall that
every smooth principal $G$-bundle arises as a smooth manifold $P$
with a free smooth left $G$-action.

What is a characteristic class? The first thing to say is that it
is a choice, for every principal $G$-bundle $P$, of a class $w(P)$
in the cohomology of the base space, $B=P/G$, which is {\it
functorial} with respect to bundle morphisms. To be precise: if
there is a morphism of principal $G$-bundles $P$ and $P'$, which
means that there is a $G$-equivariant map $f:P \rightarrow P'$,
then it should be true that
$f_B^*(w(P')) = w(P)$,
where $f_B$ denotes the map induced by $f$ between the base spaces
$f_B : B\rightarrow B'$.

Chern-Weil theory constructs such systems of classes in de Rham
cohomology in the following way.
\begin{enumerate} \item{It takes an {\it arbitrarily chosen}
connection form $A\in\Omega^1(P)\otimes \mathfrak{g}$.} \item{It
uses it to construct closed forms in
$\Omega(P)_{\mathrm{basic}}\subset \Omega(P)$. This is the ring of
{\it basic} differential forms on $P$. It is most simply
introduced as the isomorphic image under $\pi^*$ of $\Omega(B)$.
It is a crucial fact in the theory that this basic subring is
selected as the kernel of a certain Lie algebra of differential
operators arising from the generating vector fields of the
$G$-action. The main task of the theory is to canonically
construct forms from $A$ lying in the kernel of this Lie algebra.
(See Section \ref{lengthyaside} for a more detailed discussion).}
\item{It then shows that the cohomology classes in $H_{dR}(B)$ so
constructed do not depend on the choice of connection and do
indeed have the required functoriality with repsect to bundle
morphisms.}
\end{enumerate}

The heart of the Chern-Weil construction is an algebraic structure
called {\bf the Weil complex}. The Weil complex is based on the
graded algebra $S({\mathfrak{g}^*}\oplus{\mathfrak{g}}^*)$
(that is, the graded symmetric algebra on two copies of
${\mathfrak{g}}^*$, the dual of $\mathfrak{g}$, one copy in
grade 1 and one copy in grade 2). The Weil complex is this graded
algebra equipped with a certain graded Lie algebra of graded
derivations. This Lie algebra of derivations selects a key
subalgebra
$S({\mathfrak{g}}^*\oplus{\mathfrak{g}}^*)_\mathrm{basic}$,
the {\it basic} subalgebra.

What is the construction then? Consider some principal $G$-bundle
$P$. Given the fixed data of a {\it connection form}
$A\in\Omega^1(P)\otimes\mathfrak{g}$ on $P$, Chern-Weil theory
constructs a map
\[
c_A : S({\mathfrak{g}}^*\oplus{\mathfrak{g}}^*)\rightarrow
\Omega(P)
\]
by the correspondence
\[
g_1g_2\ldots g_n \overline{g}_1\ldots\overline{g}_m\mapsto
g_1(A)\wedge \ldots\wedge g_n(A) \wedge d\overline{g}_1(A)\wedge
\ldots \wedge d\overline{g}_m(A).
\]
It is immediate in the theory that this construction specializes
to a map of differential graded algebras from the basic subalgebra
of the Weil complex into the ring of basic forms on $P$:
\[
S({\mathfrak{g}}^*\oplus{\mathfrak{g}}^*)_\mathrm{basic}
\rightarrow \Omega(P)_{\mathrm{basic}} \cong \Omega(B),
\]
so that, passing to cohomology, we get a map:
\[
H^*\left(S({\mathfrak{g}}^*\oplus{\mathfrak{g}}^*)_\mathrm{basic}\right)
\rightarrow H_{dR}^*(B).
\]

The fundamental theorem of Chern-Weil theory is that this
construction does indeed give a canonical system of rings of
characteristic classes. We can summarize Chern-Weil theory with
the following commutative diagram:
\[
\xymatrix{
 & H^*(\Omega(P)_{\mathrm{basic}}) \ar@{=}[r] & H^*_{dR}(B) \\
 H^*(S({\mathfrak{g}}^*\oplus{\mathfrak{g}}^*)_\mathrm{basic}) \ar[ur]^{H(c_A)} \ar[dr]_{H(c_{A'})}& \\
  & H^*(\Omega(P')_{\mathrm{basic}}) \ar@{=}[r] & H^*_{dR}(B') \ar[uu]_{f_B^*}  }
 \]

So why should this theory have anything to say about the Duflo
isomorphism, which seems a purely algebraic fact? The connection
arises from the key computation that
\[
H(S({\mathfrak{g}}^*\oplus{\mathfrak{g}}^*)_\mathrm{basic}) =
(S\mathfrak{g})^{\mathfrak{g}}.
\]
In words: the basic cohomology of the Weil complex, which is the
universal ring of characteristic classes (consider its place in
the above diagram), is precisely the ring of
$\mathfrak{g}$-invariants in $S\mathfrak{g}$, a ring involved in
the Duflo isomorphism. (We remark that this formulation assumes
that $\mathfrak{g}$ has a symmetric, non-degenerate, invariant
inner product.)

\subsection{A Weil complex for diagrams}
This paper will start with a combinatorial reconstruction of that key
computation
\[
H(S({\mathfrak{g}}^*\oplus{\mathfrak{g}}^*)_\mathrm{basic}) =
(S\mathfrak{g})^{\mathfrak{g}}.
\]
Recall that when we abstract Wheeling from the Duflo isomorphism,
we replace $(S\mathfrak{g})^{\mathfrak{g}}$ by $\Bspace$, the
space of symmetric Jacobi diagrams. An evaluation map (the obvious, minor variation on the ``weight system" construction, \cite{BarNatan})
connects the two:
$\text{Eval}_{\mathfrak{g}} : \Bspace \rightarrow
(S\mathfrak{g})^{\mathfrak{g}}$.

We'll begin by replacing the ``something" in the
following equation \[ H((\text{Something})_\mathrm{basic}) =
\Bspace,
\]
by an appropriate construction: $\Wspace$, the commutative Weil
complex for diagrams. The complex $\Wspace$ is a combinatorial
version of the usual Weil complex
$S({\mathfrak{g}}^*\oplus{\mathfrak{g}}^*)$. Because the usual Weil
complex is built from two copies of ${\mathfrak{g}}^*$, one in
grade 1 and one in grade 2, $\Wspace$ will be built from
diagrams with two types of legs: legs of grade 1 and legs of grade
2. In other words, $\Wspace$ will be generated by the ``Weil diagrams"
introduced in Section \ref{sketch}.

The full details of the construction of $\Wspace$ will be presented in
Section \ref{commweil}. $\Wspace$ will be constructed to be an $\iota$-complex, so that we may take its basic cohomology.


The final ingredient required in the explication of the equation
$H^i(\Wspace_\mathrm{basic}) \cong \Bspace^i$ is the map
\[
\Upsilon : \Bspace \rightarrow \Wspace
\]
which actually performs the isomorphism. This map $\Upsilon$
occupies a key role in our presentation. It will be introduced in
Section \ref{basiccompute}, where the computation that
$H(\Wspace_\mathrm{basic}) = \Bspace$ will be concluded. Its
definition on some symmetric Jacobi diagram $w$ is quite simple:
choose an ordering of the legs of $w$, then expand each leg
according to the following rule:
\[
\begin{array}{ccccl}
\Upsilon :\ \ \ \ \
\raisebox{-2.6ex}{\scalebox{0.28}{\includegraphics{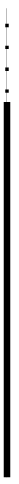}}}\ \ &
\mapsto &
\raisebox{-3ex}{\scalebox{0.27}{\includegraphics{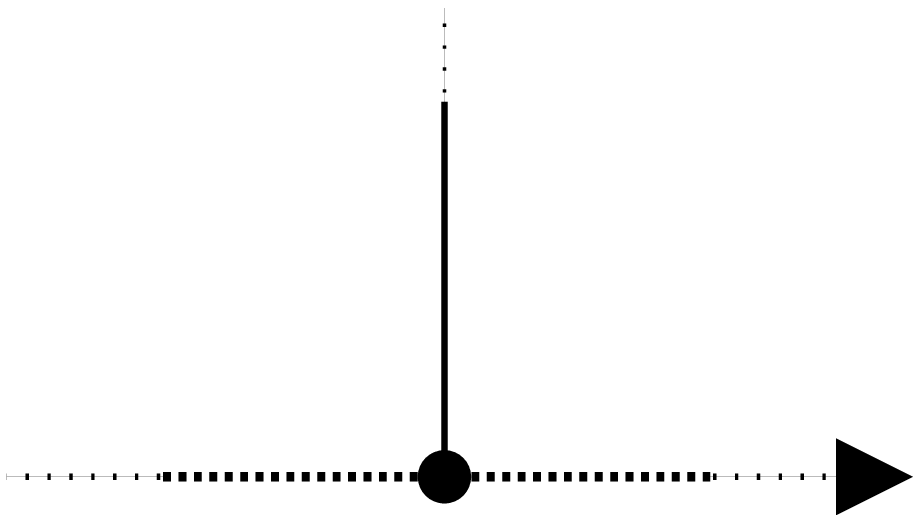}}}
-\frac{1}{2}
\raisebox{-3ex}{\scalebox{0.27}{\includegraphics{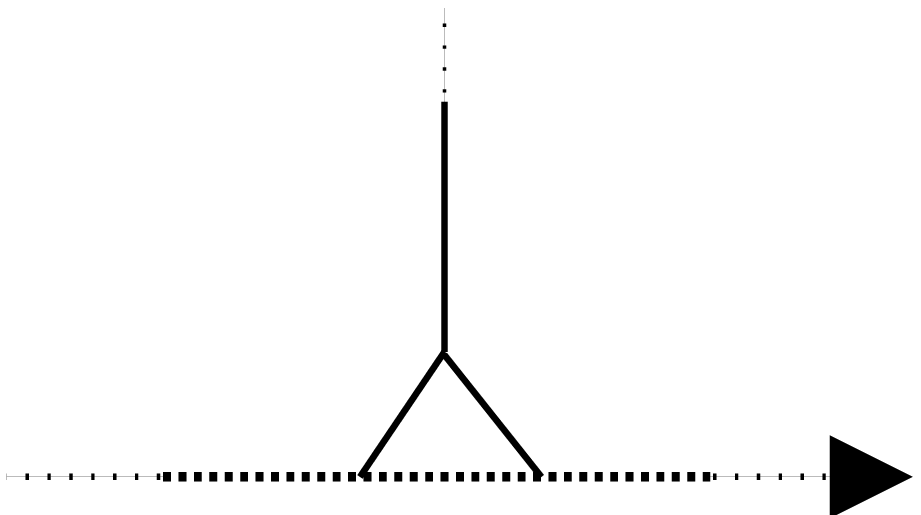}}}\ .
\end{array}
\]
Sometimes we'll refer to $\Upsilon$ as the {\it hair-splitting map}.

\subsection{Some motivation for a non-commutative Weil complex}
As the author learnt it from E. Meinrenken, this construction can
be motivated with the following line of thought: to begin, note
that there is {\it almost} a much more natural way to build these
characteristic classes in {\bf de Rham cohomology}.  The more
natural way would exploit $EG$, the classifying bundle of
principal $G$-bundles, and its associated base space $EG/G=BG$.
Recall the key property of $EG$ that the bundle pullback operation
performs a bijection, for some space $X$:
\[
\left\{ \begin{array}{l} \text{Isomorphism classes of} \\
\text{principal $G$-bundles over $X$}
\end{array}
\right\} \leftrightarrow [X:BG].
\]
So, proceeding optimistically, we could define the ring of
characteristic classes associated to some smooth principal
$G$-bundle $P$ by the map:
\[
\text{Class}_P^* : H_{dR}(BG)\rightarrow H_{dR}(B),
\]
for some choice of classifying map Class$_P:B\rightarrow BG$
giving $P$ as its corresponding pullback bundle. This ring would
be well-defined because the map Class$_P$ is unique up to
homotopy. Furthermore, this arrangement automatically gives the
required commutative diagrams:
\[
\xymatrix{
 & H_{dR}(B) \\
 H_{dR}(BG) \ar[ur]^{\text{Class}_P^*} \ar[dr]_{\text{Class}_{P'}^*}& \\
  & H_{dR}(B') \ar[uu]_{f_B^*} }
\]
The problem with this picture, of course, is that the space $BG$
is not in general a manifold, so there is no ring of differential
forms $\Omega(BG)$ with which to work. The Weil complex can be
seen as a kind of stand in for the ring of differential forms on
the classifying bundle $EG$: $ S(\mathfrak{g}^*\oplus
\mathfrak{g}^*) = ``\Omega(EG)".$

But this picture is too nice to let it go so easily. It turns out
that there does exist a point of view in which it is reasonable to
discuss ``differential forms" on $EG$: the simplicial space
construction of $EG$ (see Section 2 of \cite{MP} for a discussion
of the required simplicial techniques). $EG$ can be constructed as
the geometric realization of a certain simplicial manifold:
\[
EG = \left(\coprod_{n=0}^{\infty}\Delta^n\times G^n\right)/\sim.
\]
It follows, with some work, from a theorem of Moore's (see Section
4 of \cite{MP}) that the real cohomology of such a space can, in
fact, be computed as the total cohomology of an appropriate double
complex built from differential forms {\it on the pieces} of the
simplicial manifold (in the case at hand they are $G^n$). (Moore's
theorem is a generalization of the familiar theorem that equates
singular and simplicial cohomology). To be precise, let $C^n =
\oplus_{i=0}^n \Omega^i(G^{n-i})$. Then:
\[
H(EG,\mathbb{R}) \cong H(C^\bullet).
\]

So, with this beautiful idea in hand, we may feel that we are
ready to build a Chern-Weil theory which involves the classifying
spaces, as it should. There is a final hurdle to be cleared
though: the natural product structure on $\oplus_n\oplus_{i=0}^n
\Omega^i(G^{n-i})$ that is provided by Moore's theorem {\bf is not
graded commutative}, unlike the wedge product on the usual ring of
differential forms $\Omega(X)$. Well, we didn't really need
commutativity anyway, so we follow Alekseev and Meinrenken in
passing to the {\bf non-commutative Weil complex}. This has the
same definition as the usual Weil complex, but without introducing
commutativity (so that it is based on a {\it tensor} algebra, not
a {\it symmetric} algebra). It has natural maps to both the usual,
commutative, Weil complex, as well as the ``de Rham complex" of
$EG$:
\[
\xymatrix{ & \oplus_n\oplus_{i=0}^{n} \Omega^{i}(G^{n-i}) \\
T({\mathfrak{g}}^*\oplus {\mathfrak{g}}^*) \ar[ur]_{\text{\hspace{3cm}The usual Chern-Weil map based on a certain ``connection" on $EG$.}} \ar[dr]^{\text{\hspace{1.15cm}The canonical quotient map.}}& \\
& S({\mathfrak{g}}^*\oplus {\mathfrak{g}}^*)}
\]
The Alekseev-Meinrenken proof of the Duflo isomorphism arises from
an algebraic study of the relationship between the two Weil
complexes, $T({\mathfrak{g}}^*\oplus {\mathfrak{g}}^*)$ and
$S({\mathfrak{g}}^*\oplus {\mathfrak{g}}^*)$, exploiting properties of
$T({\mathfrak{g}}^*\oplus {\mathfrak{g}}^*)$ that derive from its place as the universal ring in a characteristic class theory.
\subsection{The Non-commutative Weil complex for Diagrams.}

So in Section \ref{noncommsect} we define $\ncw$, a
non-commutative Weil complex for diagrams, in the way that will be
obvious at this point of the development. The interplay between
$\Wspace$ and $\ncw$ is at the heart of this theory. There are two
natural maps between these $\iota$-complexes: \[ \xymatrix{
\mathcal{W} \ar@/^/[r]^{\chi_\Wspace}&
\widetilde{\mathcal{W}}\ar@/^/[l]^\tau}.
\]
The map $\chi_\Wspace$ is a key map: the graded averaging map. Its
action on some Weil diagram $w$ is to take the average of all
diagrams you can get by rearranging the legs of $w$ (each
multiplied by a sign appropriate to the rearrangement). The map
$\tau$ is the basic quotient map (corresponding to the
introduction of the commutativity relations).

{\bf Caution}: Our combinatorial analogue of the  ``non-commutative Weil complex" defined by the work \cite{AleksMein} is actually a quotient of this complex $\ncw$ by a certain system of relations. The precise corresponding structure we use, $\widehat{W}_{\mathrm{F}}$, will be discussed later, in Section \ref{theamncwc}. This complex $\ncw$ is the analog of a deeper (and simpler) structure that was discovered in the later work \cite{AleksMein2}.

\subsection{Homological Wheeling}
When we compose $\Upsilon$, the ``Hair-splitting map", with
$\chi_\Wspace$, the graded averaging map, we get a map which takes
elements of $\Bspace$ to basic cohomology classes in $\ncw$:
\[
\Bspace \stackrel{\Upsilon}{\longrightarrow}\Wbase
\stackrel{\chi_\Wspace}{\longrightarrow} \ncw_{\text{basic}}.
\]
So if we take two elements $v$ and $w$ of $\Bspace$, then
$\left(\chi_\Wspace\circ \Upsilon\right)(v\sqcup w)$
represents a basic cohomology class in $\ncw_{\text{basic}}$.

But there is another way to build a basic cohomology class in
$\ncw_{\text{basic}}$ from $v$ and $w$, using the juxtaposition
product $\#$ on $\ncw$. This is a product on basic cohomology
arising from the following product defined on the generating
diagrams:
\[
\raisebox{-3ex}{\scalebox{0.25}{\includegraphics{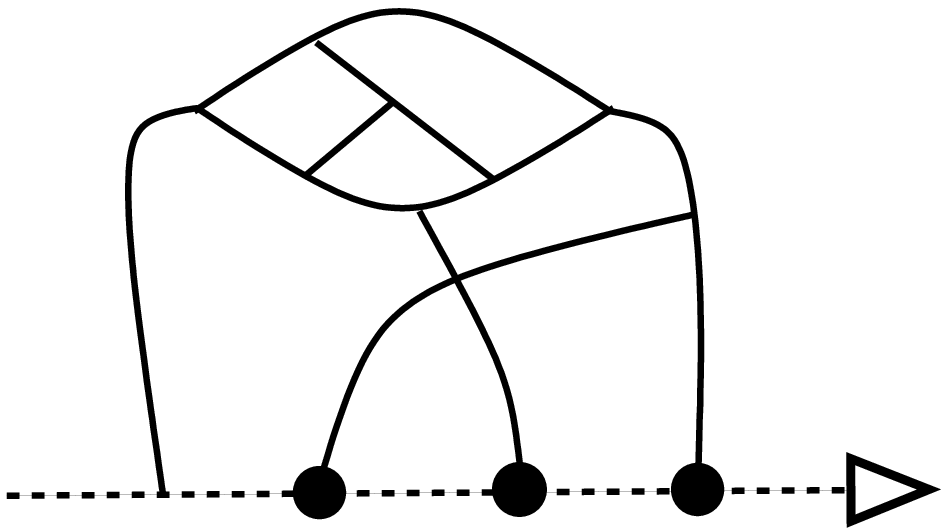}}} \#
\raisebox{-3ex}{\scalebox{0.25}{\includegraphics{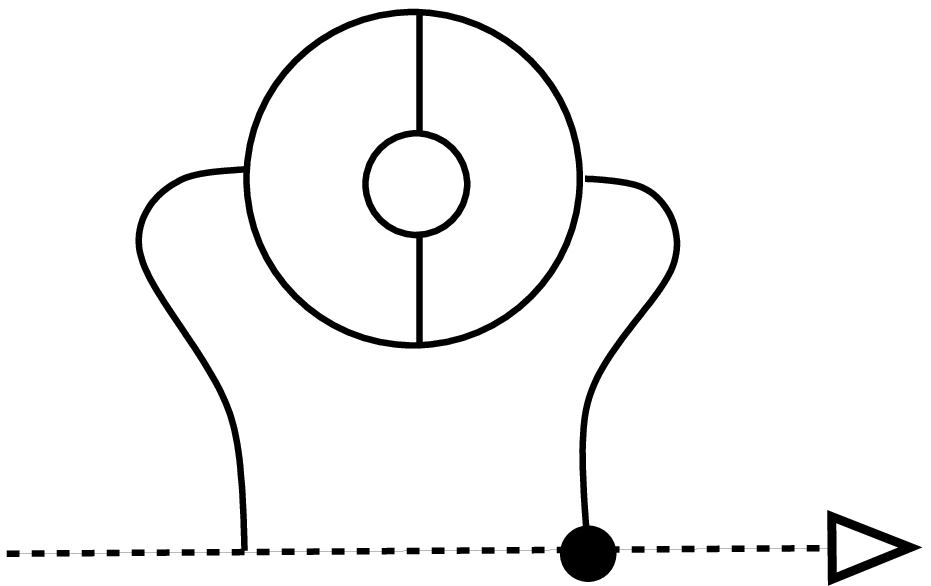}}} \ \
=\ \raisebox{-3ex}{\scalebox{0.25}{\includegraphics{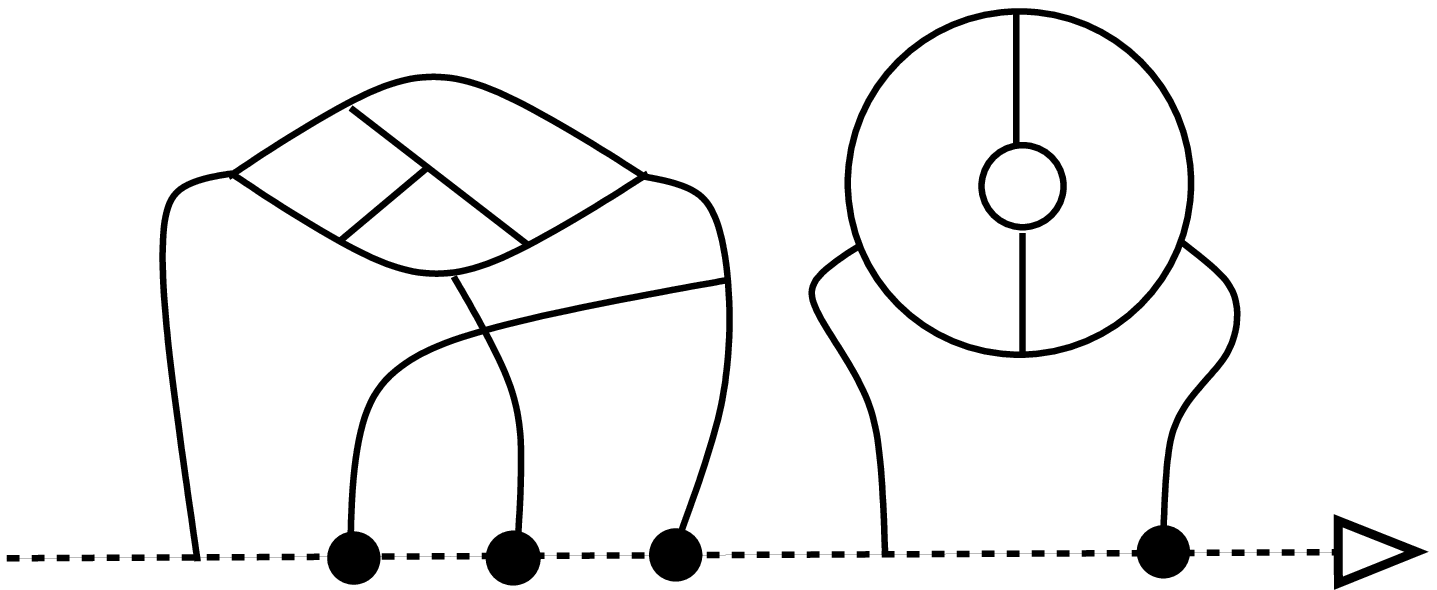}}}\ .
\]

\begin{thm}[Homological Wheeling]
Let $v$ and $w$ be elements of $\Bspace$. Then the two elements of
$\ncw_{\text{basic}}$:
\begin{itemize}
\item{$\left(\chi_\Wspace\circ \Upsilon\right)(v\sqcup w)$} \item{
$\left(\chi_\Wspace\circ \Upsilon\right)(v)\#
\left(\chi_\Wspace\circ \Upsilon\right)(w)$}
\end{itemize}
represent the same basic cohomology class.
\end{thm}

The technical fact which underlies Homological Wheeling is that
the two maps from the basic subcomplex of $\ncw$ to itself:
\[
\xymatrix{ \ncw_{\mathrm{basic}}
\ar@/^/[r]^{\chi_\Wspace\circ\tau} \ar@/_/[r]^{\mathrm{id}} &
\ncw_{\mathrm{basic}} }
\]
are chain homotopic. The homotopy is constructed in Section
\ref{homotopyconstruct}. In words: if $z\in \ncw$ represents some
basic cohomology class, then its {\it graded symmetrization}
represents the same basic cohomology class. In the case at hand,
this fact implies that:
\[
\left(\chi_\Wspace\circ \Upsilon\right)(v)\#
\left(\chi_\Wspace\circ \Upsilon\right)(w)
\]
and
\[
\left(\chi_{\Wspace}\circ\tau\right)\left(\left(\chi_\Wspace\circ
\Upsilon\right)(v)\# \left(\chi_\Wspace\circ
\Upsilon\right)(w)\right) = \left(\chi_\Wspace\circ
\Upsilon\right)(v\sqcup w)
\]
represent the same basic cohomology class, as required.

Before commencing on the formalities, let's compare Wheeling and
Homological Wheeling:
\begin{center}
\begin{tabular}{|p{5cm}|p{5cm}|} \hline
    Wheeling & Homological Wheeling \\ \hline \hline
    The Wheeling map $\partial_\Omega$: Gluing the Duflo element $\Omega$ into legs in all possible ways.
        & The Hair-splitting map $\Upsilon$: Gluing forks into
        legs in all possible ways. \\ \hline
    The averaging map $\chi_\Bspace$. & The graded averaging map
    $\chi_\Wspace$. \\ \hline
    Wheeling says: \mbox{$\left(\chi_\Bspace\circ \partial_\Omega\right)(v\sqcup
    w)$} and $\left(\chi_\Bspace\circ \partial_\Omega\right)(v) \#
    \left(\chi_\Bspace\circ \partial_\Omega\right)(w)$ are equal in $\Aspace$. &
    HW says:
    \mbox{$\left(\chi_\Wspace\circ \Upsilon\right)(v\sqcup w)$} and \mbox{$\left(\chi_\Wspace\circ \Upsilon\right)(v)\#
\left(\chi_\Wspace\circ \Upsilon\right)(w)$} represent the same
basic cohomology class in $\ncw$. \\ \hline
\end{tabular}
\end{center}

\

Homological Wheeling clarifies Wheeling in a number of different
ways:
\begin{itemize}
\item{The map $\Upsilon$ is considerably simpler and less
mysterious than $\partial_\Omega$.} \item{HW has a completely
transparent and combinatorial proof. The underlying mechanisms of
the proof are quite accessible to study and computation.}
\end{itemize}
The cost of HW (at least from the combinatorial point of view) is
that transferring statements in $\ncw$ to statements in $\Aspace$
involves a considerable amount of combinatorial work. In fact, the
bulk of this theory, from Section \ref{gettingwheeling} on into the companion work
\cite{K}, is
taken up with such computation.

\subsection{Acknowledgements}
The author has benefitted greatly from the assistance and
expertise of Dror Bar-Natan and Eckhard Meinrenken, and many
others at the University of Toronto. He also thanks Iain Moffatt
for a useful observation on the first version.

\section{$\Aspace$, $\Bspace$ and the averaging
map.}\label{spacesintro} Apart from some comments about how we
grade diagrams in this work, the material in this section is
standard. The algebras $\Aspace$ and $\Bspace$ are built from
certain graph-theoretic objects called Jacobi diagrams: $\Aspace$
from {\it ordered} Jacobi diagrams and $\Bspace$ from {\it
symmetric} Jacobi diagrams. Here is an example of a symmetric
Jacobi diagram:
\[ \scalebox{0.25}{\includegraphics{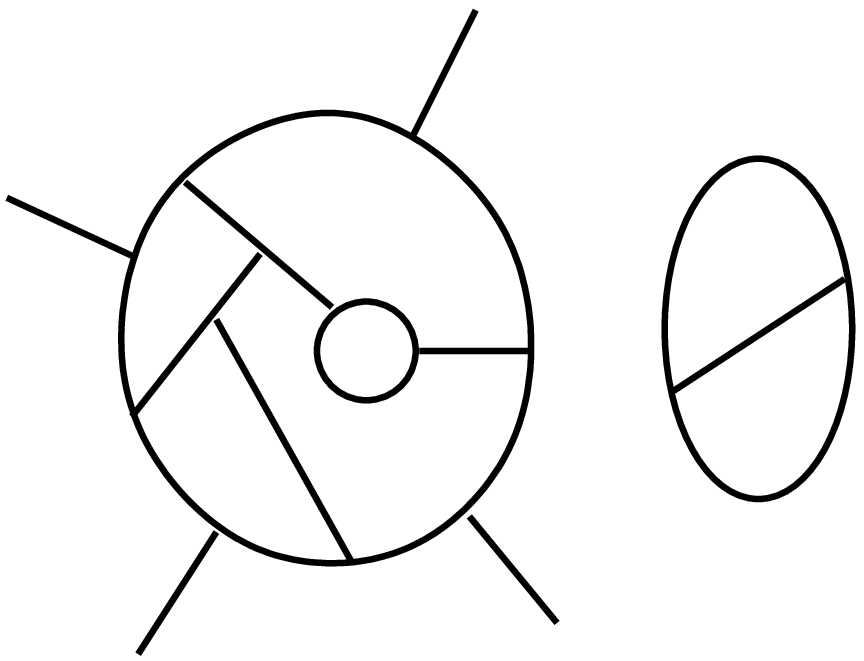}}
\]
What is it? It is first of all a graph with vertices of degree 1 and degree 3, with unoriented edges, and possibly with parallel edges and `loops' (edges both of whose endpoints coincide). In addition, each trivalent vertex with three distinct incoming  edges has an {\it orientation} (a cyclic ordering of its
incident edges). To read the orientation of a vertex from a
drawing of a symmetric Jacobi diagram (as above) simply take the
counter-clockwise ordering determined by the drawing. Two
symmetric Jacobi diagrams are {\bf isomorphic} if there is an
isomorphism of their underlying graphs which respects the
orientations at the trivalent vertices.

\begin{defn}
The space $\Bspace$ is defined to be the rational vector space
obtained by quotienting the vector space of formal finite rational
linear combinations of isomorphism classes of symmetric Jacobi
diagrams by its subspace generated by the anti-symmetry
($\ASreln$) and Jacobi ($\IHXreln$) relations:
\[
\begin{array}{lp{0.1cm}l}
\ASreln: & &
\raisebox{-4ex}{\scalebox{0.25}{\includegraphics{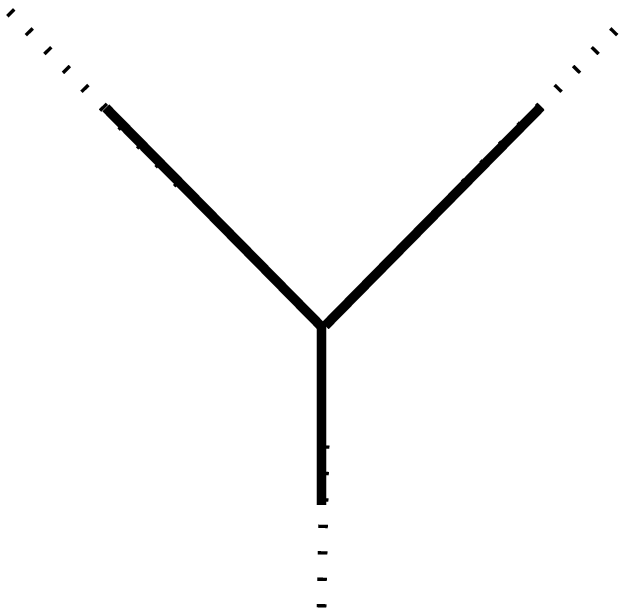}}}+\raisebox{-4ex}{\scalebox{0.25}{\includegraphics{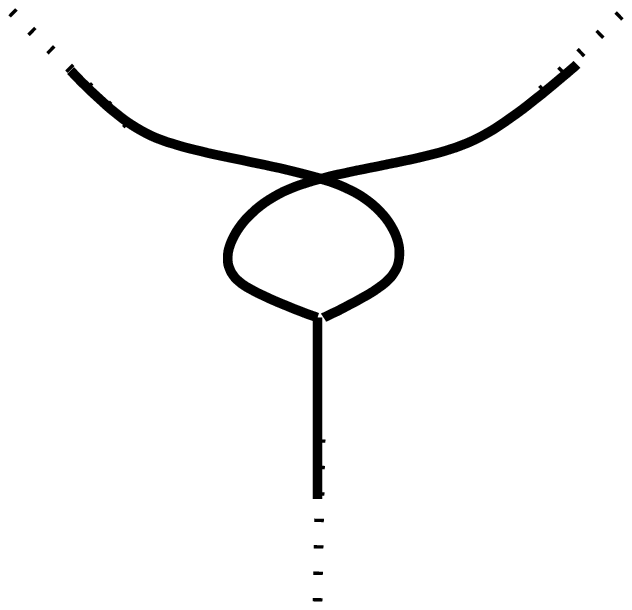}}}\
=\ 0 . \\[3ex]
\IHXreln: & &
\raisebox{-4ex}{\scalebox{0.25}{\includegraphics{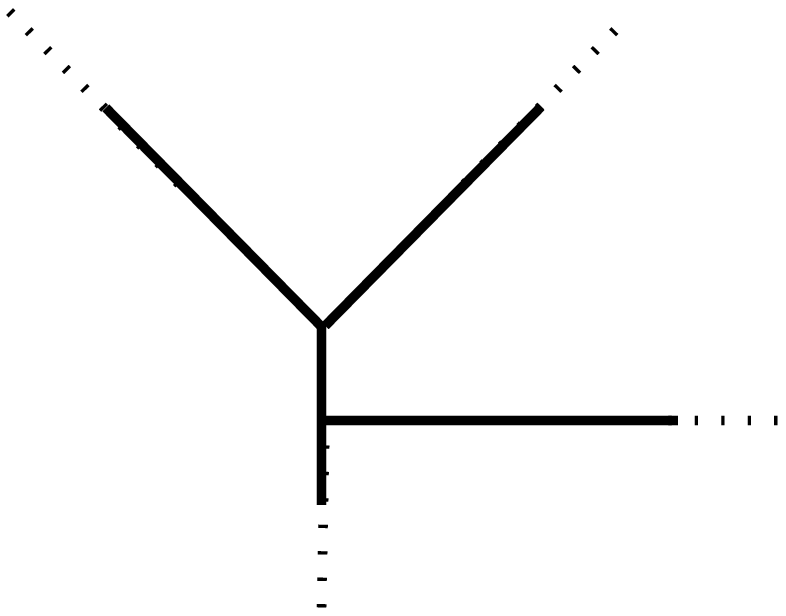}}}-
\raisebox{-4ex}{\scalebox{0.25}{\includegraphics{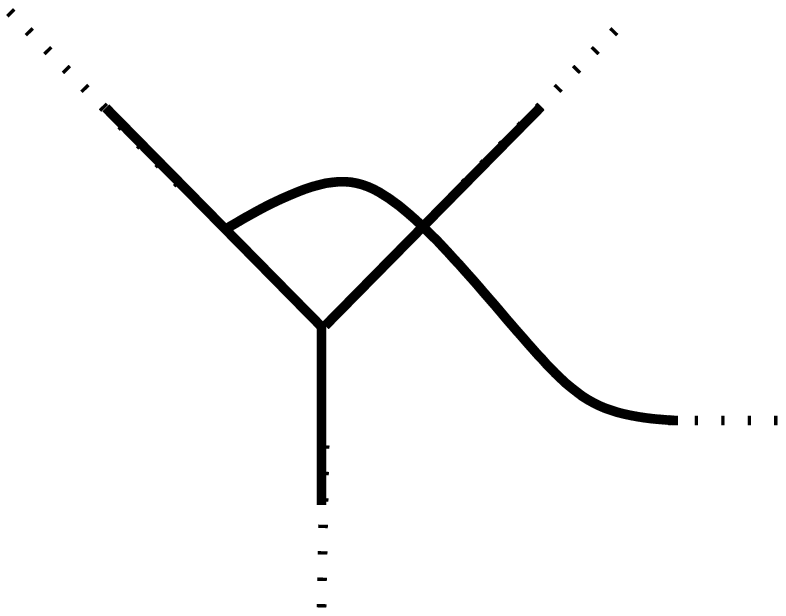}}}-
\raisebox{-4ex}{\scalebox{0.25}{\includegraphics{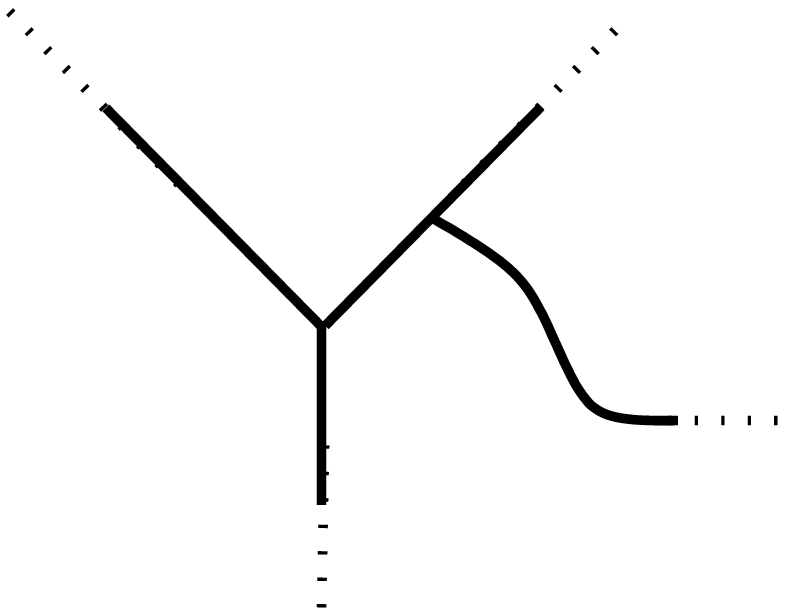}}}\ =\ 0 .
\end{array}
\]
\end{defn}
We remark that the relation which sets a diagram featuring a loop to be zero will be considered an instance of an AS relation.

Next let's consider the space $\Aspace$. An {\bf ordered Jacobi diagram}
is the same as a symmetric Jacobi diagram but has an extra piece
of structure: the set of vertices of degree 1 has been ordered.
\begin{defn}
The space $\Aspace$ is defined to be the rational vector space
obtained by quotienting the vector space of formal finite rational
linear combinations of isomorphism classes of ordered Jacobi
diagrams by its subspace generated by the anti-symmetry
($\ASreln$), Jacobi ($\IHXreln$) and permutation ($\STUreln$)
relations:
\[
\begin{array}{lp{0.1cm}l}
\STUreln: & &
\raisebox{-4ex}{\scalebox{0.23}{\includegraphics{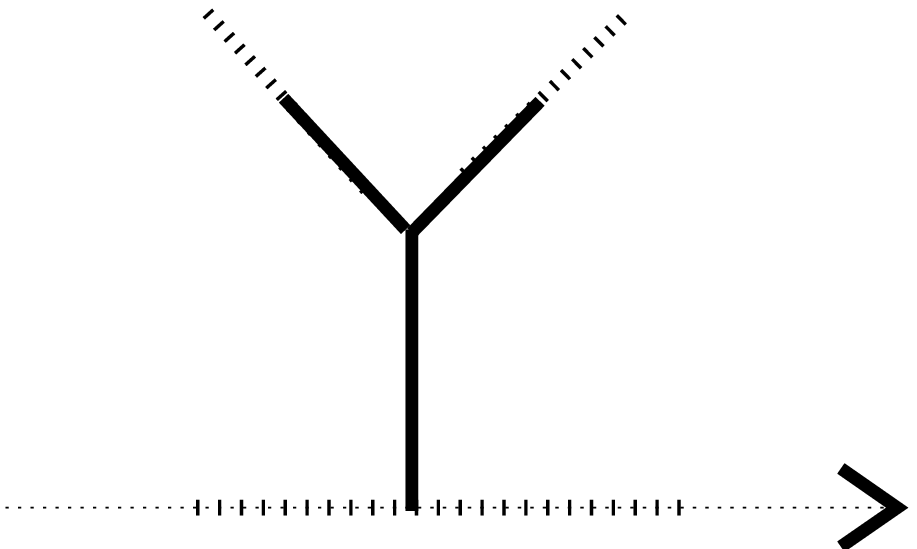}}}\ -\
\raisebox{-4ex}{\scalebox{0.23}{\includegraphics{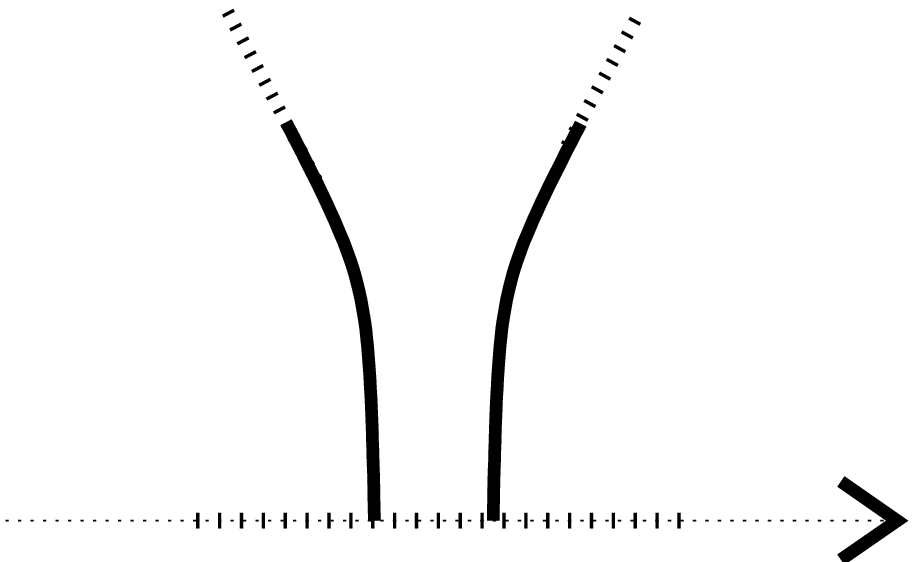}}}\ -\
\raisebox{-4ex}{\scalebox{0.23}{\includegraphics{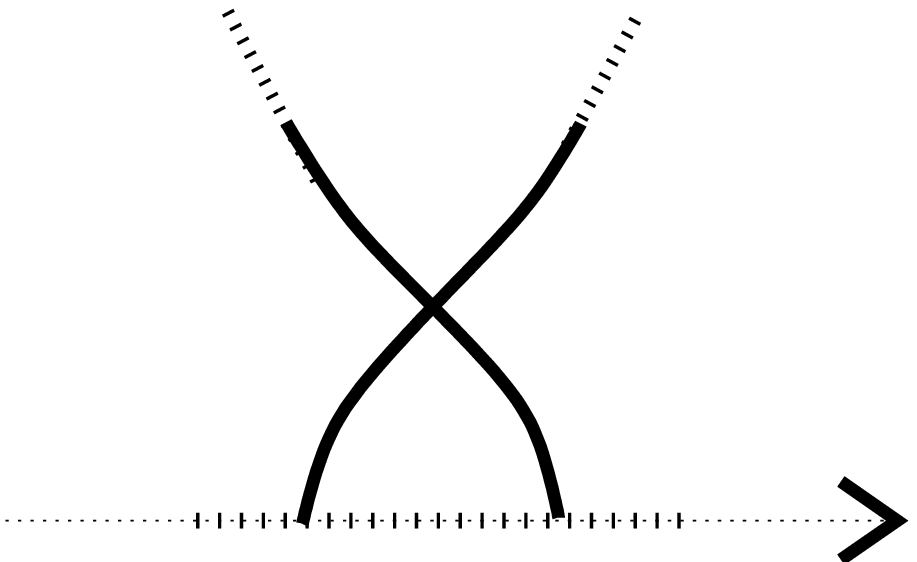}}} \ =\ 0\ .
\end{array}
\]
\end{defn}

There is a natural map between these vector spaces, the {\bf
averaging map}\footnote{Remark: in this work we require
averaging maps between a variety of different spaces. To keep the
notation logical we will record the domain space as a subscript on
the symbol $\chi$.} $\chi_\Bspace$, which was already recalled in Section
\ref{sketch}. The proof of the following statement can be
read in \cite{BarNatan}:
\begin{fpbw}
The averaging map is a vector space isomorphism.
\end{fpbw}

These spaces have natural product structures. The product on
$\Aspace$ (we'll refer to it as the ``juxtaposition product")
looks like
\[
\raisebox{-2.5ex}{\scalebox{0.24}{\includegraphics{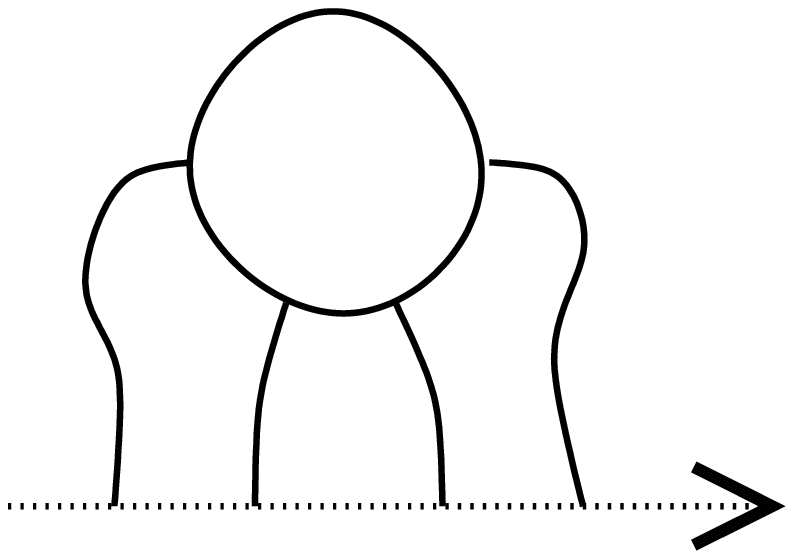}}}\
\#\
\raisebox{-2.5ex}{\scalebox{0.24}{\includegraphics{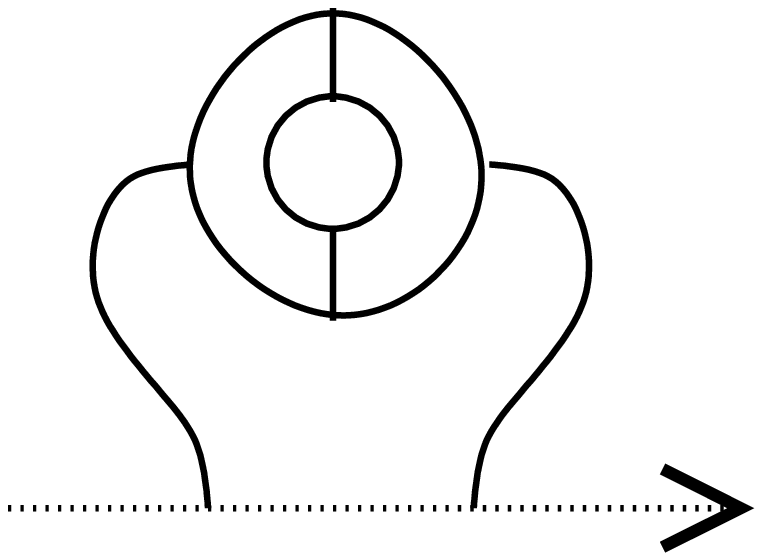}}} =
\raisebox{-2.5ex}{\scalebox{0.24}{\includegraphics{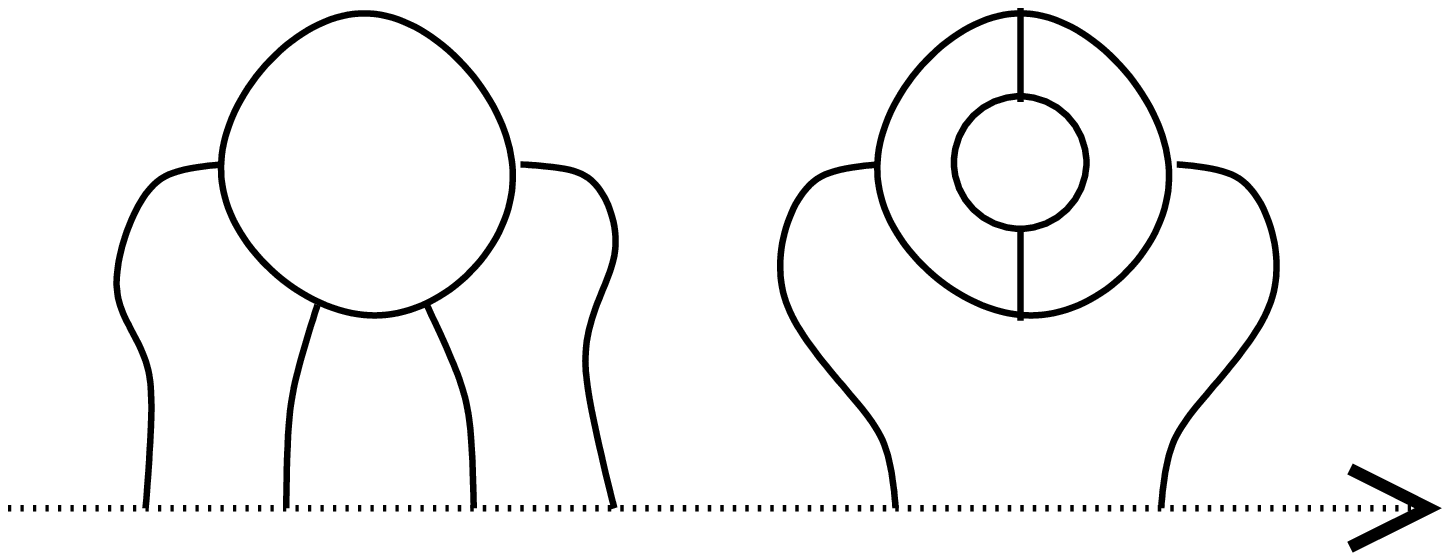}}}\ .
\]
The product on $\Bspace$ (the ``disjoint union" product) looks
like
\[
\raisebox{-2ex}{\scalebox{0.28}{\includegraphics{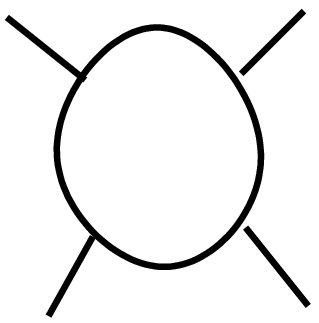}}}\
\sqcup\
\raisebox{-1.75ex}{\scalebox{0.28}{\includegraphics{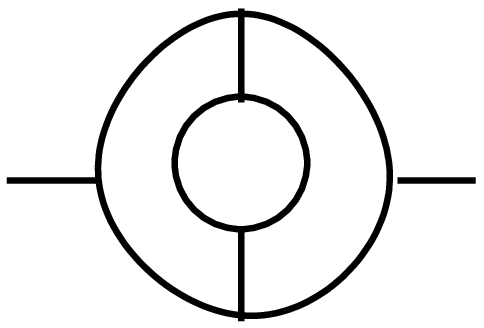}}}\ =\
\raisebox{-4ex}{\scalebox{0.23}{\includegraphics{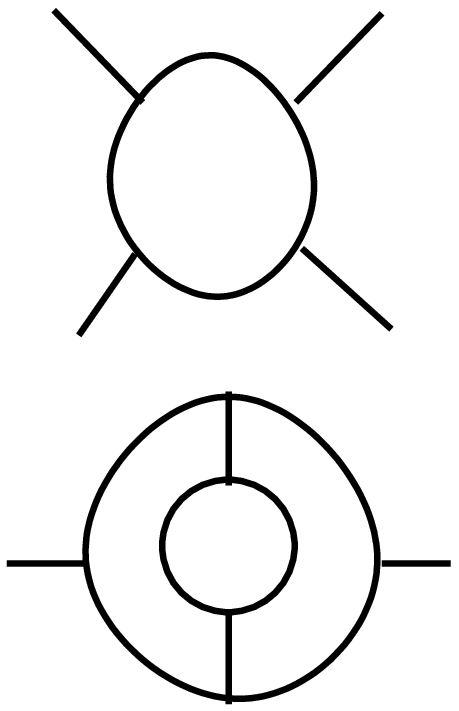}}}\ \ .
\]
Famously, while the averaging map {\it is} a vector space
isomorphism it is {\it not} an isomorphism of algebras. But if the
averaging map is {\it preceded} by a certain vector space
isomorphism of $\Bspace$, $\partial_\Omega : \Bspace \to \Bspace$,
then one {\it does} obtain an isomorphism of algebras:
$(\chi_\Bspace\circ\partial_\Omega)(a\sqcup b) =
(\chi_\Bspace\circ\partial_\Omega)(a) \#
(\chi_\Bspace\circ\partial_\Omega)(b).$
This, of course, is the Wheeling theorem.
Our next task, then, is to recall the isomorphism
$\partial_\Omega$.

\subsection{Operating with diagrams on diagrams.}

Let $D_1$ and $D_2$ be two symmetric Jacobi diagrams. The notation
$\partial_{D_1}(D_2)$ will denote the result of {\it operating}
with $D_1$ on $D_2$. That is,  $\partial_{D_1}(D_2)$ is defined to
be the sum of all the possible symmetric Jacobi diagrams that can
be obtained by gluing {\it all} of the legs of $D_1$ to {\it some}
of the legs of $D_2$.

Let's say that again with more combinatorial precision. Let $D_1$
have $n$ legs, and number them in some way. Let $D_2$ have $m$
legs, and number them in some way. For some injection
$\sigma: \{1,\ldots,n\} \mapsto \{1,\ldots,m\}$
let $D_\sigma$ denote the diagram that results when the legs of
$D_1$ are joined to the legs of $D_2$ according to the map
$\sigma$. For example, if
\[
D_1 =
\raisebox{-1.6ex}{\scalebox{0.25}{\includegraphics{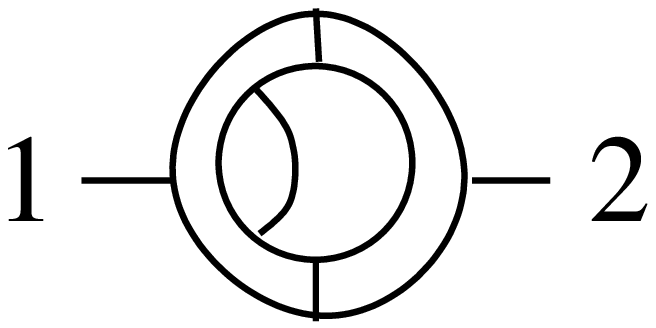}}}\ \
\text{and}\ \ D_2 =
\raisebox{-2.9ex}{\scalebox{0.25}{\includegraphics{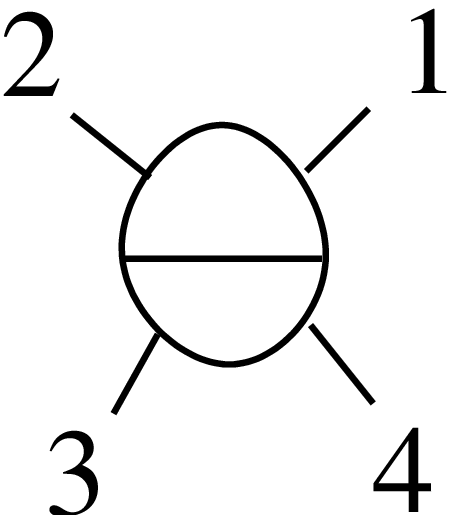}}}\ \
\text{then}\ \ D_{\left(\substack{1 2 \\ 3 2}\right)} =
\raisebox{-6.5ex}{\scalebox{0.25}{\includegraphics{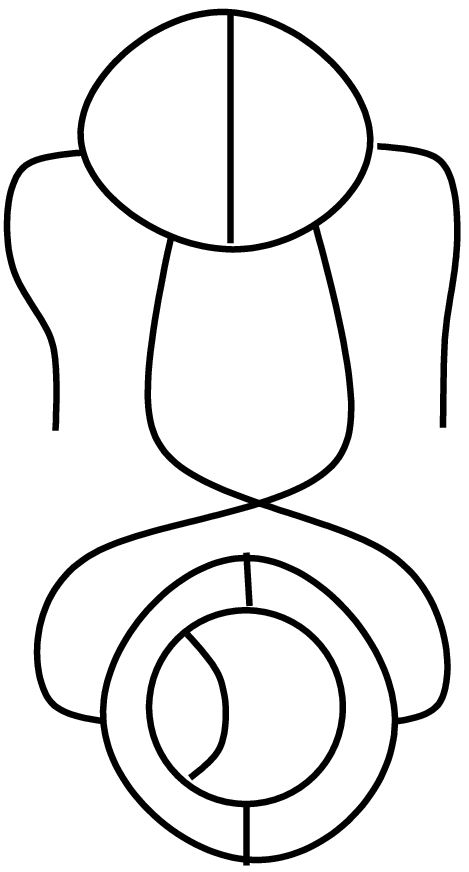}}}\ \
.
\]
In these terms the diagram operation can be defined:
\[
\partial_{D_1}(D_2) = \sum_{\substack{\text{injections}\
\sigma \\
\{1,\ldots,n\}\to\{1,\ldots,m\}} } D_\sigma\ .
\]
(Note that if $n>m$ then the result is zero.) We can extend this
operation to the case when $D_1$ and $D_2$ are finite linear
combinations of symmetric Jacobi diagrams in an obvious way.

The wheeling isomorphism $\partial_\Omega$ is obtained by
operating with a certain remarkable {\it power series} of
diagrams, $\Omega$. Before introducing this power series we must
say a few words about gradings.

\subsection{Gradings and $\Omega$}

It is traditional to grade the spaces $\Aspace$ and $\Bspace$ by
half the number of vertices in the diagram. This grading has no
use in this paper, so, to avoid complication, we will ignore it.

On the other hand the {\it leg-grade} of a diagram will play a
crucial role.
\begin{defn}
Define the {\bf leg-grade} of a symmetric Jacobi diagram to be
twice the number of legs of the diagram.
\end{defn}
The reason for the weight of $2$-per-leg will become clear very
soon. If $\Bspace^n$ denotes the subspace of leg-grade $n$ then we
have a direct-sum decomposition $
\Bspace = \bigoplus_{n=0}^{\infty} \Bspace^n$.

When we refer to a {\bf power series $\Gamma$ of symmetric Jacobi
diagrams} we are referring to an element of the direct product
vector space $\Gamma \in \prod_{n=0}^{\infty} \Bspace^n$.
Note that if $\Gamma$ is such a power series, then we get a
well-defined map $\partial_\Gamma : \Bspace \to \Bspace$
because all but finitely many of the terms of $\Gamma$ contribute
zero.

To introduce $\Omega$ we will employ a
notation that gets used in disparate works in the
literature. If we draw a diagram, orient an edge of the diagram at
some point, and label that point with a power series in some
formal variable $a$, then:
\begin{multline*}
\raisebox{-6ex}{\scalebox{0.21}{\includegraphics{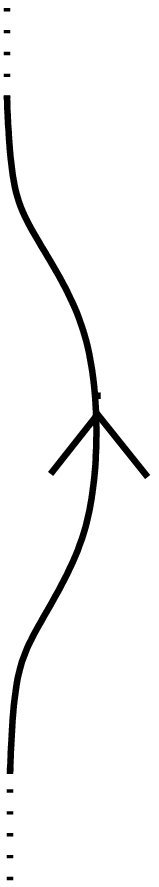}}} \
c_0 +  c_1a + c_2a^2 + c_3a^3 + \ldots  \\  \ \ \mbox{denotes}\ \
\ \ \ \ c_0\
\raisebox{-6ex}{\scalebox{0.21}{\includegraphics{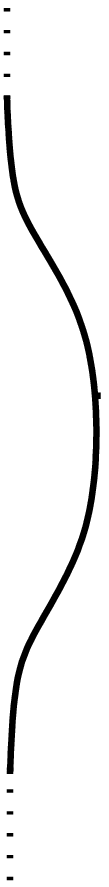}}}\
\ +\ \  c_1\
\raisebox{-6ex}{\scalebox{0.21}{\includegraphics{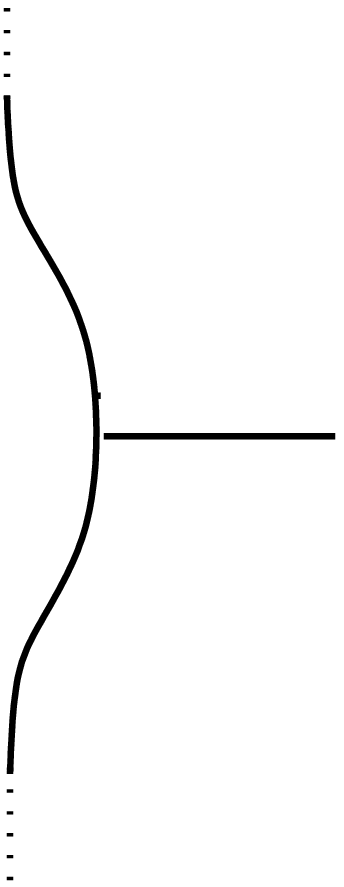}}} \
\ +\ \  c_2\
\raisebox{-6ex}{\scalebox{0.21}{\includegraphics{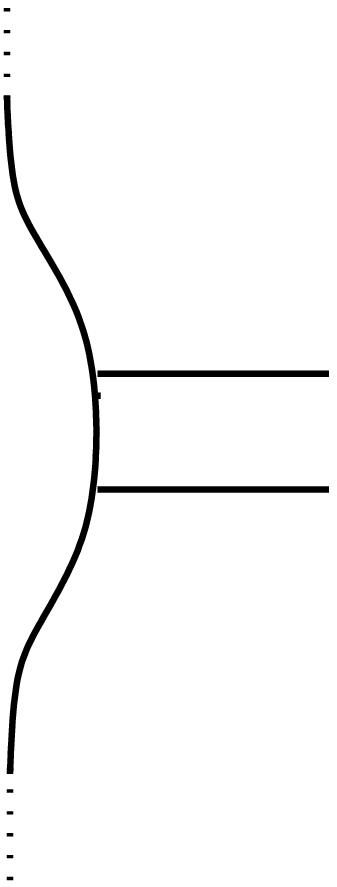}}} \
\ +\ \ c_3\
\raisebox{-6ex}{\scalebox{0.21}{\includegraphics{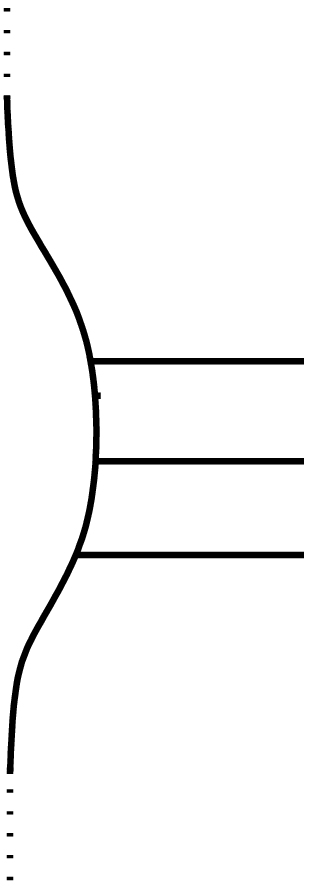}}} \
\ +\ \ldots\ \ \ \ .
\end{multline*}
\begin{defn}
The Wheels element, $\Omega$, is the formal power series of
symmetric Jacobi diagrams defined by the expression
\[
\Omega = \exp_{\sqcup}\left(\frac{1}{2}\
\raisebox{-3.4ex}{\scalebox{0.28}{\includegraphics{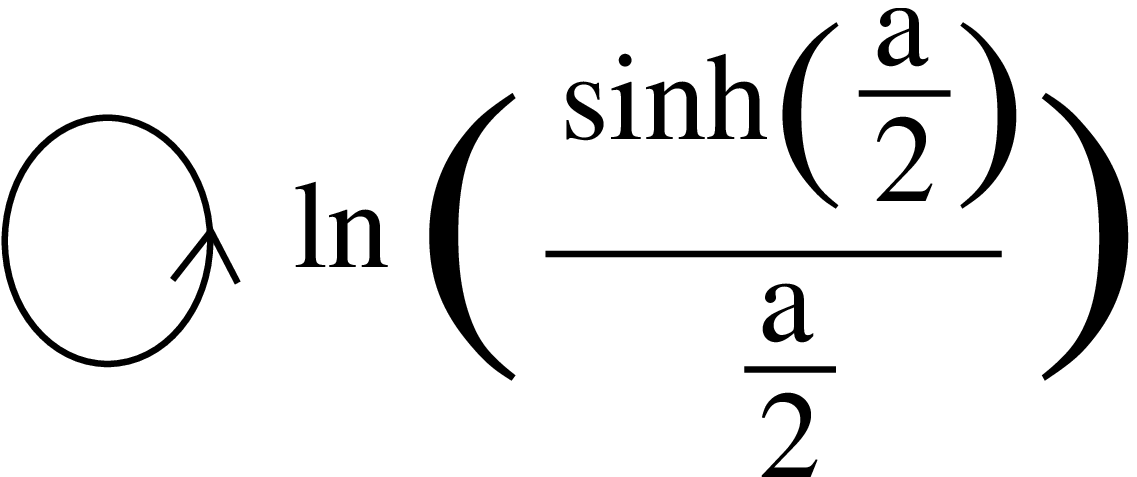}}}
\right)\ \in\ \prod_{n=0}^{\infty}\Bspace^n\ .
\]
\end{defn}
The combinatorial mechanisms that produce this element are the
main subject of the companion work \cite{K}.

\section{The commutative Weil complex for
diagrams}\label{commweil}
The purpose of this section is the introduction of a certain
cochain complex $\Wbase$. The next section will describe a certain
map of cochain complexes (with the spaces $\Bspace^i$ assembled
into a complex with zero differential),
$\Upsilon^i_{\mathrm{basic}} : \Bspace^i \rightarrow \Wbasei$,
which gives isomorphisms in cohomology\footnote{A remark
concerning our notation: if $f$ is a chain map between cochain
complexes then $H(f)$ (or sometimes $H^i(f)$) will denote the
induced map on cohomology.}:
\[
H^i(\Upsilon_{\mathrm{basic}})\,:\,\Bspace^i\cong H^i(\Bspace)
\stackrel{\cong}{\longrightarrow} H^i(\Wbase).
\]


Our goal is a clean presentation of the combinatorial
constructions. These constructions may appear quite mysterious at
first glance, however, so to provide some context to the theory,
we'll close the section by constructing a ``characteristic
class-valued evaluation map". Given the data of a compact Lie group
$G$, a smooth principal $G$-bundle $P$ with corresponding base $B=P\backslash G$, and a connection form $A$
on $P$, we'll construct a map of complexes
$\text{Eval}^{(G,P,A)} : \Wspace \rightarrow \Omega(P)$
(this map {\it evaluates} Weil diagrams in $\Omega(P)$)
which yields a map
\[
H\left(\text{Eval}^{(G,P,A)}_{\mathrm{basic}}\right) :
\Bspace\cong H(\Wbase) \rightarrow
H(\Omega(P)_{\mathrm{basic}})\cong H_{dR}(B).
\]

\subsection{Weil diagrams}\label{commweildiags}
The distinguishing feature of Weil diagrams is that their legs are
{\it graded}. Every leg of a Weil diagram is either of grade 1 or
of grade 2. The grade 2 legs are distinguished by drawing them with a fat dot.
The {\bf leg-grade} of a Weil diagram is the total grade of the
legs. Thus, the following example has leg-grade 6:
\begin{center}
\raisebox{-1.4ex}{\scalebox{0.3}{\includegraphics{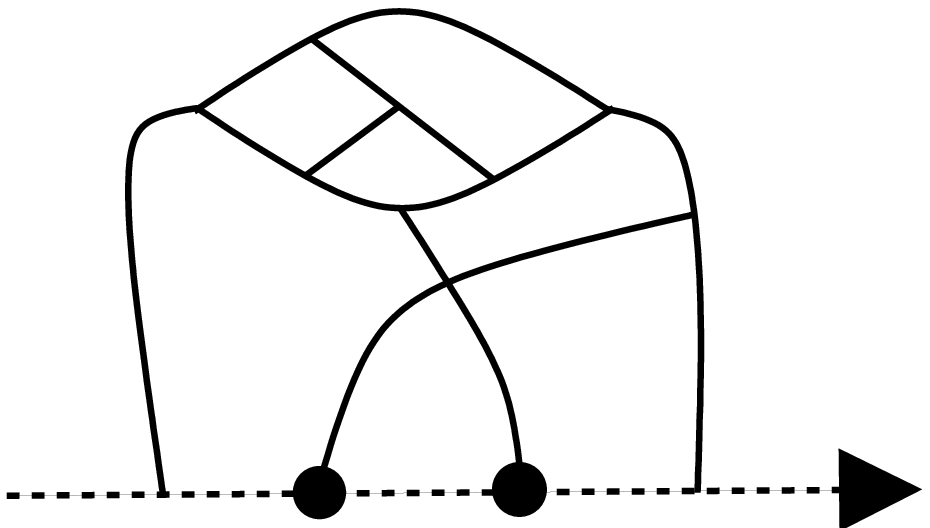} }}\
\end{center}

The precise definition of {\bf Weil diagram} is unsurprising:
it is a graph with vertices of degree 1 and 3, where each degree 3
vertex is oriented, and the set of degree 1 vertices $U$ is
ordered and weighted by a grading $U\to
\{1,2\}$.

Occasionally a ``Weil
diagram with a special degree 1 vertex" is used, which is just a Weil
diagram except that exactly one of the degree 1 vertices is
neither included in the total ordering nor assigned a leg-grading.
\subsection{Permuting the legs of Weil diagrams}
When constructing spaces from these diagrams, we'll always employ
the familiar AS and IHX relations. The numerous spaces we employ
will differ from each other according to how they treat the legs
of the diagrams. As there are really quite a number of different
sets of relations applied to legs in this work we'll remind the
reader of which space we are in by drawing the arrow on the
orienting line in a different style depending on which space the
diagram is in.

To begin: consider relations which say that we are allowed to move
the legs around freely, as long as when we transpose an adjacent
pair of legs we introduce the sign $(-1)^{xy}$, where $x$ and $y$
are the grades of the two legs involved:
\[
\begin{array}{lrcl}
{\mbox{Perm}}_1:\hspace{4ex} &
\raisebox{-2.75ex}{\scalebox{0.27}{\includegraphics{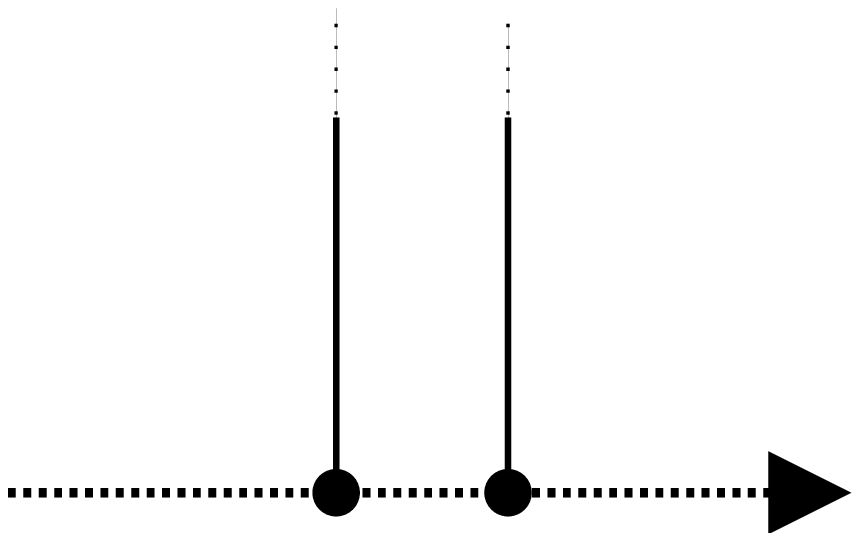}}}\ -\
\raisebox{-2.75ex}{\scalebox{0.27}{\includegraphics{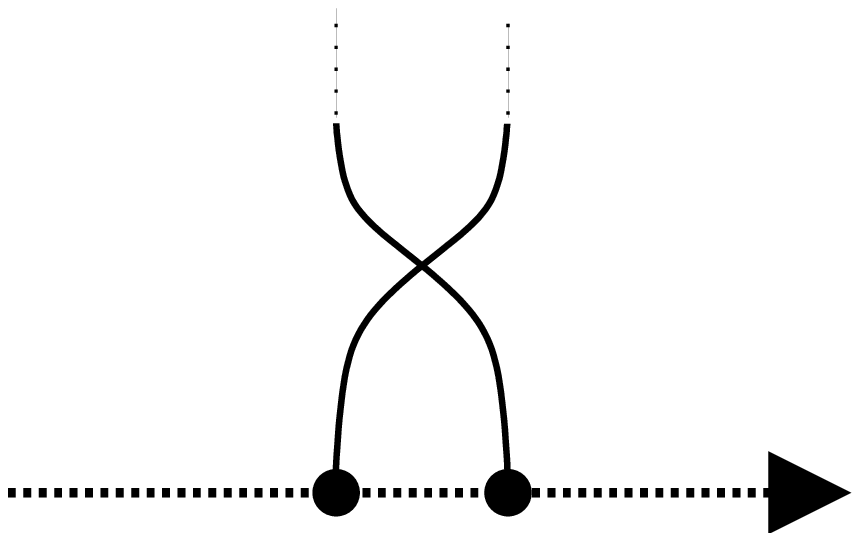}}} & = & 0 \\
{\mbox{Perm}}_2:\hspace{4ex} &
\raisebox{-2.75ex}[1\height]{\scalebox{0.27}{\includegraphics{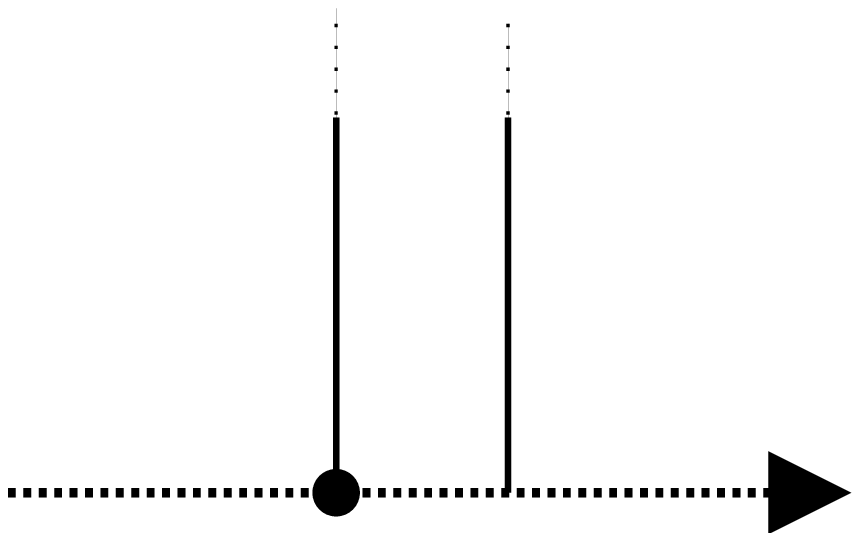}}}
\ -\
\raisebox{-2.75ex}{\scalebox{0.27}{\includegraphics{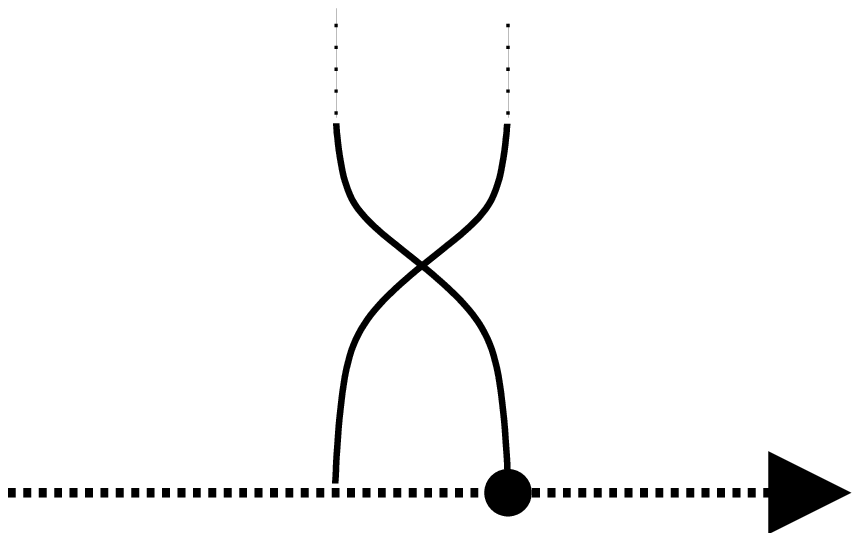}}} & = & 0 \\
{\mbox{Perm}}_3:\hspace{4ex} &
\raisebox{-2.75ex}[1\height]{\scalebox{0.27}{\includegraphics{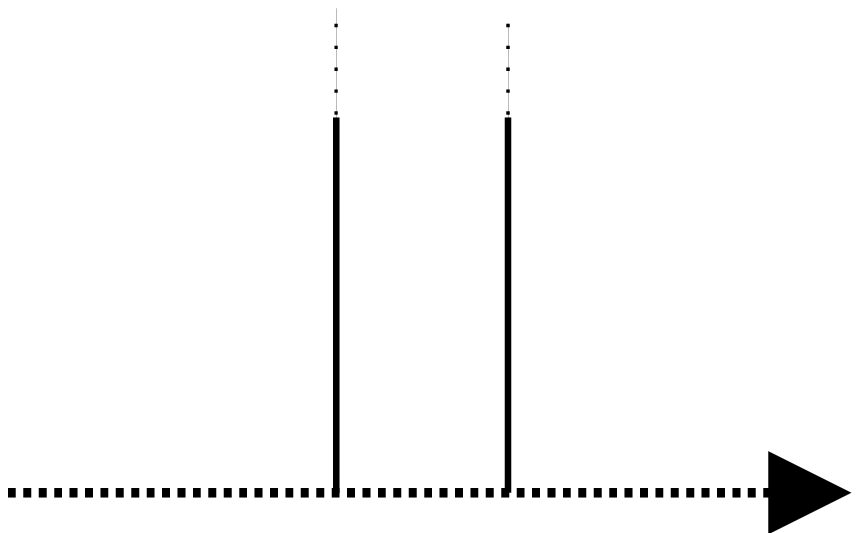}}}
\ +\ \raisebox{-2.75ex}{\scalebox{0.27}{\includegraphics{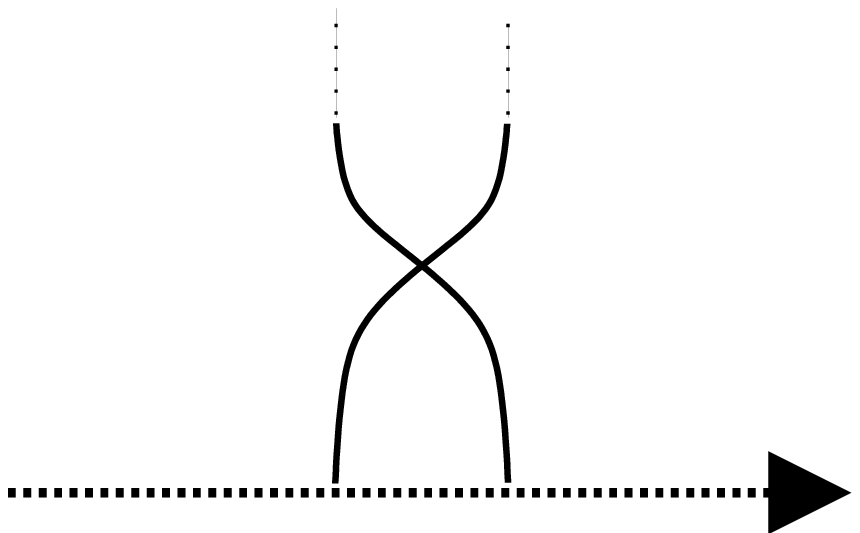}}}
& = & 0
\end{array}
\]

\begin{defn}
The vector space $\Wspace^i$ is
defined to be the quotient of the vector space of formal finite
$\mathbb{Q}$-linear combinations of isomorphism classes of Weil
diagrams of leg-grade $i$ by the subspace generated by AS, IHX and
Perm relations. If $i<0$ then $\Wspace^i$ is understood to be zero.
\end{defn}

\subsubsection{A clarification concerning $\Wspace^0$.}
We regard the empty diagram as a Weil diagram. Thus, for example,
\[
\frac{1}{2} -
\frac{7}{6}\,\raisebox{-2ex}{\scalebox{0.24}{\includegraphics{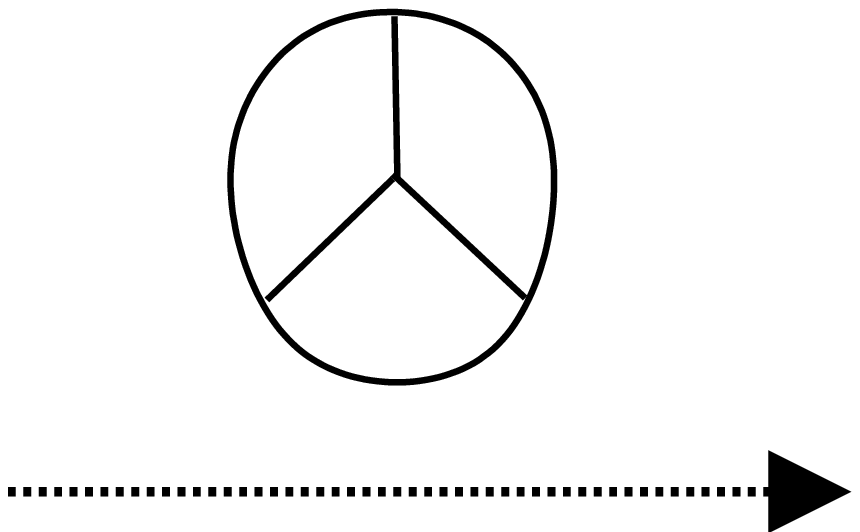}}}
\]
is an element of $\Wspace^0$.

To assemble these vector spaces into a complex we first need to
learn how to operate on them with ``differential operators".

\subsection{Formal linear differential operators.}\label{Fldo}

A formal linear differential operator of grade $j$ is a sequence
$\{\mathcal{D}^0, \mathcal{D}^1, \ldots \}$ of linear maps
$\mathcal{D}^i: \mathcal{W}^i \rightarrow \mathcal{W}^{i+j}$
defined by a pair of substitution
rules. The substitution rules are specified by the following data:
$X$, a Weil diagram of leg-grade $j+1$ with a special degree 1
vertex, and $Y$, a Weil diagram of leg-grade $j+2$ with a special
degree 1 vertex. More generally $X$ and $Y$ can be formal finite combinations of such diagrams.

Given $w$, a Weil diagram of leg-grade $i$, the evaluation of
$\mathcal{D}^i(w)$ begins by placing a copy of the word
``$\mathcal{D}($" on the far left-hand end of the base vector and
a copy of the symbol ``$)$" on the far right-hand end of the base
vector.
 The word ``$\mathcal{D}($" is then pushed towards the
``$)$" by the repeated application of the two substitution rules:
\begin{eqnarray*}
\raisebox{-3ex}{\scalebox{0.3}{\includegraphics{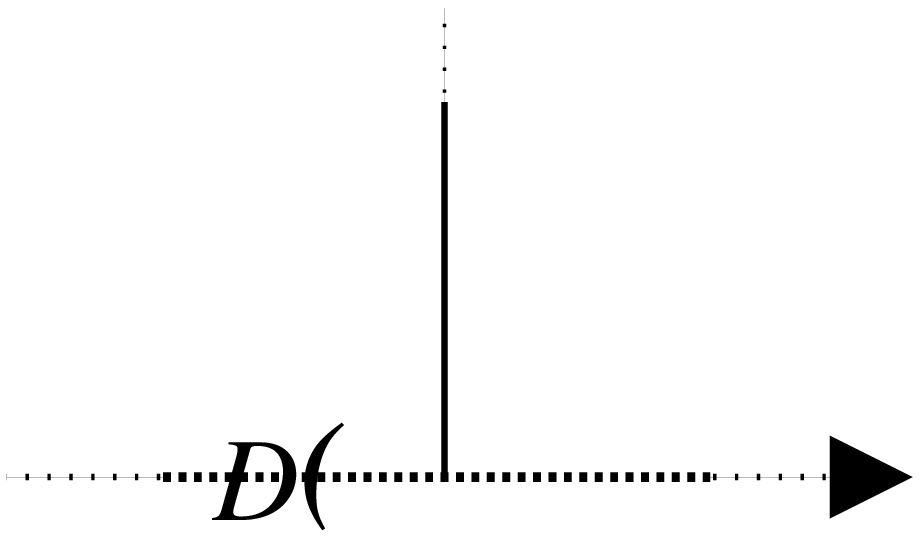}}} &
\leadsto &
\raisebox{-3ex}{\scalebox{0.3}{\includegraphics{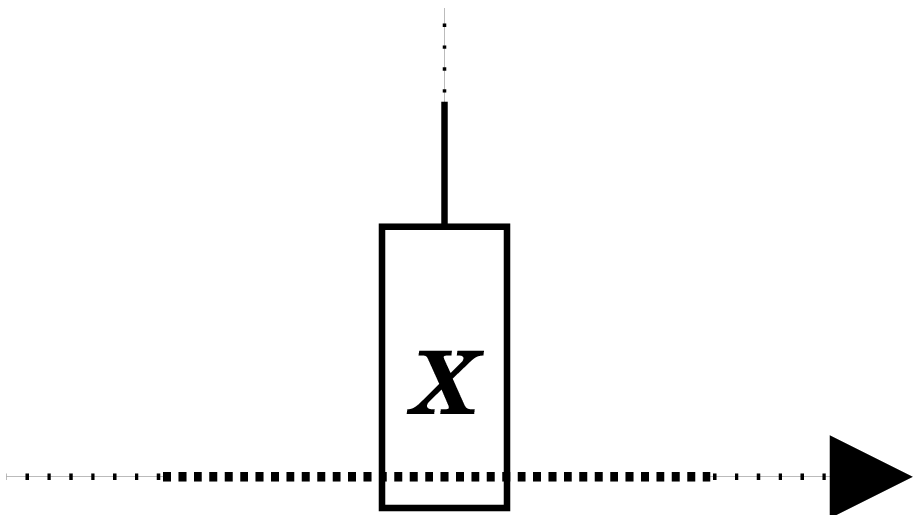}}}\
+(-1)^{1j}\
\raisebox{-3ex}{\scalebox{0.3}{\includegraphics{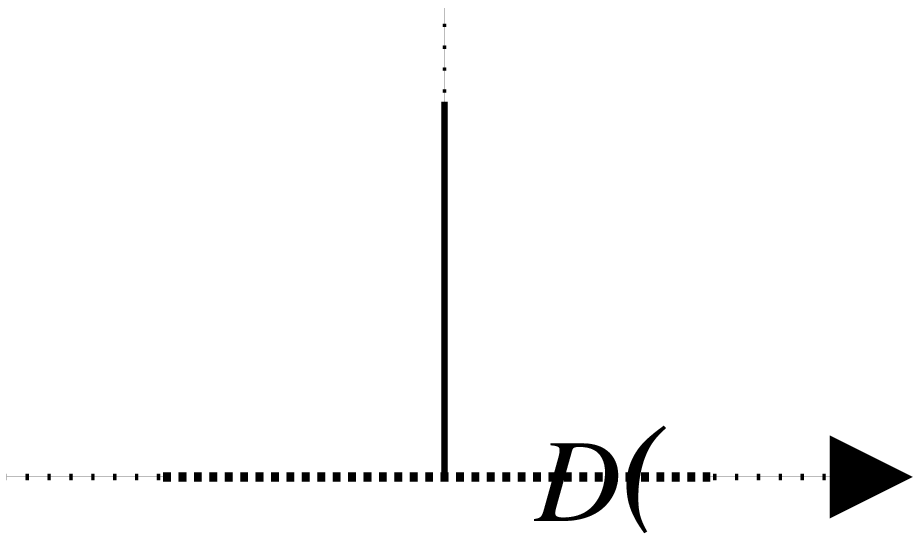}}}\
,
\\[0.15cm]
\raisebox{-3ex}{\scalebox{0.3}{\includegraphics{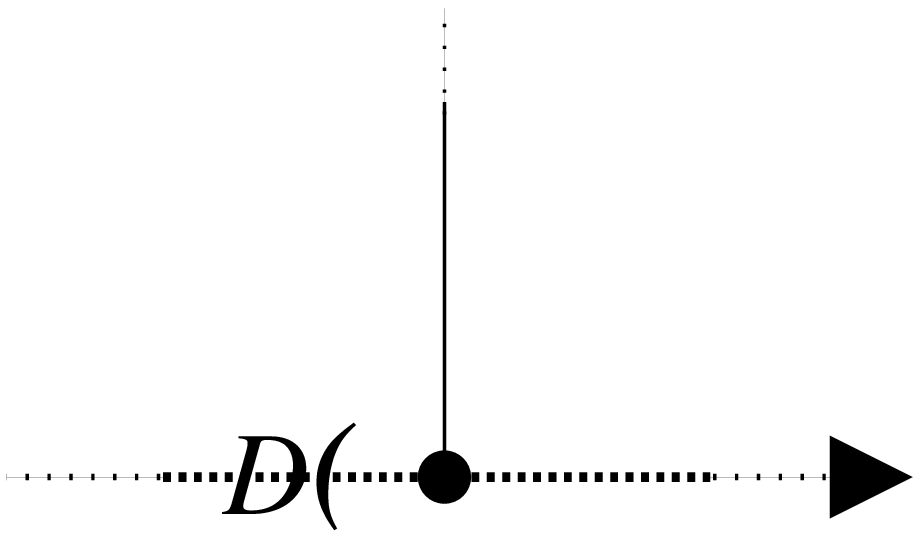}}} &
\leadsto &
\raisebox{-3ex}{\scalebox{0.3}{\includegraphics{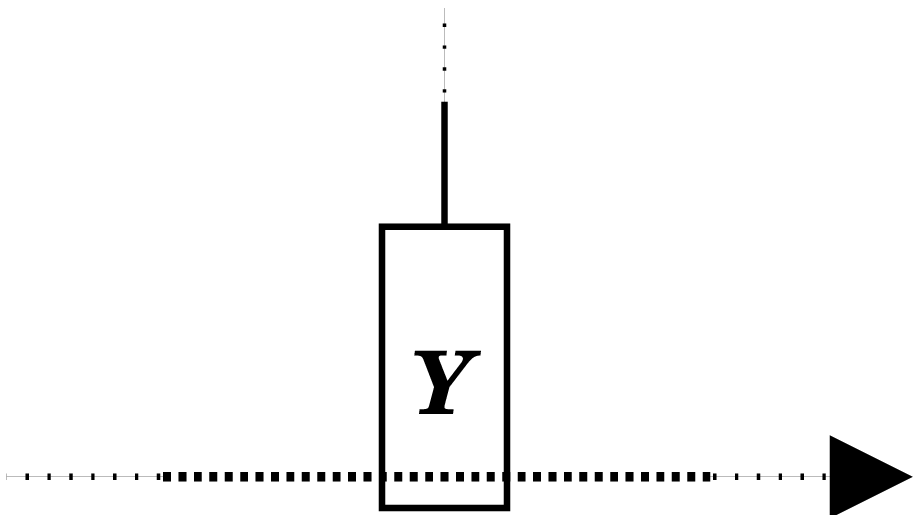}}}\
+(-1)^{2j}\
\raisebox{-3ex}{\scalebox{0.3}{\includegraphics{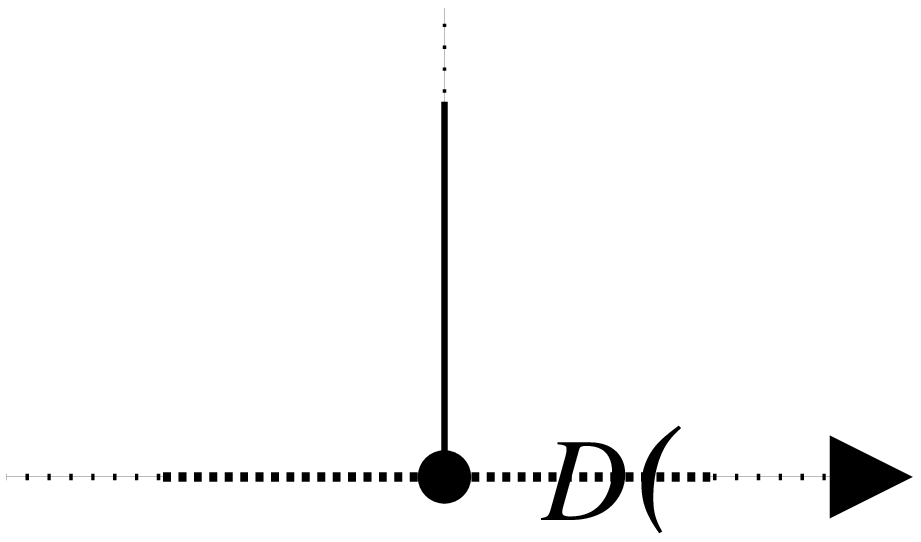}}}\ .\\[0.15cm]
\end{eqnarray*}
When the word ``$\mathcal{D}($" reaches the symbol ``$)$" then
that diagram is eliminated and the procedure terminates.
In some situations the
closing bracket appears earlier along the orienting
line.
In the typical situation, when the closing bracket appears
on the far-right of the orienting line, we'll omit the brackets
altogether. The reader can check that this operation respects
the {Perm} relations.

\subsection{Example: the differential.}\label{diff} 
\begin{figure}
\caption{An example of the calculation of $d$. \label{dexamp}}
\parbox{12cm}{
\begin{eqnarray*}
\lefteqn{\raisebox{-3ex}{\scalebox{0.3}{\includegraphics{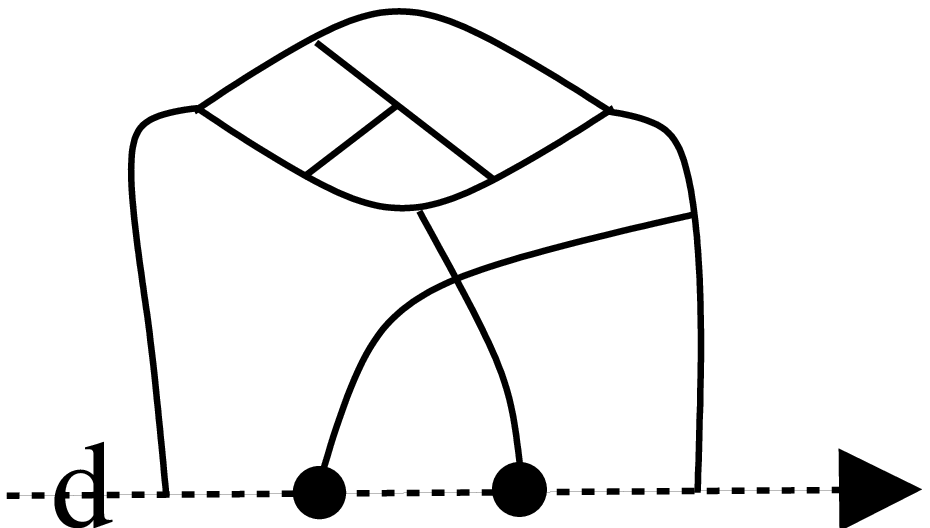}}}  \leadsto
\raisebox{-3ex}{\scalebox{0.3}{\includegraphics{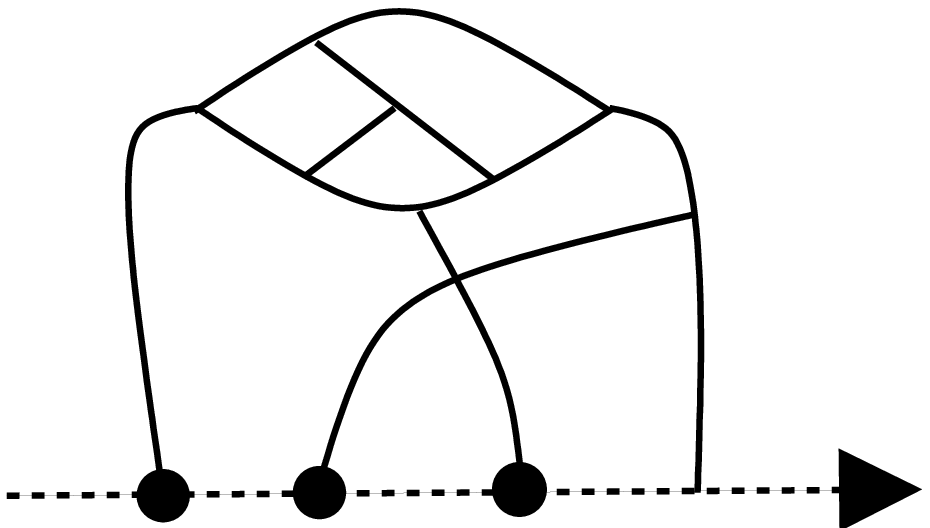}}}\ -\
\raisebox{-3ex}{\scalebox{0.3}{\includegraphics{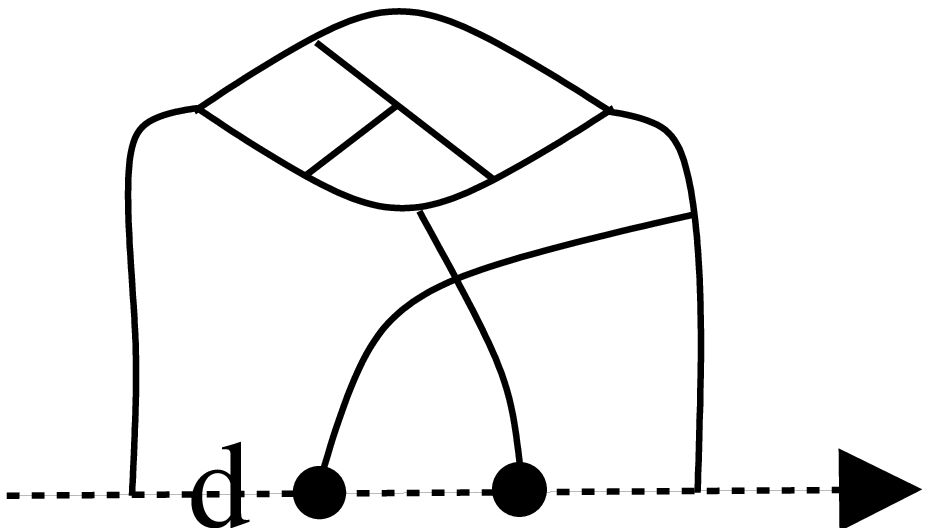}}}}
\\
& \leadsto &
\raisebox{-3ex}{\scalebox{0.3}{\includegraphics{diffexampAA}}}\ -\
\raisebox{-3ex}{\scalebox{0.3}{\includegraphics{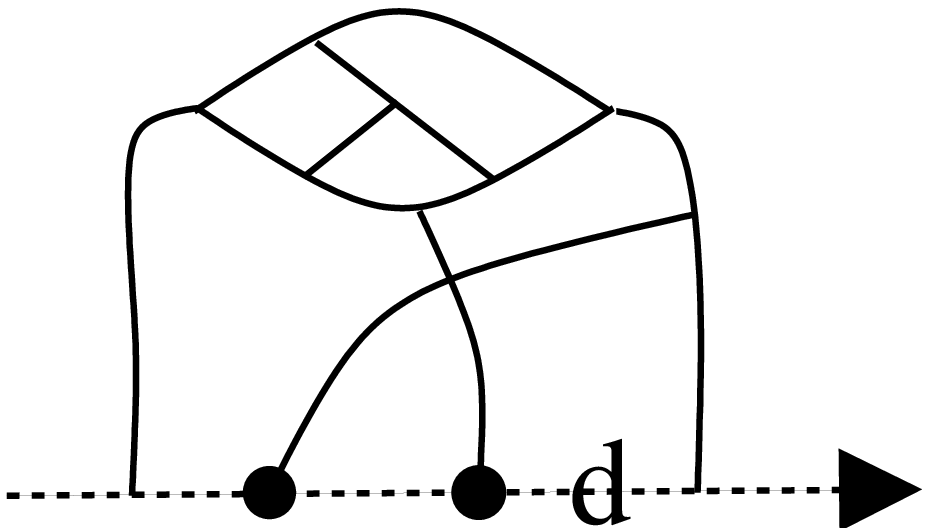}}}\ \\
& \leadsto &
\raisebox{-3ex}{\scalebox{0.3}{\includegraphics{diffexampAA}}}\ -\
\raisebox{-3ex}{\scalebox{0.3}{\includegraphics{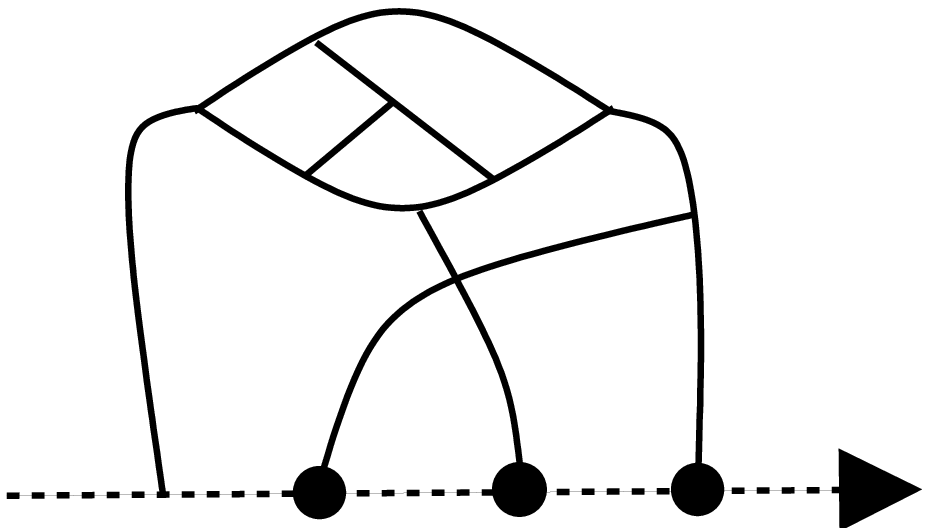}}}\ +
\raisebox{-3ex}{\scalebox{0.3}{\includegraphics{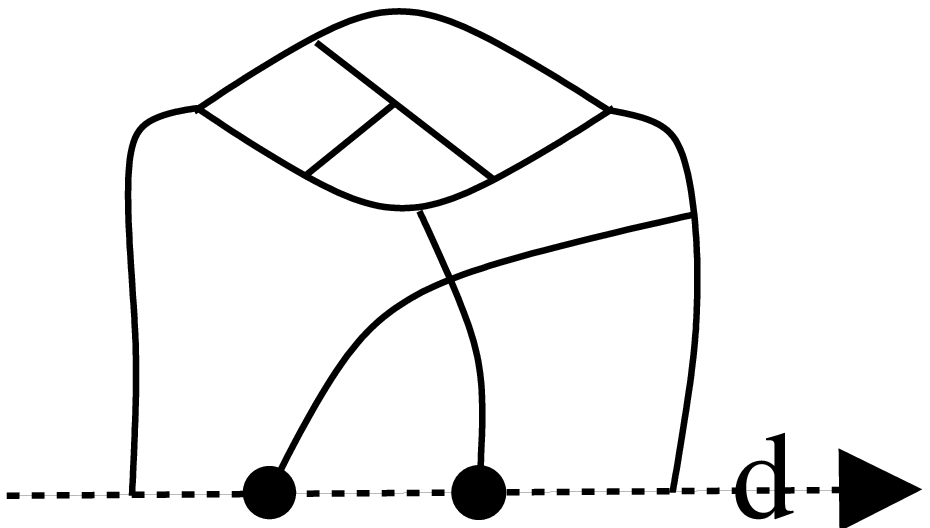}}}\ \ .
\end{eqnarray*}
Thus:
\[
d^{6}\left(
\raisebox{-3ex}{\scalebox{0.3}{\includegraphics{diffexamp}}}
\right) =
\raisebox{-3ex}{\scalebox{0.3}{\includegraphics{diffexampAA}}} \
-\ \raisebox{-3ex}{\scalebox{0.3}{\includegraphics{diffexampAB}}}\
\ .
\]}
\underline{\hspace{7cm}}
\end{figure}

Let us illustrate all this by an important
specific example. Define the differential $d^i: \mathcal{W}^i
\rightarrow \mathcal{W}^{i+1}$ to be the formal linear
differential operator of degree $+1$ corresponding to the
following substitution rules:
\[
\begin{array}{rccl}
\raisebox{-3ex}{\scalebox{0.3}{\includegraphics{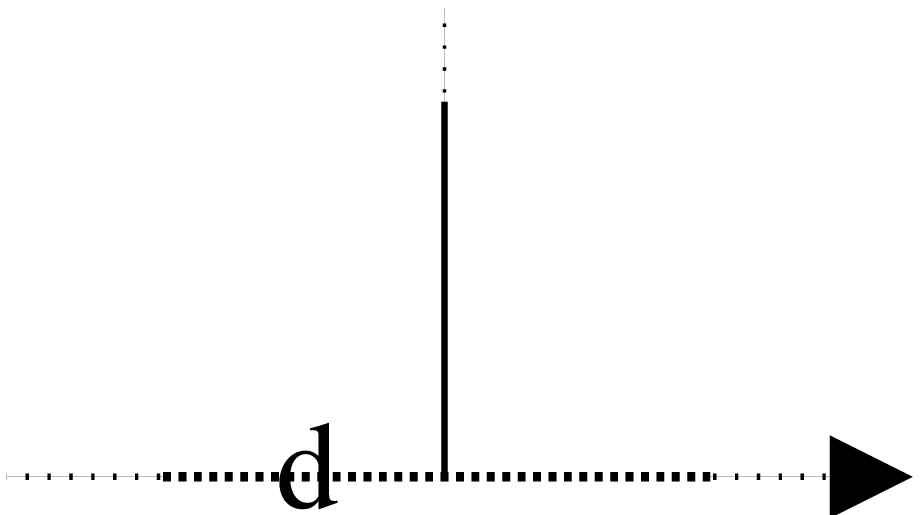}}} &
\leadsto &
\raisebox{-3ex}{\scalebox{0.3}{\includegraphics{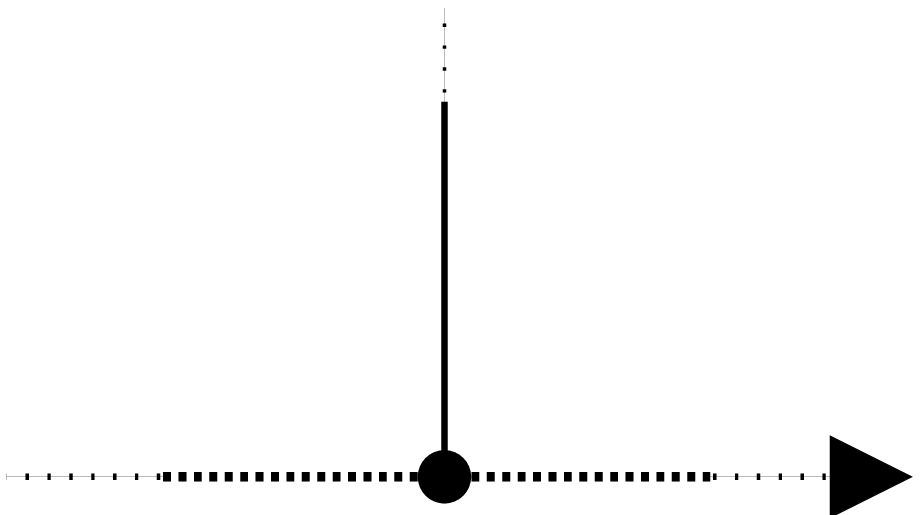}}}\ & -\
\raisebox{-3ex}{\scalebox{0.3}{\includegraphics{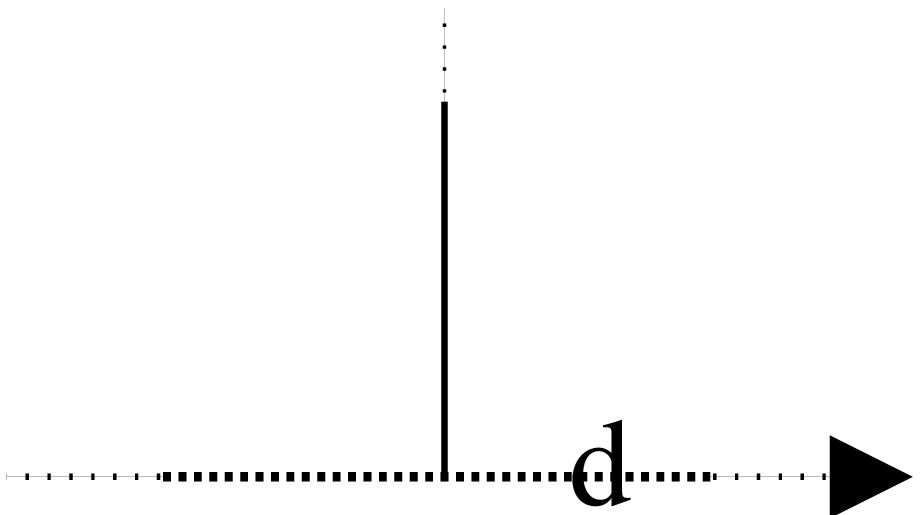}}}\ ,\
\\[0.5cm]
\raisebox{-3ex}{\scalebox{0.3}{\includegraphics{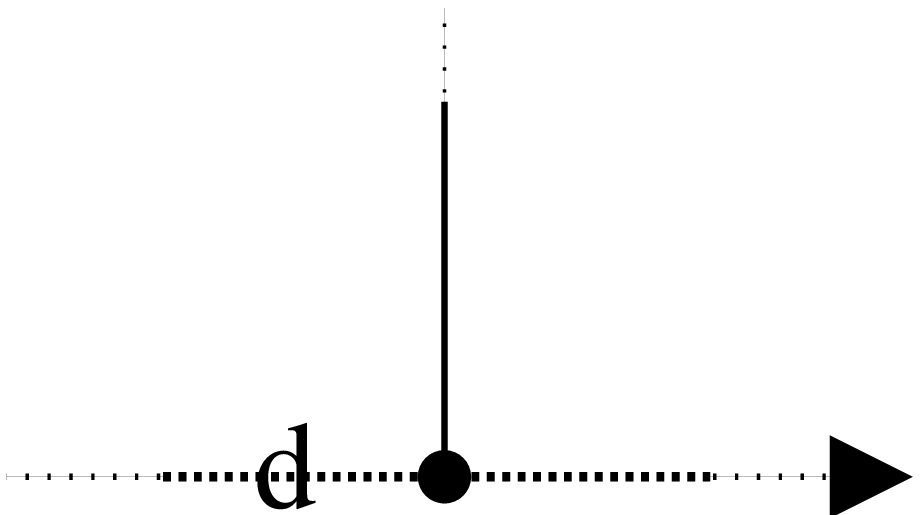}}} &
\leadsto & \makebox[0.3cm][l]{0} & +\
\raisebox{-3ex}{\scalebox{0.3}{\includegraphics{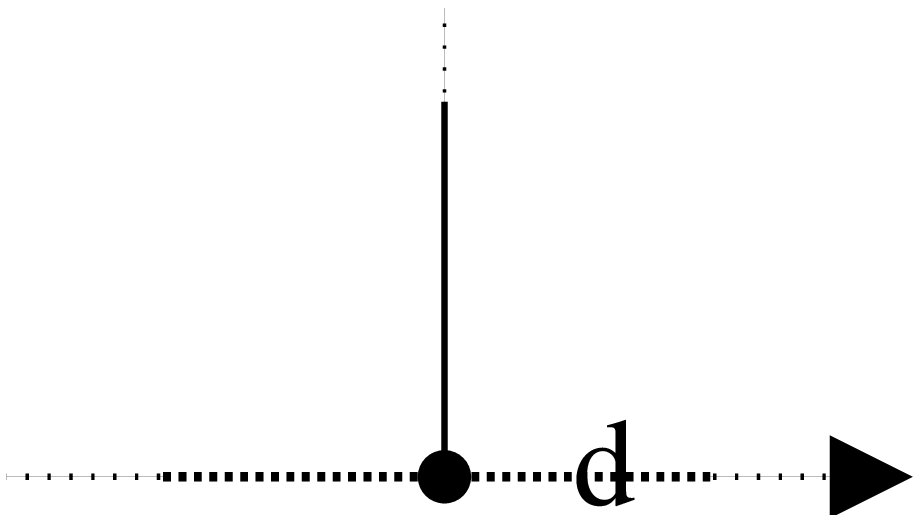}}}\ .\\[0.5cm]
\end{array}
\]
Figure \ref{dexamp} gives an example of the calculation of the
value of $d$.
Before showing that this map actually is a differential (i.e. that
$d\circ d=0$), we'll record a useful lemma.

\subsection{The Lie algebra of formal linear differential
operators.}\label{liediff}
Let $F$ and $G$ be a pair of a pair of formal
linear differential operators, of grades $|F|$ and
$|G|$. Their graded commutator is defined by
$[F,G] = F\circ G - (-1)^{|F||G|} G\circ F$.
The following lemma (whose proof is a formal exercise in the definitions) says that $[F,G]$ is {\it again} a formal linear differential operator.
\begin{lem}\label{commlem}
The graded commutator $[F,G]$ is the formal linear differential operator of
grade $|F|+|G|$ associated to the substitution rules
\[
\raisebox{-3ex}[0.75\height]{\scalebox{0.3}{\includegraphics{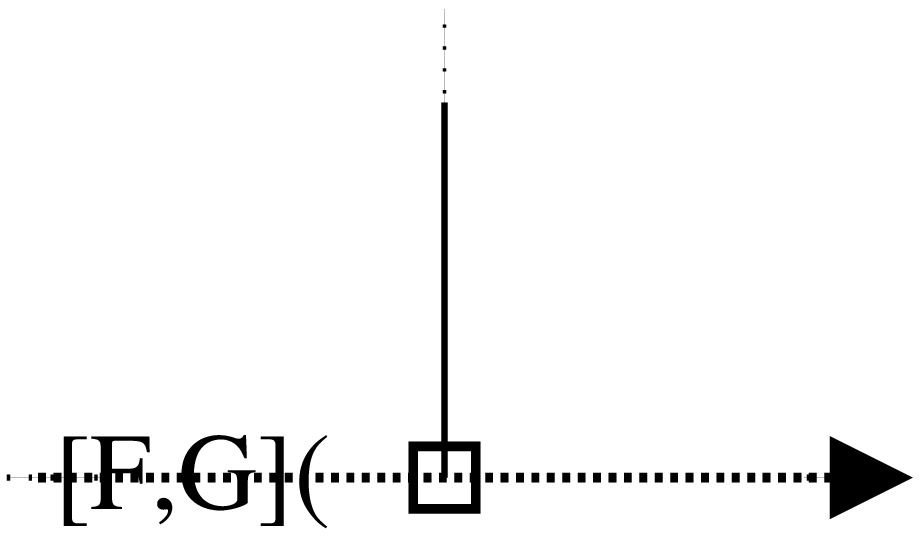}}}
\leadsto
\raisebox{-3ex}[0.5\height]{\scalebox{0.3}{\includegraphics{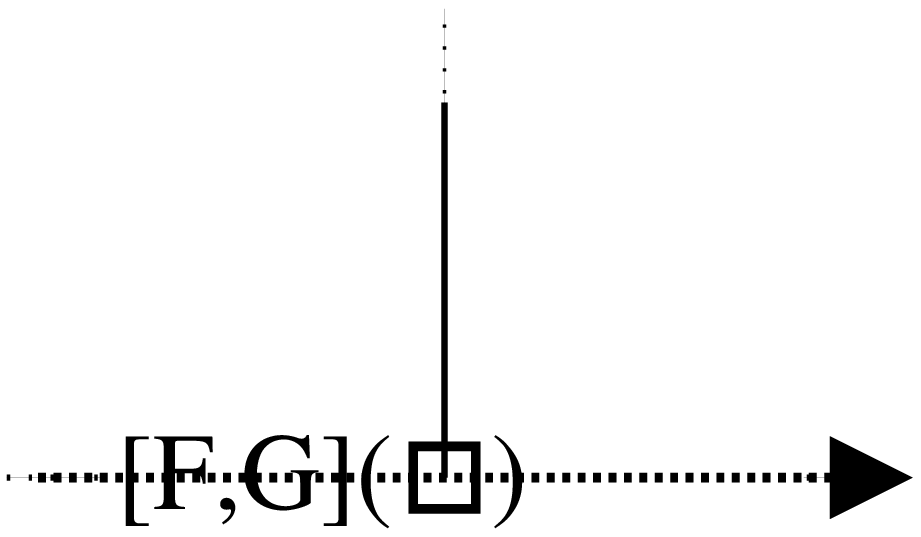}}}
+(-1)^{(|F|+|G|)l}\
\raisebox{-3ex}[0.5\height]{\scalebox{0.3}{\includegraphics{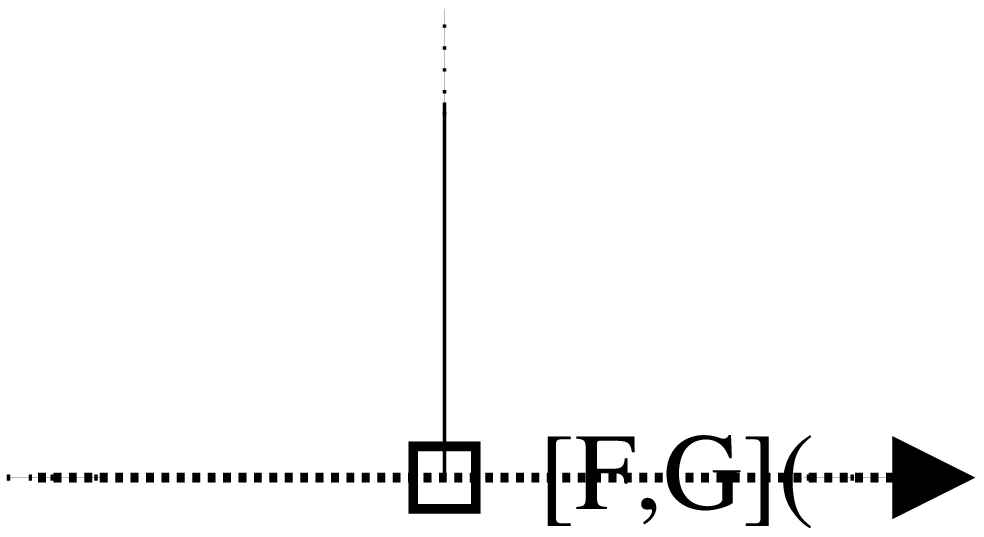}}}\ , \\
\]
where $l$ denotes the grade of the leg in the box.
\end{lem}


\

\subsection{The complex $\mathcal{W}$.}

The system
$
\mathcal{W}: 0 \rightarrow \mathcal{W}^0
\stackrel{d^0}{\longrightarrow} \mathcal{W}^1
\stackrel{d^1}{\longrightarrow} \mathcal{W}^2 \longrightarrow
\ldots
$
forms a differential complex. This is because $d\circ
d=\frac{1}{2}[d,d]$, and Lemma \ref{commlem} implies that $[d,d] =
0$.

\begin{prop}
\[
H^i(\Wspace) = \left\{
\begin{array}{lcl}
\Wspace^0 & & \mbox{if\, $i=0$,} \\[0.15cm]
0 & & \mbox{otherwise.}
\end{array}
\right.
\]
\end{prop}
\begin{proof}
Let $t$ be the formal linear differential operator of grade $-1$
defined by the following substitution rule for grade 2 legs
\[
\begin{array}{rccl}
\raisebox{-3ex}{\scalebox{0.27}{\includegraphics{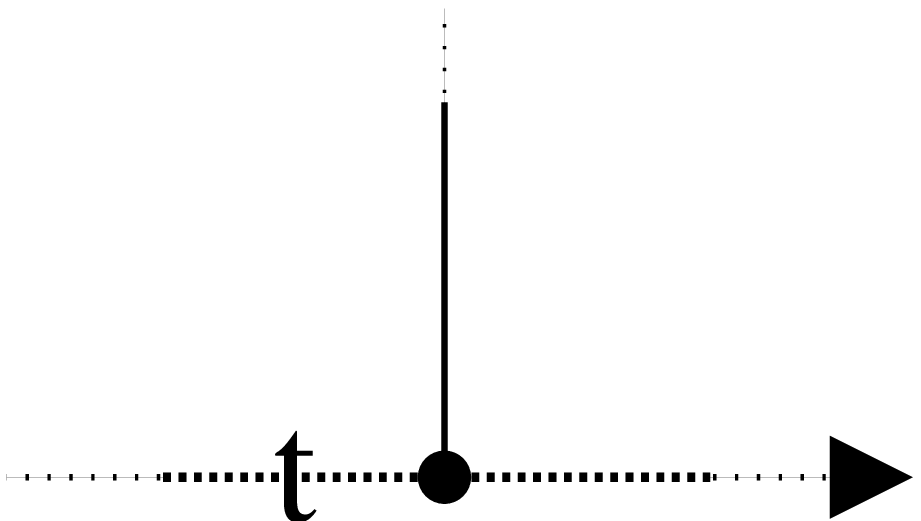}}} &
\leadsto &
\raisebox{-3ex}{\scalebox{0.27}{\includegraphics{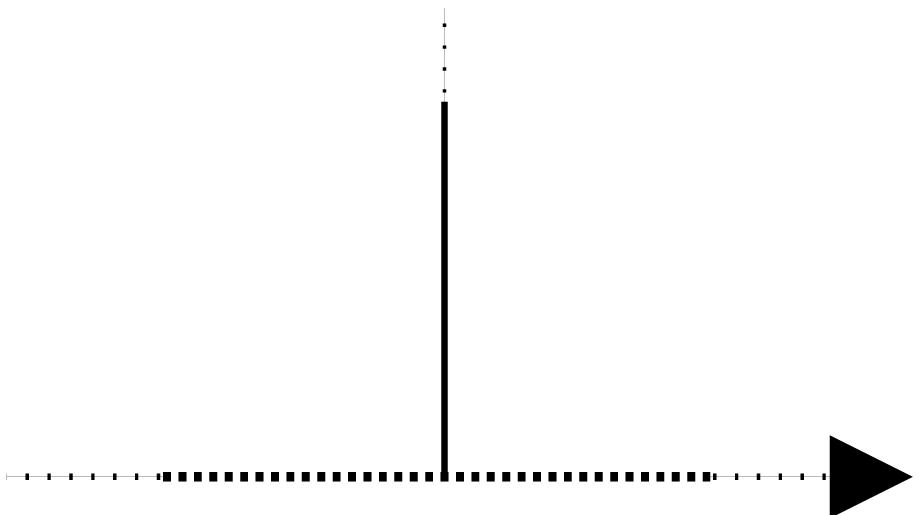}}}\ & +
\raisebox{-3ex}{\scalebox{0.27}{\includegraphics{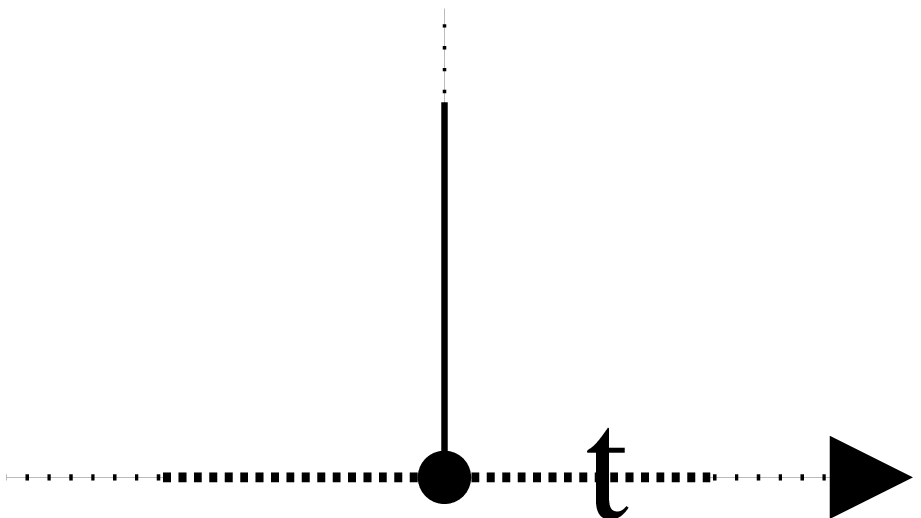}}} \\[0.75cm]
\end{array}
\]
together with the rule that says that $t$ graded-commutes through grade 1 legs:
\[
\begin{array}{rccl}
\raisebox{-3ex}{\scalebox{0.27}{\includegraphics{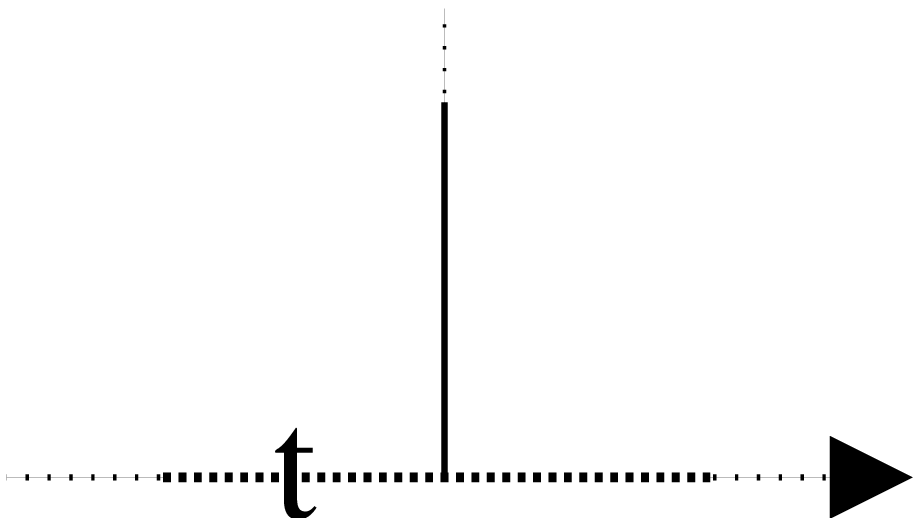}}} &
\leadsto & \makebox[0.3cm][l]{0} & -\
\raisebox{-3ex}{\scalebox{0.27}{\includegraphics{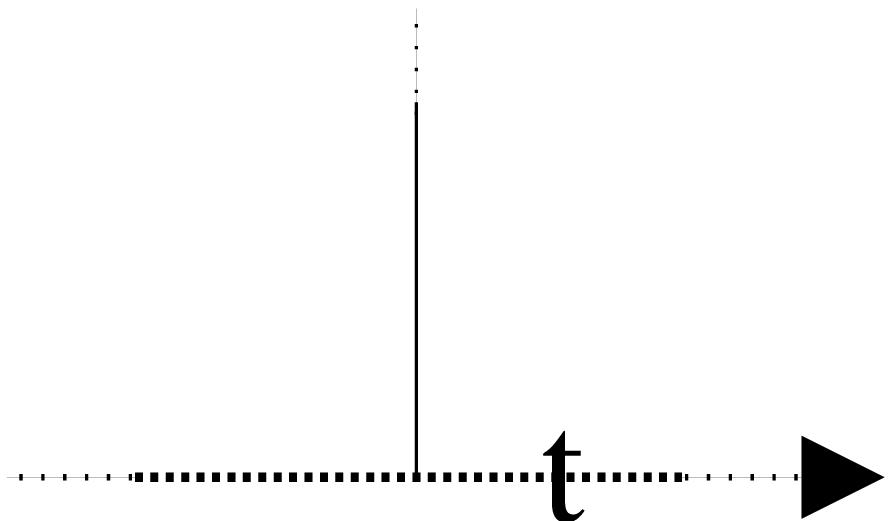}}}
\end{array}
\]
It follows from Lemma \ref{commlem} that
$[t,d](D) = \left(\mbox{$\#$ of legs of $D$}\right) D,$ where $D$ is an arbitrary Weil diagram.
So for every $i>0$ define a map $s^i:\mathcal{W}^i\rightarrow
\mathcal{W}^{i-1}$ by
\[
s^i(D) = \frac{1}{\left(\mbox{$\#$ of legs of $D$}\right)}t^i(D).
\]
This is the required contracting homotopy:
$d^{i-1}\circ s^i + s^{i+1}\circ d^i = \mathrm{id}^i$, when $i\geq 1$.
\end{proof}

So the complex ${\mathcal W}$ has trivial cohomology. There is,
however, a  subcomplex of ${\mathcal W}$, the {\it basic}
subcomplex $\Wbase$, whose cohomology spaces are canonically
isomorphic to the $\{\mathcal{B}^i\}$ (as will be discussed in Section \ref{basiccompute}).
To cut out the subcomplex $\Wbase \subset \Wspace$, we need a
simple structure we'll call an $\iota$-complex.

\subsection{$\iota$-complexes.}\label{ics}

An $\iota$-complex is a pair of complexes together with a grade
$-1$ map $\iota$ between them. Our first $\iota$-complex, to be
defined presently, looks like this:
\[
\xymatrix{ 0 \ar[r] & \mathcal{W}^0 \ar[r]^{d^0}
\ar[dl]_{\iota^0}& \mathcal{W}^1
\ar[r]^{d^1} \ar[dl]_{\iota^1} & \mathcal{W}^2 \ar[dl]_{\iota^2} \ar[r]^{d^2} & \ldots \\
0 \ar[r] & \mathcal{W}^0_\iota \ar[r]^{d^0_\iota} &
\mathcal{W}^1_\iota \ar[r]^{d^1_\iota}
& \mathcal{W}^2_\iota \ar[r]^{d^2_\iota} & \ldots \\
}
\]
To be precise: the equations $d^{i+1}\circ d^{i} = 0$,
$d^{i+1}_\iota\circ d^{i}_\iota=0$ and $\iota^{i+1}\circ d^i =
-d^{i-1}_\iota\circ \iota^{i}$ hold in the above system. (This
last equation has a $(-1)$ because $d$ is grade $1$ and $\iota$ is
grade $-1$.)

\begin{defn}
The {\bf basic subcomplex} $\mathcal{K}_{\mathrm{basic}}$
associated to an $\iota$-complex \\
$(\mathcal{K},\mathcal{K}_\iota,\iota)$ is defined by setting
$\mathcal{K}_{\mathrm{basic}}^i=\mathrm{Ker}(\iota^i)$ and by
defining the differential to be the restriction of the
differential on $\mathcal{K}$. (The rule that $\iota$ commutes
with $d$ implies that this actually is a subcomplex.\,)
\end{defn}

\begin{defn}
A map $f:(\mathcal{K},\mathcal{K}_\iota,\iota_{\mathcal{K}})
\rightarrow (\mathcal{L},\mathcal{L}_\iota,\iota_{\mathcal{L}})$
between a pair of $\iota$-complexes is just a pair of sequences of
maps $(\,\{f^i: \mathcal{K}^i\rightarrow
\mathcal{L}^i\},\{f^i_\iota:\mathcal{K}^i_\iota \rightarrow
\mathcal{L}^i_\iota\}\,)$ commuting with the $d$'s and the
$\iota$'s. (Commuting up to sign as determined by the
understanding that $d$ and $d_\iota$ are regarded as grade $1$
maps, $f$ and $f_\iota$ as grade $0$ maps, and $\iota$ as a grade
$-1$ map.)

A map $f:(\mathcal{K},\mathcal{K}_\iota,\iota_{\mathcal{K}})
\rightarrow (\mathcal{L},\mathcal{L}_\iota,\iota_{\mathcal{L}})$
between $\iota$-complexes restricts to a chain map between their
basic subcomplexes. The restriction will be denoted
\[
f_\mathrm{basic}:\mathcal{K}_{\mathrm{basic}} \rightarrow
\mathcal{L}_{\mathrm{basic}}.
\]
\end{defn}

\subsection{The $\iota$-complex
\((\mathcal{W},\mathcal{W}_\iota,\iota)\).}\label{wics} The space
$\mathcal{W}^{i}_\iota$ is defined in exactly the same way as the
space $\mathcal{W}^{i}$, except that it is based on Weil diagrams
with precisely one special degree 1 vertex (in diagrams the
special vertex will be labelled by the symbol $\iota$). For
example, the following diagram is a generator of
${\mathcal{W}_{\iota}^{6}}$:
\[
\raisebox{-5ex}[0.5\height]{\scalebox{0.32}{\includegraphics{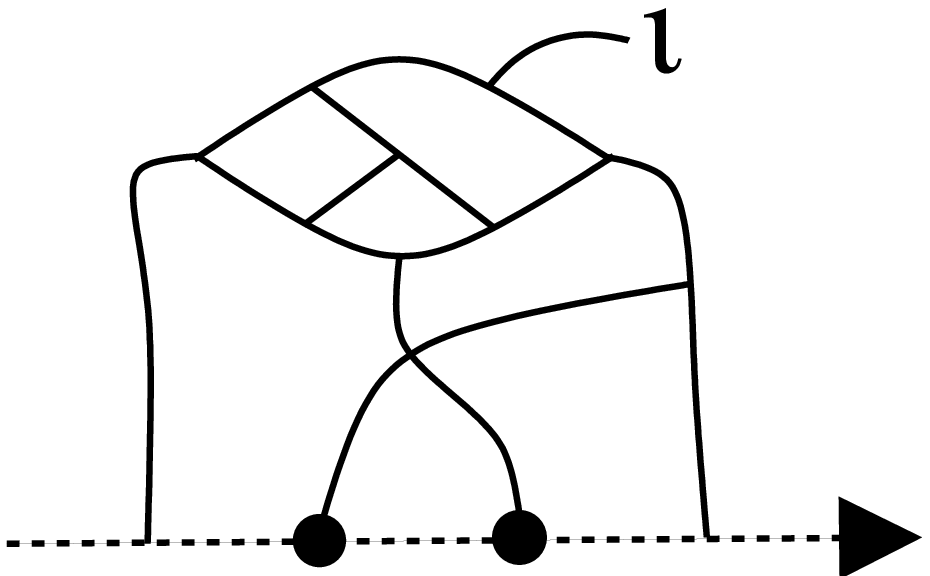}}}
\]
The map $\iota$ is the formal linear differential operator of
grade $-1$ defined by the substitution rules:
\[
\begin{array}{rccl}
\raisebox{-3ex}{\scalebox{0.3}{\includegraphics{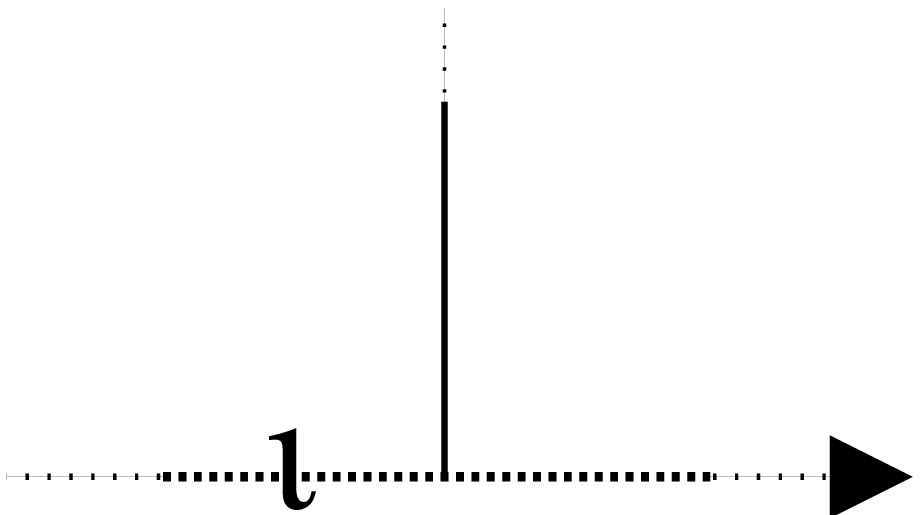}}} & \leadsto
& \raisebox{-3ex}{\scalebox{0.3}{\includegraphics{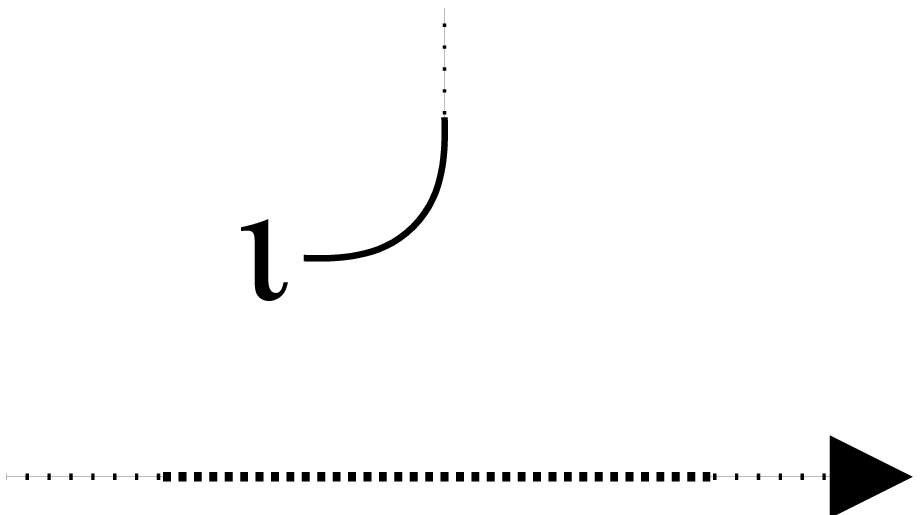}}} -\
\raisebox{-3ex}{\scalebox{0.3}{\includegraphics{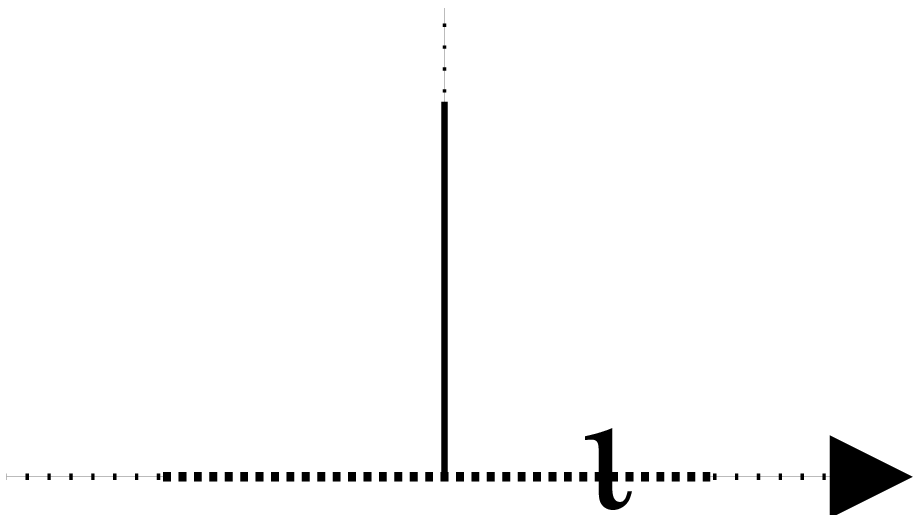}}}\ \ ,
\\[0.6cm]
\raisebox{-3ex}{\scalebox{0.3}{\includegraphics{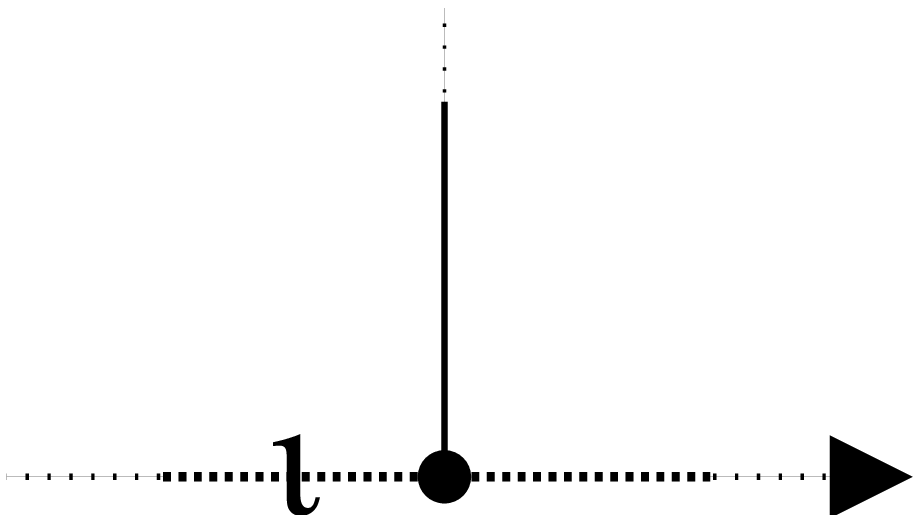}}} & \leadsto
& \raisebox{-3ex}{\scalebox{0.3}{\includegraphics{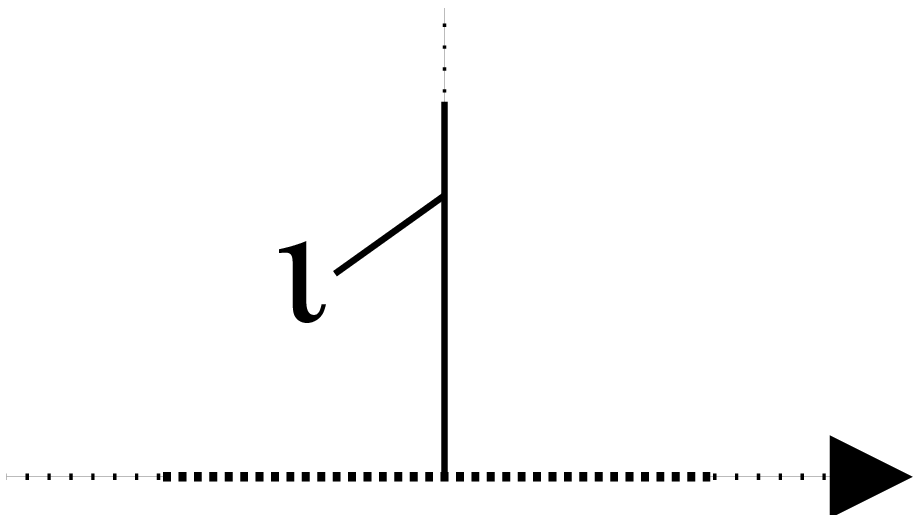}}} +\
\raisebox{-3ex}{\scalebox{0.3}{\includegraphics{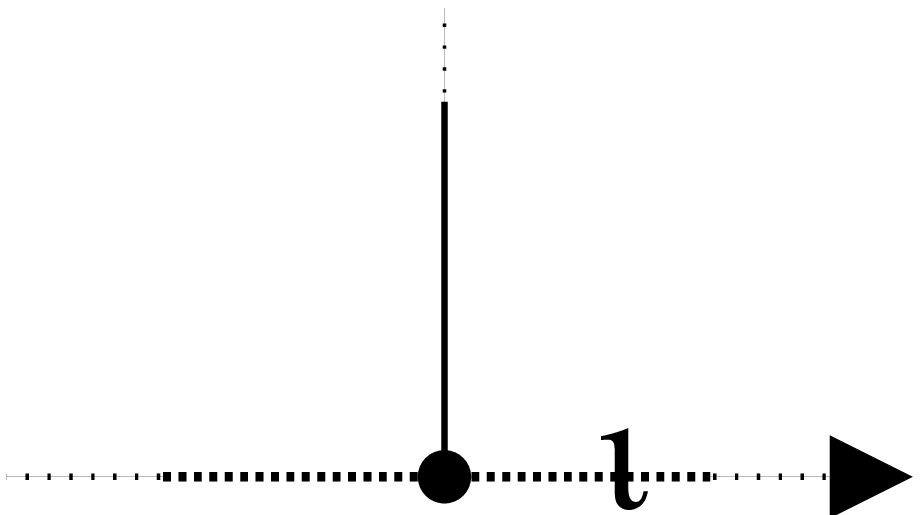}}}\ \
.\\[0.6cm]
\end{array}
\]
\begin{prop}\label{wisanidga}
With these definitions, $(\mathcal{W},\mathcal{W}_\iota,\iota)$
forms an $\iota$-complex.
\end{prop}
\begin{proof}
The only remaining thing to verify is that $[\iota,d]=0$. Lemma
$\ref{commlem}$ tells us that $[\iota,d]$ is the formal linear
differential operator associated to the following substitution
rule (the leg in the box can be of either type):
\[
\begin{array}{rccl}
\raisebox{-3ex}{\scalebox{0.32}{\includegraphics{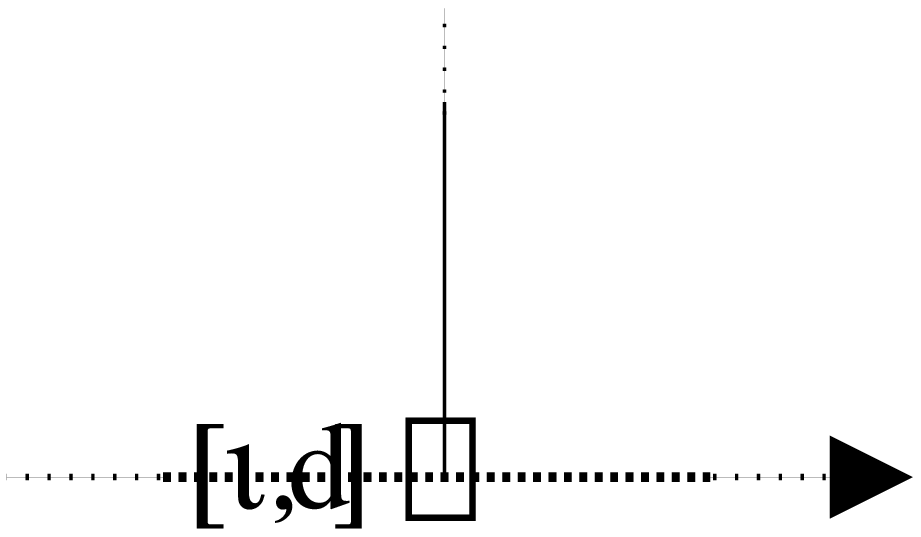}}} &
\leadsto &
\raisebox{-3ex}{\scalebox{0.32}{\includegraphics{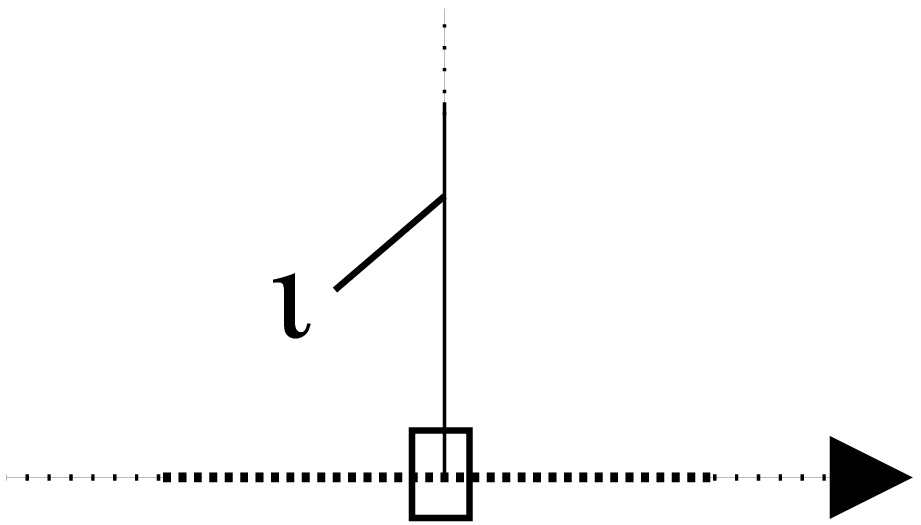}}} +
\raisebox{-3ex}{\scalebox{0.32}{\includegraphics{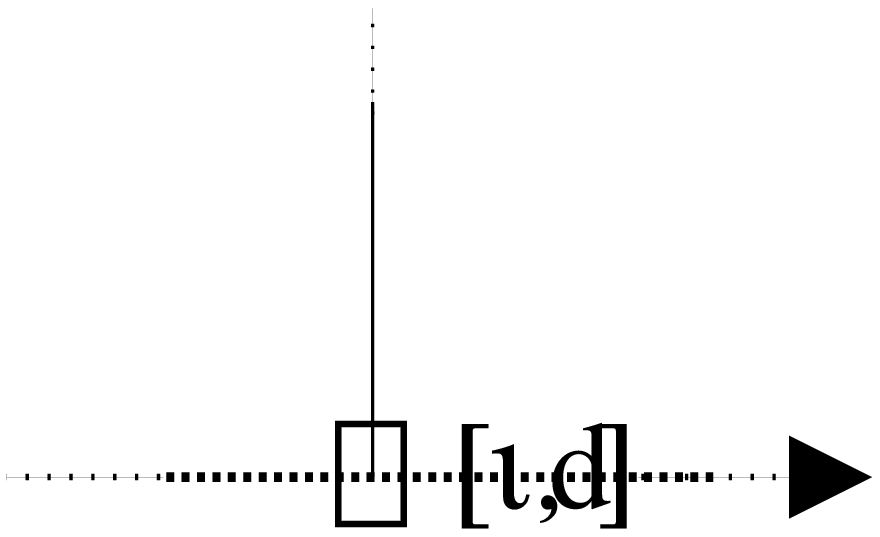}}}
\end{array}
\]
This operator is precisely zero on ${\mathcal{W}}$. To see why
this is true we'll examine its action on a generic Weil diagram.
To begin:
\begin{eqnarray*}
\raisebox{-3ex}{\scalebox{0.28}{\includegraphics{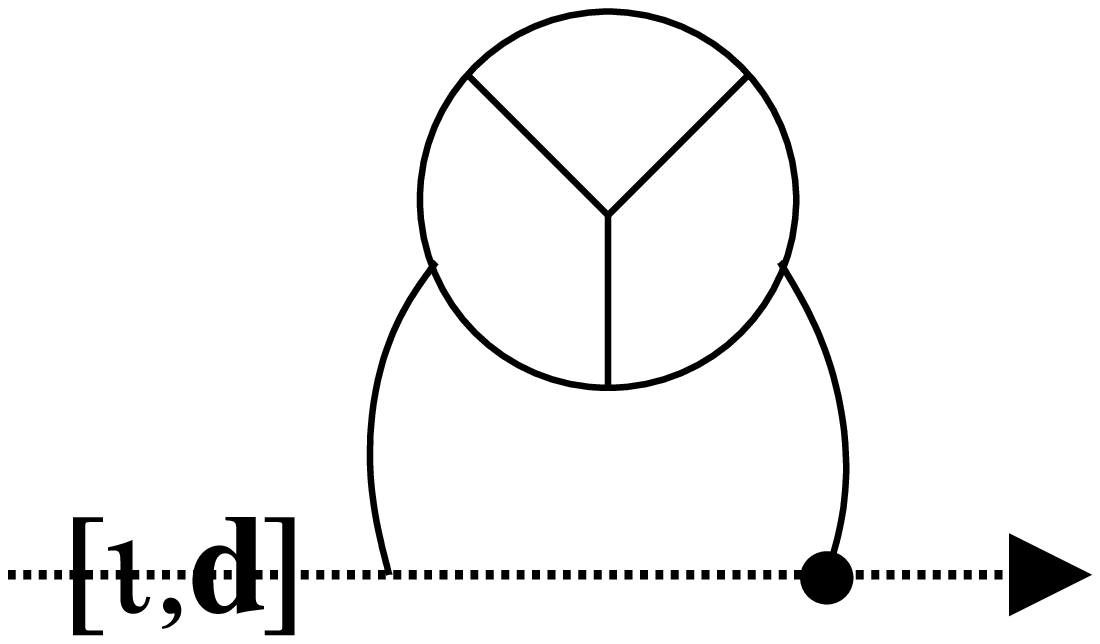}}} & = &
\raisebox{-3ex}{\scalebox{0.28}{\includegraphics{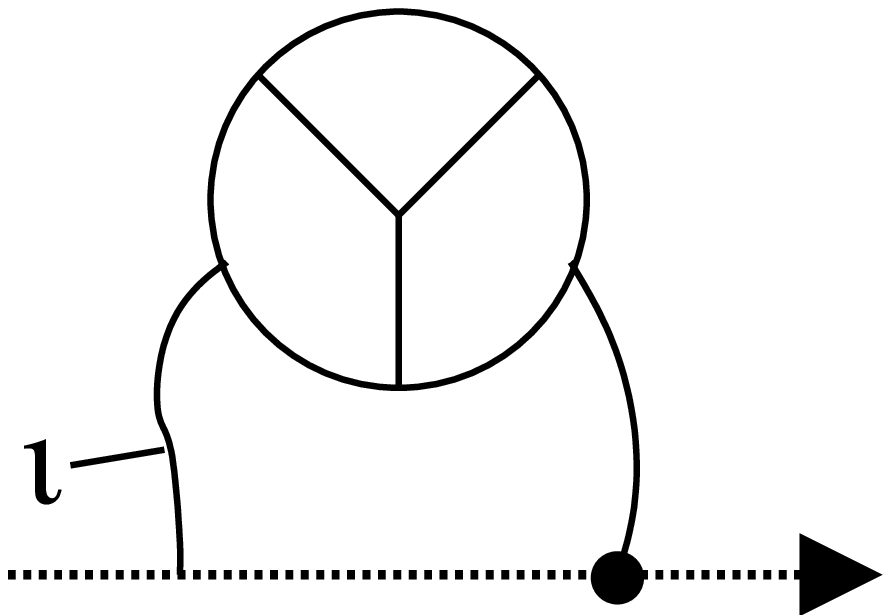}}}\
+\ \raisebox{-3ex}{\scalebox{0.28}{\includegraphics{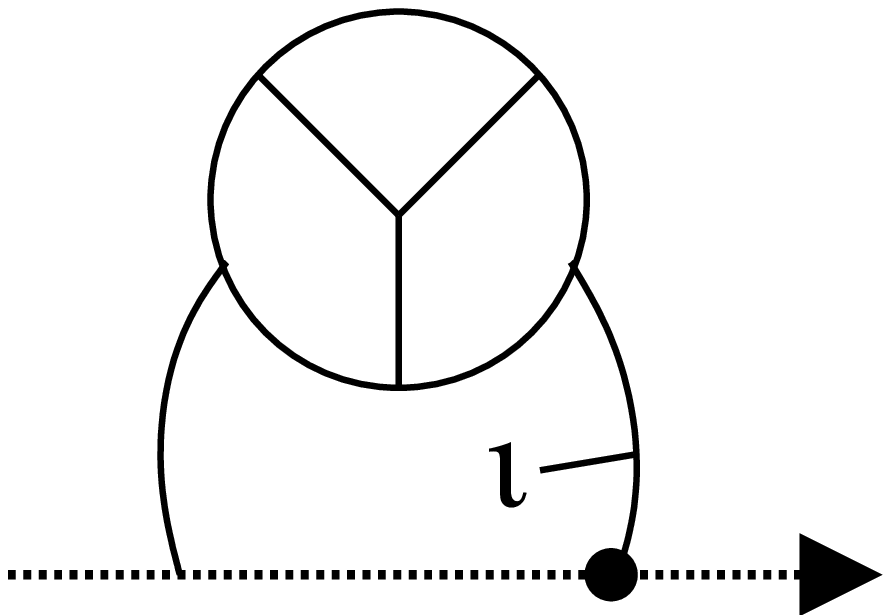}}} \\
& = & \raisebox{-3ex}{\scalebox{0.28}{\includegraphics{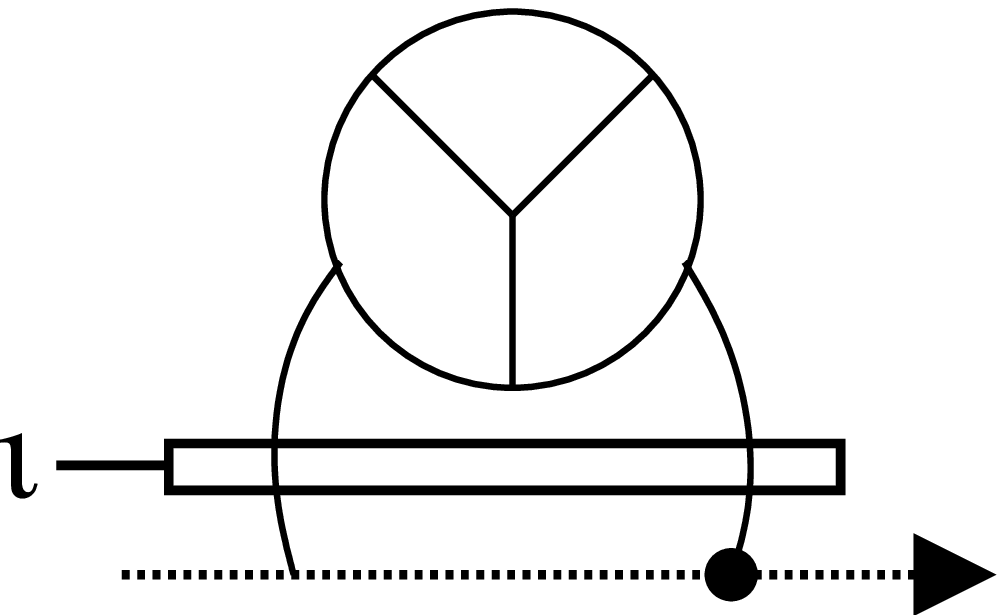}}}\
\ .
\end{eqnarray*}
The box above represents the sum over all ways of placing the end
of the $\iota$-labelled edge on one of the edges going through the
box. Then we can just use the {AS} and {{IHX}} relations to
``sweep" the box up and off the diagram:
\begin{eqnarray*}
&\raisebox{-3ex}{\scalebox{0.26}{\includegraphics{sweepD}}}&
 \stackrel{\mathrm{(IHX)}}{=}\ \
\raisebox{-3ex}{\scalebox{0.26}{\includegraphics{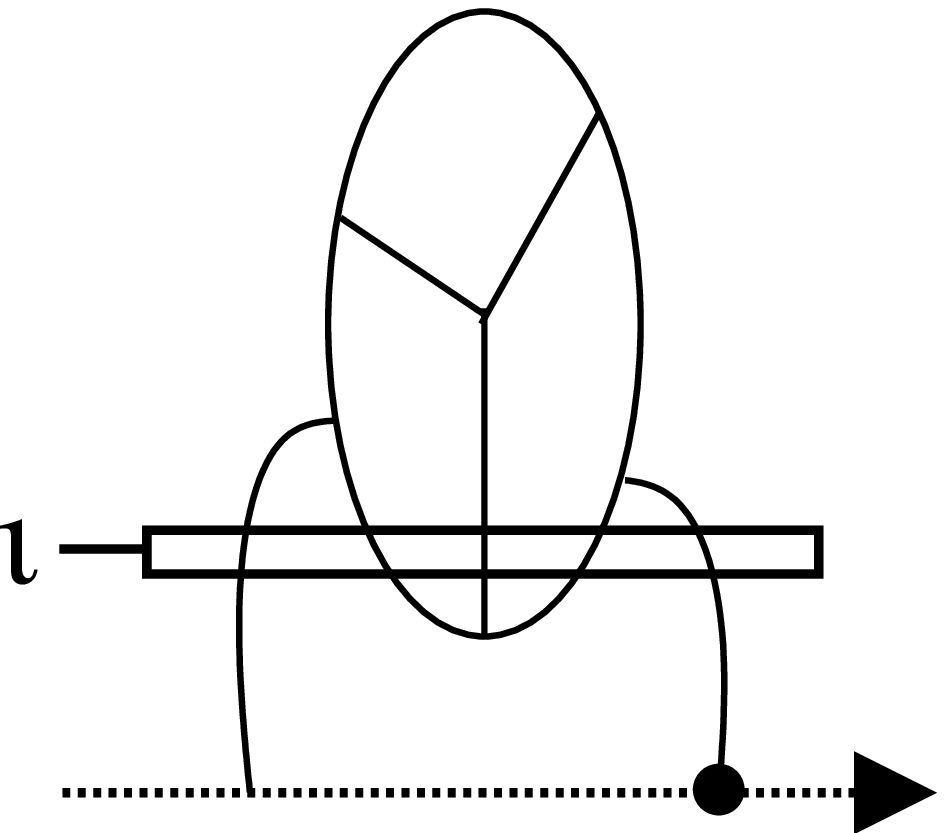}}} \ \
\stackrel{\mathrm{(IHX)}}{=}\ \
\raisebox{-3ex}{\scalebox{0.26}{\includegraphics{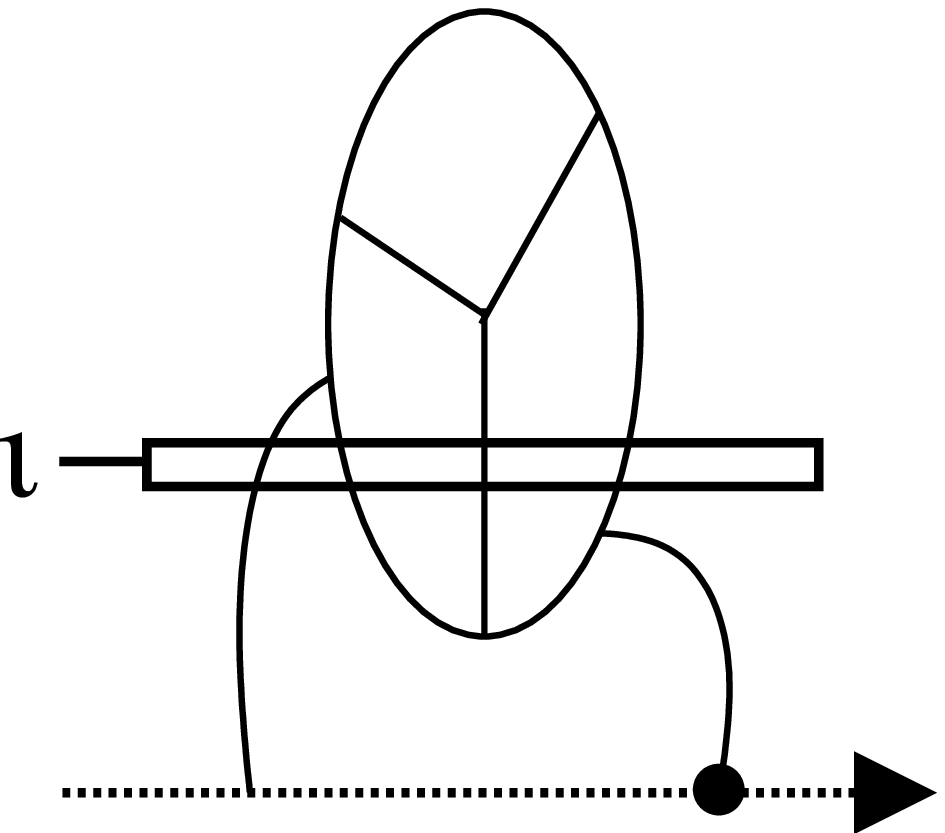}}}
\\
& & \stackrel{\mathrm{(IHX)}}{=}\ \
 \ldots
 \ \ \stackrel{\mathrm{(IHX)}}{=}\ \
\raisebox{-3ex}{\scalebox{0.26}{\includegraphics{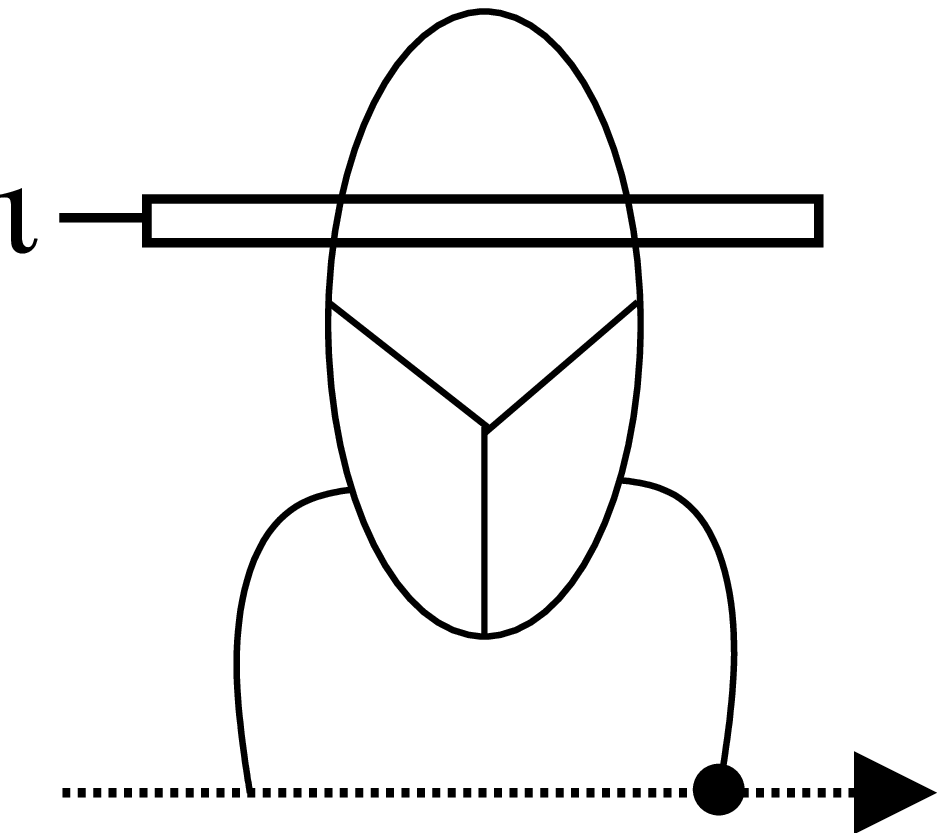}}} \ \ \stackrel{\mathrm{(AS)}}{=}\ \ \ 0.
\end{eqnarray*}
\end{proof}
So the triple $(\Wspace,\Wspace_\iota,\iota)$ forms an
$\iota$-complex. Thus we can pass to the basic subcomplex
$\Wbase$. The computation of the corresponding basic cohomology is
the subject of Section \ref{basiccompute}.

\subsection{The characteristic class-valued evaluation
map.}\label{lengthyaside}

The most commonly asked question that arises when this formalism
is described is: what is the $\iota$ for?

As we discussed in the introduction, these structures start as a formal abstraction of expressions inside Chern-Weil theory.
In this section we'll make the relationship more precise, explaining how to `evaluate' Weil diagrams inside
the de Rham complex on a principal bundle, so as to build characteristic classes.

So consider a compact Lie group $G$ and a smooth manifold $P$ with a
free, smooth left $G$-action. It turns out that the space of
orbits $P/G$, which we'll denote $B$, may be uniquely equipped
with the structure of a smooth manifold so that $P$, $B$ and the
projection map \mbox{$\pi:P\rightarrow B$} form a smooth principal
$G$-bundle.
Indeed, it turns out that every smooth principal $G$-bundle arises
in this fashion.

The Chern-Weil theory of characteristic classes begins by
examining the corresponding pull-back map of differential forms
$\pi^* : \Omega(B) \rightarrow \Omega(P)$
and asking: is this map an injection of graded differential algebras? If so, can we characterize
the image of this map? Yes, and yes, are the answers.

The characterization requires the {\it generating vector fields}
of the $G$-action on $P$. Let $\mathfrak{g}$ denote the Lie
algebra of $G$. Given a vector $\xi\in \mathfrak{g}$, the vector
of the corresponding generating vector field $\xi_P$ at a point
$p$ of $P$ is given by the derivative:
\[
\left(\xi_P
f\right)(p)=\left.\frac{\mathrm{d}}{\mathrm{d}t}f\left(
\mathrm{exp}\left(-t\xi\right).p\right)\right|_{t=0}.
\]
There are two graded differential operators on $\Omega(P)$ that
are naturally associated to a generating vector field $\xi_P$:
\begin{itemize}
\item{ $\mathcal{L}_\xi  :  \Omega^\bullet(P) \rightarrow
\Omega^{\bullet}(P),$ the Lie derivative along $\xi_P$,}\item{
$\iota_\xi : \Omega^\bullet(P) \rightarrow \Omega^{\bullet-1}(P),$
the corresponding contraction operator. Given $\tau$, a 1-form on
$P$, $\iota_\xi$ simply evaluates $\tau$ on the vector field
$\xi_P$. In other words, $\left(\iota_\xi\tau\right)(p) =
\tau_p((\xi_P)_p)$.}
\end{itemize}

\begin{fact} The graded differential operators $\mathcal{L}_\xi$, $\iota_\xi$ and $d$ satisfy the following
graded commutation relations:
\begin{equation}\label{gradcomms}
\begin{array}{lll}
[\mathcal{L}_\xi,\mathcal{L}_{\xi'}]=\mathcal{L}_{[\xi,\xi']} &
[\mathcal{L}_\xi,\iota_{\xi'}]=\iota_{[\xi,\xi']} &
[\mathcal{L}_\xi,\mathrm{d}]=0 \\

[\iota_\xi,\mathrm{d}] = \mathcal{L}_\xi &
[\iota_\xi,\iota_{\xi'}]=0 & \\
\end{array}
\end{equation}
\end{fact}

Consider, now, some differential form $\omega$ on $B$, and
consider its pull-back $\pi^*(\omega)$ to $P$. Our first
observation is that any such pull-back has the property that it is
annihilated by every $\mathcal{L}_\xi$ and also by every
$\iota_\xi$. Before noting the simple reasons for this fact, let's
encode it in a definition:

\begin{defn}
A differential form on $P$ is said to  be {\bf basic} if it is
annihilated by every $\mathcal{L}_\xi$ and also by every
$\iota_\xi$. Let $\Omega(P)_{\mathrm{basic}}$ denote the ring of
basic forms.
\end{defn}

\begin{prop}
The pull-back map $\pi^*$ maps into the ring of basic forms:
\[
\pi^*(\Omega(B)) \subset \Omega(P)_{\mathrm{basic}}.
\]
\end{prop}
\begin{proof}
(1st point: $\mathcal{L}_\xi\left(\pi^*(\omega)\right)=0$.) A
direct implementation of the definitions gives:
\begin{eqnarray*}
\left(\mathcal{L}_\xi \left(\pi^*(\omega)\right)\right)_p & = &
\frac{\mathrm{d}}{\mathrm{d}t}\left(
\left(\text{exp}(-t\xi)^*\circ
\pi^*\right)(\omega)_p\right) \\
& = & \frac{\mathrm{d}}{\mathrm{d}t}\left( \left(\pi \circ
\text{exp}(-t\xi)
\right)^*(\omega)_p\right) \\
& = & \frac{\mathrm{d}}{\mathrm{d}t}\left(\pi^*(\omega)_p\right),
\\[0.1cm]
& = & 0,
\end{eqnarray*}
where to get the penultimate line we have used the fact that the
projection map is invariant. (2nd point:
$\iota_\xi\left(\pi^*(\omega)\right)=0$.) This is for the simple
reason that the projection map $\pi_*$ sends any vector from a
generating vector field to zero.
\end{proof}
Actually, quite a bit more is true. The following improvement
follows from the local triviality of a principal bundle.
\begin{fact}
The pull back map gives an isomorphism of differential graded
algebras \[ \pi^* : \Omega(B)
\stackrel{\cong}{\longrightarrow}\Omega(P)_{\mathrm{basic}}.\]
\end{fact}

The next step in Chern-Weil theory is to consider {\it connection
forms} on $P$, as we will presently recall. Chern-Weil theory
asks: how can we build basic forms on $P$ from the components of a
connection form?

A {\bf connection form} is an element $A\in
\mathfrak{g}\otimes\Omega^1(P)$ satisfying the two properties
that:
\begin{itemize}
\item[($\star$)]{$\mathcal{L}_\xi(A)=\text{ad}_\xi(A)$, where
ad$_\xi(A)$ denotes the adjoint action of $\xi$ on the first
factor of $\mathfrak{g}\otimes \Omega^1(P)$.}
\item[($\star\star$)]{$\iota_\xi(A) = \xi.$}
\end{itemize}
Connection forms are equivalent to {\it connections} on $P$. Recall: Letting $V_pP$ (the ``vertical
subspace") denote the subspace of $T_pP$ spanned by the generating
vector fields at $p$, a {\bf connection} is a smooth,
$G$-equivariant decomposition $T_pP = V_pP\oplus H_pP$ (where
$H_pP$ is called the ``horizontal subspace"). To get the
corresponding connection form from a connection: observe that the
decomposition
\[
T_pP = V_pP \oplus H_pP \cong \mathfrak{g} \oplus H_pP
\]
gives, at every point $p$ of $P$, a linear map $T_pP \rightarrow
\mathfrak{g}$. That the element of
$\mathfrak{g}\otimes\Omega^1(P)$ that this map corresponds to
satisfies the two properties above is a very instructive exercise.

\subsubsection{Evaluating diagrams in $\Omega(P)$.}

With these preliminaries in hand we can now turn to constructing
a {\it chain map}
\[
\text{Eval}^{(G,P,A)} : \Wspace \rightarrow \Omega(P),
\]
depending on the initial data of a compact Lie group $G$, a
principal $G$-bundle $P$, and a connection form $A$ on $P$. We
will then indicate that this map restricts to a map between the
basic subcomplexes
$\text{Eval}^{(G,P,A)}_{\text{basic}} : \Wspace_{\mathrm{basic}}
\rightarrow \Omega(P)_{\mathrm{basic}}$,
so that we end up with a map
\[
H\left(\text{Eval}^{(G,P,A)}_{\mathrm{basic}}\right) :
H(\Wbase)\rightarrow H(\Omega(P)_{\text{basic}}) \cong H(B).
\]
We will refer to this map as the {\bf characteristic class-valued
evaluation map}. The computation of $H(\Wbase)$, the ``universal ring" in this characteristic class theory, is the subject of
the next section. (The answer is $\Bspace$.)

The map $\text{Eval}^{(G,P,A)}$ will be constructed
as a state-sum in the familiar way.
For completeness, we'll describe the method in detail here. To
begin, let's fix some notation to refer to three important tensors
associated to the initial data.

({\it 1. The inner product}.) The Lie algebra $\mathfrak{g}$, being
the Lie algebra of a compact Lie group, comes equipped with an
inner product which is:
\begin{enumerate}
\item{Symmetric, so that $\langle v,w\rangle=\langle w,v\rangle$;}
\item{Invariant, so that $\langle[u,v],w\rangle + \langle v ,
[u,w]\rangle=0$;} \item{Non-degenerate, so that the map from
$\mathfrak{g}$ to ${\mathfrak{g}}^*$ given by $v\mapsto \langle
v,.\,\rangle$ is an isomorphism. Let
$\varpi:{\mathfrak{g}}^*\rightarrow\mathfrak{g}$ denote the
inverse map.}
\end{enumerate}

({\it 2. The Casimir element}.) The identity map from $\mathfrak{g}$
to itself can be viewed as an element of
${\mathfrak{g}}^*\otimes\mathfrak{g}$. (In coordinates:
introduce a basis $\{t_a\}_{a=1..\text{dim}\,\mathfrak{g}}$. Let
$\{{t}^*_a\}$ denote the set of dual vectors defined by the rule
${t}^*_a(t_b)=\delta_{ab}$. Then the identity map is $\sum_a
{t}^*_a\otimes t_a$.) If we take the tensor representing the
identity map and apply $\varpi$ to the first factor (what
physicists call ``lowering the index") then we get a tensor in
$\mathfrak{g}\otimes\mathfrak{g}$ which we'll refer to as the {\it
Casimir} element. Write this element as a sum $\sum_{a}s_a\otimes
t_a.$

({\it 3. The connection form}.) The
connection form $A$ is an element of $\mathfrak{g}\otimes
\Omega^1(P)$. Write it as a sum: $A = \sum_i   r_i \otimes \omega_i$.

We can now construct the state-sum. Consider some Weil
diagram. To begin, chop the diagram into pieces by cutting every
edge at its midpoint. For example:
\[
\raisebox{-4ex}{\scalebox{0.22}{\includegraphics{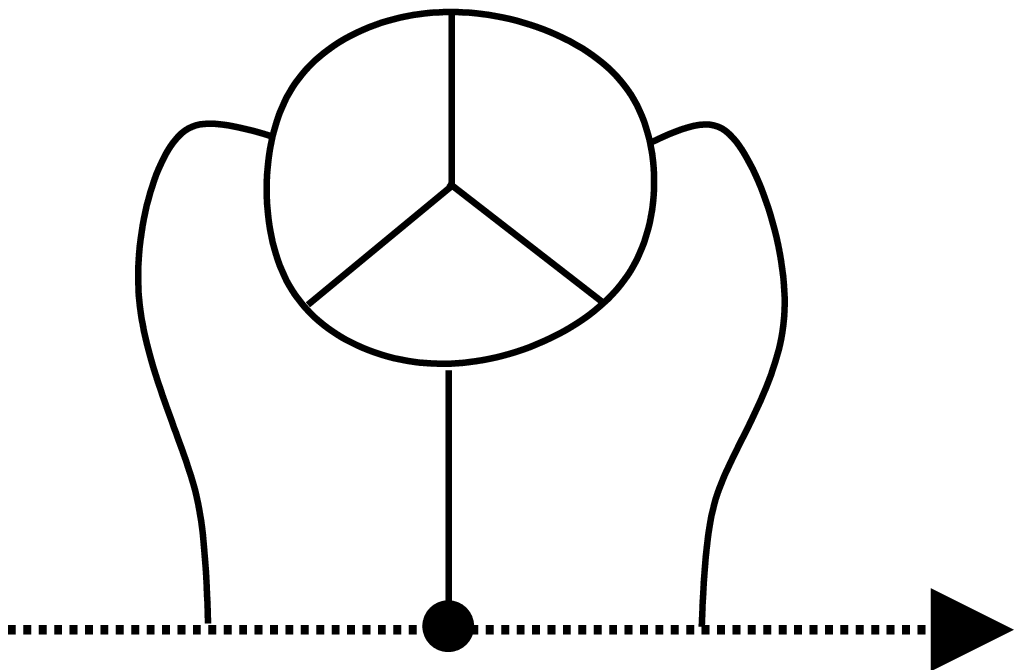}}} \
\mapsto\
\raisebox{-4ex}{\scalebox{0.22}{\includegraphics{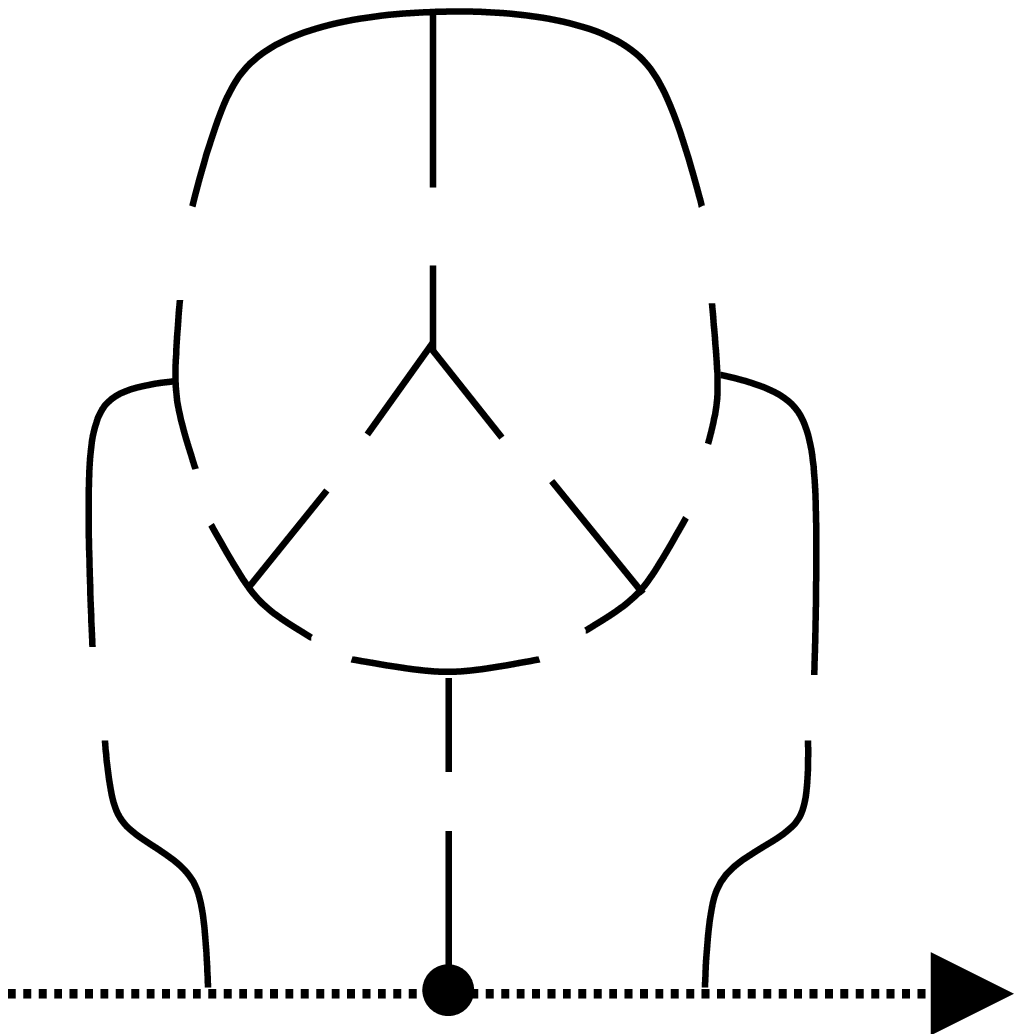}}}
\]
The next step is to decorate the chopped-up diagram with tensors.
Every piece of the diagram is decorated with a tensor according to
what kind of piece it is. The rules are as follows, for trivalent vertices, leg-grade 1 legs, and
leg-grade 2 legs, respectively. The relevant tensors are expressed using the notations just discussed.
\[
\sum_{a,b}\
\raisebox{-4ex}{\scalebox{0.2}{\includegraphics{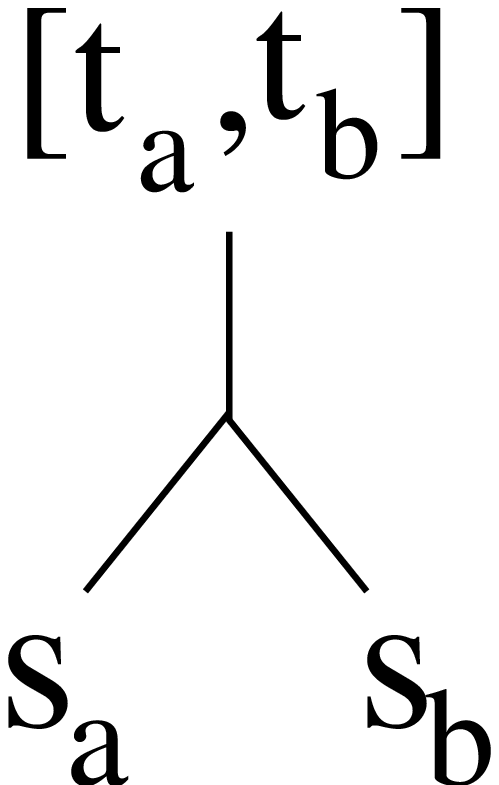}}}\ \ \ \ \ \ \  \ \ \ \
\sum_i
\raisebox{-5ex}{\scalebox{0.2}{\includegraphics{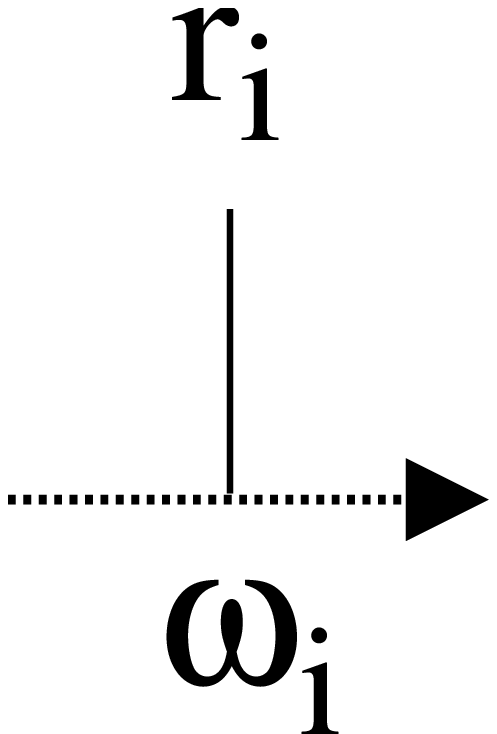}}}\ \ \ \ \ \ \ \ \ \ \
\sum_i
\raisebox{-5ex}{\scalebox{0.2}{\includegraphics{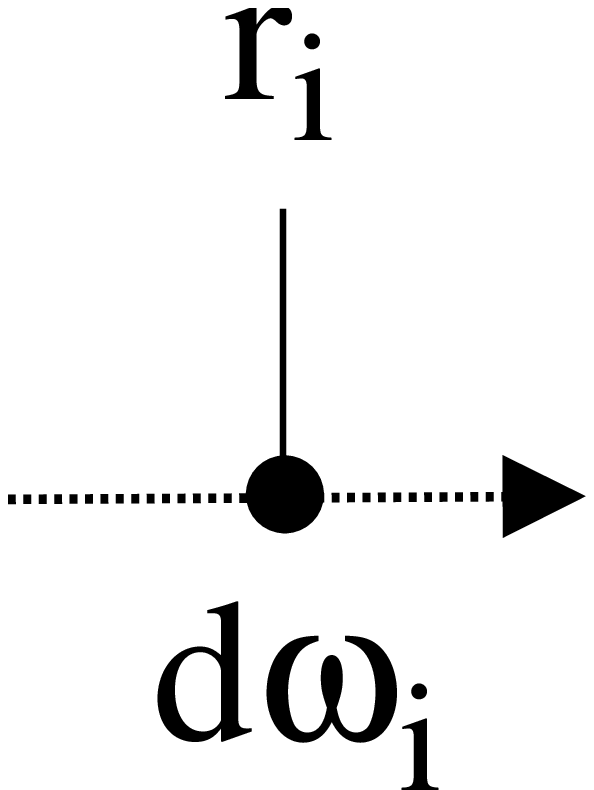}}}\ \ \
\]

\begin{figure}
\caption{Constructing the state sum.\label{evalexamp}}
\[
\label{statesumexamp}
\sum_{a,\ldots,r}\raisebox{-12ex}{\scalebox{0.26}{\includegraphics{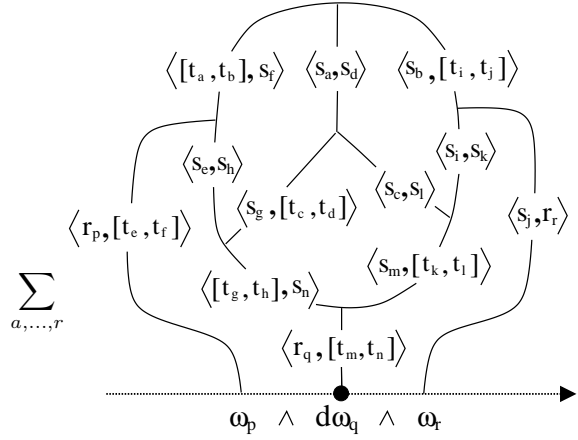}}}\
\ .
\]
\end{figure}

To finish the construction: use the inner product to pair up every
pair of elements that face each other where you broke the edges;
then take the wedge product of the differential forms that appear
along the orienting line at the bottom. An example is exhibited in Figure \ref{evalexamp}.

\begin{prop}
The map $\mathrm{Eval}^{(G,P,A)} : \Wspace \rightarrow
\Omega(P)$ is well-defined. In other words, it respects the
anti-symmetry, Jacobi, and leg permutation relations, and is a
chain map.
\end{prop}

The map respects the anti-symmetry and Jacobi relations for
usual reasons. That the map respects the leg permutation
relations obviously follows from the grading of the differential
forms involved in the state-sum. It is a chain map by construction.

\begin{prop}
The map $\mathrm{Eval}^{(G,P,A)}$ restricts to a map
\[
\mathrm{Eval}^{(G,P,A)}_{\mathrm{basic}} : \Wbase \rightarrow
\Omega(P)_\mathrm{basic}
\]
between the basic subcomplexes.
\end{prop}

We'll leave the proof as an instructive exercise. The reader needs to check that a form in the image of this map is annihilated by every
$\mathcal{L}_\xi$ and by every $\iota_{\xi}$. In the case of $\mathcal{L}_\xi$: first observe the effect that the degree 0 differential operator $\mathcal{L}_\xi$ has on an expression like shown in Figure \ref{evalexamp} above. Then use the defining property $(\star)$ of a connection to convert that expression an appropriate diagrammatic expression. Finish with a sweeping argument. For the case of $\iota_{\xi}$: use the defining property $(\star\star)$ of a connection to show that the result of the action by $\iota_{\xi}$ must factor through the diagrammatic $\iota$.

\section{Introducing the map $\Upsilon$.}
\label{basiccompute} \label{changebasistoF} In this section we'll
see how to map $\Bspace$, the familiar space of symmetric Jacobi
diagrams, into $\Wspace$, the just-introduced commutative Weil
complex for diagrams, in order to get isomorphisms in basic
cohomology:
\[
H^i(\Upsilon_{\mathrm{basic}})\,:\,\Bspace^i\cong
H^i(\Bspace_\text{basic}) \stackrel{\cong}{\longrightarrow}
H^i(\Wbase).
\]
We'll introduce $\Upsilon$ as the composition of
two other maps of $\iota$-complexes: \[ (\Bspace,0,0)
\xrightarrow{\ \phi_\Bspace\ }
(\mathcal{W}_{\mathrm{F}},\mathcal{W}_{\mathrm{F}\iota},\iota)
\xrightarrow{\baseFtobull} (\Wspace,\Wspace_\iota,\iota).
\]
Formally speaking, this section will be rather routine. We'll
introduce the $\iota$-complex
$(\mathcal{W}_{\mathrm{F}},\mathcal{W}_{\mathrm{F}\iota},\iota)$.
Then we'll introduce the maps $\baseFtobull$ and $\phi_\Bspace$\footnote{There is a $\Bspace$ suffix in this notation to logically distinguish this map from a later $\Aspace$ version $\phi_\Aspace$.}. Along the
way we'll explain that they both induce isomorphisms in basic
cohomology.

\subsection{A change of basis: The $\iota$-complex
\((\mathcal{W_\mathrm{F}},\mathcal{W}_{\mathrm{F}\iota},\iota)\).}

This should be regarded as the $\iota$-complex
$(\mathcal{W},\mathcal{W}_\iota,\iota)$ viewed under a different
basis.
This complex employs Weil diagrams
with two types of leg: grade 2 legs labelled by the symbol {\bf
F}, and grade 1 legs with no decoration. These legs obey the usual
permutation rules appropriate to this grading. For example:
\[
\raisebox{-3ex}{\scalebox{0.30}{\includegraphics{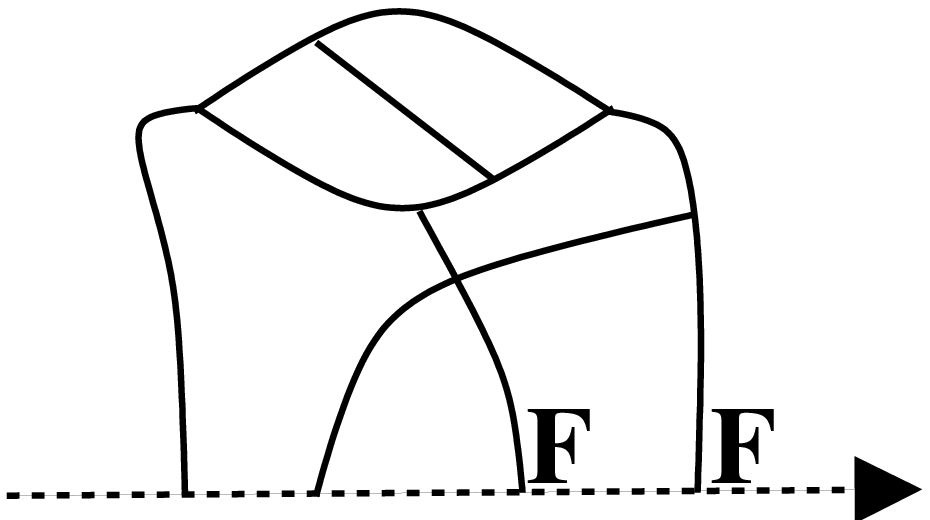}}} =\ \ \
-\raisebox{-3ex}{\scalebox{0.30}{\includegraphics{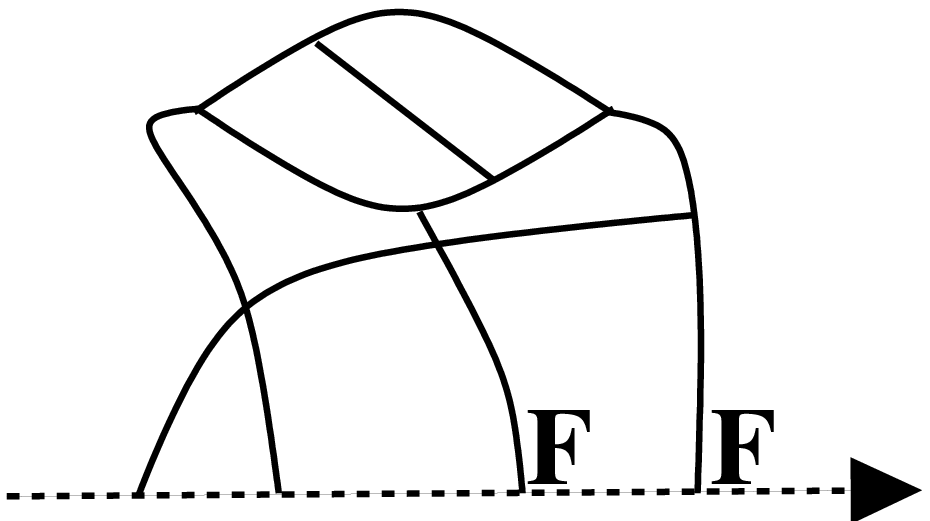}}} =\ \ \
-\raisebox{-3ex}{\scalebox{0.30}{\includegraphics{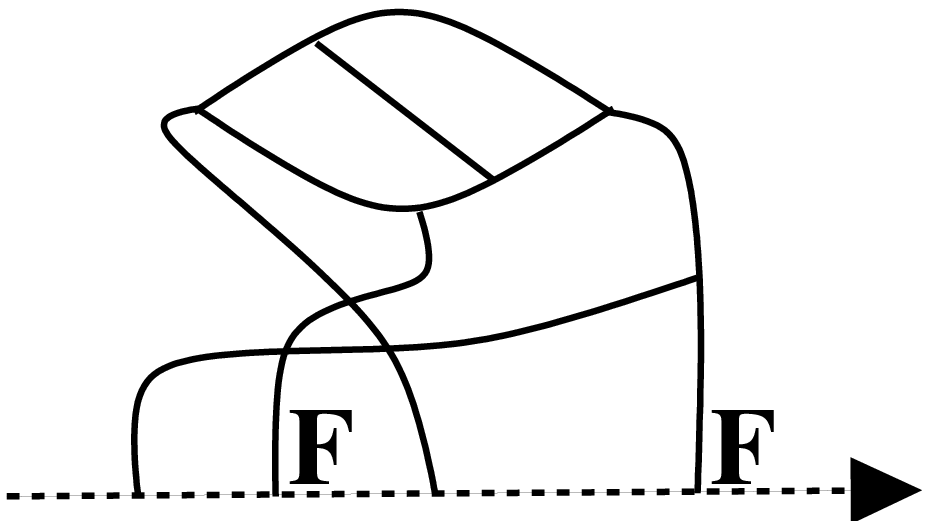}}}
\]
The spaces $\mathcal{W}_{\mathrm{F}}^i$ and
$\mathcal{W}_{\mathrm{F}\iota}^i$ are defined in the obvious way.
Indeed, the $\iota$-complexes
$(\mathcal{W},\mathcal{W}_\iota,\iota)$ and
$(\mathcal{W}_{\mathrm{F}},\mathcal{W}_{\mathrm{F}\iota},\iota)$
differ from each other only in the definition of the maps $\iota$
and $d$.
In this case, the map $\iota$ is defined to be the formal linear
differential operator defined by the substitution rules
\[
\begin{array}{ccccc}
\raisebox{-3ex}{\scalebox{0.25}{\includegraphics{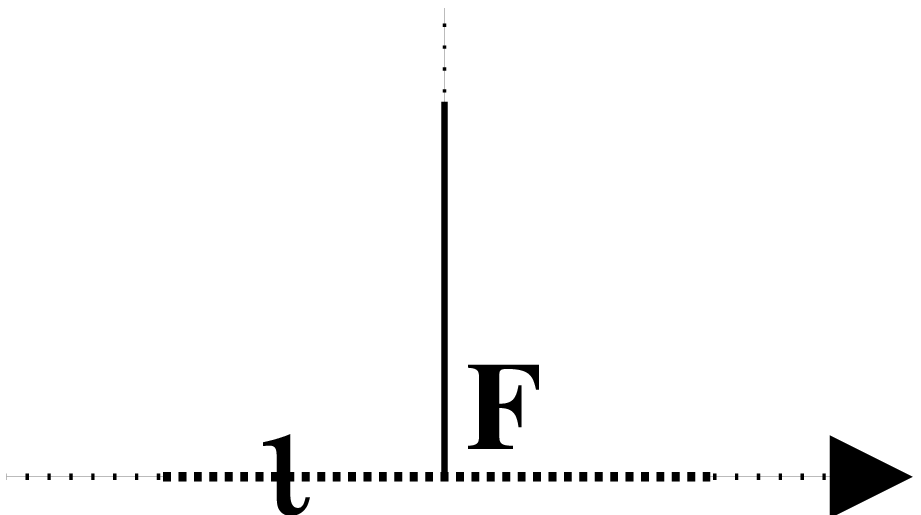}}}
 & \leadsto & 0 & + &
\raisebox{-3ex}{\scalebox{0.25}{\includegraphics{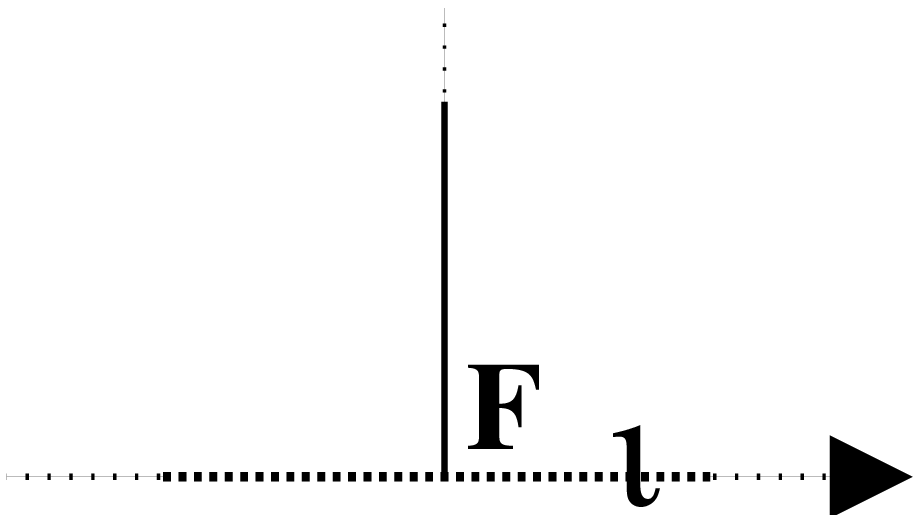}}} \\[0.75cm]
\raisebox{-3ex}{\scalebox{0.25}{\includegraphics{iopa}}} &
\leadsto &
\raisebox{-3ex}{\scalebox{0.25}{\includegraphics{iopb}}} & - &
\raisebox{-3ex}{\scalebox{0.25}{\includegraphics{iopc}}}
\end{array}
\]
This simplification in $\iota$ comes at the cost of a more
complicated differential:
\begin{eqnarray*}
\raisebox{-3ex}{\scalebox{0.25}{\includegraphics{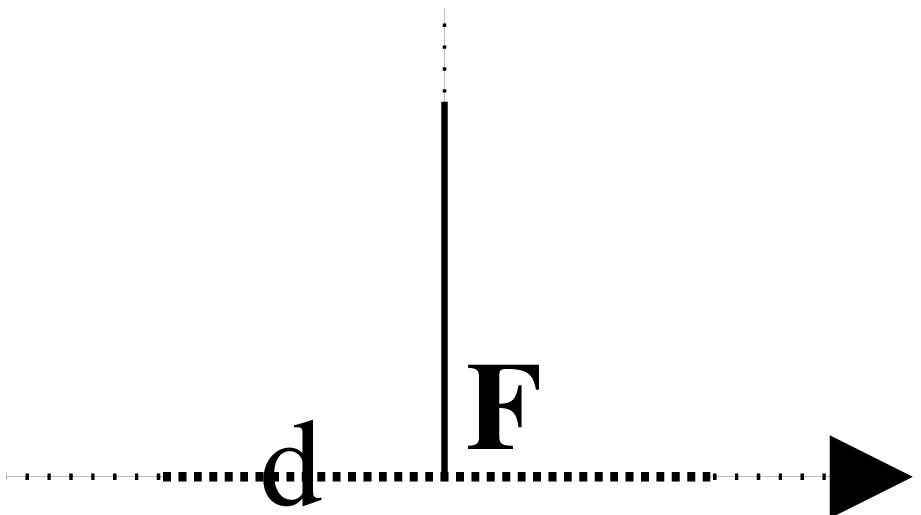}}}
 & \leadsto &
\raisebox{-3ex}{\scalebox{0.25}{\includegraphics{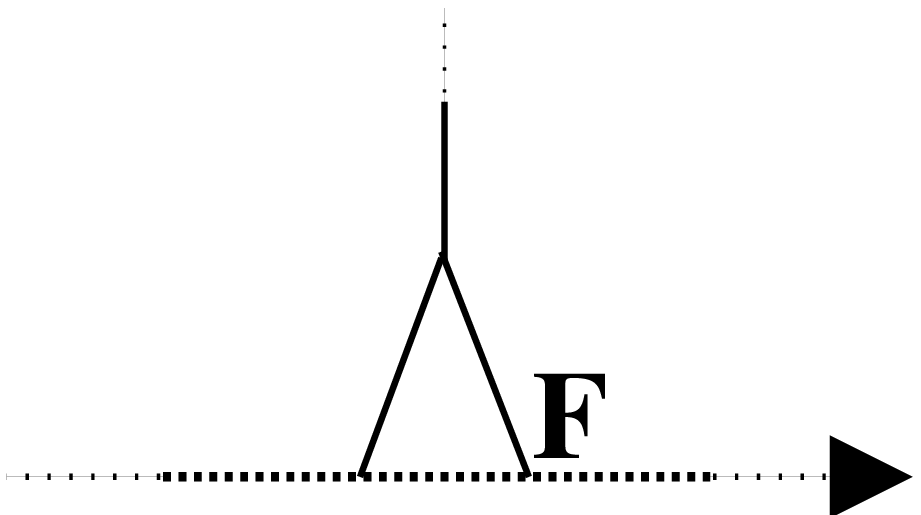}}}
+\raisebox{-3ex}{\scalebox{0.25}{\includegraphics{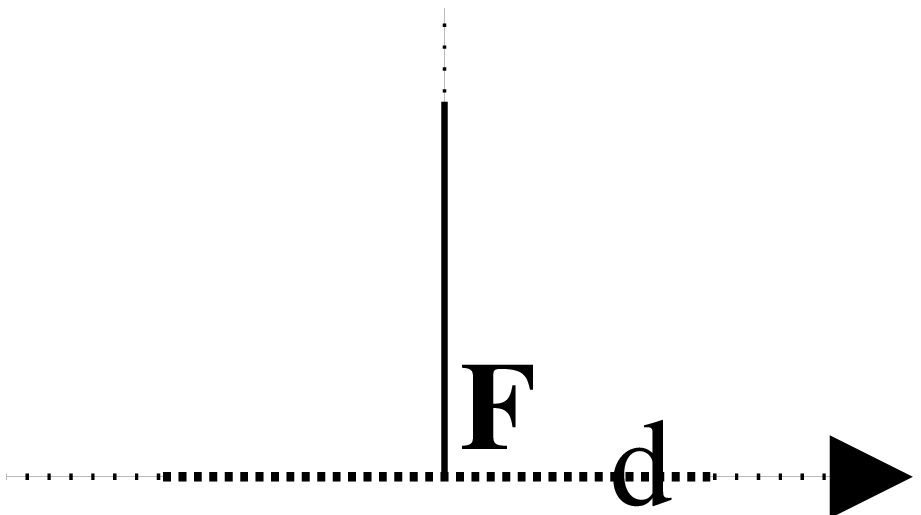}}} \\
\raisebox{-3ex}{\scalebox{0.25}{\includegraphics{diffopB}}} &
\leadsto &
\left(\,\raisebox{-3ex}{\scalebox{0.25}{\includegraphics{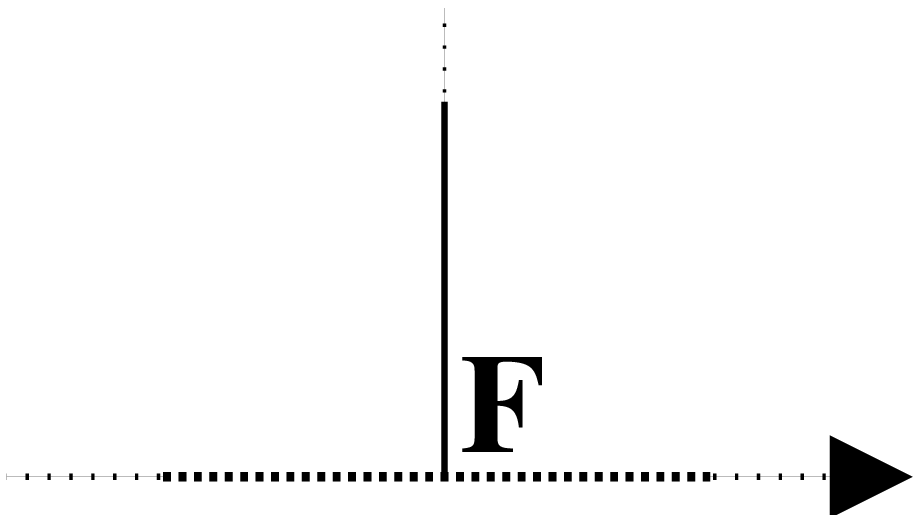}}}
+\frac{1}{2}
\raisebox{-3ex}{\scalebox{0.25}{\includegraphics{curvC}}}\,\right)
- \raisebox{-3ex}{\scalebox{0.25}{\includegraphics{diffopB2}}}
\end{eqnarray*}
With these definitions,
$(\mathcal{W}_{\mathrm{F}},\mathcal{W}_{\mathrm{F}\iota},\iota)$
forms an $\iota$-complex. (This is a quick calculation based on
Lemma \ref{commlem} together with a sweeping argument.)

\begin{defn}
Define maps $\baseFtobull^i: \mathcal{W}^i_{\mathrm{F}}
\rightarrow \mathcal{W}^i$ and
$B_{\mathrm{F}\rightarrow\bullet,\iota}^i:
\mathcal{W}^i_{\mathrm{F}\iota} \rightarrow \mathcal{W}^i_\iota$
by expanding F-labelled legs as follows:
\[
\begin{array}{rccl}
\raisebox{-3ex}{\scalebox{0.25}{\includegraphics{curvA}}} &
\mapsto &
\raisebox{-3ex}{\scalebox{0.25}{\includegraphics{curvB}}}
-\frac{1}{2}
\raisebox{-3ex}{\scalebox{0.25}{\includegraphics{curvC}}}\ .
\end{array}
\]
\end{defn}
For example:
%
\begin{multline*}
\baseFtobull^4\left(\raisebox{-3ex}{\scalebox{0.22}{\includegraphics{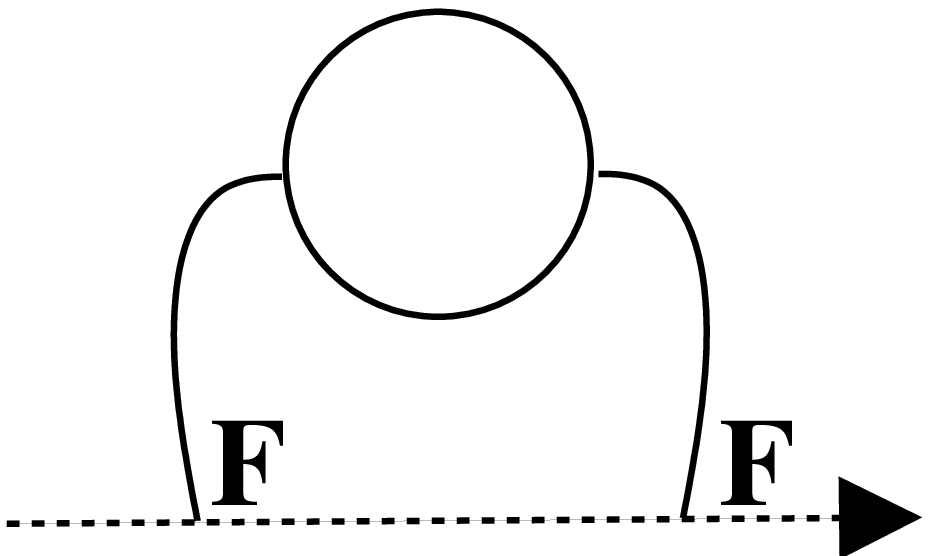}}}\right)
 =  \raisebox{-3ex}{\scalebox{0.22}{\includegraphics{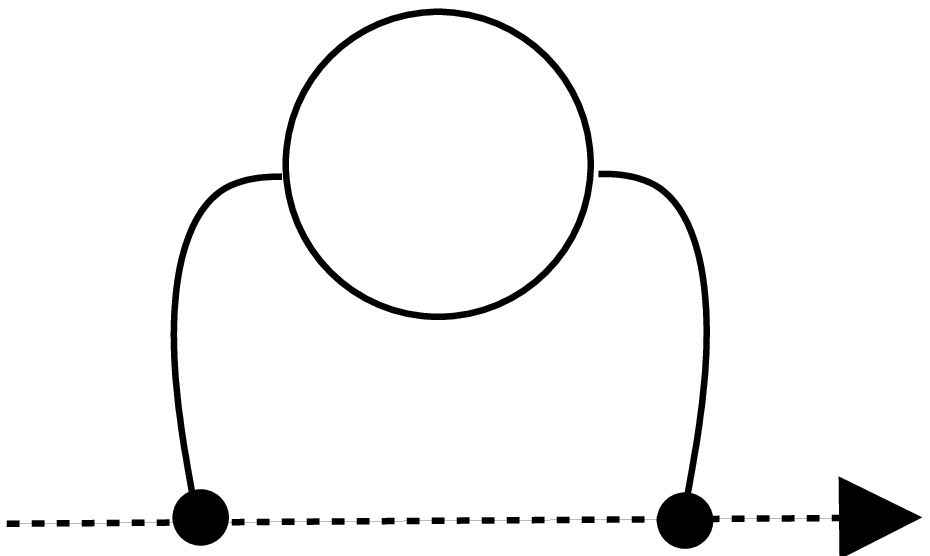}}}
-\frac{1}{2}
\raisebox{-3ex}{\scalebox{0.22}{\includegraphics{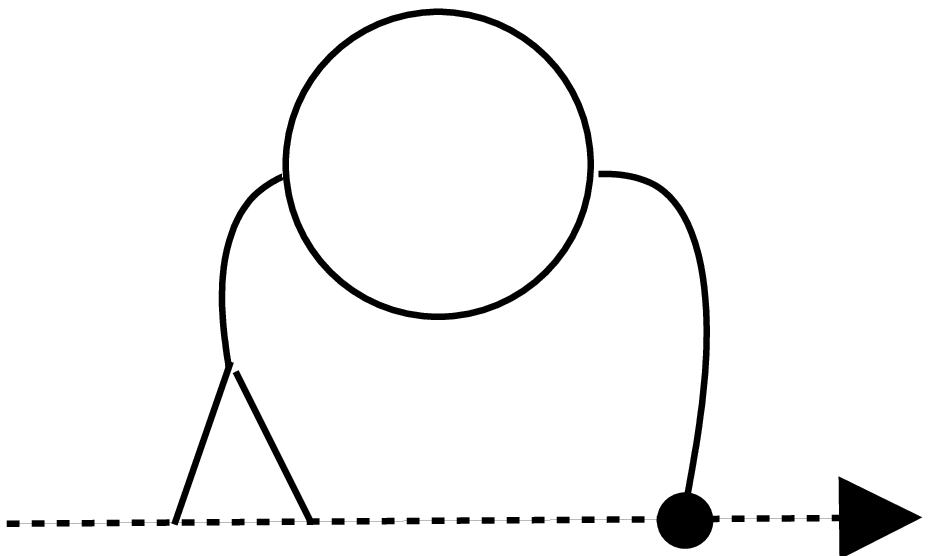}}} \\
   -\frac{1}{2}
\raisebox{-3ex}{\scalebox{0.22}{\includegraphics{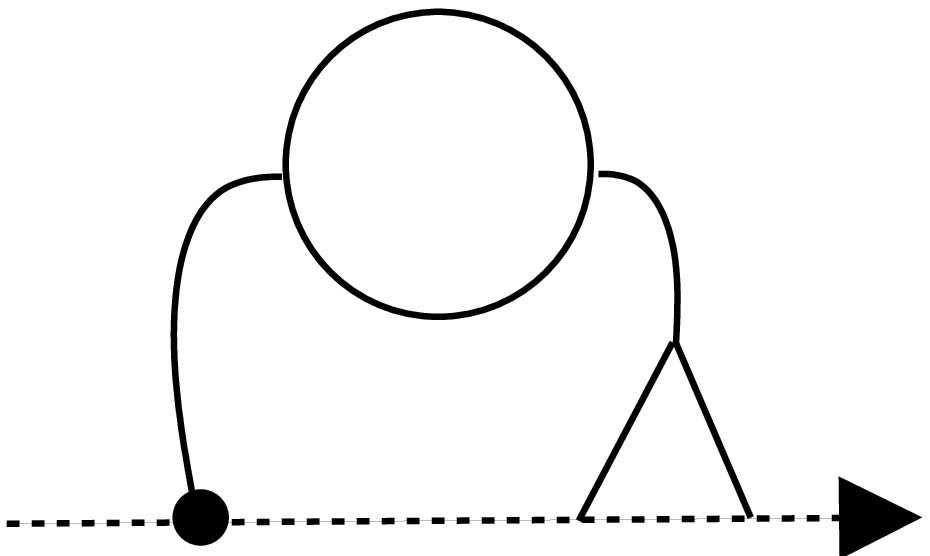}}} \
+\frac{1}{4}
\raisebox{-3ex}{\scalebox{0.22}{\includegraphics{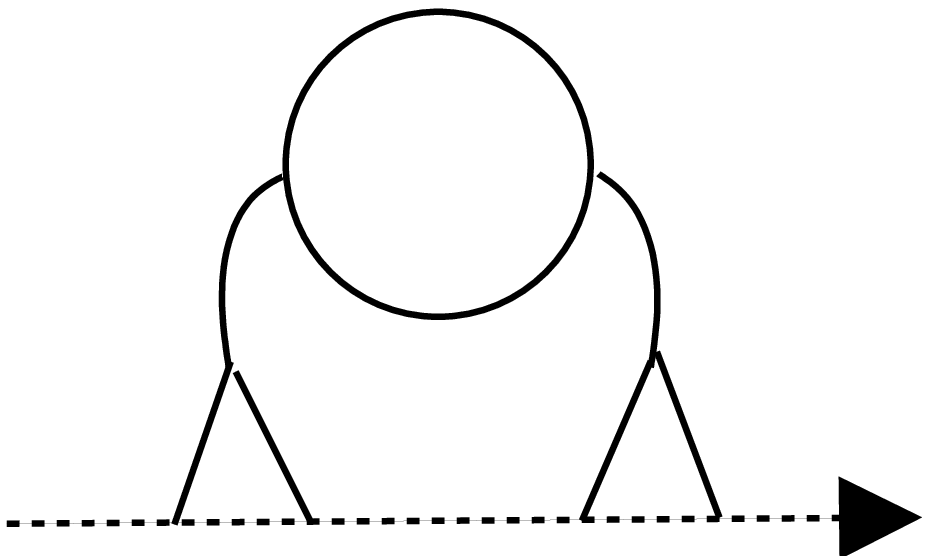}}}.
\end{multline*}

\begin{prop}
The `change of basis' map $\baseFtobull$ is a map of
$\iota$-complexes.
\end{prop}
\begin{proof} Let's begin by showing that the equation $d\circ
\baseFtobull = \baseFtobull \circ d$ holds. First consider how $d$
operates on the expansion of an ${\bf F}$-labelled leg:
\begin{eqnarray*}
\lefteqn{\raisebox{-3ex}{\scalebox{0.25}{\includegraphics{dfA}}}
 =  \raisebox{-3ex}{\scalebox{0.25}{\includegraphics{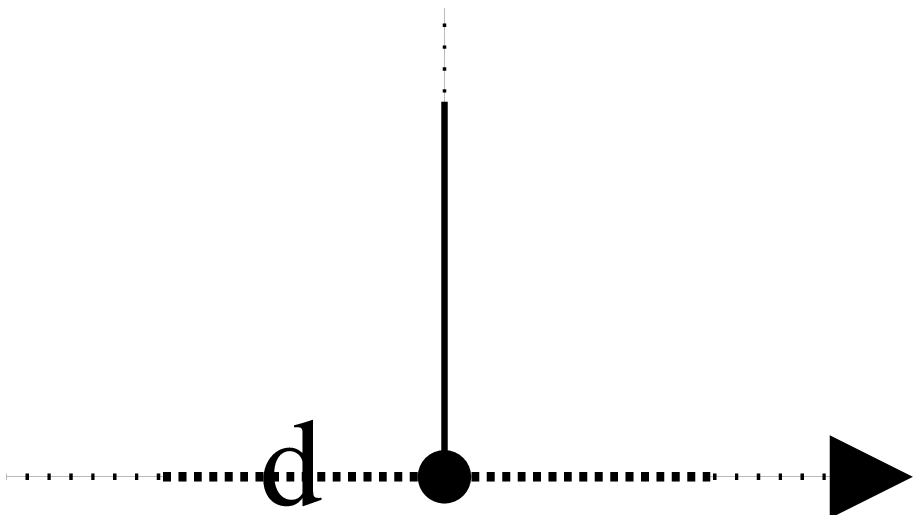}}}
-\frac{1}{2}
\raisebox{-3ex}{\scalebox{0.25}{\includegraphics{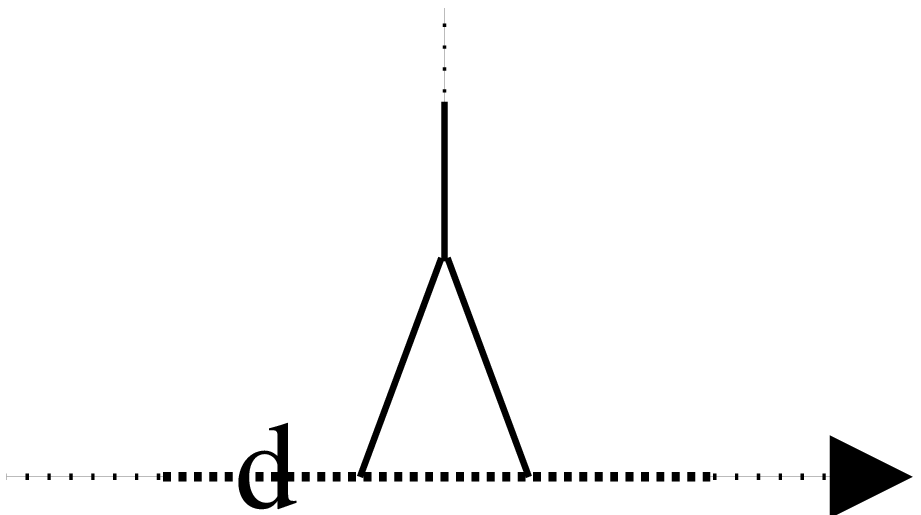}}}} \\[0.1cm]
& \leadsto & 0 +
\raisebox{-3ex}{\scalebox{0.25}{\includegraphics{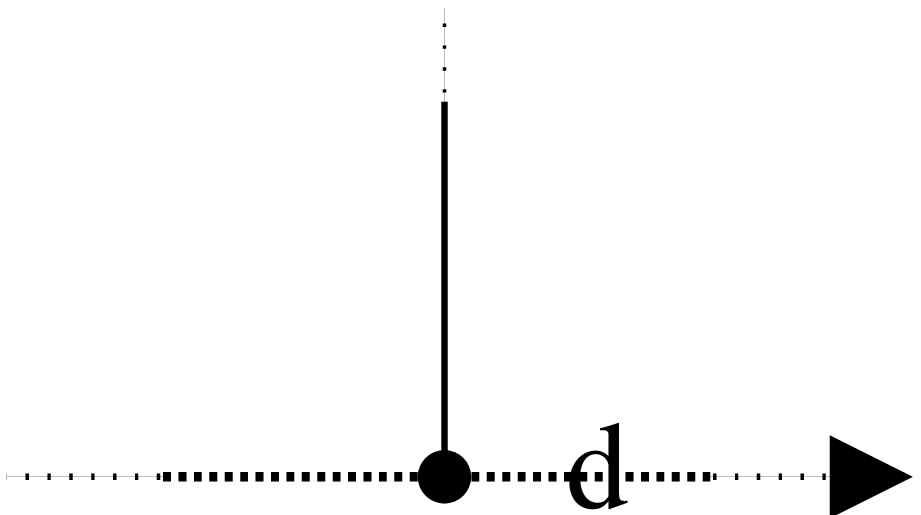}}}
 -\frac{1}{2}
\raisebox{-3ex}{\scalebox{0.25}{\includegraphics{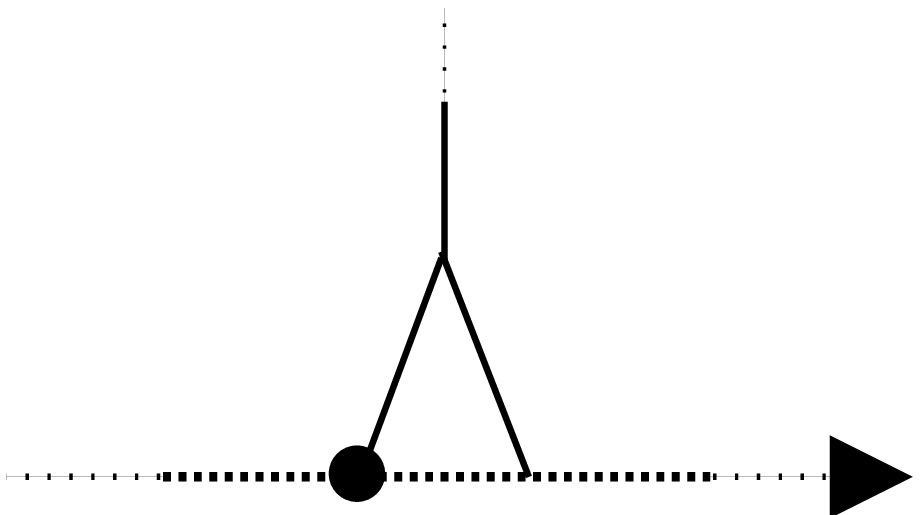}}}
+\frac{1}{2}
\raisebox{-3ex}{\scalebox{0.25}{\includegraphics{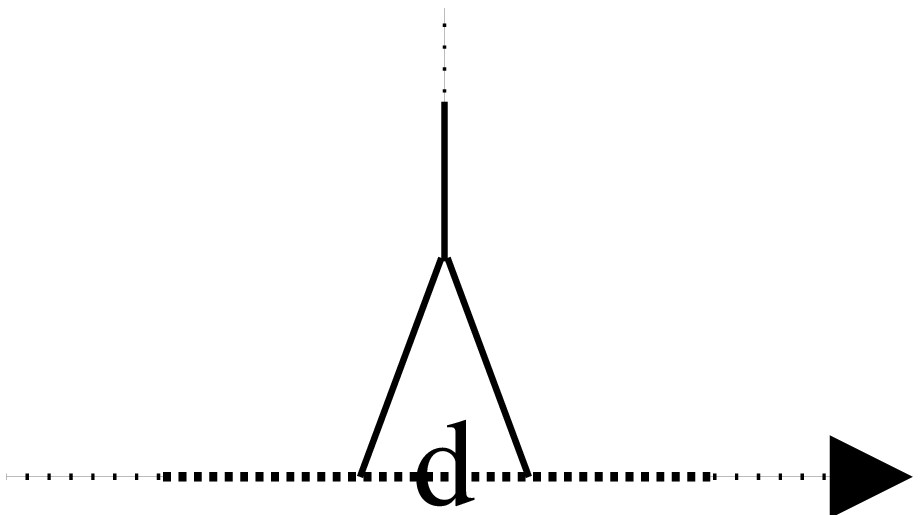}}} \\[0.1cm]
& \leadsto &
\raisebox{-3ex}{\scalebox{0.25}{\includegraphics{dfK}}}
 -\frac{1}{2}
\raisebox{-3ex}{\scalebox{0.25}{\includegraphics{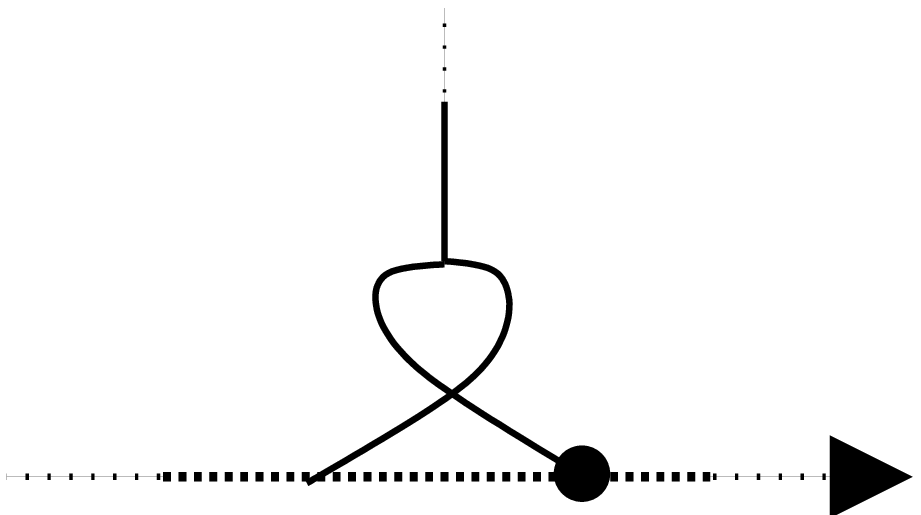}}}
+\frac{1}{2}
\raisebox{-3ex}{\scalebox{0.25}{\includegraphics{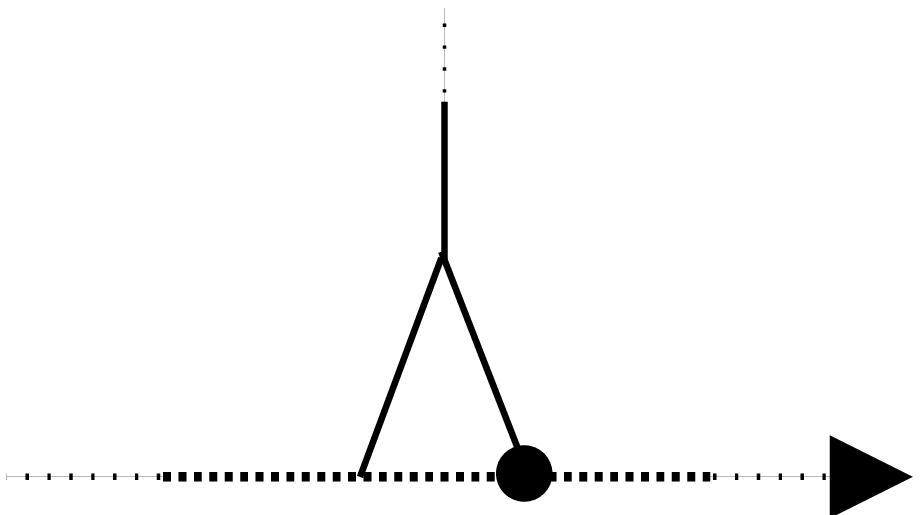}}}
-\frac{1}{2}
\raisebox{-3ex}{\scalebox{0.25}{\includegraphics{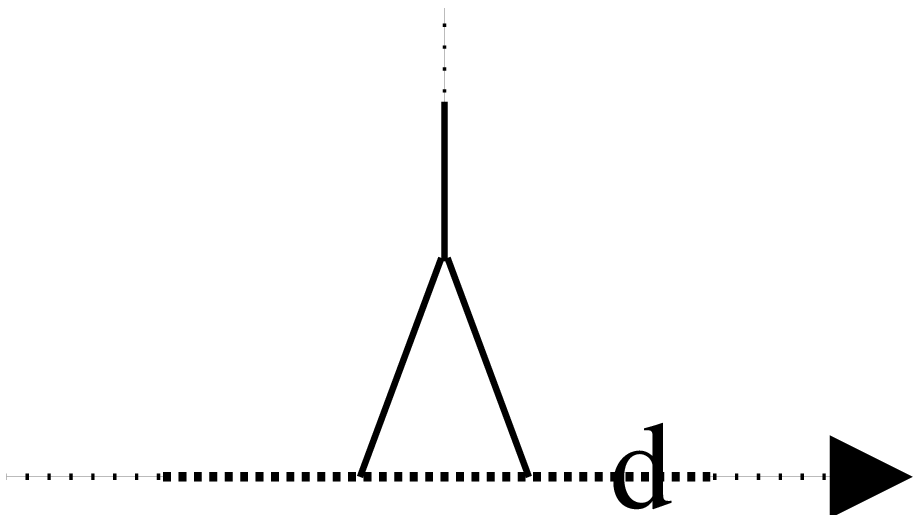}}} \\[0.1cm]
& = & \raisebox{-3ex}{\scalebox{0.25}{\includegraphics{dfF}}} +
\raisebox{-3ex}{\scalebox{0.25}{\includegraphics{dfJ}}} \\[0.1cm]
& = & \raisebox{-3ex}{\scalebox{0.25}{\includegraphics{dfX}}}
+\frac{1}{2}\raisebox{-3ex}{\scalebox{0.25}{\includegraphics{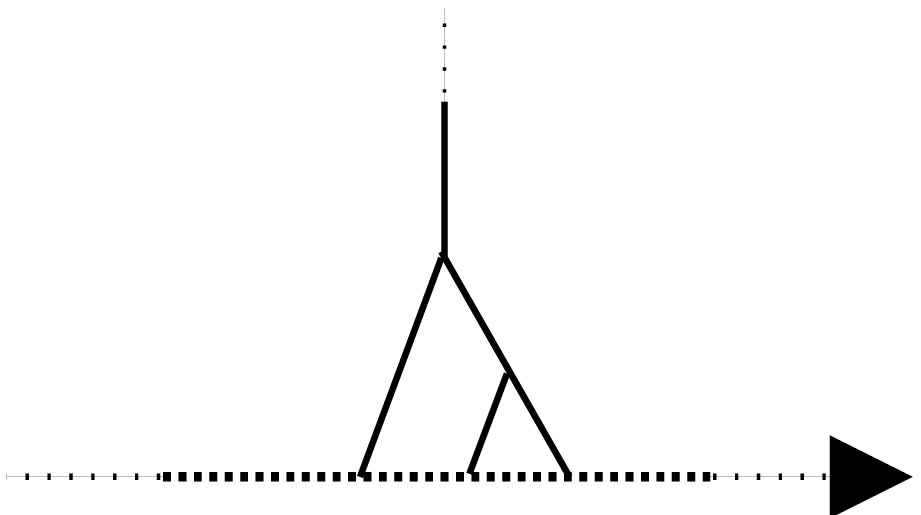}}}
 +
\raisebox{-3ex}{\scalebox{0.25}{\includegraphics{dfJ}}}\ .\\[0.1cm]
\end{eqnarray*}
The required equation holds if the second term above is zero. The
following argument provides that fact:
\begin{multline*}
\raisebox{-3ex}{\scalebox{0.25}{\includegraphics{dfY}}}\ \
\stackrel{(\mathrm{IHX})}{=}\ \
\raisebox{-3ex}{\scalebox{0.25}{\includegraphics{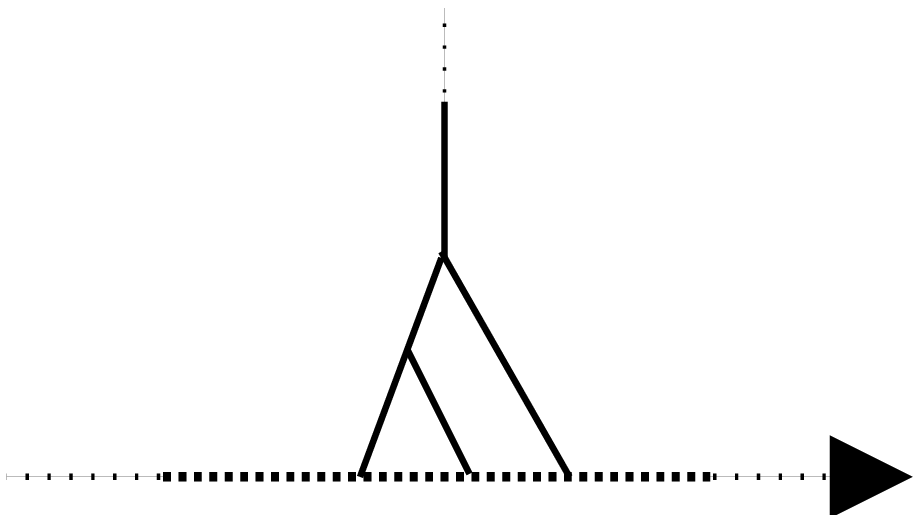}}} +
\raisebox{-3ex}{\scalebox{0.25}{\includegraphics{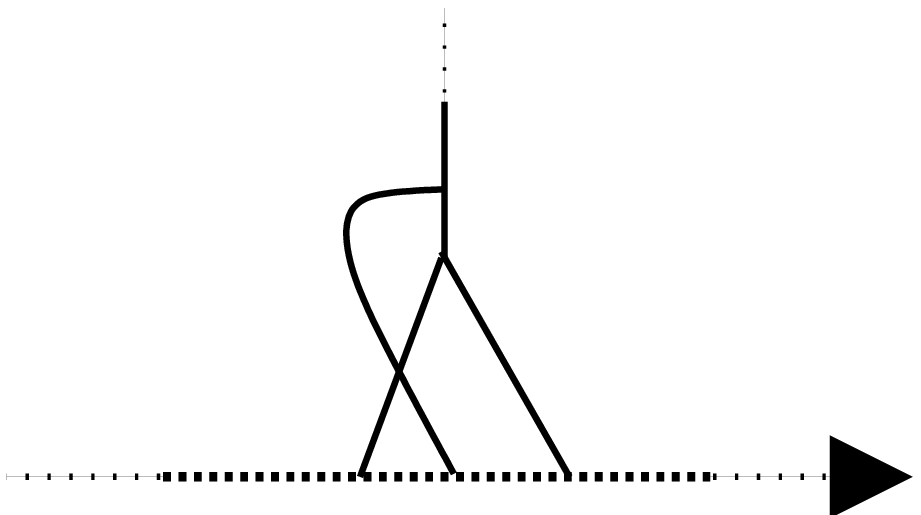}}} \\[0.1cm]
 \stackrel{(\mathrm{Perm})}{=}
\raisebox{-3ex}{\scalebox{0.25}{\includegraphics{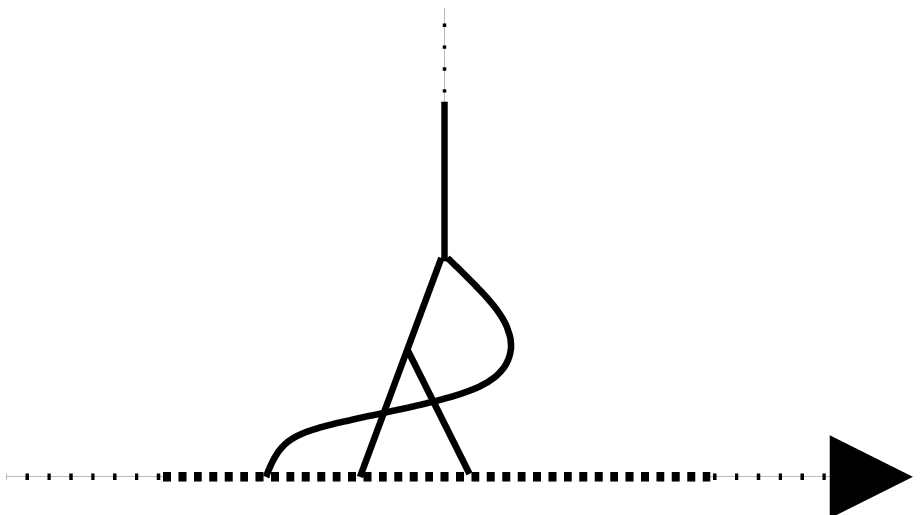}}} -
\raisebox{-3ex}{\scalebox{0.25}{\includegraphics{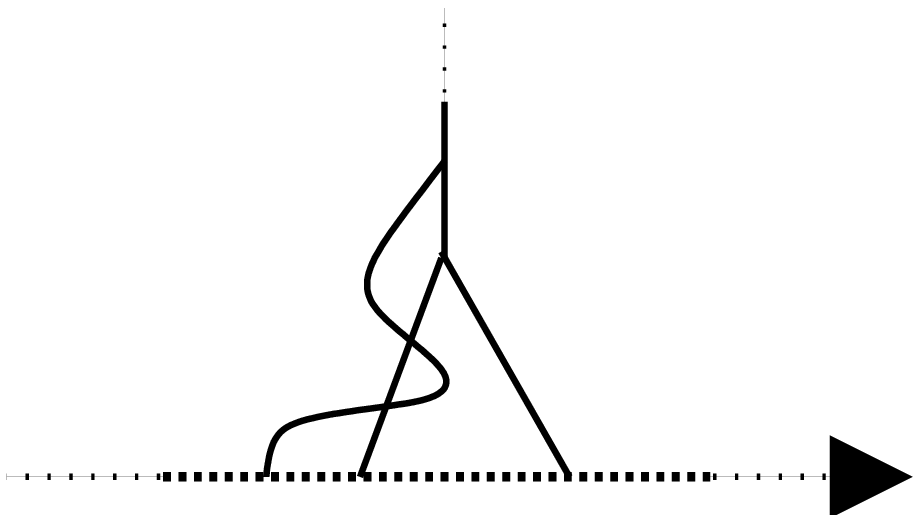}}}\ \  \stackrel{(\mathrm{AS})}{=}\ \
-2\raisebox{-3ex}{\scalebox{0.25}{\includegraphics{dfY}}}.
\end{multline*}
The case of the action of $d$ on the grade $1$ legs is treated
straightforwardly. The equation $\iota\circ \baseFtobull =
\baseFtobull \circ \iota$ is also straightforward.
\end{proof}
\begin{prop} \label{basisiso}The maps
$\baseFtobull^i$ and $B_{\mathrm{F}\rightarrow\bullet,\iota}^i$
are vector space isomorphisms, and so yield isomorphisms in basic
cohomology:
\[
H^i\left((\baseFtobull)_{\mathrm{basic}}\right) :
H^i((\mathcal{W}_{\mathrm{F}})_\mathrm{basic})
\stackrel{\cong}{\longrightarrow}
H^i(\mathcal{W}_\mathrm{basic}).
\]
\end{prop}
\begin{proof} The inverse $\basebulltoF^i:
\mathcal{W}^i \rightarrow \mathcal{W}^i_{\mathrm{F}}$ is given by
replacing legs as follows:
\[
\begin{array}{rccl}
\raisebox{-3ex}{\scalebox{0.24}{\includegraphics{curvB}}} &
\mapsto &
\raisebox{-3ex}{\scalebox{0.24}{\includegraphics{curvA}}}
+\frac{1}{2}
\raisebox{-3ex}{\scalebox{0.24}{\includegraphics{curvC}}}
\end{array}
\]
\end{proof}

\subsection{The calculation of the basic cohomology of
\((\Wspace,\Wspace_\iota,\iota)\)} \label{basicsubcomplexintro} We
now turn to the following factor in the definition of $\Upsilon$:
\[ \underbrace{(\Bspace,0,0)
\xrightarrow{\ \phi_\Bspace\ }
(\mathcal{W}_{\mathrm{F}},\mathcal{W}_{\mathrm{F}\iota},\iota)}
\xrightarrow{\baseFtobull} (\Wspace,\Wspace_\iota,\iota).
\]
We begin by assembling the spaces ${\Bspace^i}$ into an
$\iota$-complex in the following way (recalling that symmetric
Jacobi diagrams are graded by 2-per-leg):
\[
\xymatrix{ 0 \ar[r] & \mathcal{B}^0 \ar[r]^{0} \ar[dl]_{0}& 0
\ar[r]^{0} \ar[dl]_{0}&\mathcal{B}^2 \ar[r]^{0} \ar[dl]_{0} & 0
\ar[r]^{0} \ar[dl]_{0}&
 \mathcal{B}^4 \ar[dl]_{0} \ar[r]^{0} & 0 \ar[r]^{0} \ar[dl]_{0}&
 \ldots \\
0 \ar[r]^0 &0 \ar[r]^0 &0 \ar[r]^0 & 0 \ar[r]^{0} & 0 \ar[r]^{0}
& 0 \ar[r]^{0} & \ldots \\
}
\]
Note the obvious fact that
\[
H^i(\Bspace_\mathrm{basic}) = \left\{
\begin{array}{cc}
\Bspace^{i} & \mbox{if $i$ is even,} \\[0.15cm]
0 & \mbox{if $i$ is odd.}
\end{array}
\right.
\]


Define a map $\phi_\Bspace^i: \Bspace^i \rightarrow \Wspace_\mathrm{F}^i$
on some symmetric Jacobi diagram $D$ by simply choosing some
ordering of its legs and labelling each of them with an {\bf F}.
For example:
\[
\phi_\Bspace^8\left(\raisebox{-3.5ex}{\scalebox{0.26}{\includegraphics{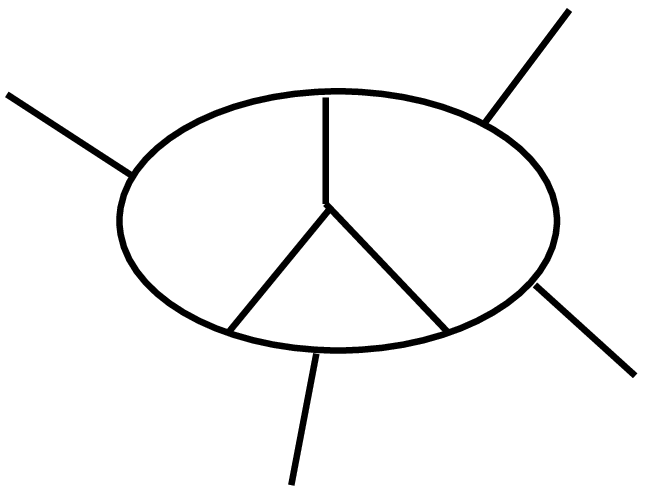}}}\right)
\ =\
\raisebox{-5.5ex}{\scalebox{0.26}{\includegraphics{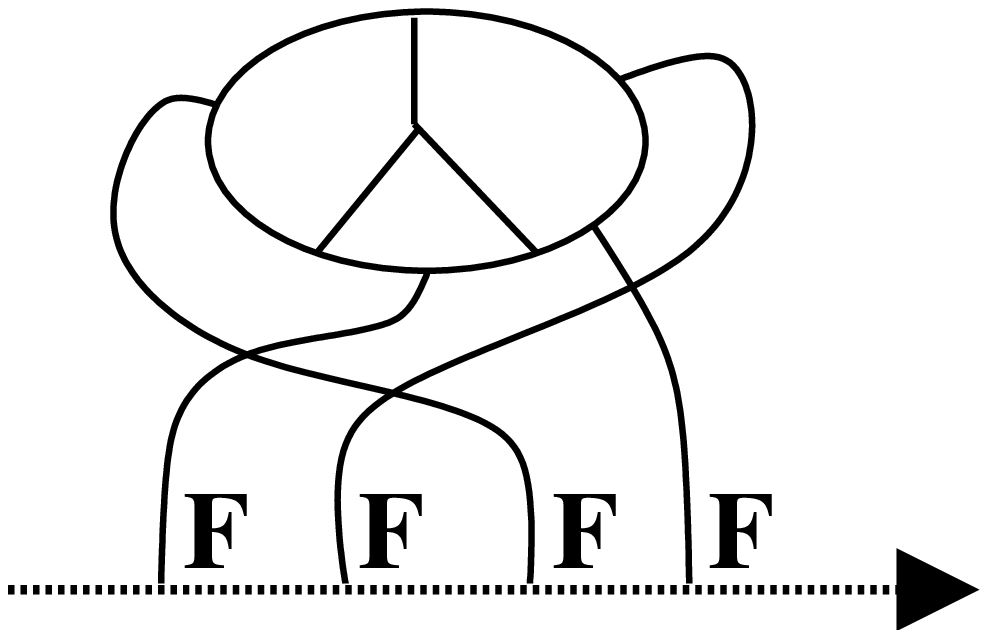}}}\ \ .
\]
\begin{prop}\label{phiisaniotamap}
The map $\phi_\Bspace$ is a map of $\iota$-complexes.
\end{prop}
\begin{proof}
The only thing we need to show is that if $D$ is a symmetric
Jacobi diagram then $\iota(\phi_\Bspace(D))=0$ and $d(\phi_\Bspace(D))=0$. The
first equation is clear (given the defined operation of $\iota$ on
curvature legs). The following example will show why the second
equation holds:
\begin{eqnarray*}
\lefteqn{ \raisebox{-3ex}{\scalebox{0.28}{\includegraphics{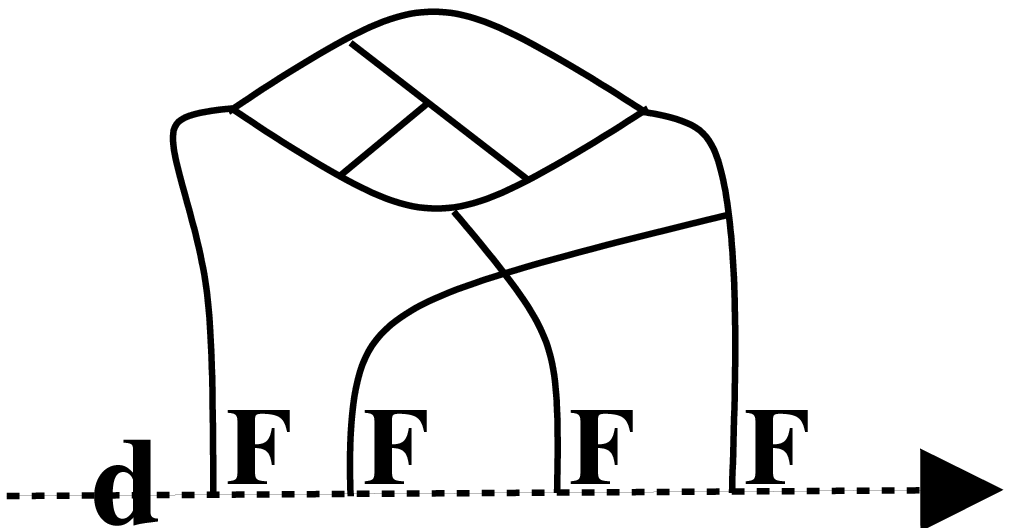}}} } \\
& \leadsto &
\raisebox{-3ex}{\scalebox{0.28}{\includegraphics{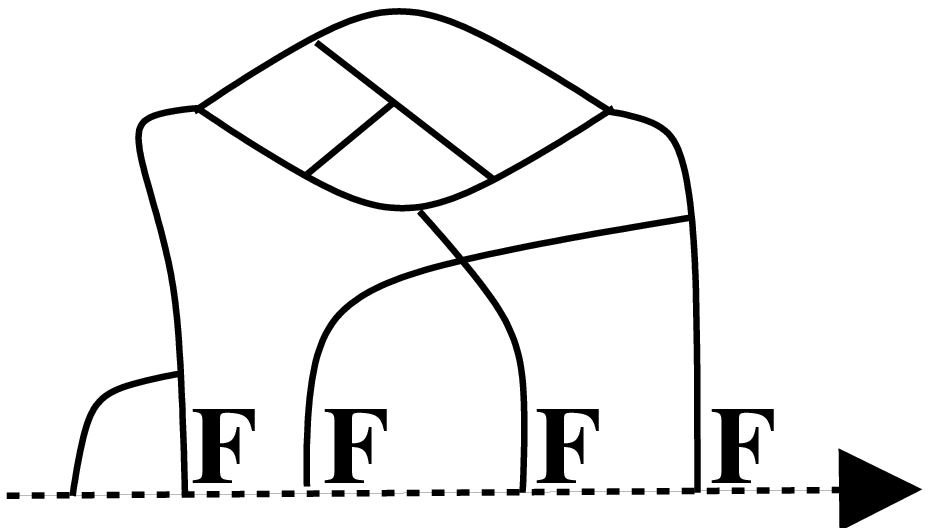}}}
+\raisebox{-3ex}{\scalebox{0.28}{\includegraphics{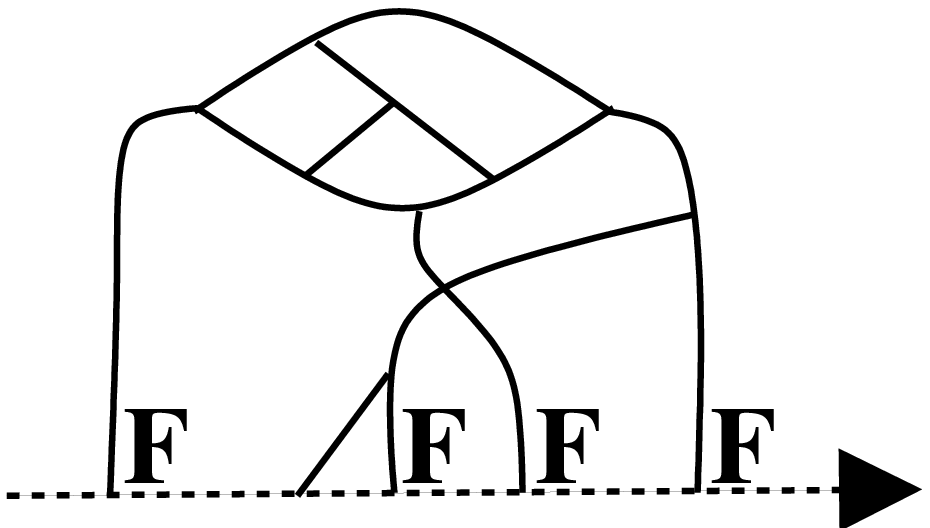}}}
+\raisebox{-3ex}{\scalebox{0.28}{\includegraphics{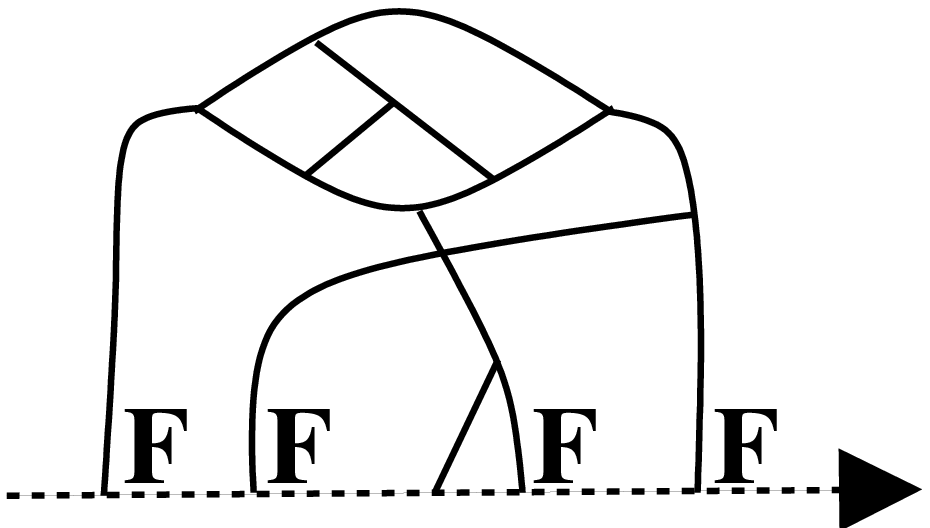}}}
+\raisebox{-3ex}{\scalebox{0.28}{\includegraphics{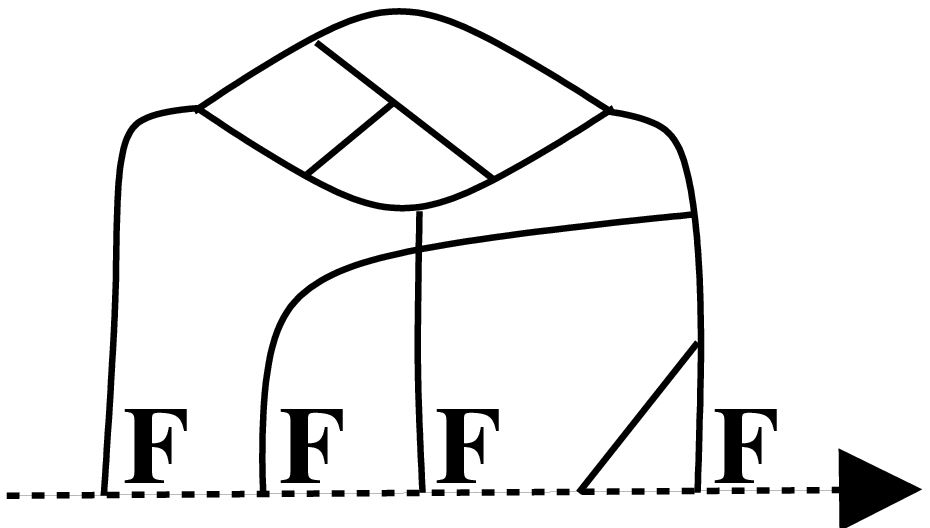}}} \\
& = & \raisebox{-3ex}{\scalebox{0.28}{\includegraphics{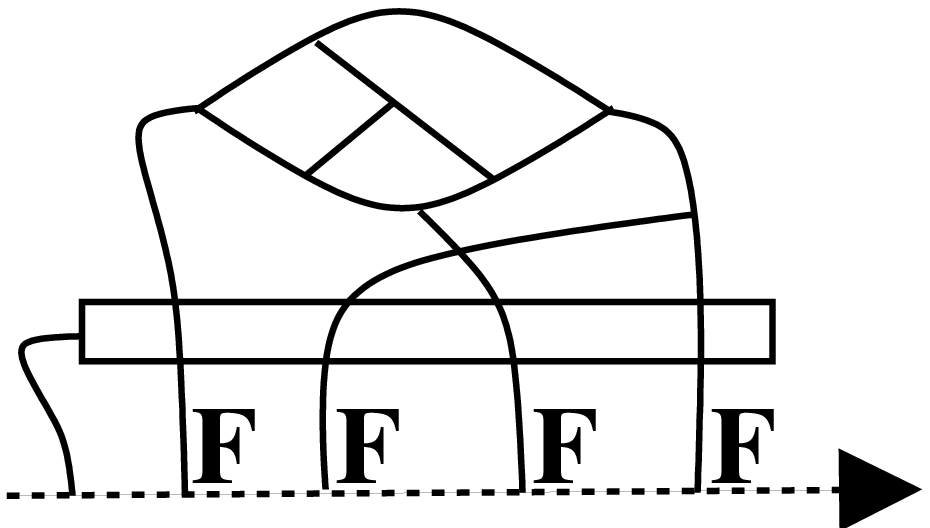}}}\ \
=\ \ 0.
\end{eqnarray*}
\end{proof}
\begin{prop}\label{phiiso} The map $\phi_\Bspace$ yields isomorphisms on
basic cohomology:
\[
H^i((\phi_\Bspace)_\mathrm{basic}):
H^i(\Bspace_\mathrm{basic})\stackrel{\cong}{\longrightarrow}
H^i((\Wspace_\mathrm{F})_\mathrm{basic}).
\]
\end{prop}
\begin{proof} (In this proof the symbol $\iota$ refers to the $\iota$ map in the complex $\Wspace_{\mathrm{F}}$.) It suffices to show that:
\[ \left\{
\begin{array}{lp{1cm}l}
\mbox{$\phi_\Bspace^j:\Bspace^j \stackrel{\cong}{\longrightarrow}
\mathrm{ker}\left(\iota^j\right)$} & & \mbox{if $j$ is even,} \\[0.05cm]
\mathrm{ker}\left(\iota^j\right)=0 & & \mbox{if $j$ is odd.}
\end{array}\right.
\]
To determine ker$\left(\iota^j\right)$: observe that there is a
direct-sum decomposition
\[
\Wspace^i_\mathrm{F} = \oplus_{j=0}^{i}\Wspace^{i,j}_\mathrm{F}
\]
where $\Wspace^{i,j}_\mathrm{F}$ is the subspace generated by Weil
diagrams with exactly $j$ grade 1 legs.
This direct-sum
decomposition exists because there are no relations involving
diagrams with varying numbers of such legs. Note that
$\Wspace^{i,0}$ is obviously isomorphic to $\Bspace^i$.

The proposition follows immediately from the claim that
ker$(\iota^i)=\Wspace_\mathrm{F}^{i,0}$. Note that it is immediate
from the action of $\iota$ on {\bf F}-labelled legs that
ker$(\iota^i)\supset\Wspace_\mathrm{F}^{i,0}$. To prove that
ker$(\iota^i)\subset\Wspace_\mathrm{F}^{i,0},$ let
$\hat{\iota}^i:\Wspace_{\mathrm{F}\iota}^{i-1}\rightarrow
\Wspace_{\mathrm{F}}^i$ be the map which takes the special
$\iota$-labelled leg and places it at the far left hand end of the
orienting line. For example: \[ \hat{\iota}^6\left(
\raisebox{-4ex}{\scalebox{0.25}{\includegraphics{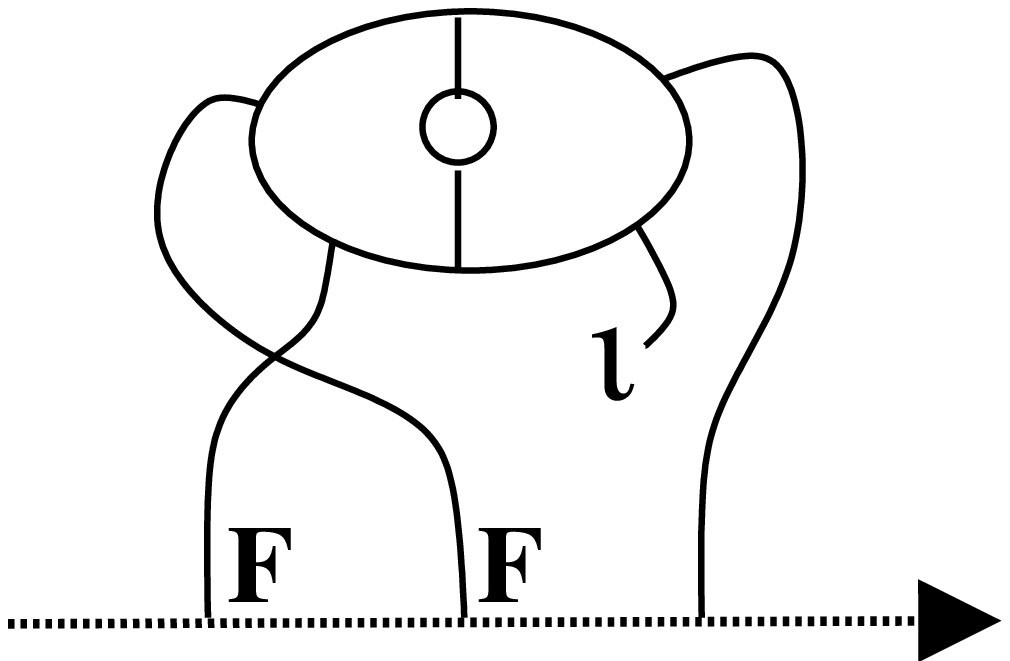}}}
\right)\ \ =\ \
\raisebox{-4ex}{\scalebox{0.25}{\includegraphics{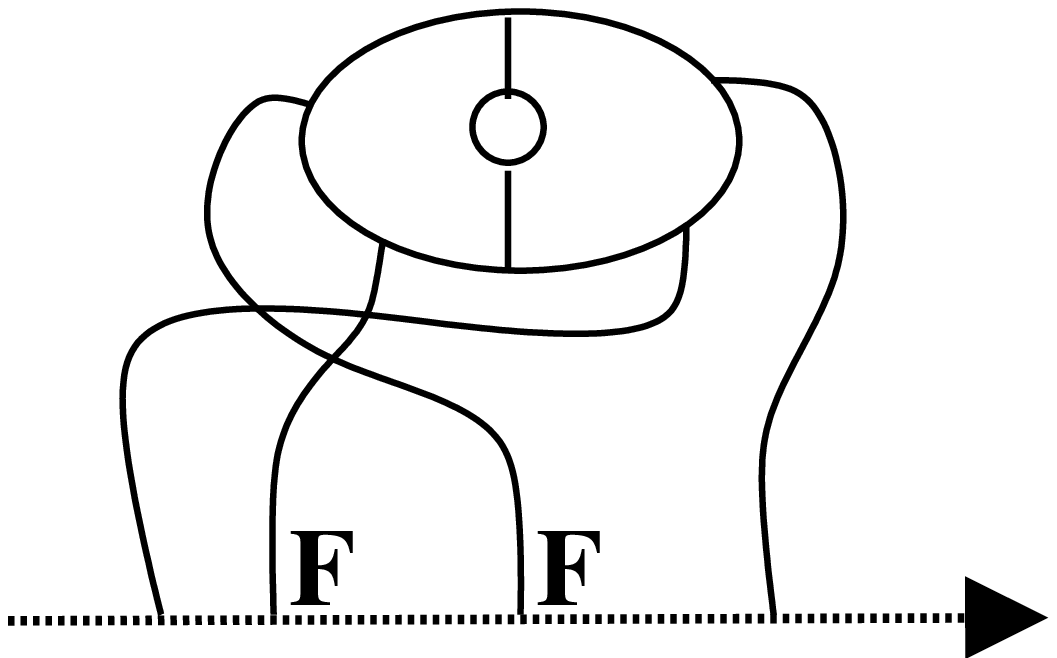}}}
\]
Observe that $(\hat{\iota}^i\circ\iota^i)(D) = jD$, if $D\in
\Wspace^{i,j}_\mathrm{F}$. Thus
ker$(\iota^i)\subset\mathrm{ker}(\hat{\iota}^i\circ\iota^i)=\Wspace^{i,0}$,
establishing the claim.
\end{proof}

To summarize:
If we define maps $\Upsilon^{i}$ by the composition $\Upsilon^i = \baseFtobull^i \circ \phi_\Bspace^i$,
then Propositions \ref{basisiso} and  \ref{phiiso} give us the
following theorem.
\begin{thm}
The induced maps
$H^i( \Upsilon_\mathrm{basic} ) : \Bspace^i \cong
H^i(\Bspace_\mathrm{basic})\rightarrow H^i(
\Wspace_\mathrm{basic})$
are isomorphisms.
\end{thm}


\section{The non-commutative Weil complex.}\label{noncommsect}

In this section we embed the commutative Weil complex for diagrams
into a larger $\iota$-complex, the {\it non-commutative} Weil
complex for diagrams:
\[
\chi_\Wspace : (\Wspace,\Wspace_\iota,\iota) \to
(\ncw,\ncw_\iota,\iota).
\]
This larger complex is built in the same way as $\Wspace$ but
without introducing the permutation relations; the embedding
$\chi_\Wspace$ is the graded averaging map. The key technical
theorem regarding this $\iota$-complex says, in essence, that
every cocycle $z$ in $\ncw_{\mathrm{basic}}$ has its corresponding
basic cohomology class $[z]$ represented by its graded
symmetrisation $\left(\chi_\Wspace\circ\tau\right)(z)$.

\begin{defn}
Define vector spaces $\ncw^i$ in exactly the same way as the
vector spaces $\mathcal{W}^i$ but without introducing the
permutation relations.
\end{defn}
For the purposes of clarity, we'll always draw generators of
$\ncw$, which we'll refer to as {\bf non-commutative Weil
diagrams}, with a certain arrow head on their orienting lines.
Observe: \[
\raisebox{-3ex}{\scalebox{0.25}{\includegraphics{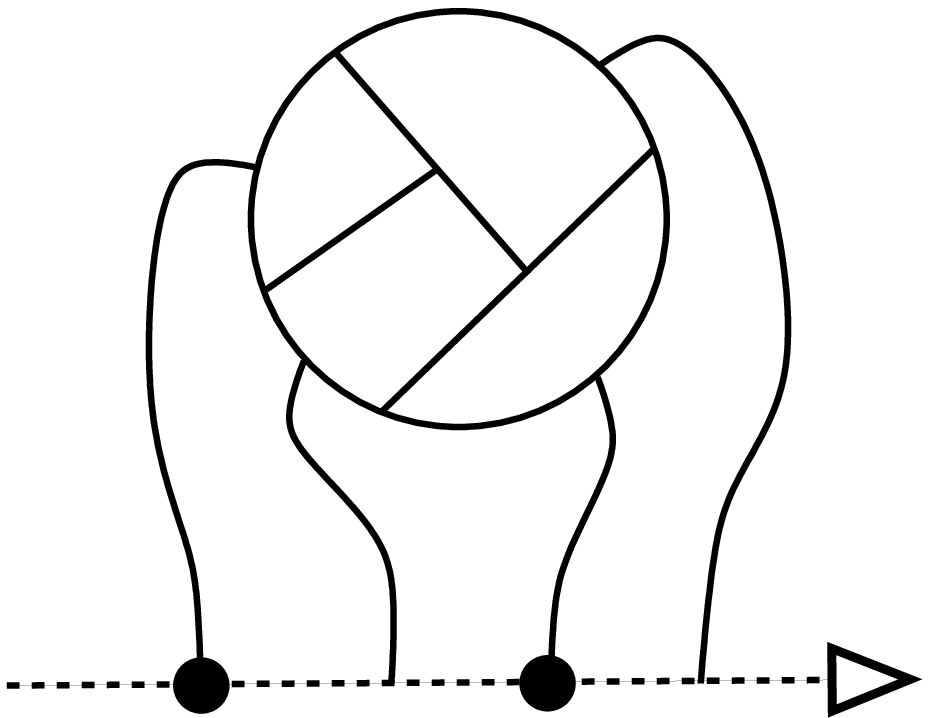}}} \
\ \ \neq\ \ \
\raisebox{-3ex}{\scalebox{0.25}{\includegraphics{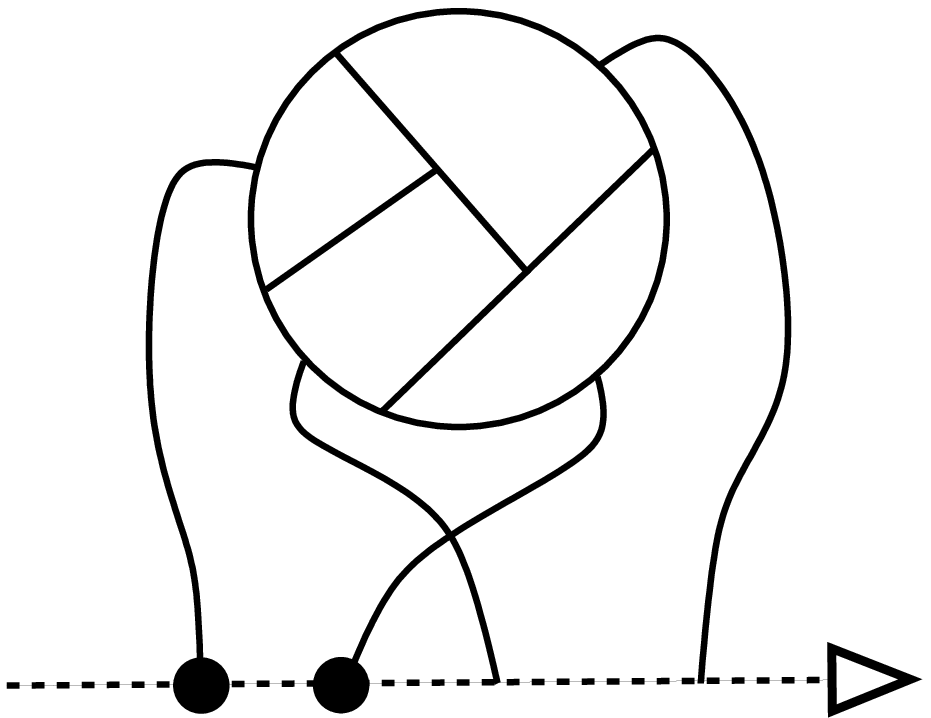}}}
\ \ \ \neq\ \ \ -\
\raisebox{-3ex}{\scalebox{0.25}{\includegraphics{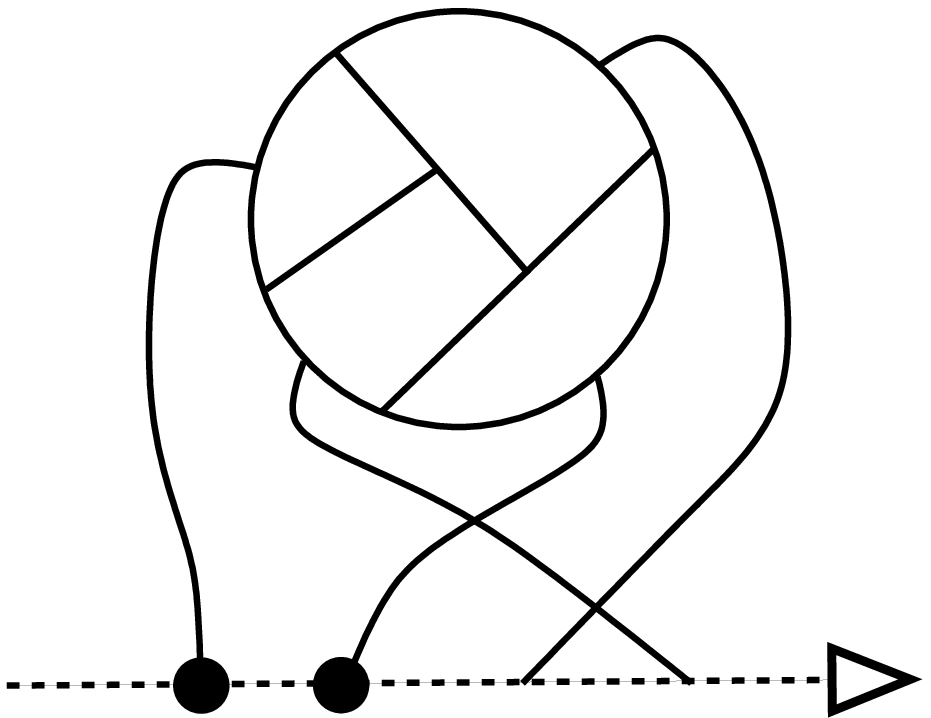}}}\
\ \ \ \ \mbox{in\ \ $\widetilde{\mathcal{W}}^6$.}
\]
Define $d$, $\ncw_\iota$, and $\iota$ in the obvious way. These
definitions form an $\iota$-complex (the calculation that
$[d,\iota]=0$ is exactly the same here as it is for $\Wspace$, Lemma \ref{wisanidga}).

The basic subcomplex here seems much harder to construct
explicitly than it is in the commutative case. It is a startling
and central fact, however, that the two $\iota$-complexes have
canonically isomorphic basic cohomology spaces. We'll show this by
constructing a chain equivalence of $\iota$-complexes:
\[
\xymatrix{ \mathcal{W} \ar@/^/[r]^{\chi_\Wspace}&
\widetilde{\mathcal{W}}\ar@/^/[l]^\tau}.
\]
\begin{defn}
Define the {\bf graded averaging maps}
\[
\chi_\Wspace^i : \mathcal{W}^i \rightarrow
\widetilde{\mathcal{W}}^i\ \ \ \mbox{and}\ \ \
\chi_{\Wspace,\iota}^i : \mathcal{W}^i_\iota \rightarrow
\widetilde{\mathcal{W}}^i_\iota
\]
by declaring the value of $\chi_\Wspace$ (and similarly for
$\chi_{\Wspace,\iota}$) on some commutative Weil diagram to be the
average of the (signed) orderings of its legs. For example:
\begin{eqnarray*}
\chi_\Wspace^4\left(
\raisebox{-3ex}{\scalebox{0.22}{\includegraphics{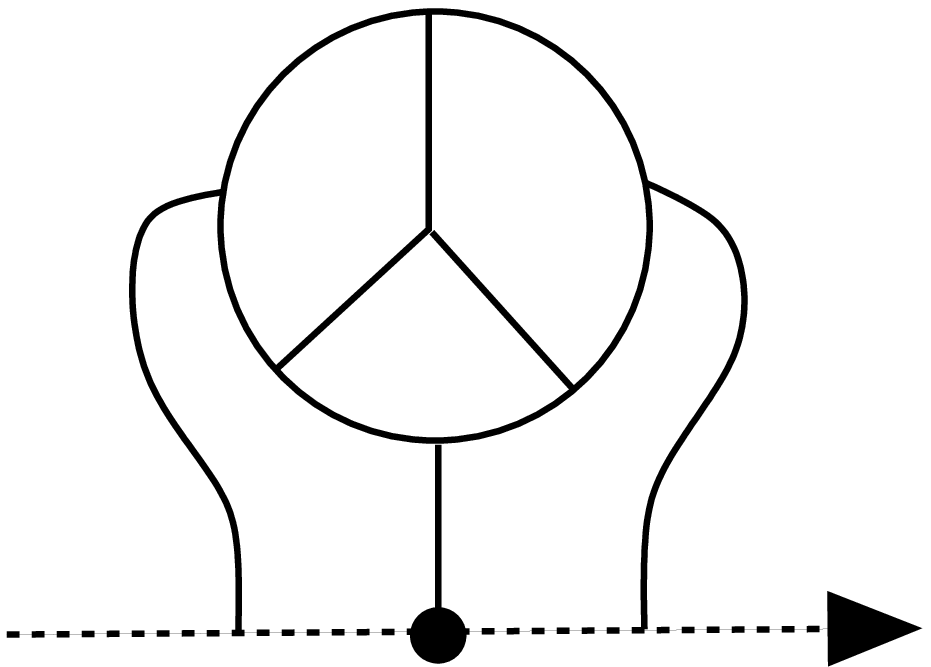}}}
\right) & = & \frac{1}{6!}\left(
\raisebox{-3ex}{\scalebox{0.22}{\includegraphics{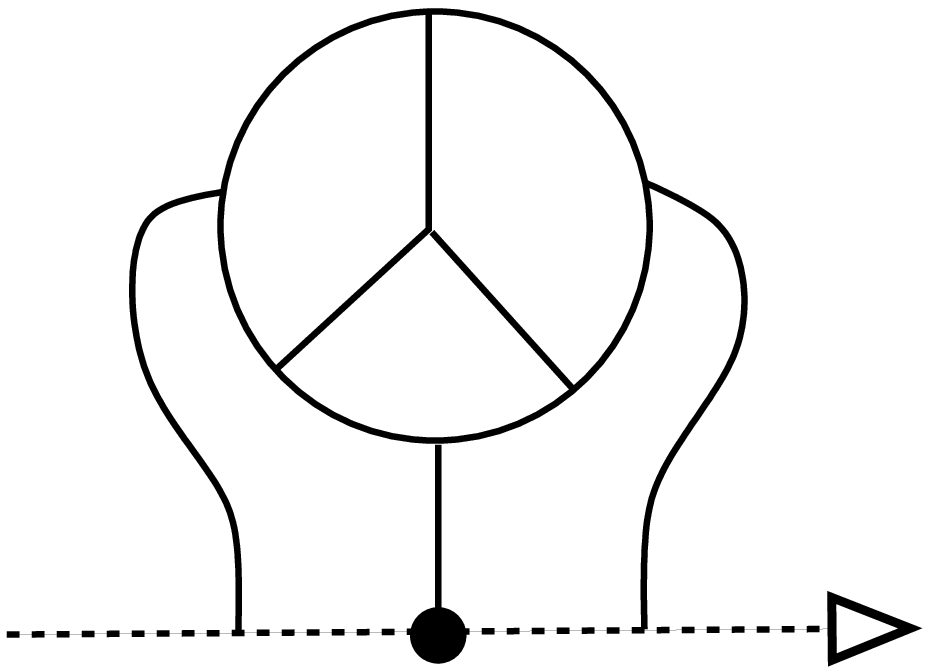}}}
+
\raisebox{-3ex}{\scalebox{0.22}{\includegraphics{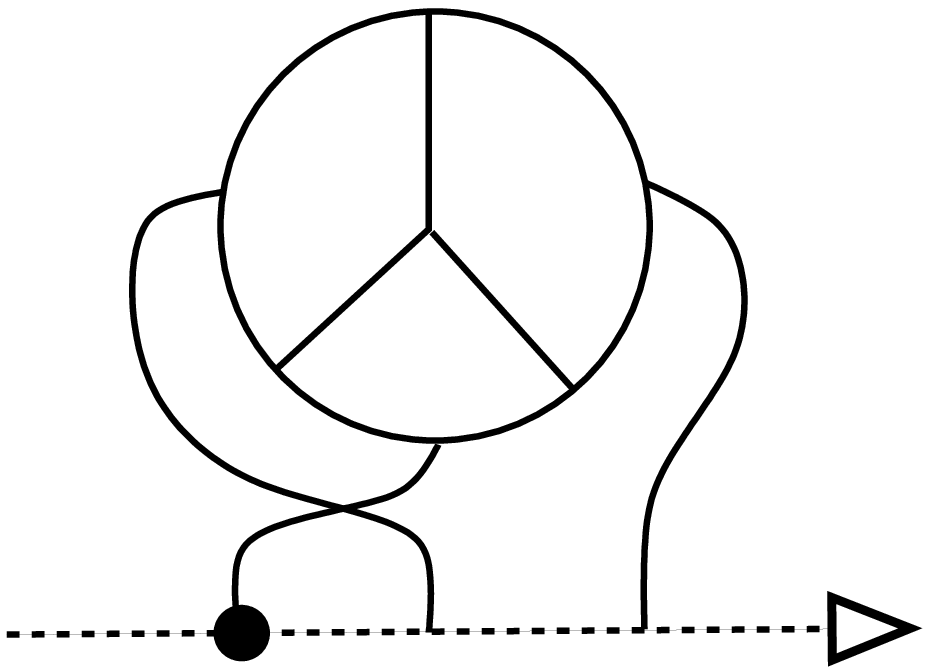}}}
-
\raisebox{-3ex}{\scalebox{0.22}{\includegraphics{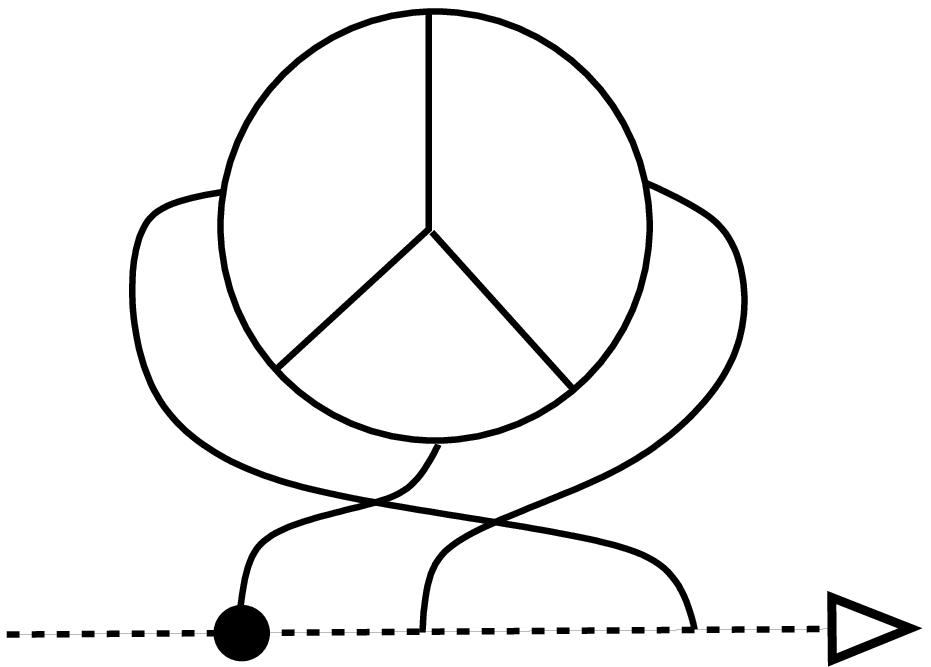}}}
\right.
\\ & & \left.
\ \ \ \,-
\raisebox{-3ex}{\scalebox{0.22}{\includegraphics{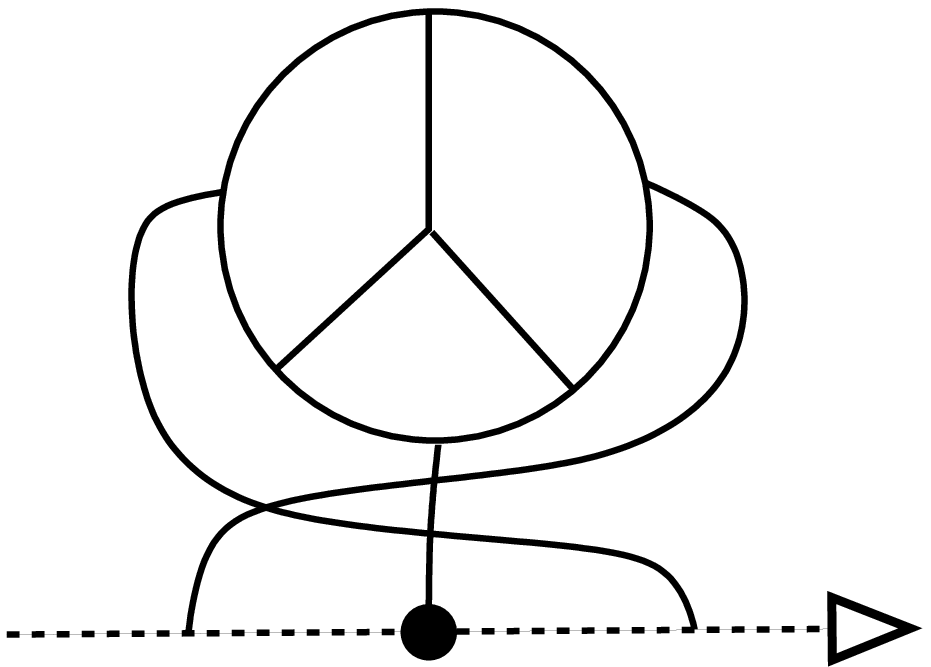}}}
-
\raisebox{-3ex}{\scalebox{0.22}{\includegraphics{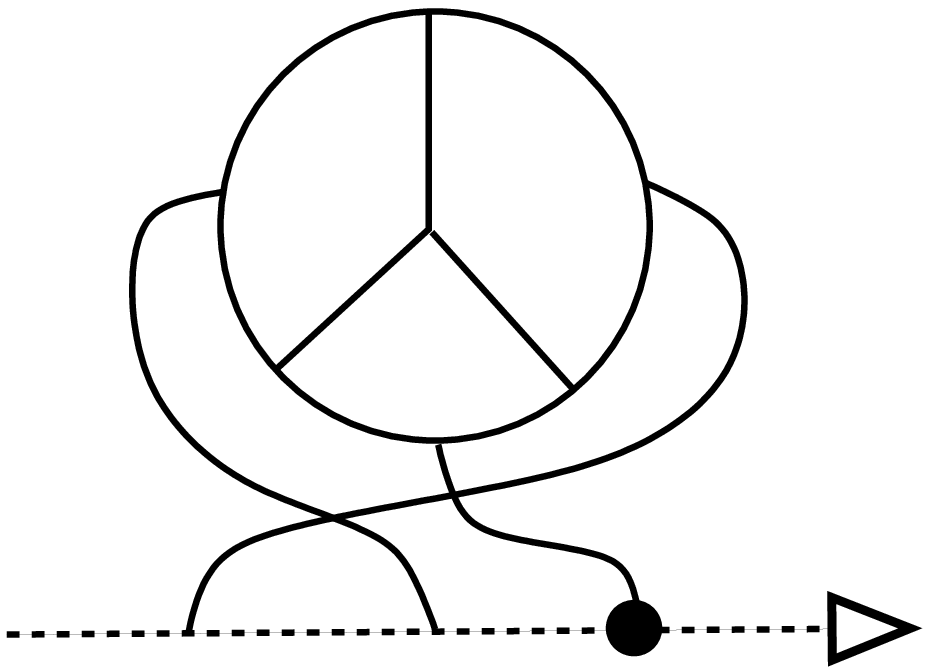}}}
+
\raisebox{-3ex}{\scalebox{0.22}{\includegraphics{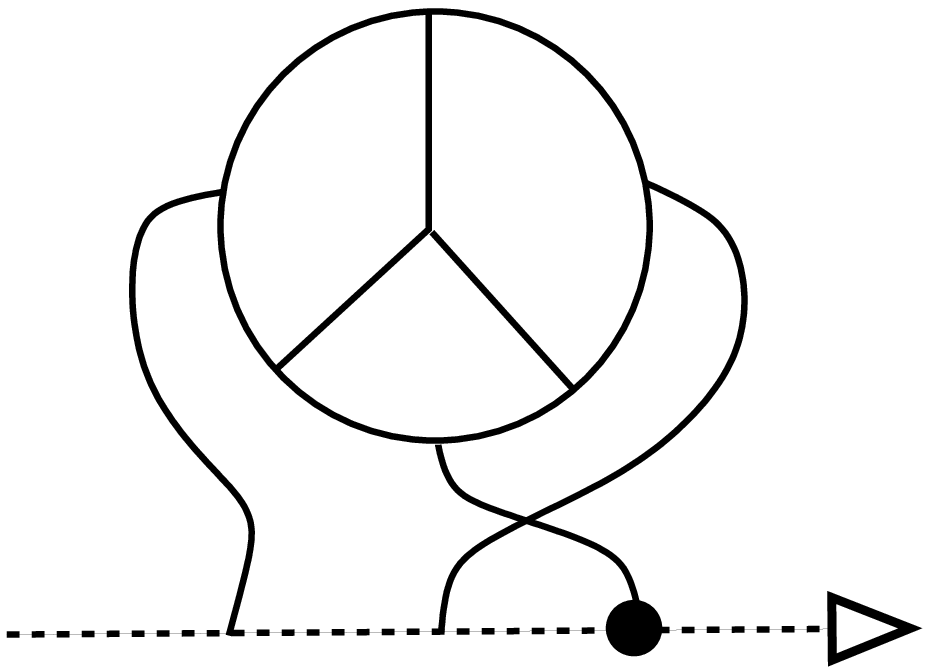}}}\right).
\end{eqnarray*}
\end{defn}

\begin{defn}
Define the {\bf forget-the-ordering} maps
\[
\tau^i : \ncw^i \rightarrow \Wspace^i\ \ \ \ \mbox{and}\ \ \ \
\tau^i_\iota : \ncw^i_\iota \rightarrow \Wspace^i_\iota
\]
by simply obtaining a commutative Weil diagram from a given
non-commutative Weil diagram by forgetting the ordering of the
legs. For example:
\[
\tau^4\left(
\raisebox{-3.5ex}{\scalebox{0.22}{\includegraphics{superaverageF}}}
\right) =
\raisebox{-3.5ex}{\scalebox{0.22}{\includegraphics{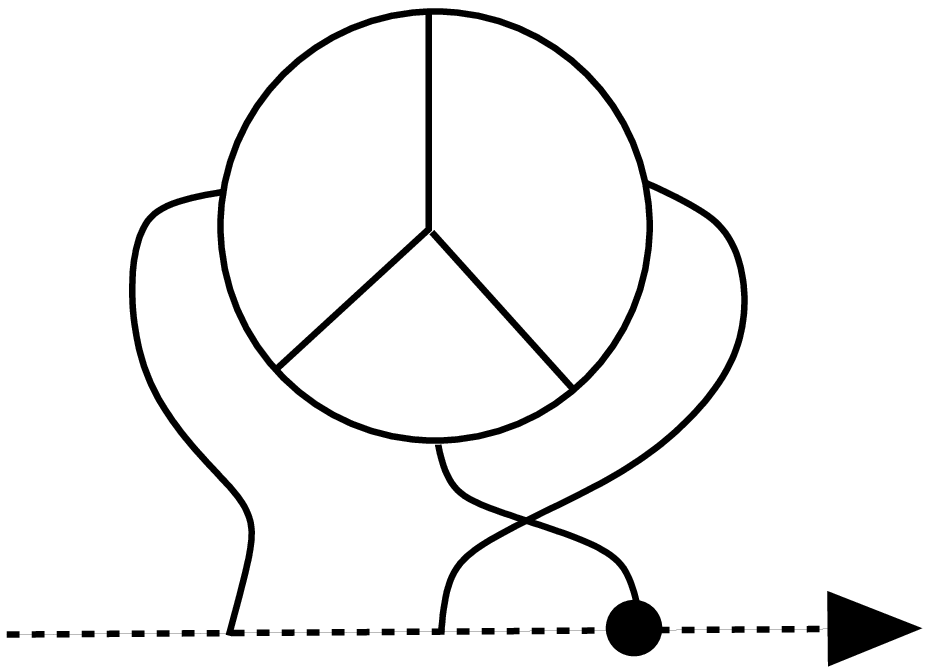}}}\
\in\ \Wspace^4.
\]
\end{defn}
Note that $\chi_\Wspace:\Wspace \rightarrow \ncw$ and $\tau: \ncw
\rightarrow \Wspace$ are both maps of $\iota$-complexes. (This is
immediate for $\tau$; for $\chi_\Wspace$ it may take a few minutes
reflection.)

The composition $\tau\circ\chi_\Wspace$ is clearly equal to the
identity $\text{id}_{\Wspace}$. Furthermore, we shall presently
declare the two maps of $\iota$-complexes $\text{id}_{\ncw}$ and
$\chi_\Wspace\circ\tau$ to be ``$\iota$-chain homotopic". This
fact is the technical heart of this work. Let us first spell out
what ``$\iota$-chain homotopic" means.

An {\bf $\iota$-chain homotopy} $s$ between two maps of
$\iota$-complexes \[ \xymatrix{
(\mathcal{K},\mathcal{K}_\iota,\iota) \ar@/^/[r]^{f} \ar@/_/[r]_g
&  (\mathcal{L},\mathcal{L}_\iota,\iota) }.
\]
is a pair of sequences of maps
\[
s^i : \mathcal{K}^i \rightarrow \mathcal{L}^{i-1}\ \ \ \text{and}\
\ \ s^i_\iota : \mathcal{K}^i_\iota \rightarrow
\mathcal{L}^{i-1}_\iota
\]
such that
\begin{enumerate}
\item{$d^{i-1}_{\mathcal{L}}\circ s^i + s^{i+1}\circ
d^i_{\mathcal{K}} = f^i - g^i$,}
\item{$d^{i-1}_{\mathcal{L},\iota}\circ s_\iota^i +
s^{i+1}_\iota\circ d_{\mathcal{K},\iota}^i = f_\iota^i -
g_\iota^i$,} \item{$s^{i-1}_\iota \circ \iota^i + \iota^{i-1}
\circ s^i = 0$.}
\end{enumerate}
It will not surprise the reader to learn that if two maps of
$\iota$-complexes are $\iota$-chain homotopic, then the induced
maps on the basic subcomplexes
\[
\xymatrix{ \mathcal{K}_{\text{basic}}
\ar@/^/[r]^{f_{\text{basic}}} \ar@/_/[r]_{g_\text{basic}} &
\mathcal{L}_{\text{basic}}}
\]
are chain homotopic.

\begin{thm}[The key technical theorem]\label{thebigtheorem}
There exists an $\iota$-chain homotopy
\[
s: \ncw\rightarrow\ncw \ \
\]
between the two maps of $\iota$-complexes $\text{id}_{\ncw}$ and
$\chi_\Wspace\circ\tau$.
\end{thm}
The construction of this $\iota$-chain homotopy is the subject of
Section \ref{homotopyconstruct}.

\section{Homological wheeling}\label{HWsection}

In this section we will state the Homological Wheeling (HW)
theorem, and point out why it is an immediate consequence of the
existence of the $\iota$-chain homotopy described in Theorem
\ref{thebigtheorem}. In Section \ref{homotopyconstruct} that chain
homotopy will be constructed, and in Section \ref{gettingwheeling}
we'll explain how the usual Wheeling theorem is recovered from HW
when certain relations are introduced.

\subsection{The statement}
First we'll make some comments concerning the ingredients of HW.
In Section \ref{basicsubcomplexintro} we showed how to take
symmetric Jacobi diagrams (i.e. generators of $\Bspace$) and
construct basic cocycles in $\Wspace$. Recall this map $\Upsilon$:
order the legs of the symmetric Jacobi diagram, label each leg
with an ${\bf F}$, then expand these legs into the usual basis.
If we continue on and compose $\Upsilon$ with the graded averaging
map $\chi_\Wspace$, then we have constructed basic cocycles in
$\ncw$ (because the graded averaging map is a map of
$\iota$-complexes). HW concerns these cocycles. Formalizing this
composition:
\begin{defn}
Let ${\bf H}^i : \mathcal{B}^i \rightarrow H^i(\ncw_\text{basic})$
denote the linear map given by the formula $
{\bf H}^i = H^i((\chi\circ\Upsilon)_\mathrm{basic})$.
\end{defn}

The other ingredient in HW is the natural graded product on the
$\iota$-complex $\ncw$. For two non-commutative diagrams $D_1$ and
$D_2$, $D_1\#D_2$ is defined by placing $D_2$ to the right of
$D_1$ on the orienting line. For example: \[
\raisebox{-3ex}{\scalebox{0.25}{\includegraphics{productA}}} \#
\raisebox{-3ex}{\scalebox{0.25}{\includegraphics{productB}}} \ \
=\ \raisebox{-3ex}{\scalebox{0.25}{\includegraphics{productC}}}\ .
\]
Observe that $d$ and $\iota$ satisfy the graded Leibniz rule with
repsect to this product, so that this product descends to a
product on the basic cohomology.

\begin{HW}
Let $v$ be an element of $\mathcal{B}^i$ and $w$ be an element of
$\mathcal{B}^j$. Then:
\[
{\bf H}^{i+j}(v\sqcup w)= {\bf H}^i(v)\#{\bf H}^j(w) \in
H^{i+j}(\ncw_{\mathrm{basic}}).
\]
\end{HW}
It is worth re-expressing this theorem in more concrete terms.
\begin{thm}[A re-expression of Homological
Wheeling]\label{reexpress} Let $v$ be an element of
$\mathcal{B}^i$ and $w$ be an element of $\mathcal{B}^j$. Then
there exists an element $x_{v,w}$ of
$\widetilde{\mathcal{W}}^{i+j-1}_\mathrm{basic}$, (that is, an
element of $\widetilde{\mathcal{W}}^{i+j-1}$ satisfying
$\iota(x_{v,w})=0$), with the property that
\[
(\chi_\Wspace\circ\Upsilon)^i(v)\#(\chi_\Wspace\circ\Upsilon)^j(w)
= (\chi_\Wspace\circ\Upsilon)^{i+j}(v\sqcup w) + d(x_{v,w})\ \in\
\ncw^{i+j}\ .
\]\\[-0.75cm]
\end{thm}

\subsection{The error term}
Here we'll make a brief comment about how the error term in the
above statement - $d(x_{v,w})$ - fits into the big picture. When
we introduce relations to recover Wheeling (as described in
Section \ref{gettingwheeling}) we'll use the knowledge that
$\iota(x_{v,w})=0$ to show that that term vanishes. This is
more-or-less the reason we have to carry information about the
action of the map $\iota$ through the calculation.

\subsection{How Theorem \ref{thebigtheorem} implies HW} Consider the
left-hand side of the equation stated in Theorem \ref{reexpress}:
$(\chi_\Wspace\circ\Upsilon)(v)\#(\chi_\Wspace\circ\Upsilon)(w)$.
If we insert this element into the equation $\text{id}_{\ncw} = \chi_\Wspace\circ\tau + s\circ d + d\circ s$,
then we learn that it is equal to
\[
\left(\chi_\Wspace\circ\tau\right)
\left((\chi_\Wspace\circ\Upsilon)(v)\#(\chi_\Wspace\circ\Upsilon)(w)\right)
+
d\left(s\left((\chi_\Wspace\circ\Upsilon)(v)\#(\chi_\Wspace\circ\Upsilon)(w)\right)\right),
\]
where we have used the fact that
$(\chi_\Wspace\circ\Upsilon)(v)\#(\chi_\Wspace\circ\Upsilon)(w)$
is a cocycle in $\ncw$.

It takes but a moment to agree that
\[
(\chi_\Wspace\circ\tau)\left(
(\chi_\Wspace\circ\Upsilon)(v)\#(\chi_\Wspace\circ\Upsilon)(w)
\right) = (\chi_\Wspace\circ\Upsilon)(v\sqcup w)\ ,
\]
and, setting
\[
x_{v,w}=s\left((\chi_\Wspace\circ\Upsilon)(v)\#(\chi_\Wspace\circ\Upsilon)(w)\right),
\]
we have the required right-hand side (noting that
$\iota(x_{v,w})=0$ because $s$ commutes with $\iota$ and
$(\chi_\Wspace\circ\Upsilon)(v)\#(\chi_\Wspace\circ\Upsilon)(w)$
is a basic element of $\ncw$).
\begin{flushright}
$\Box$
\end{flushright}

\section{The construction of the $\iota$-chain homotopy \(s\).}
\label{homotopyconstruct} The construction of the $\iota$-chain
homotopy $s$ and the verification of its properties is a 100\%
combinatorial exercize. The combinatorial work takes place in two
$\iota$-complexes $\mathcal{T}$ and $\mathcal{T}_{\mathrm{dR}}$.
The $\iota$-complex $\mathcal{T}$ is an intermediary between
$\Wspace$ and $\ncw$, with some legs behaving as if they are in
$\Wspace$ and some legs behaving as if they are in $\ncw$. The
complex $\mathcal{T}_{\mathrm{dR}}$ may be viewed as a formal
version of the result of tensoring $\mathcal{T}$ with the de Rham
complex over the interval, which will allow us to mimic a Stokes
theorem argument combinatorially.
\subsection{The $\iota$-complex $\mathcal{T}$.}
The $\iota$-complex $\mathcal{T}$ is based on diagrams of the
following sort, which will be referred to as
$\mathcal{T}$-diagrams: \[
\raisebox{-3ex}{\scalebox{0.24}{\includegraphics{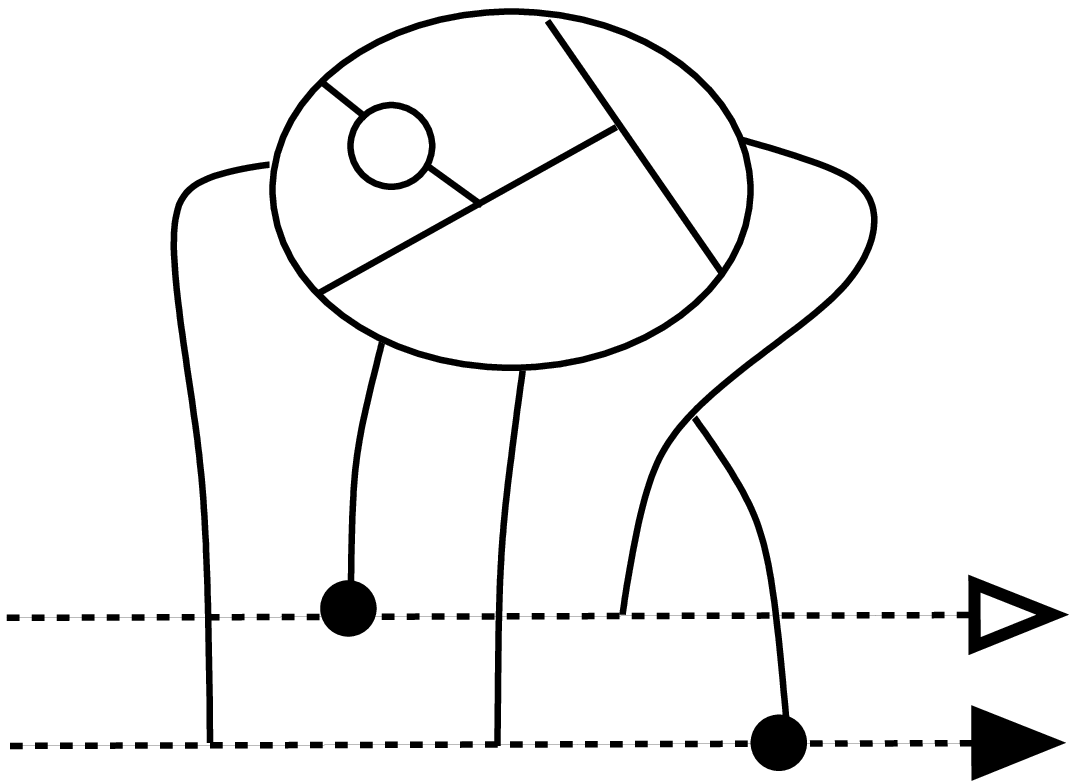}}}\ \ .
\]

In a $\mathcal{T}$-diagram two legs never lie on the same vertical
line.
Thus the set of legs of a $\mathcal{T}$-diagram is ordered. The
rules for permuting the order of legs depend on which orienting
lines they lie on. If they lie on different orienting lines, or
both lie on the bottom (the {\it commutative}) orienting line,
then they can be permuted, up to sign, in the usual way. If,
however, they both lie on the top (the {\it non-commutative})
line, then they cannot be permuted. Thus, in $\mathcal{T}^7$:
\[
\raisebox{-3ex}{\scalebox{0.22}{\includegraphics{TexampA}}} \ =\ \
-\ \raisebox{-3ex}{\scalebox{0.22}{\includegraphics{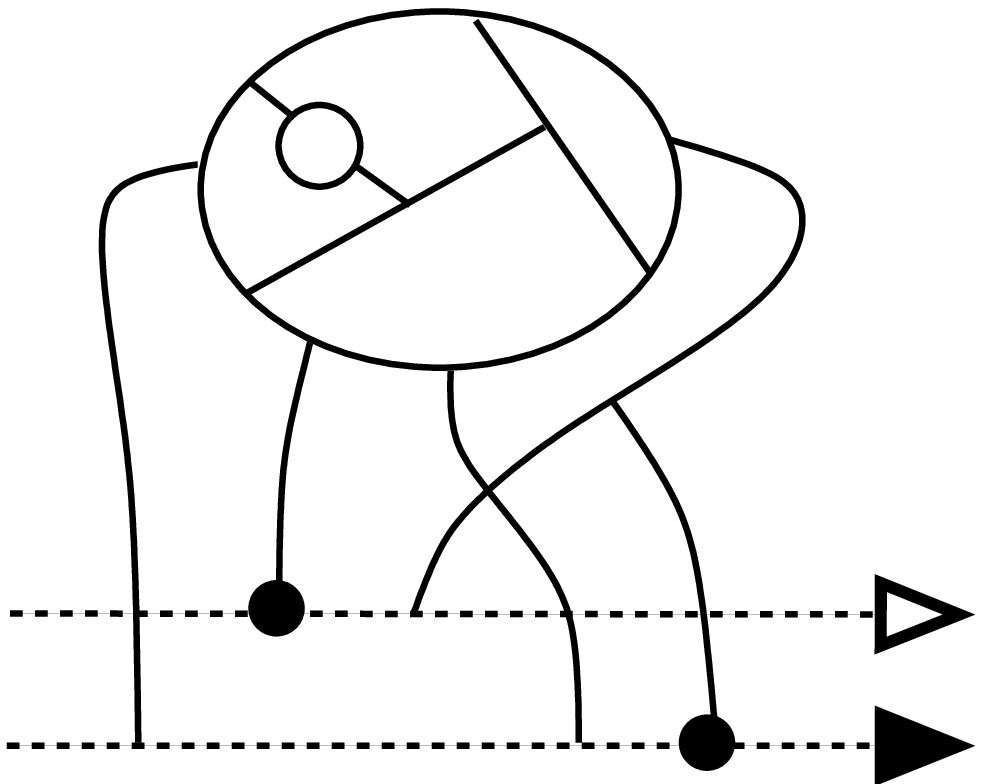}}} \
=\ \ -\
\raisebox{-3ex}{\scalebox{0.22}{\includegraphics{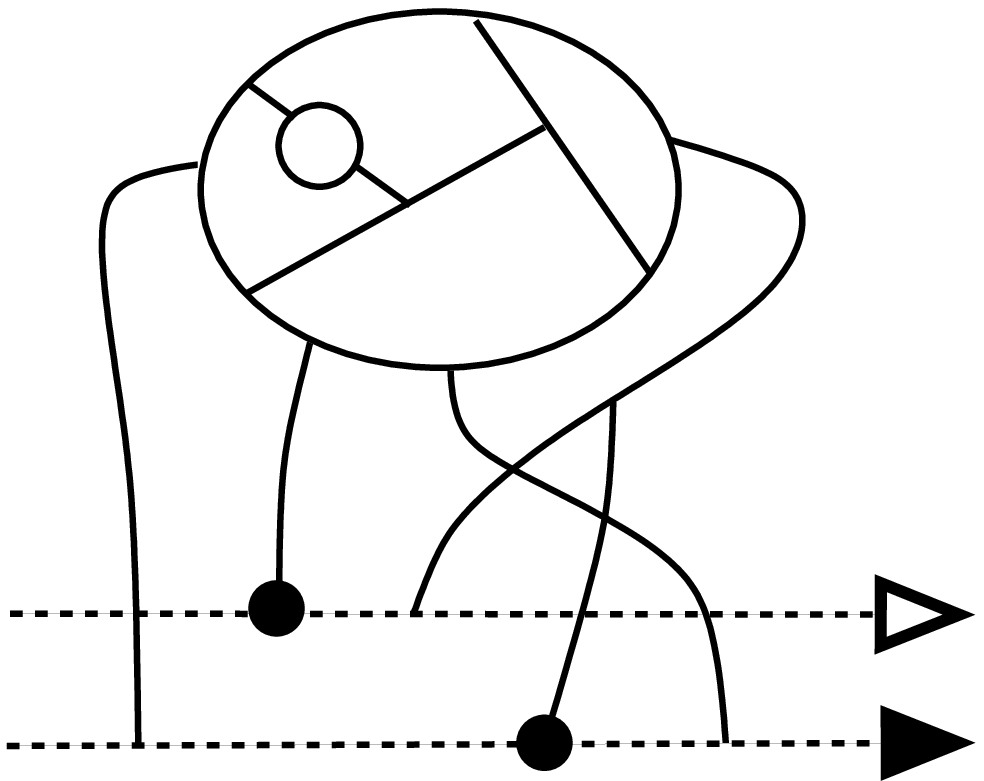}}} \
\neq\ \ -\
\raisebox{-3ex}{\scalebox{0.22}{\includegraphics{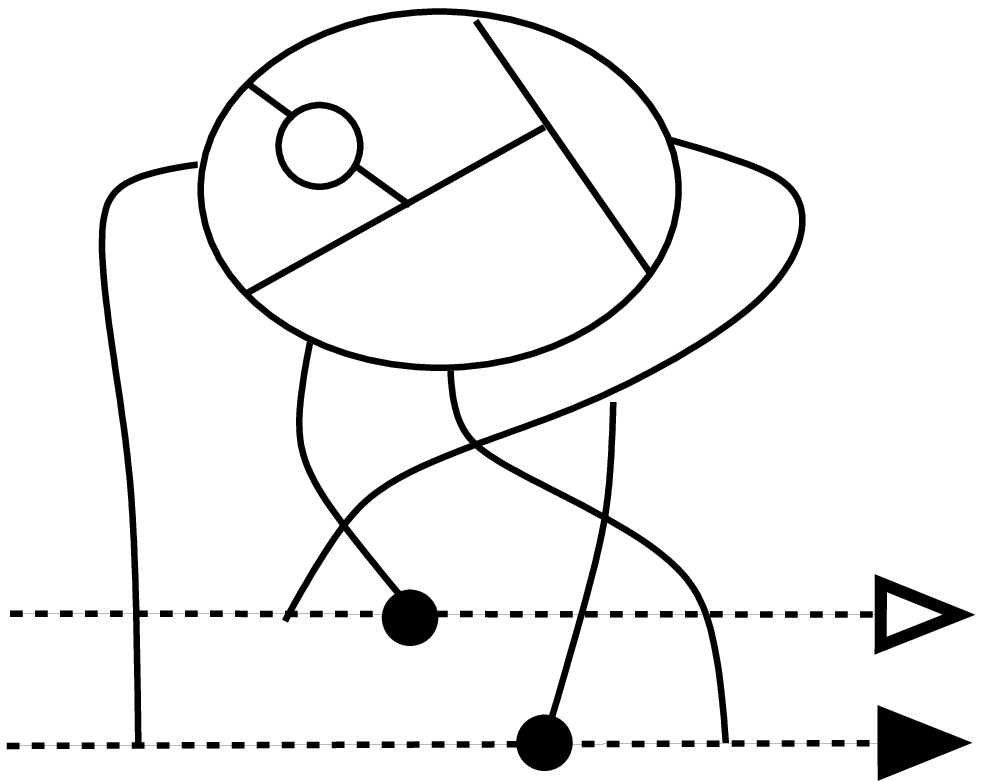}}}
\]

\begin{figure}
\caption{An example of the calculation of $d_\mathcal{T}$. \label{dTexamp}}
\parbox{12cm}{
\begin{eqnarray*}
\raisebox{-5.5ex}{\scalebox{0.22}{\includegraphics{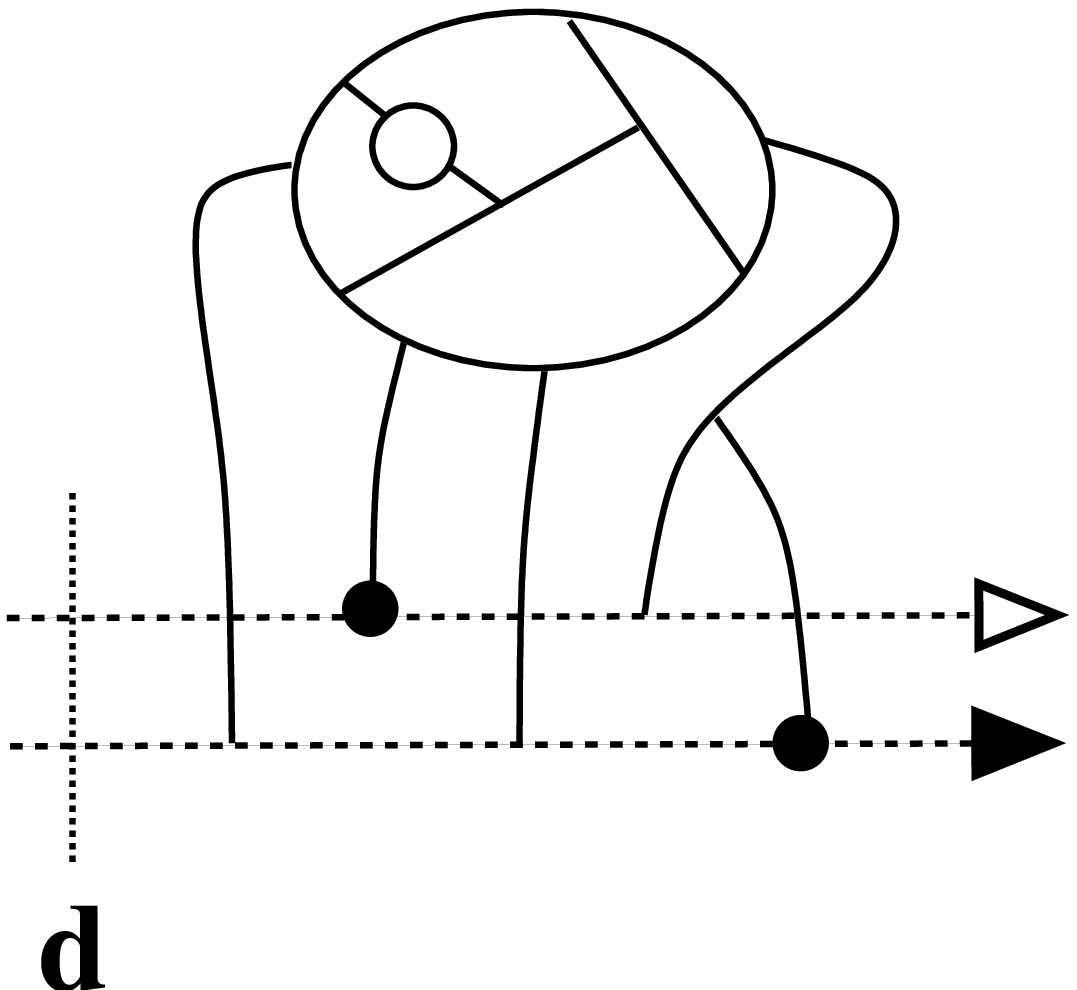}}} &
\leadsto & \raisebox{-3ex}{\scalebox{0.22}{\includegraphics{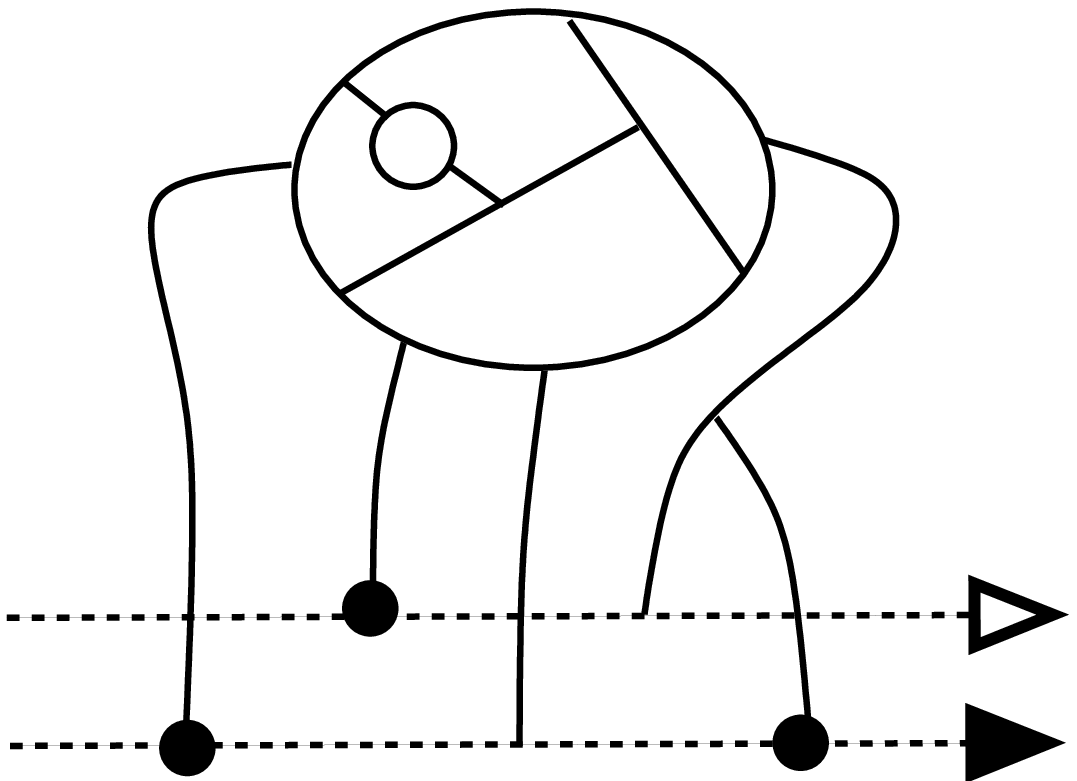}}}
\ -\ \raisebox{-5.75ex}{\scalebox{0.22}{\includegraphics{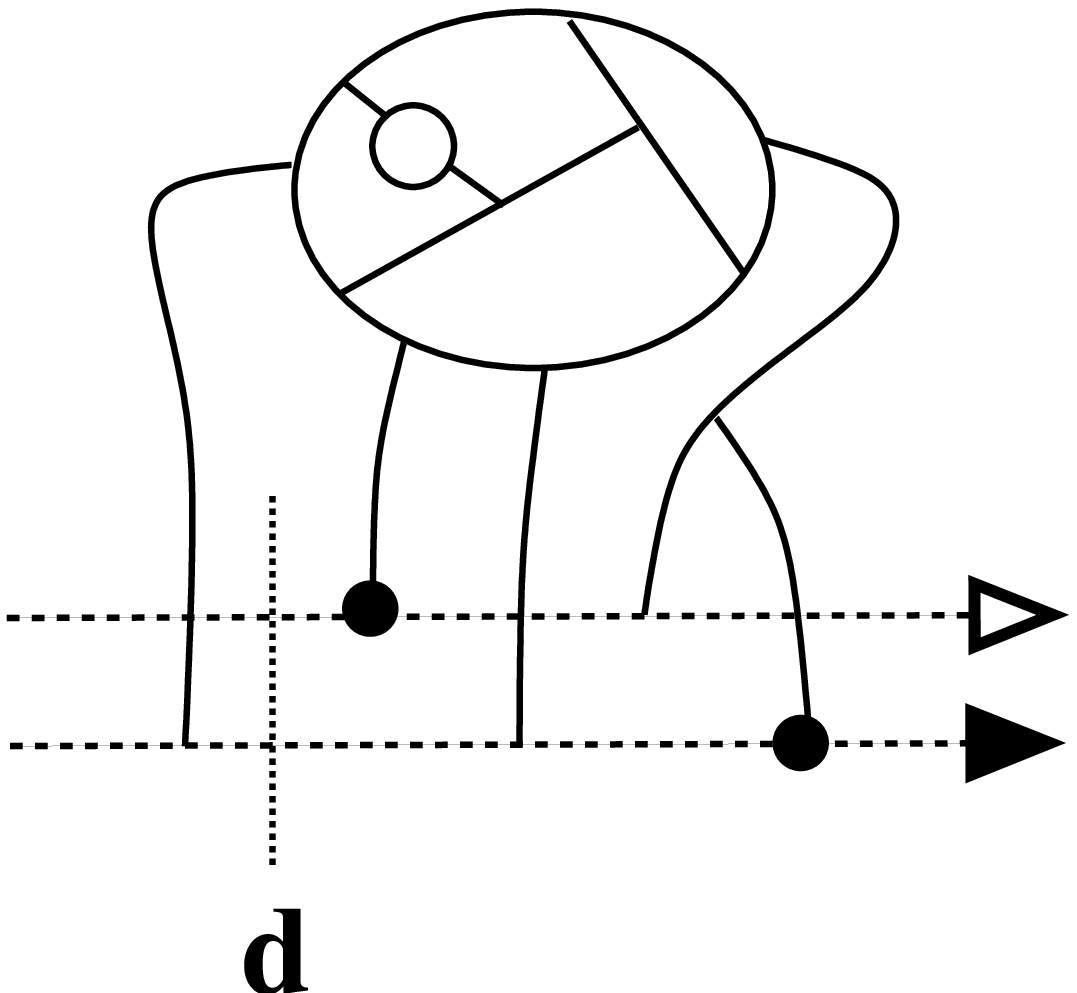}}}
\\[0.2cm]
& \leadsto &
\raisebox{-3ex}{\scalebox{0.22}{\includegraphics{TDB}}} \ -\
\raisebox{-5.75ex}{\scalebox{0.22}{\includegraphics{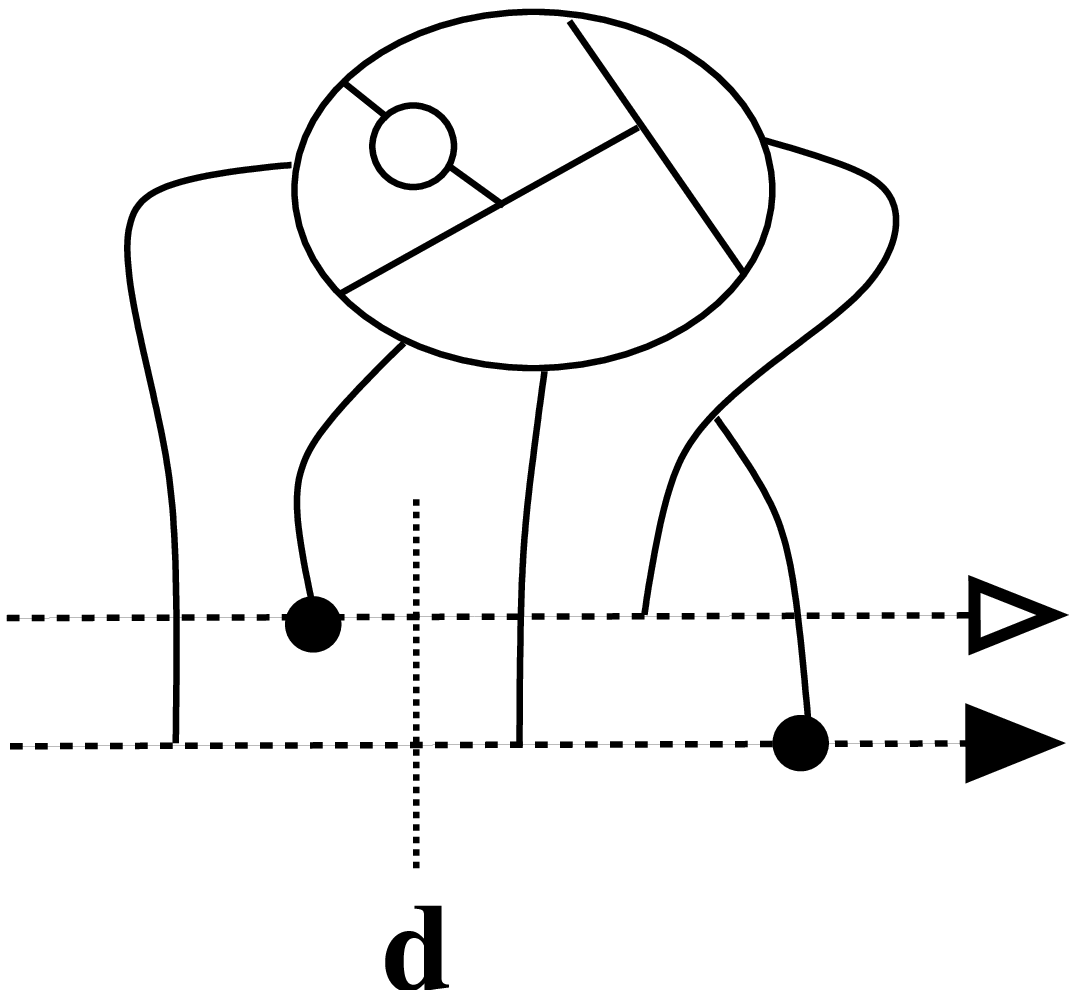}}}
\\[0.2cm]
& \leadsto &
\raisebox{-3ex}{\scalebox{0.22}{\includegraphics{TDB}}} \ -\
\raisebox{-3ex}{\scalebox{0.22}{\includegraphics{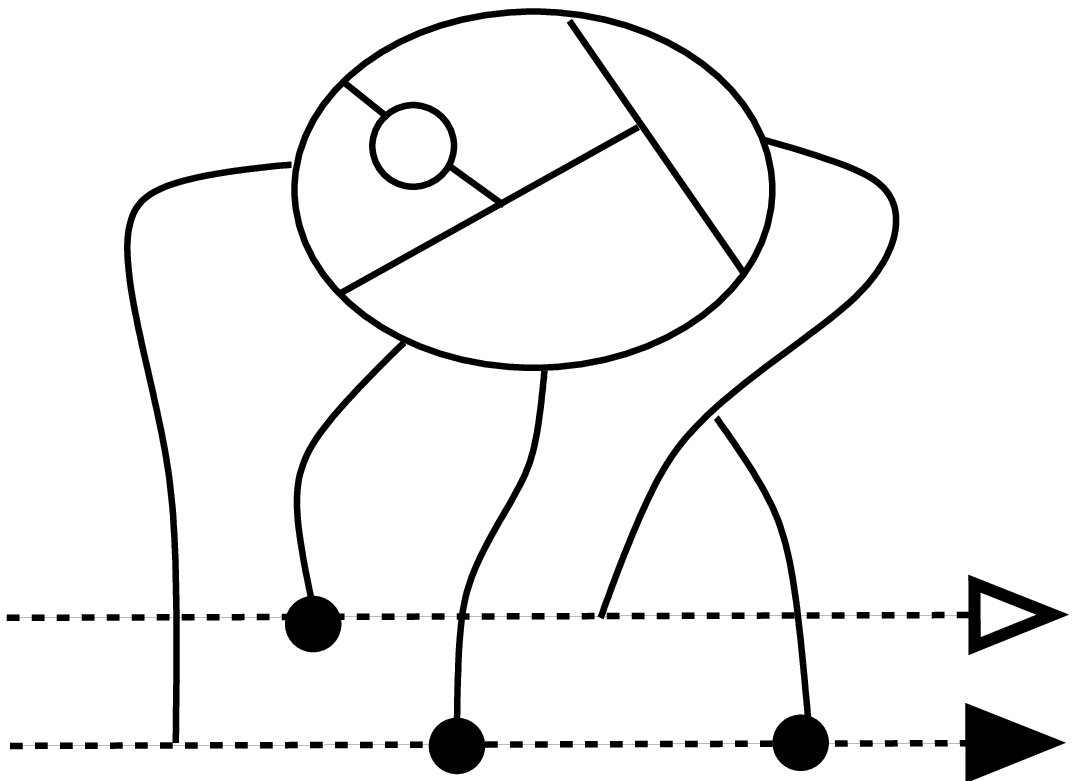}}} \ +\
\raisebox{-5.75ex}{\scalebox{0.22}{\includegraphics{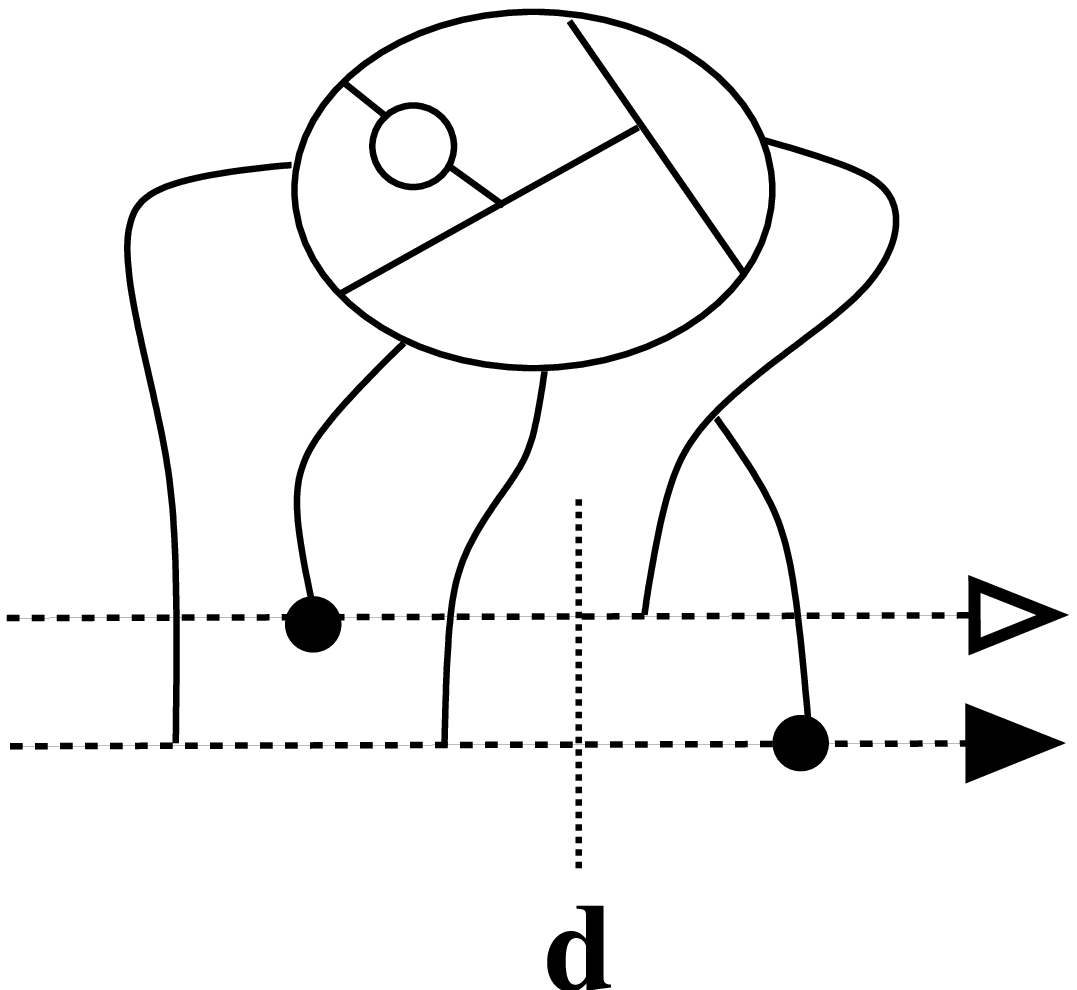}}}\ .
\end{eqnarray*}
Thus:\\[-0.5ex]
\[
d^7_\mathcal{T}\left(\raisebox{-3ex}[0.65\height]{\scalebox{0.22}{\includegraphics{TexampA}}}\right)
= \raisebox{-3ex}{\scalebox{0.22}{\includegraphics{TDB}}} \ -\
\raisebox{-3ex}{\scalebox{0.22}{\includegraphics{TDE}}} \ +\
\raisebox{-3ex}{\scalebox{0.22}{\includegraphics{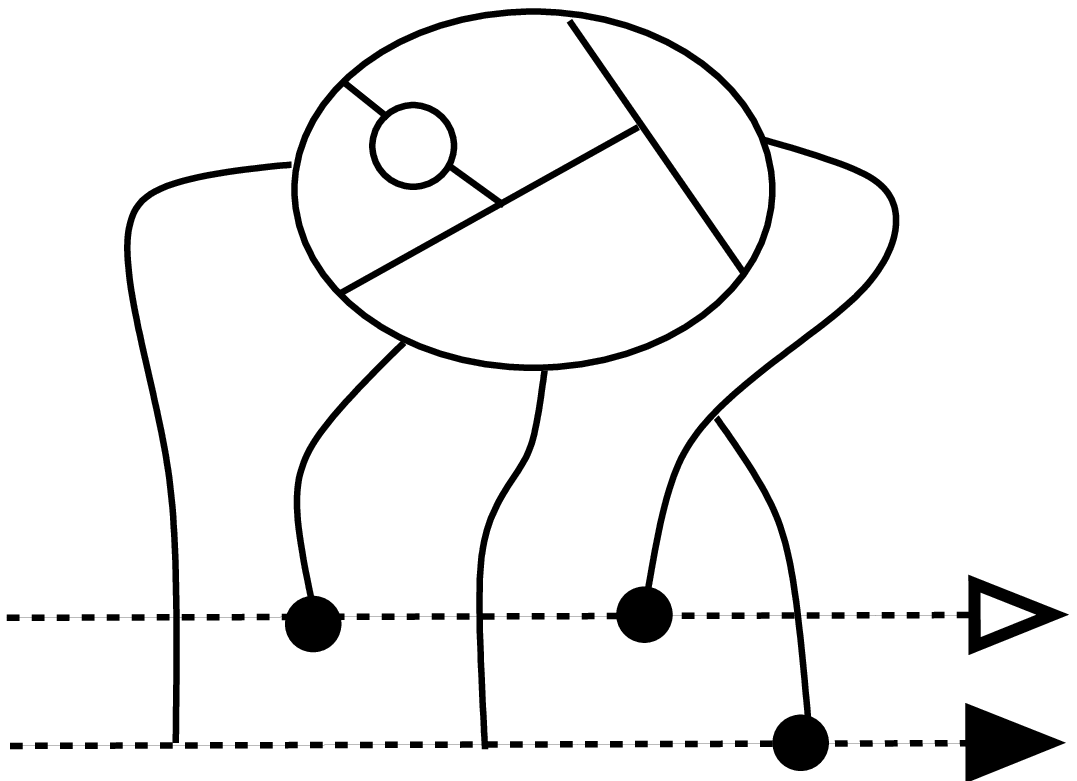}}}\ .
\]
}
\underline{\hspace{7cm}}
\end{figure}

We'll apply formal linear differential operators to these spaces.
To apply a differential operator one thinks of a vertical line
sweeping along the pair of orienting lines. Figure \ref{dTexamp} shows an example of the calculation of this differential.
The map $\iota$ is defined similarly, and with these definitions
$(\mathcal{T},\mathcal{T}_\iota,\iota_\mathcal{T})$ forms an
$\iota$-complex.

Constructions in $\mathcal{T}$ can be translated to $\ncw$ by
means of a certain map of $\iota$-complexes $
\omega : (\mathcal{T},\mathcal{T}_\iota,\iota_{\mathcal{T}}) \to
(\ncw,\ncw_\iota,\iota_{\ncw})$.
\begin{defn}\label{newomega}
The value of the map $\omega$ on some $\mathcal{T}$-diagram is
constructed in two steps.
\begin{enumerate}
\item{ Permute (with the appropriate signs) the legs of the
diagram so that all the legs lying on the commutative orienting
line lie to the right of all the legs lying on the non-commutative
orienting line. For example:
\[
\raisebox{-3ex}{\scalebox{0.22}{\includegraphics{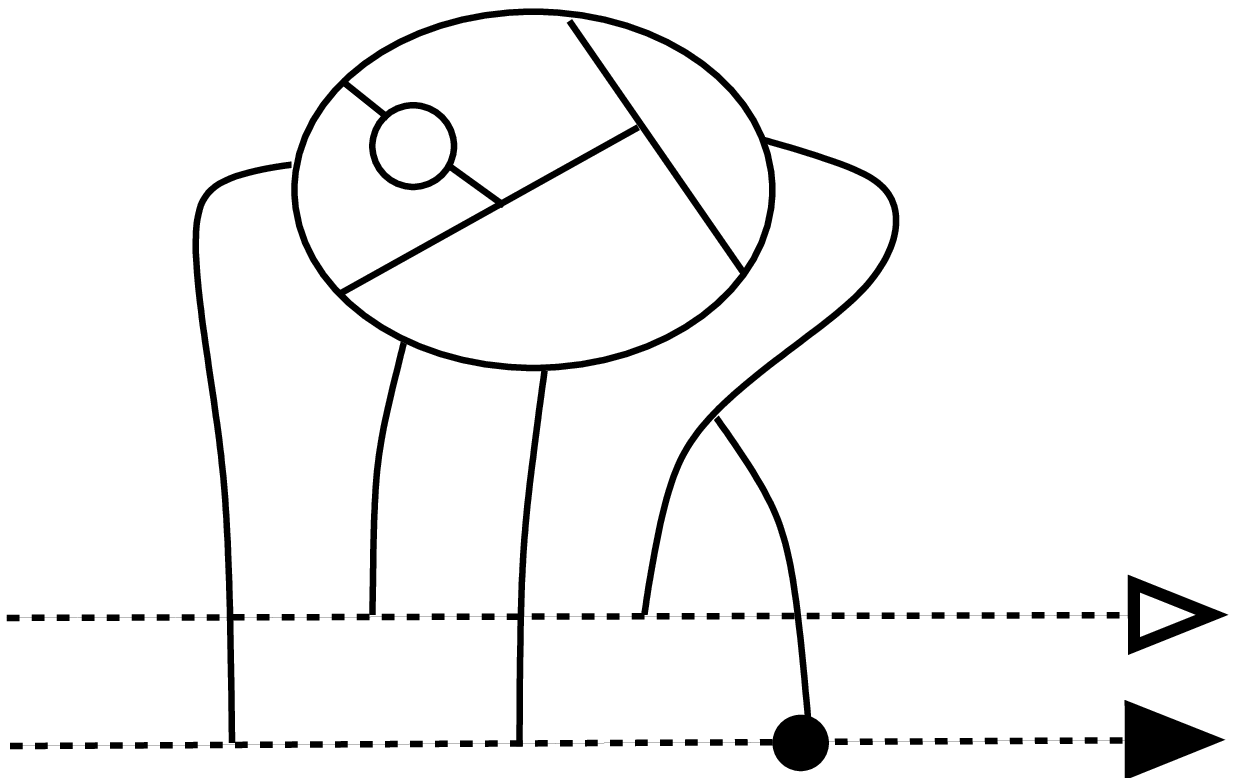}}}\ \
=\ \ -\
\raisebox{-3ex}{\scalebox{0.22}{\includegraphics{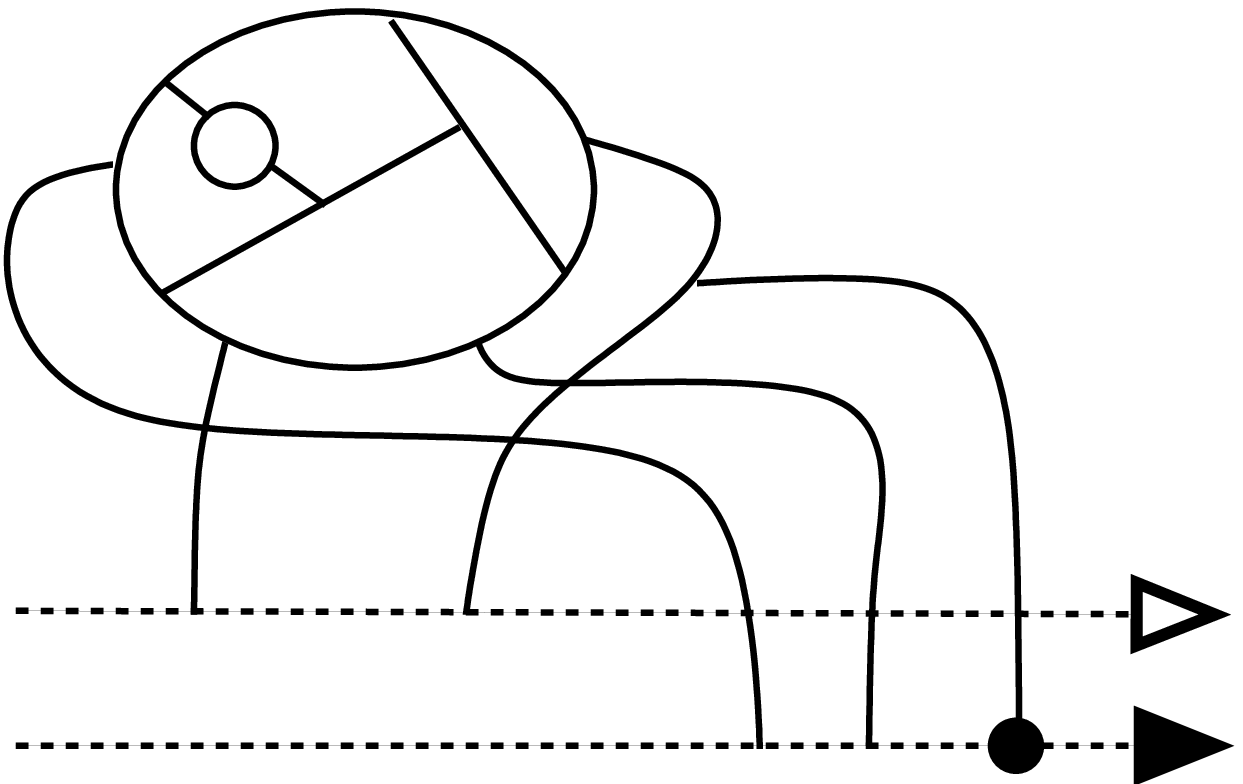}}}\ .
\]}
\item{ Now take the average of all (signed) permutations of the
legs on the commutative orienting line and adjoin the result to
the right-hand end of the non-commutative orienting line. The example is concluded in Figure \ref{omegaexamp}.}
\end{enumerate}
\begin{figure}
\caption{The second step in the definition of $\omega$. \label{omegaexamp}}
\parbox{12cm}{
\begin{eqnarray*}
\lefteqn{ -
\raisebox{-3ex}{\scalebox{0.22}{\includegraphics{omegadefnB}}}\ \
\leadsto\ \
-\raisebox{-3ex}{\scalebox{0.22}{\includegraphics{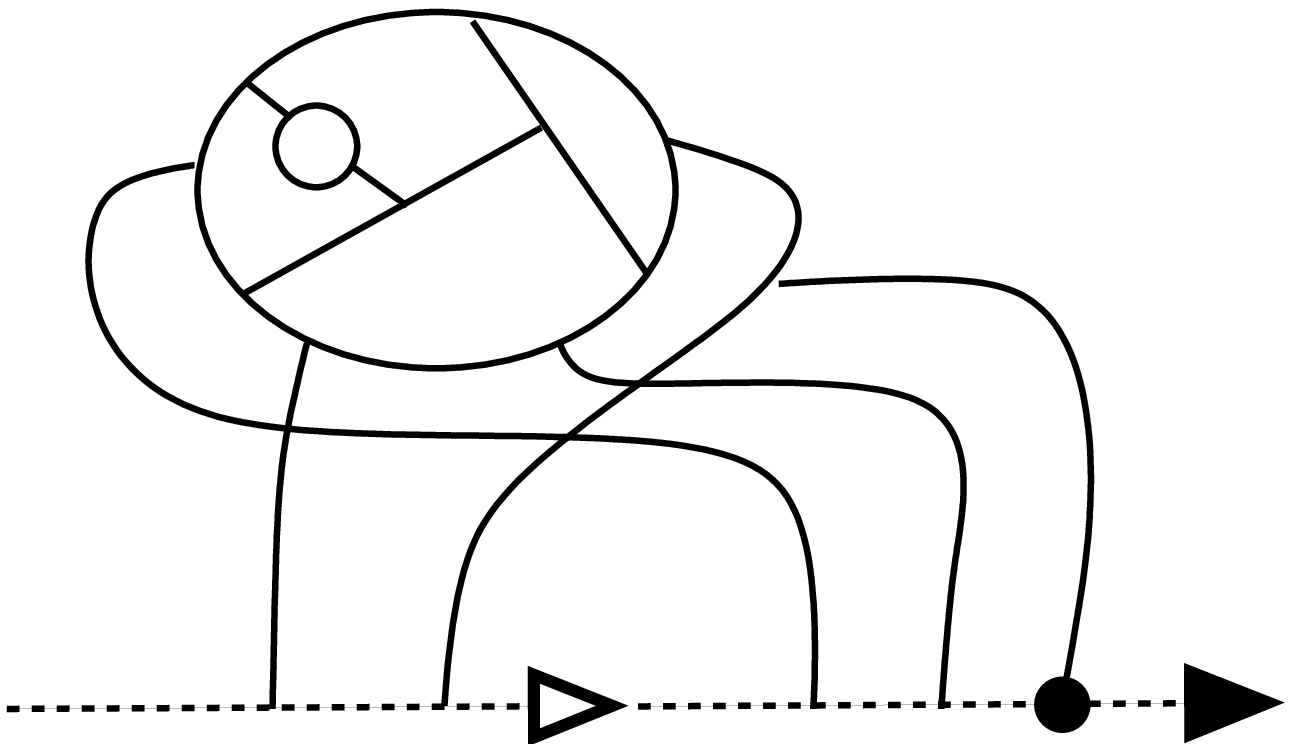}}}}
\\[0.15cm]
& \leadsto &
-\frac{1}{6}\raisebox{-3ex}{\scalebox{0.22}{\includegraphics{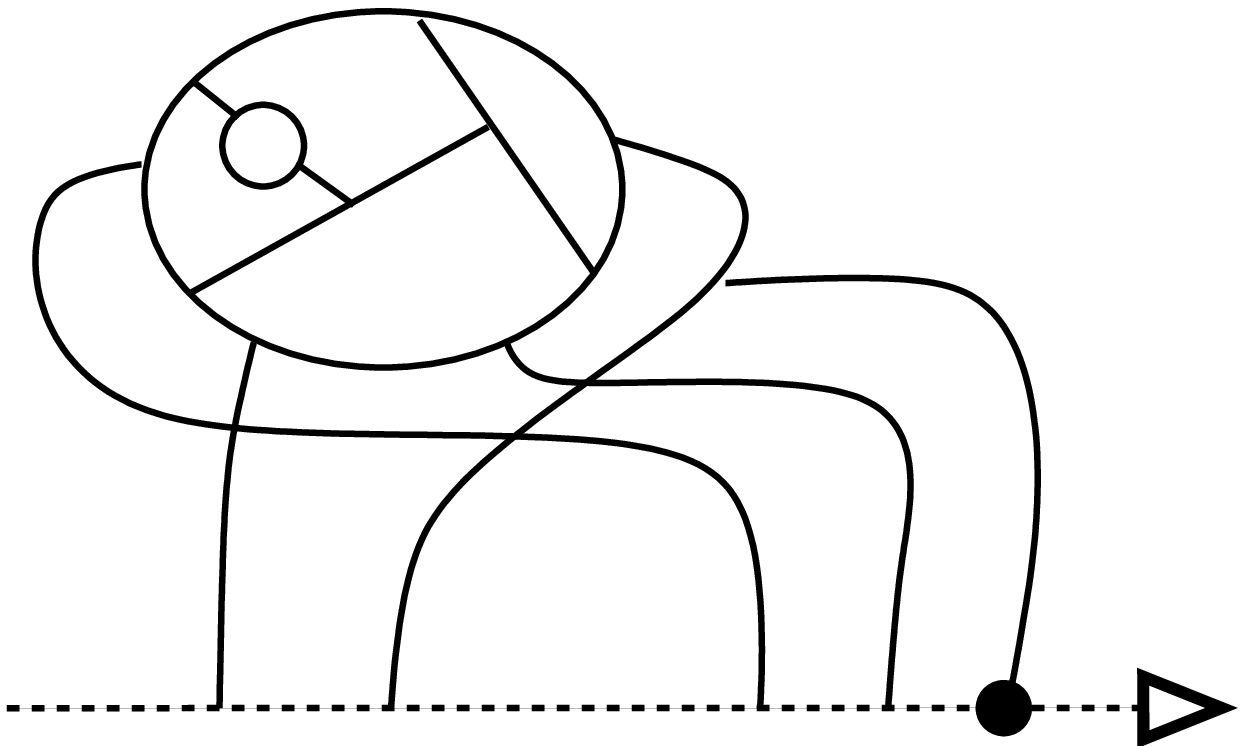}}}
-\frac{1}{6}\raisebox{-3ex}{\scalebox{0.22}{\includegraphics{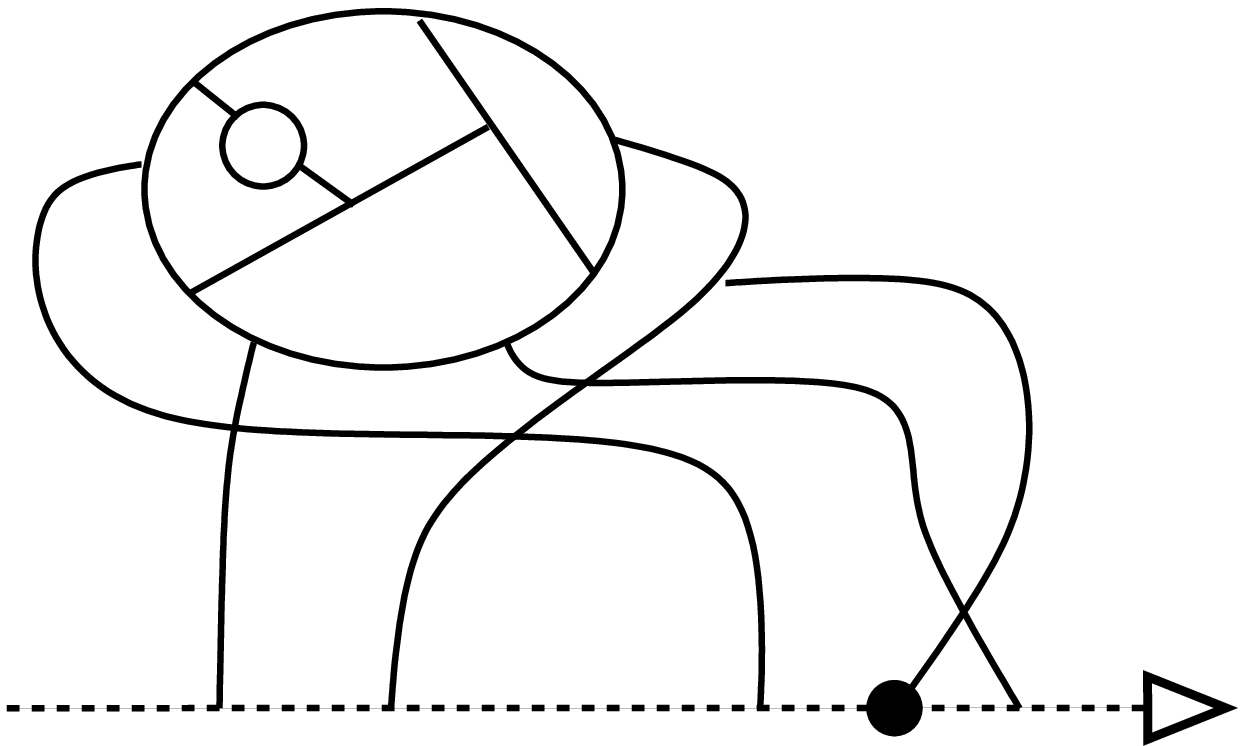}}}
-\frac{1}{6}\raisebox{-3ex}{\scalebox{0.22}{\includegraphics{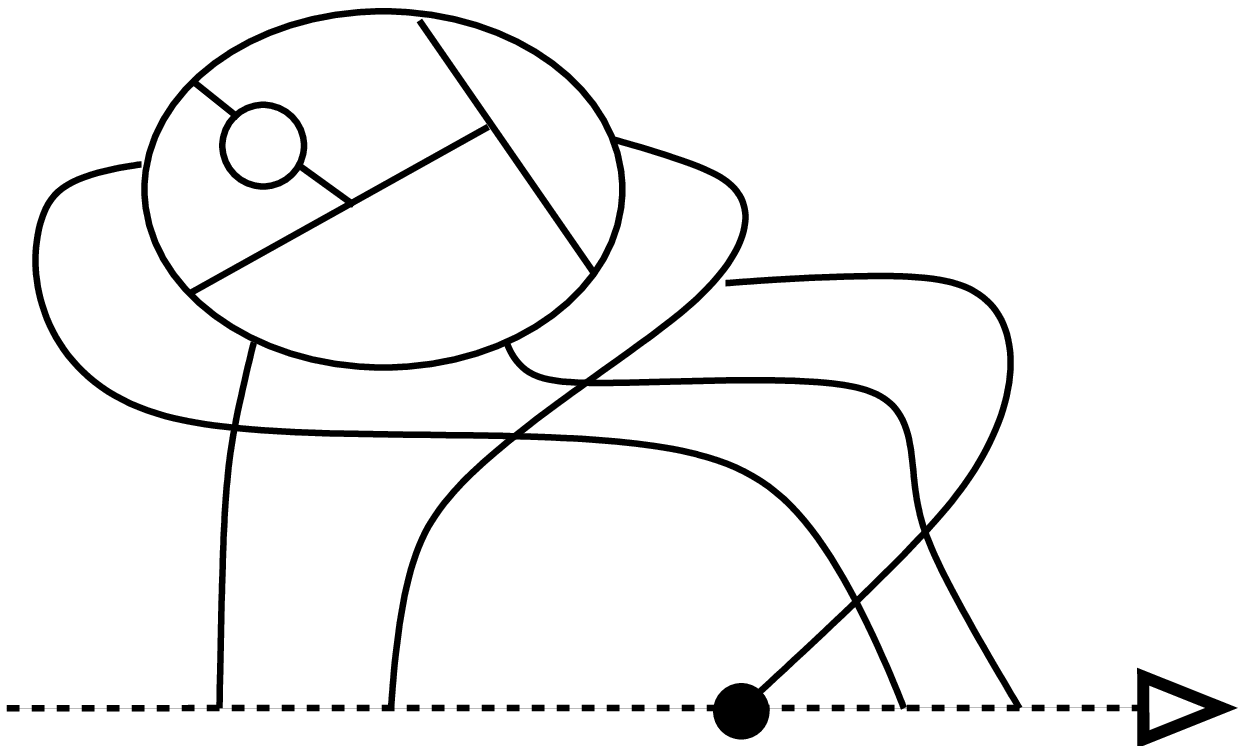}}}\\[0.15cm]
& &
+\frac{1}{6}\raisebox{-3ex}{\scalebox{0.22}{\includegraphics{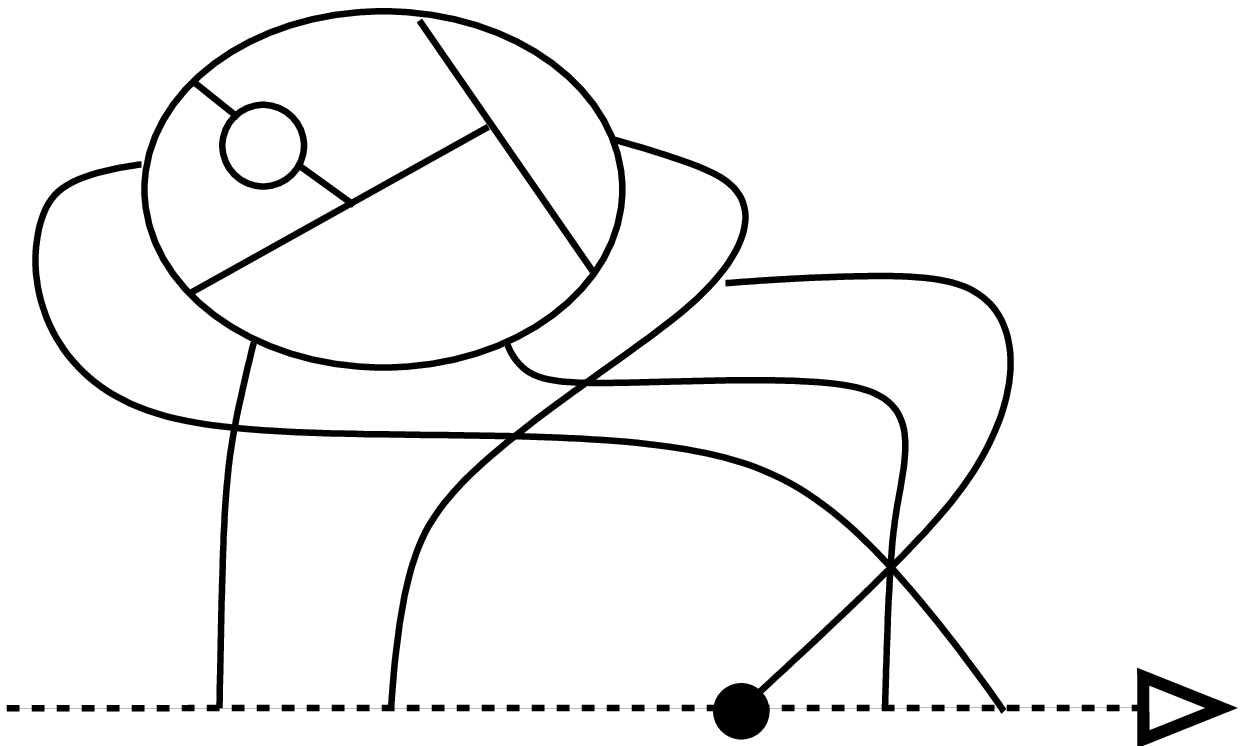}}}
+\frac{1}{6}\raisebox{-3ex}{\scalebox{0.22}{\includegraphics{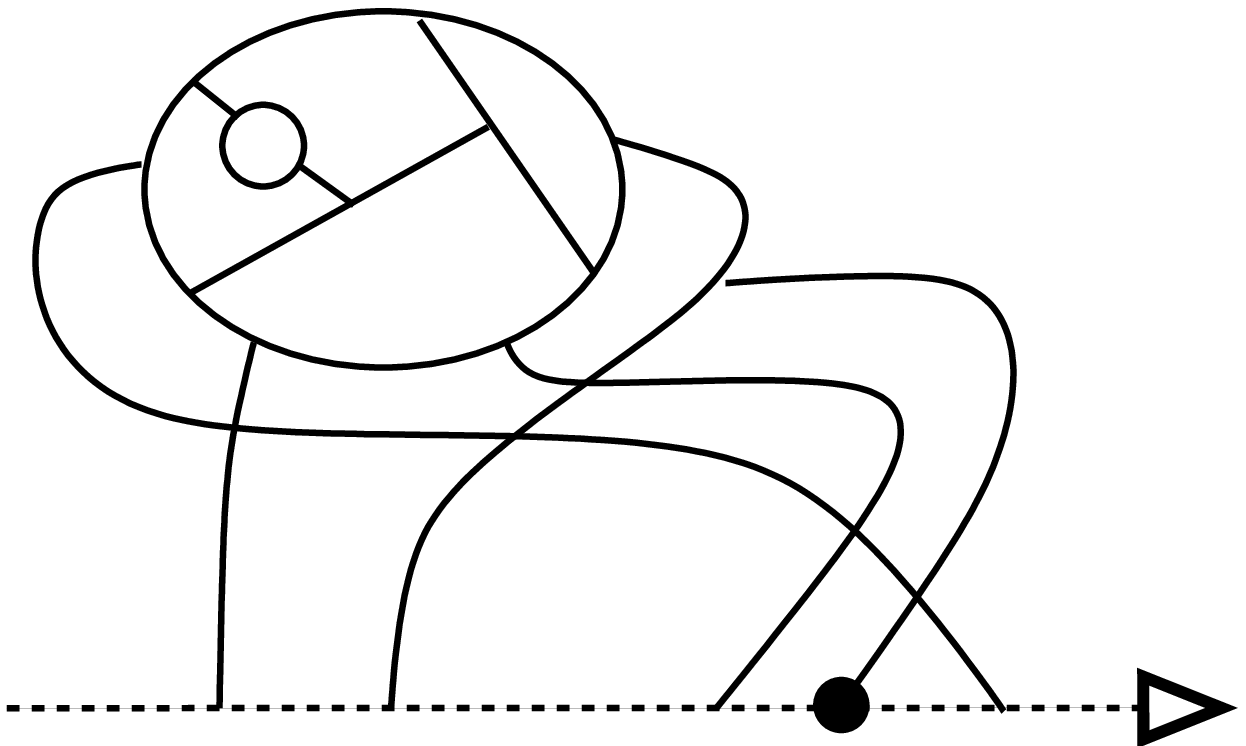}}}
+\frac{1}{6}\raisebox{-3ex}{\scalebox{0.22}{\includegraphics{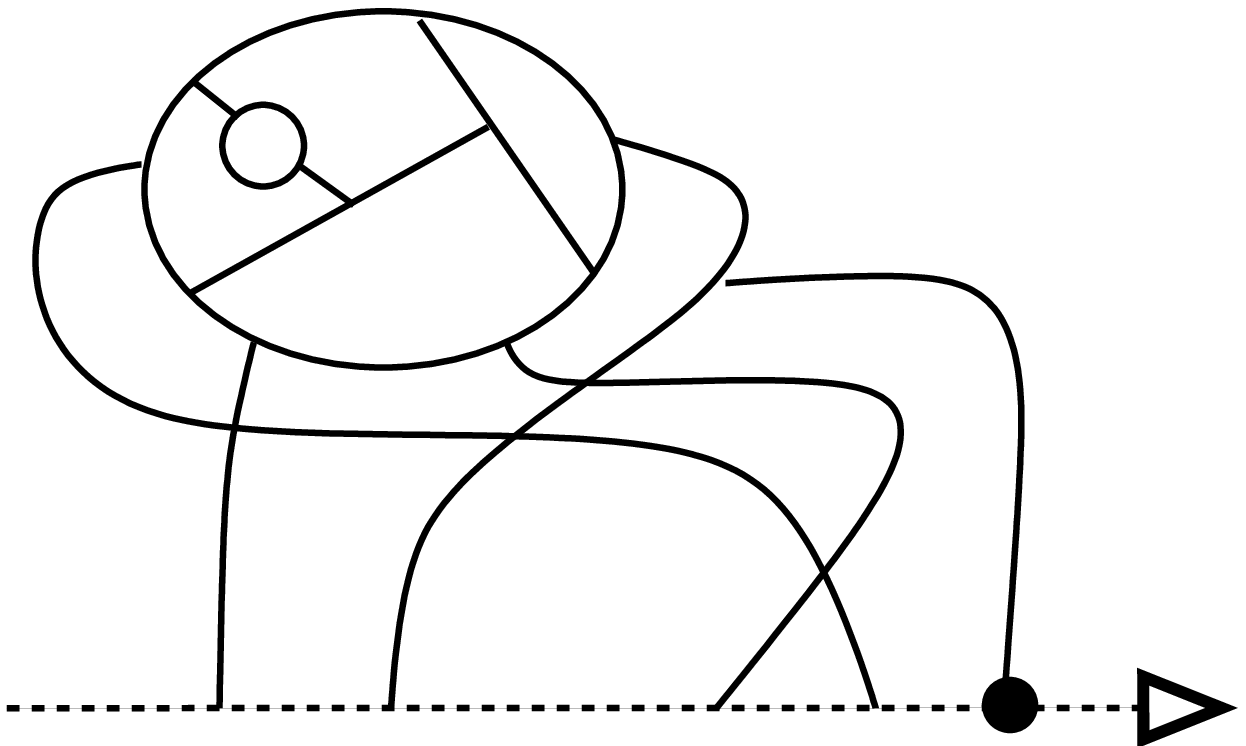}}}\ .\\[0.15cm]
\end{eqnarray*}}
\underline{\hspace{7cm}}
\end{figure}

\end{defn}
There are two natural maps of $\iota$-complexes from the
$\iota$-complex $\widetilde{\mathcal{W}}$ to the $\iota$-complex
$\mathcal{T}$. They are $i_n:\widetilde{W}\rightarrow\mathcal{T}$,
the map which puts all legs of the original diagram on the top
({\bf n}on-commutative) line, as in:
\[ i_n^6\left(
\raisebox{-3ex}{\scalebox{0.25}{\includegraphics{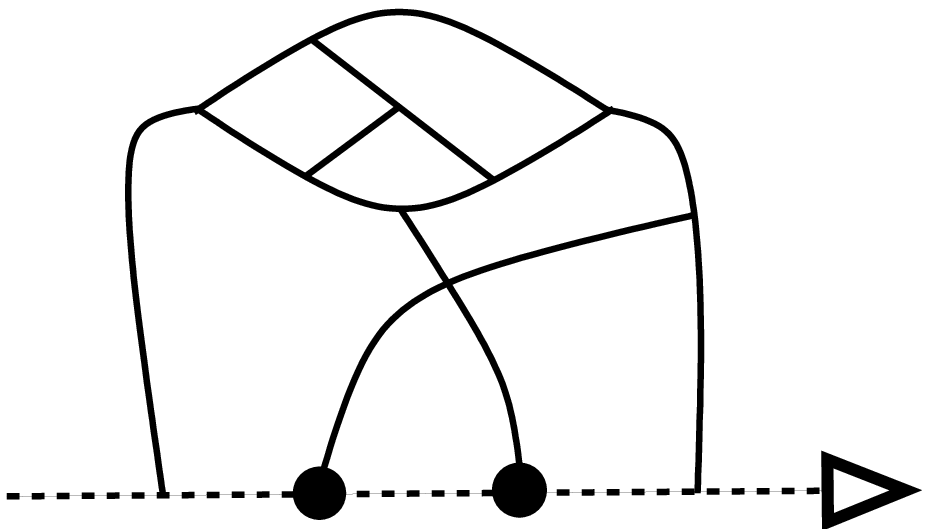}}}\right)
\ =\
\raisebox{-3ex}{\scalebox{0.25}{\includegraphics{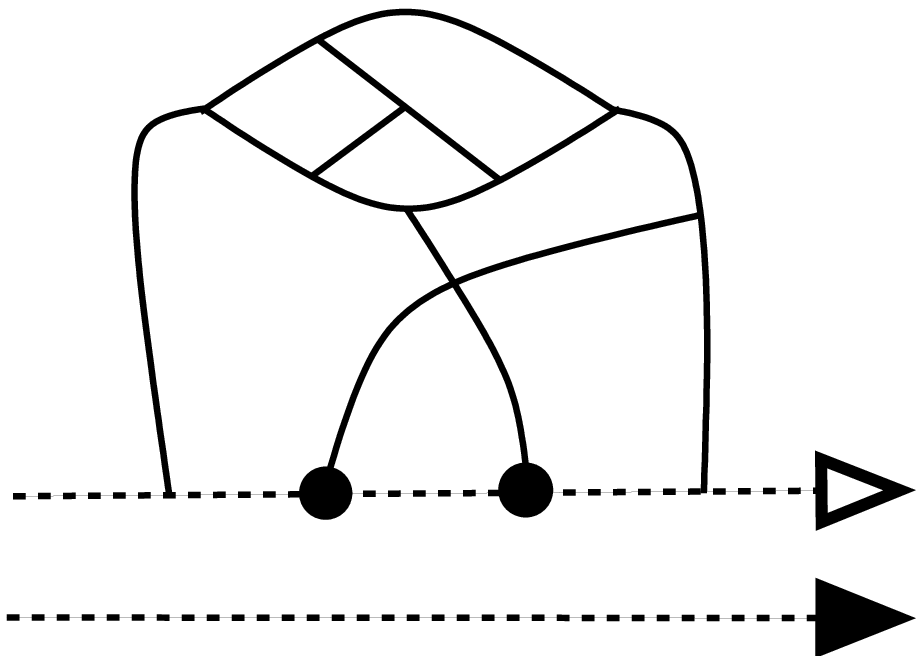}}}\ \
,
\]
and $i_c:\widetilde{W}\rightarrow\mathcal{T}$, the corresponding map which puts
all legs of the original diagram on the bottom ({\bf c}ommutative)
line.
Observe that
\[
\chi_\Wspace\circ\tau = \omega\circ i_c\ \ \ \mbox{and that}\ \ \
\mathrm{id}_{\ncw} = \omega\circ i_n.
\]

\begin{thm}\label{keytheorem}
There exists an $\iota$-chain homotopy $s_\mathcal{T}:\widetilde{\mathcal{W}}\rightarrow \mathcal{T}$
between the two maps $i_n$ and $i_c$.
\end{thm}

The construction of this $\iota$-chain homotopy is the subject of
the next two sub-sections. The key theorem, Theorem
\ref{thebigtheorem}, is a quick consequence of it. The required
$\iota$-homotopy is defined by the formula $
s = \omega\circ s_\mathcal{T}$.

{\it Proof of Theorem \ref{thebigtheorem}.} We must check that
this map satisfies the properties of an $\iota$-chain homotopy
between the maps $\chi_\Wspace\circ\tau$ and $\text{id}_{\ncw}$.
To begin, note that $s$ commutes with $\iota$. This is because
$s_\mathcal{T}$ commutes with $\iota$ (as asserted by Theorem
\ref{keytheorem}) and because $\omega$ is a map of
$\iota$-complexes. We must also check:
\begin{eqnarray*}
\mathrm{id}_{\ncw} - \chi_\Wspace\circ\tau & = & \omega \circ
\left( i_n - i_c
\right) \\
& = & \omega \circ \left( d\circ s_\mathcal{T} +
s_\mathcal{T}\circ d \right), \\
& = & d\circ\omega\circ s_\mathcal{T} + \omega\circ
s_\mathcal{T}\circ d, \\
& = & d \circ s + s\circ d.
\end{eqnarray*}
To obtain the penultimate line above we used the fact that
$\omega$ is a map of $\iota$-complexes.
\begin{flushright}
$\Box$
\end{flushright}

\subsection{The $\iota$-complex $\mathcal{T}_\mathrm{dR}$}

It remains for us to prove Theorem \ref{keytheorem}: to construct
an $\iota$-chain homotopy $s_\mathcal{T}$ between the map $i_n$
which puts all legs on the non-commutative line, and the map $i_c$
which puts all legs on the commutative line. This construction
will employ a subsidiary $\iota$-complex:
$\mathcal{T}_\mathrm{dR}$.

The $\iota$-complex $\mathcal{T}_\mathrm{dR}$ is based on diagrams
which have a third orienting line. No legs lie on this third line;
instead, this third line can be labelled by hollow discs, which
are of grade 0, and filled-in discs, which are of grade 1. These
circles can be moved about with the permutation rules appropriate
to that grading. For example, in $\mathcal{T}_{\mathrm{dR}}^7$:
\[
\raisebox{-3ex}{\scalebox{0.23}{\includegraphics{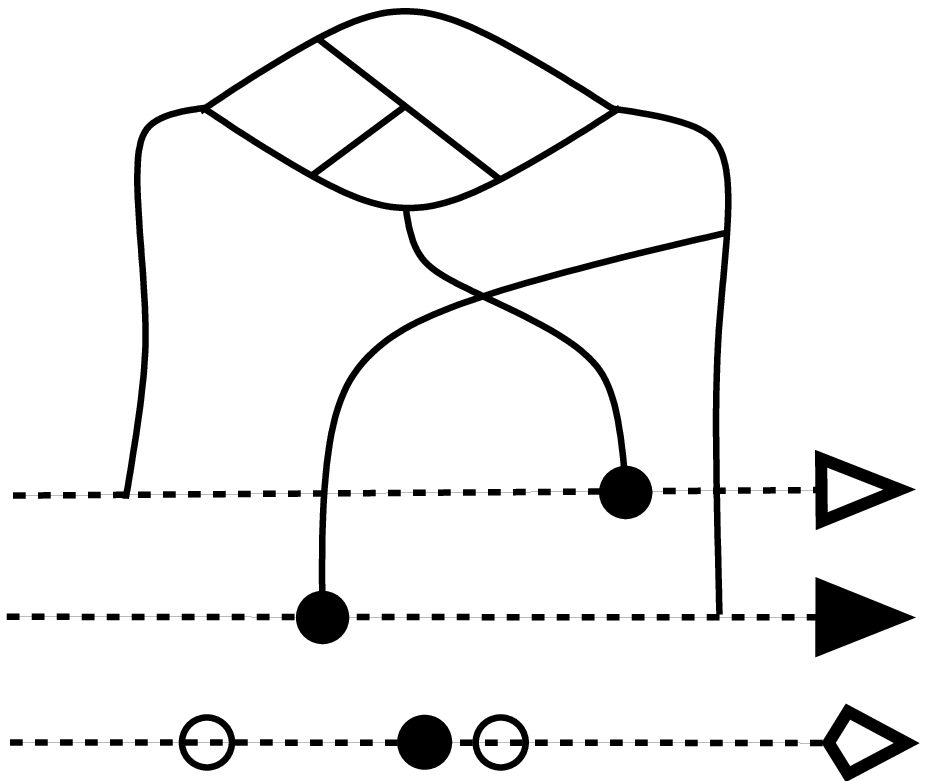}}}\
\ =\ \
\raisebox{-3ex}{\scalebox{0.23}{\includegraphics{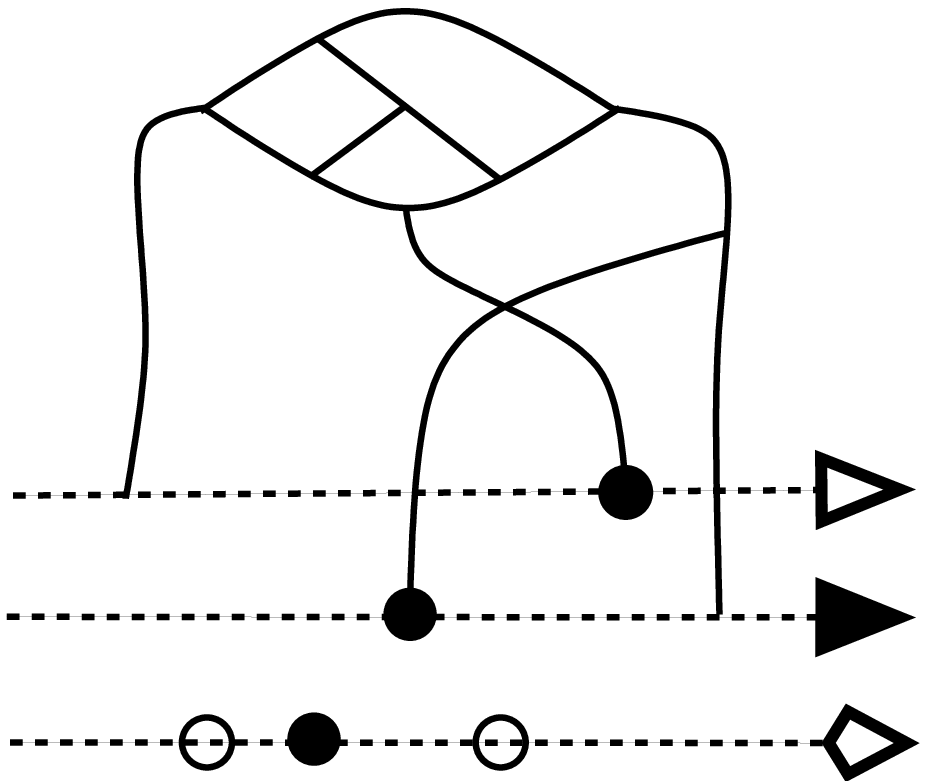}}}\
\ =\ \
\raisebox{-3ex}{\scalebox{0.23}{\includegraphics{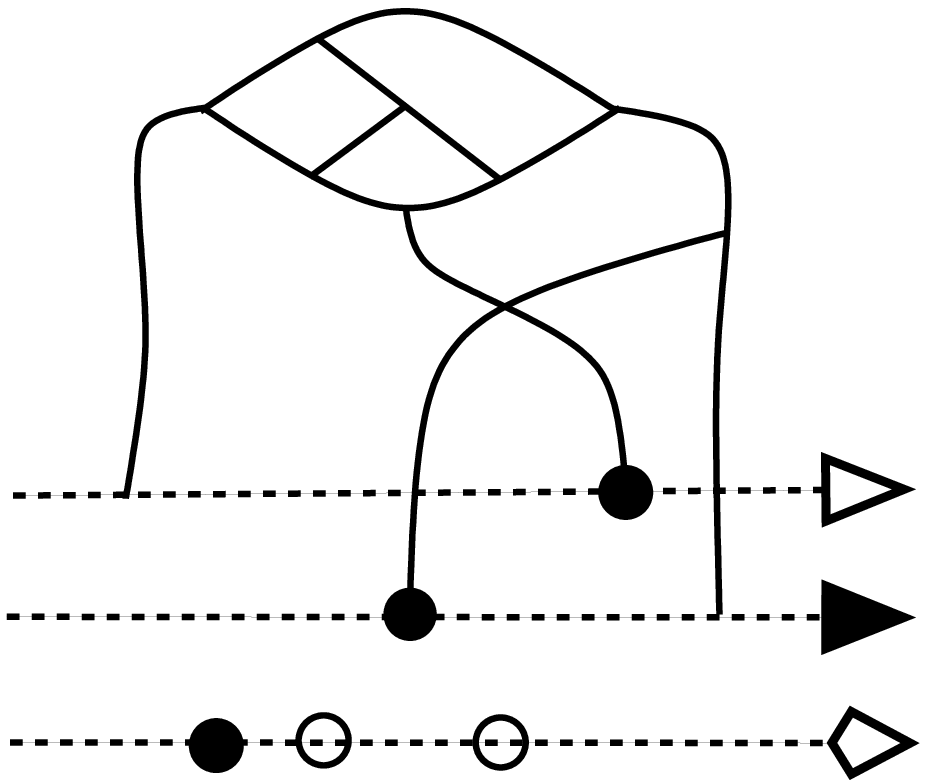}}}\
\ =\ \ -\
\raisebox{-3ex}{\scalebox{0.23}{\includegraphics{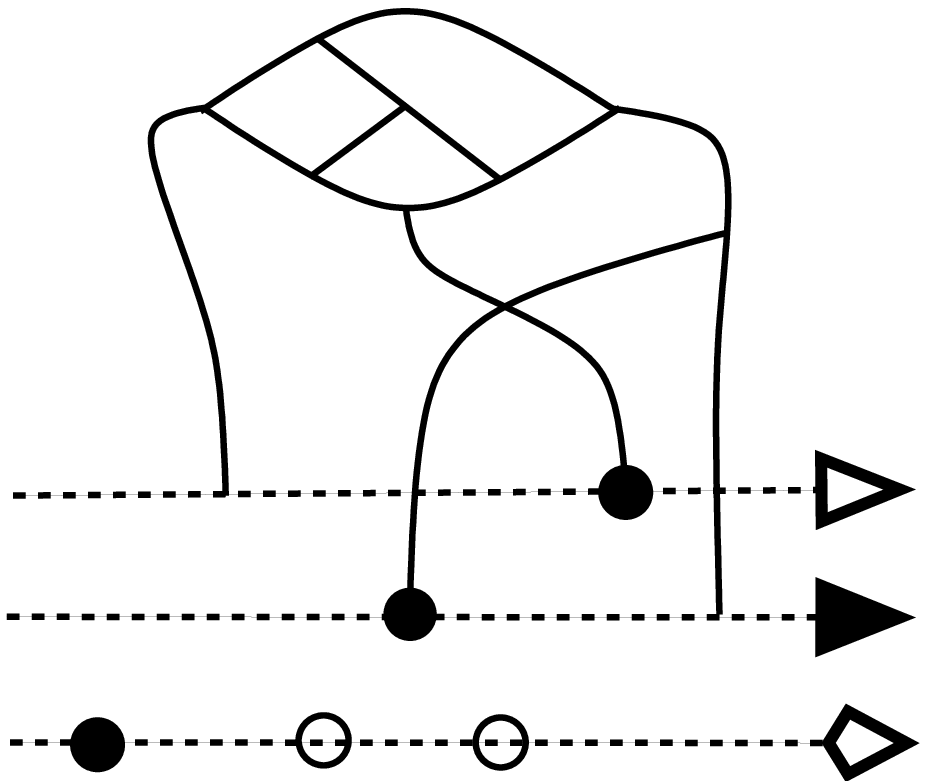}}}
\]
It is a consequence of the permutation rules that a diagram with
more than one filled-in disc is precisely zero.

It will be convenient for the discussion to come to define the
differential as a sum of two grade $1$ operators
\begin{equation}\label{derivdefn}
d_{\mathcal{T}_{\mathrm{dR}}} = d_{\mathrm{T}} + d_\bullet.
\end{equation}
The operator $d_\mathrm{T}$ is defined by the usual substitution
rules operating on the first two lines, such as
\[
\raisebox{-7ex}{\scalebox{0.23}{\includegraphics{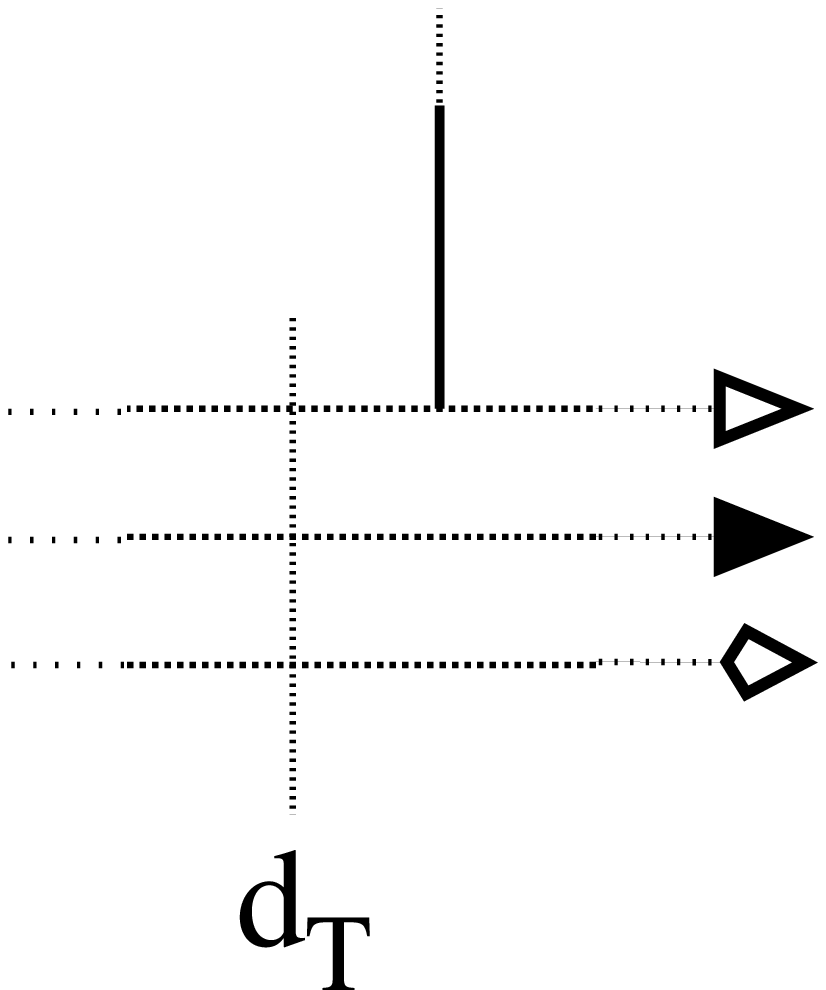}}}\
\ \leadsto\ \
\raisebox{-2.2ex}{\scalebox{0.23}{\includegraphics{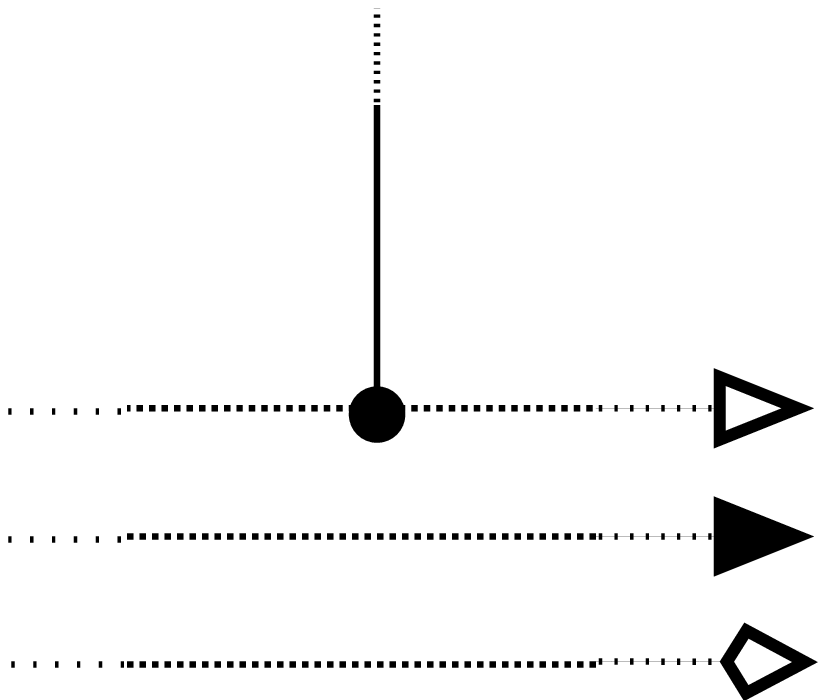}}}\
\ -\ \
\raisebox{-7ex}{\scalebox{0.23}{\includegraphics{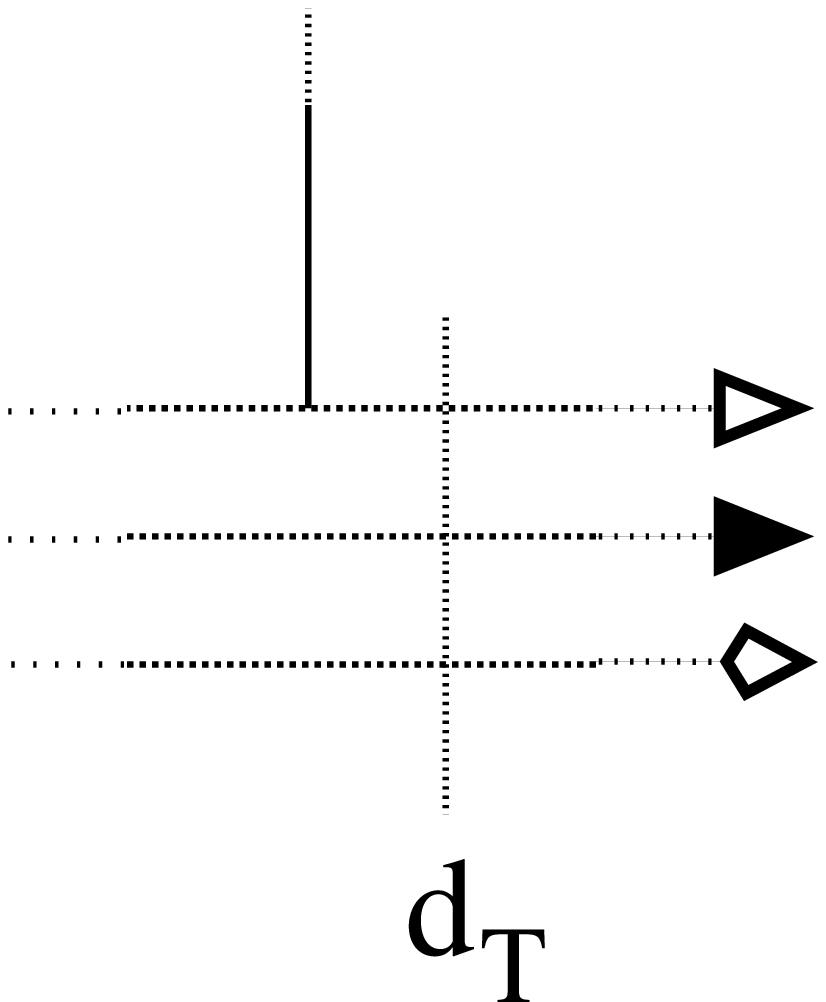}}}\
\ ,
\]
and rules that say that $d_\mathrm{T}$ (graded) commutes through
any labels on the third line:
\[
\raisebox{-7ex}{\scalebox{0.23}{\includegraphics{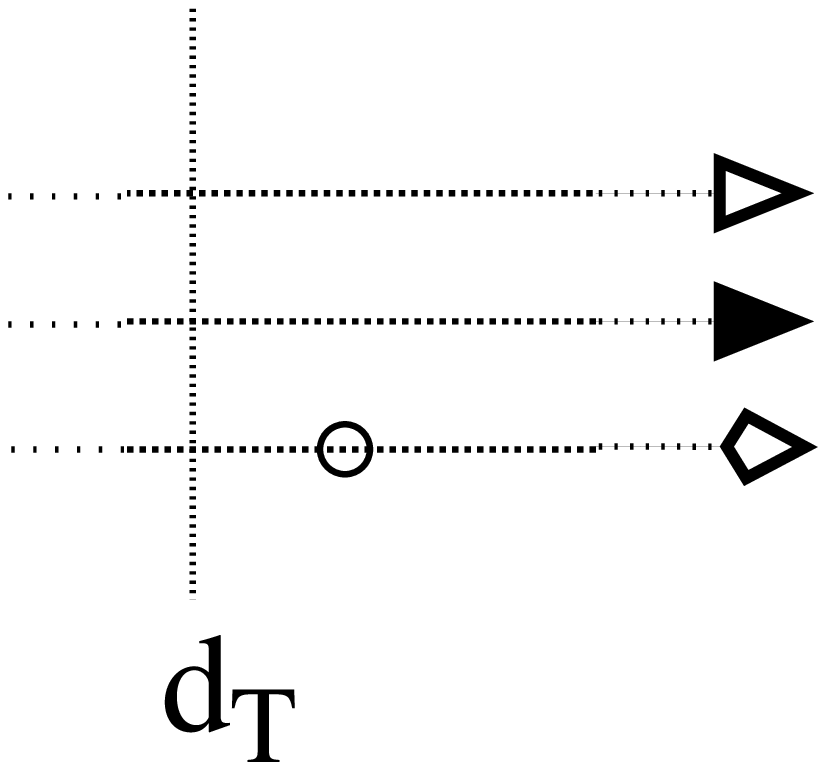}}}\,
\leadsto
\raisebox{-7ex}{\scalebox{0.23}{\includegraphics{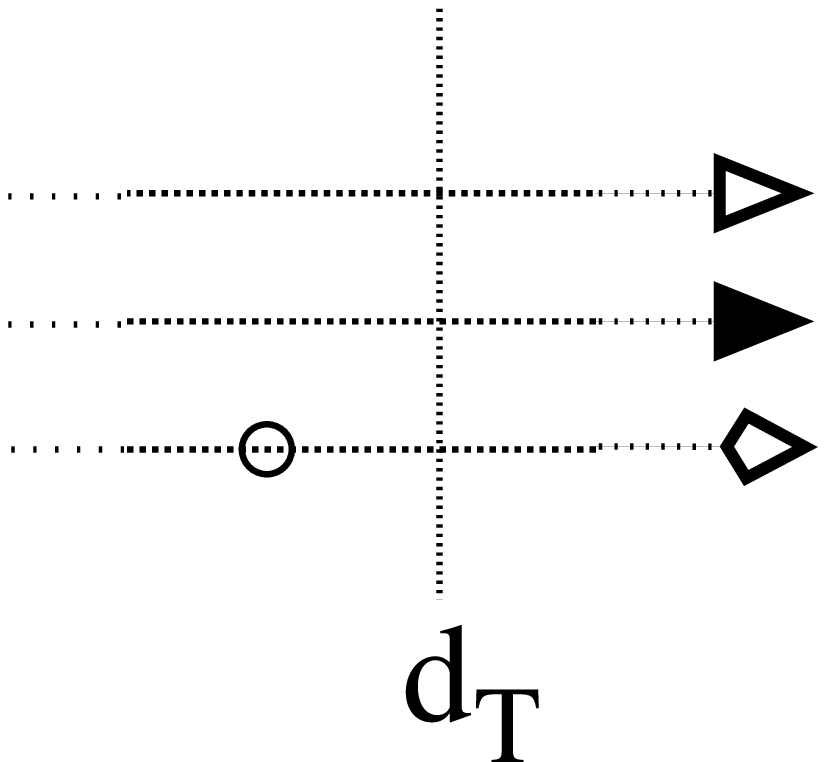}}}\
\ \ \mbox{and}\ \ \
\raisebox{-7ex}{\scalebox{0.23}{\includegraphics{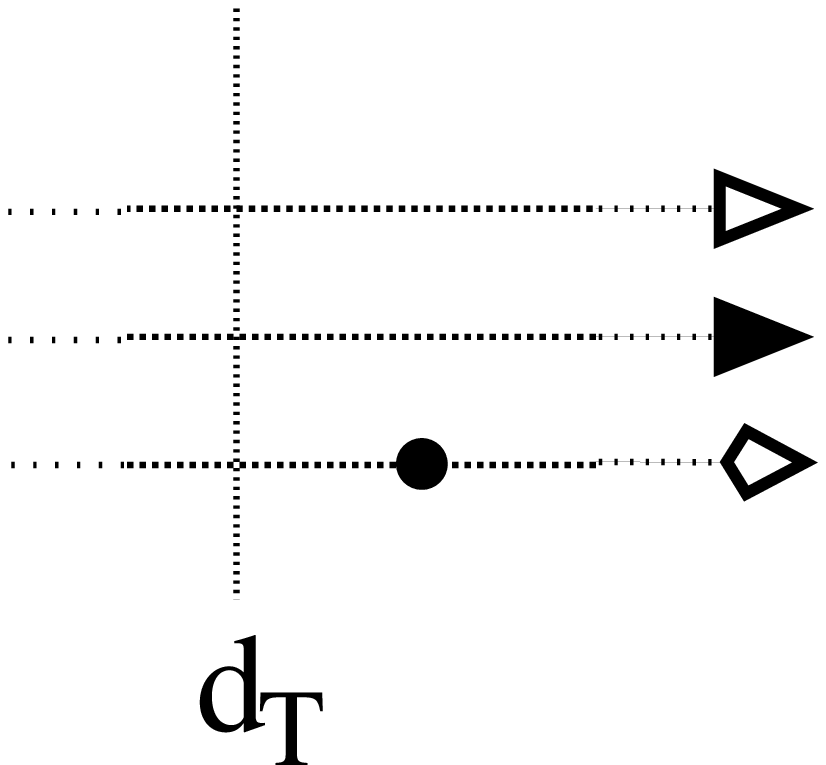}}}\,
\leadsto\ \ -\
\raisebox{-7ex}{\scalebox{0.23}{\includegraphics{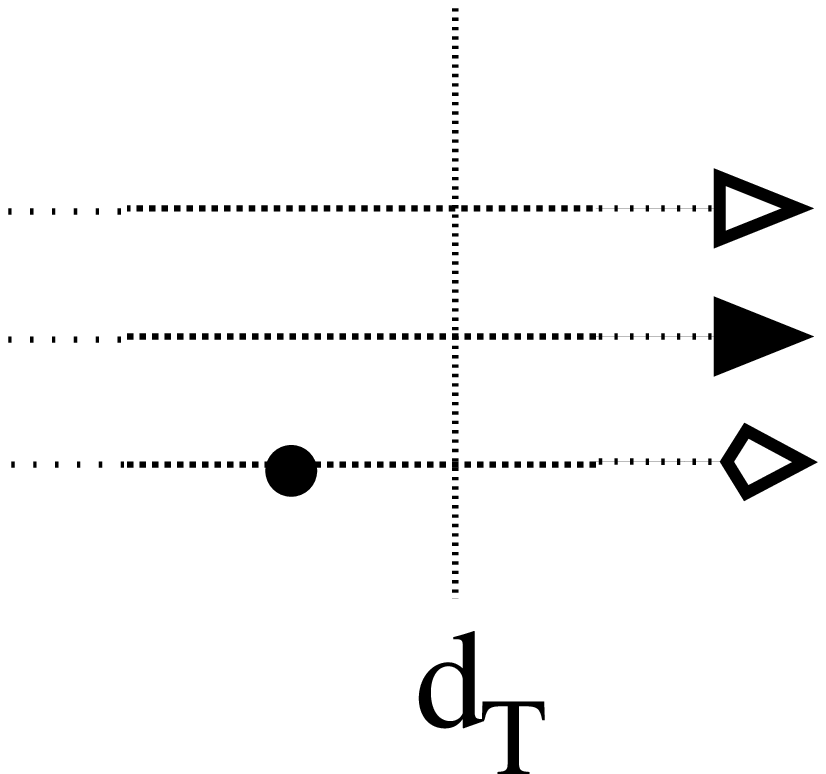}}}\
\ .
\]
On the other hand, the operator $d_\bullet$ is defined by rules
that say that it (graded) commutes through any legs on the first
two lines
together with the following rules which apply to the third line:
\[
\begin{array}{ccccc}
\raisebox{-7ex}{\scalebox{0.23}{\includegraphics{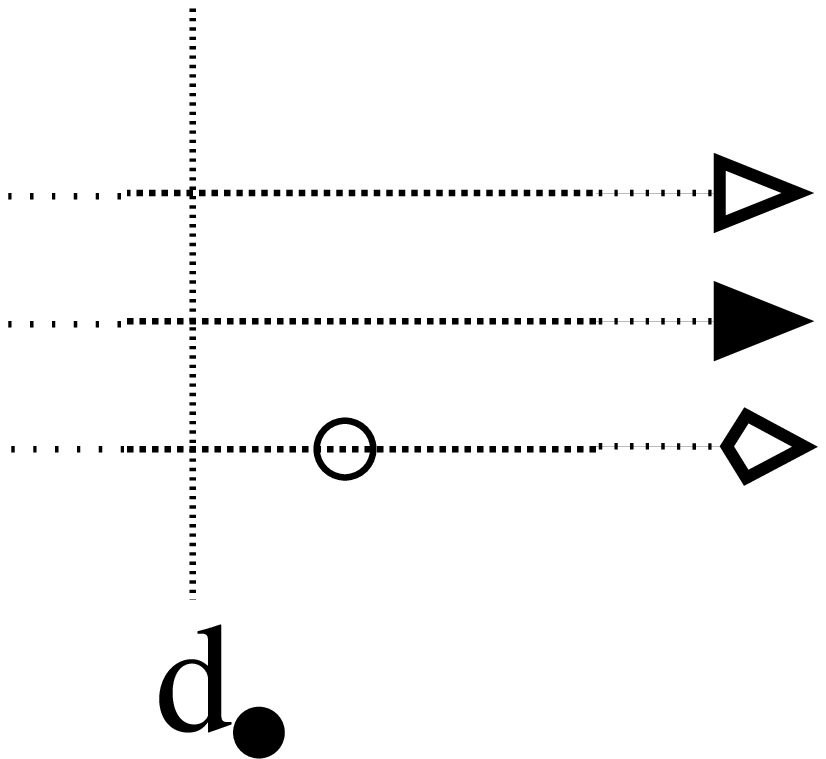}}}
& \leadsto &
\raisebox{-3ex}{\scalebox{0.23}{\includegraphics{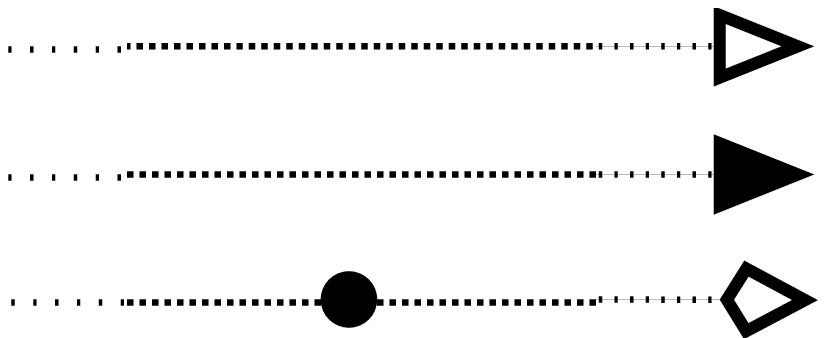}}}
& + &
\raisebox{-7ex}{\scalebox{0.23}{\includegraphics{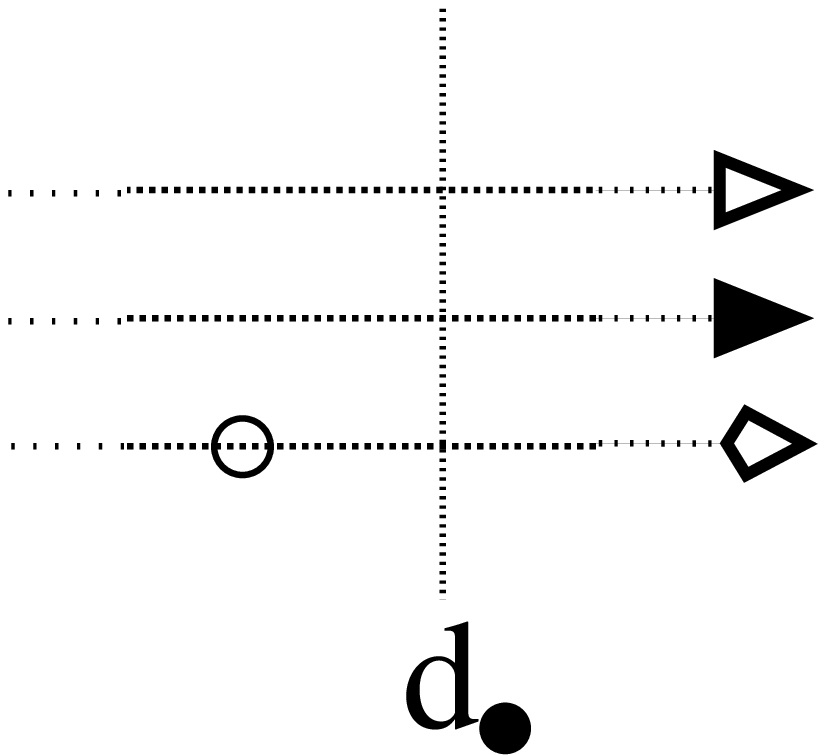}}}\
\ ,
\\[1.3cm]
\raisebox{-7ex}{\scalebox{0.23}{\includegraphics{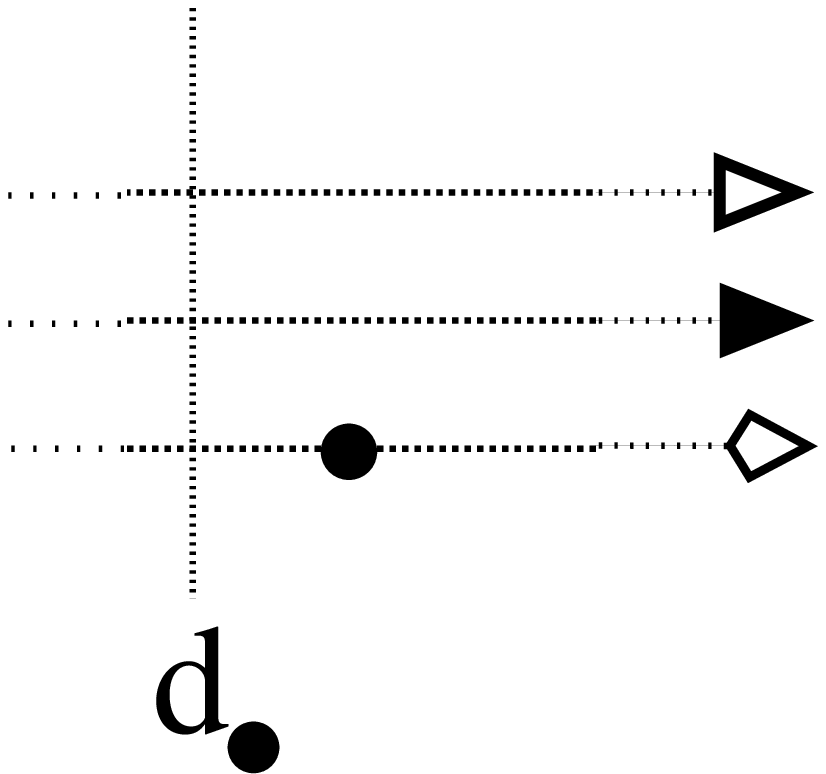}}}
& \leadsto & 0 &  - &{}\
\raisebox{-7ex}{\scalebox{0.23}{\includegraphics{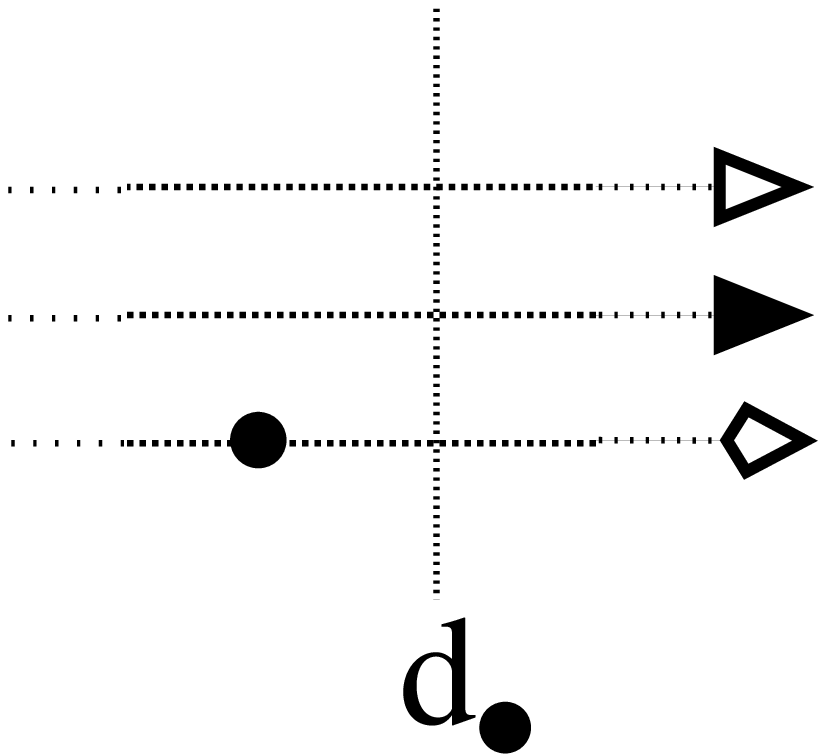}}}\
\ .
\end{array}
\]
The map $\iota$ is defined by the usual rules for $\iota$ applied to the first
two lines,
together with rules which say that $\iota$ (graded) commutes
through anything on the third line.
We invite the reader to check that, with these definitions,
$\mathcal{T}_{\mathrm{dR}}$ does indeed form an $\iota$-complex.
%
%

\subsection{The construction of the \(\iota\)-chain homotopy
\(s_\mathcal{T}\).}

What is this curious complex $\mathcal{T}_{\mathrm{dR}}$ all about
then? The empty disc can be thought of as a formal parameter (it
may be useful to mentally replace the open discs with the symbol
`$t$'), and the filled disc is its differential (`d$t$'). We chose
to use discs to emphasise the fully combinatorial nature of the
manipulations of this section.

We'll begin the construction of $s_\mathcal{T}$ by introducing the
elements of the following diagram of maps of $\iota$-complexes:
\[
\xymatrix{ & & & \mathcal{T} \\ \ncw \ar[urrr]^{i_n}
\ar[drrr]_{i_c} \ar[rr]^{\ \ \ \theta} & &
\mathcal{T}_{\mathrm{dR}}
  \ar[ur]_{\mathrm{Ev}_{\circ\rightarrow 1}}
\ar[dr]^{\mathrm{Ev}_{\circ\rightarrow 0}}
& \\
& & & \mathcal{T} }
\]
In words: we are going to $\iota$-map $\ncw$ into
$\mathcal{T}_{\mathrm{dR}}$ in such a way that when we then set
the formal variable to $0$ (this will be
$\mathrm{Ev}_{\circ\rightarrow 0}$) then all legs are pushed to
the commutative line (in other words, we get $i_c$), while when we
set the formal variable to $1$ (i.e.
$\mathrm{Ev}_{\circ\rightarrow 1}$) then all legs are pushed to
the non-commutative line (i.e. $i_n$). Later, we'll exploit a
formal version of the fundamental theorem of calculus to get an
expression for the difference $i_n-i_c$ that is the subject of
this construction. \\

\begin{figure}
\caption{The replacement rules defining the map $\theta$.\label{thetarules}}
\begin{eqnarray*}
\raisebox{-3.5ex}{\scalebox{0.23}{\includegraphics{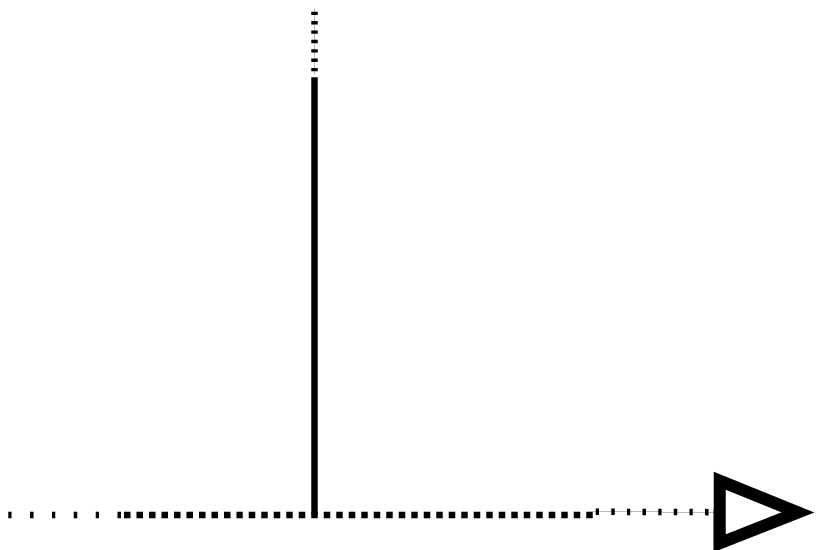}}} &
\mapsto & \ \
\raisebox{-3.5ex}{\scalebox{0.23}{\includegraphics{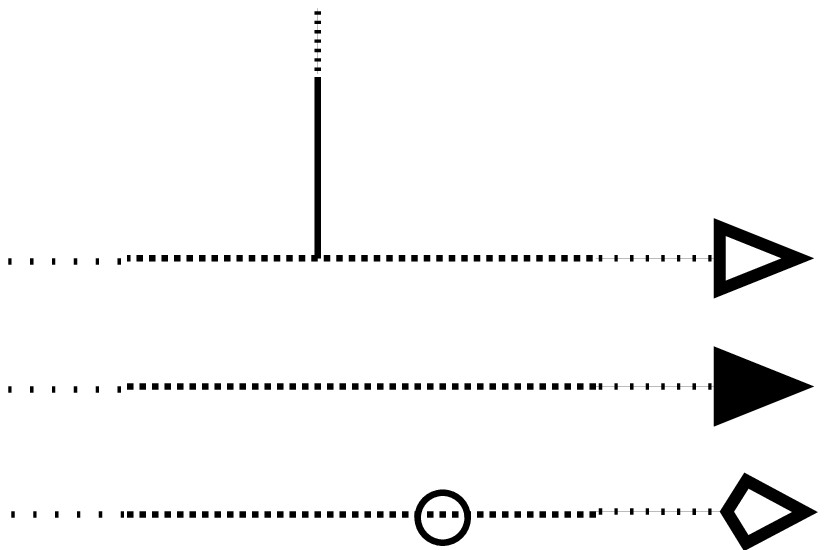}}}\
+\
\raisebox{-3.5ex}{\scalebox{0.23}{\includegraphics{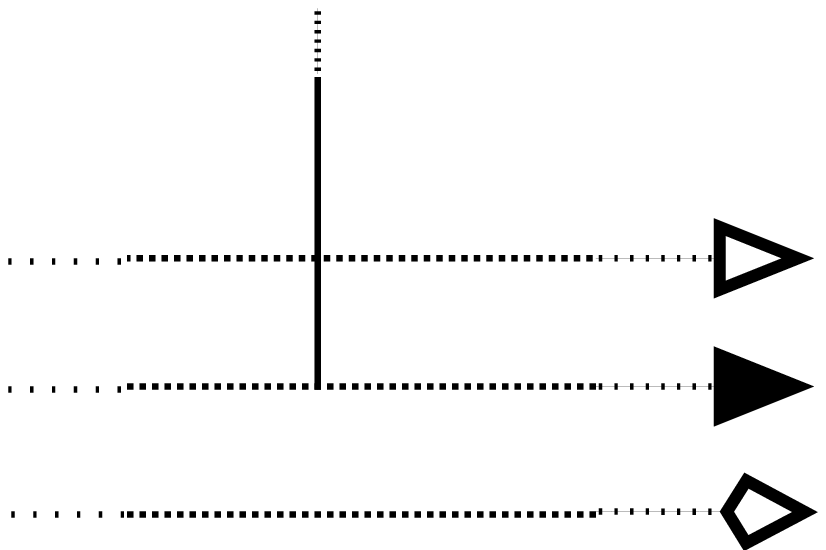}}}\
-\
\raisebox{-3.5ex}{\scalebox{0.23}{\includegraphics{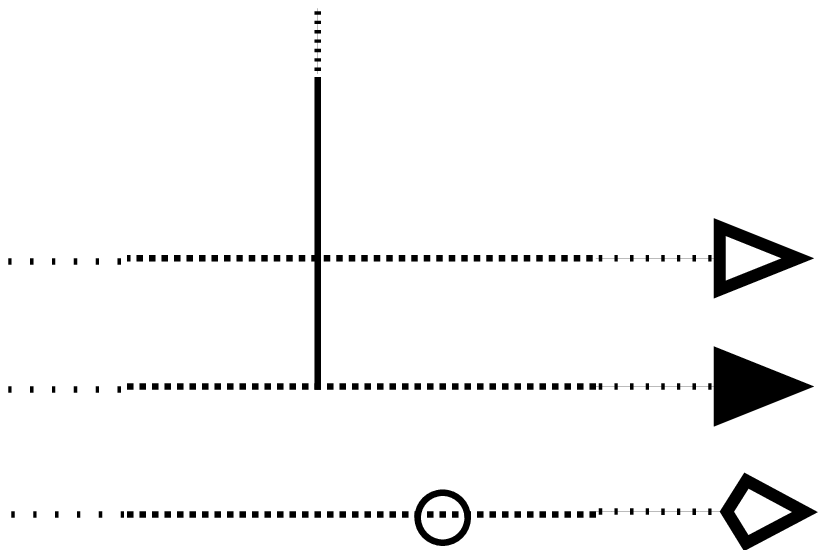}}}\
\ \ ,
\\[0.25cm]
\raisebox{-3.5ex}{\scalebox{0.23}{\includegraphics{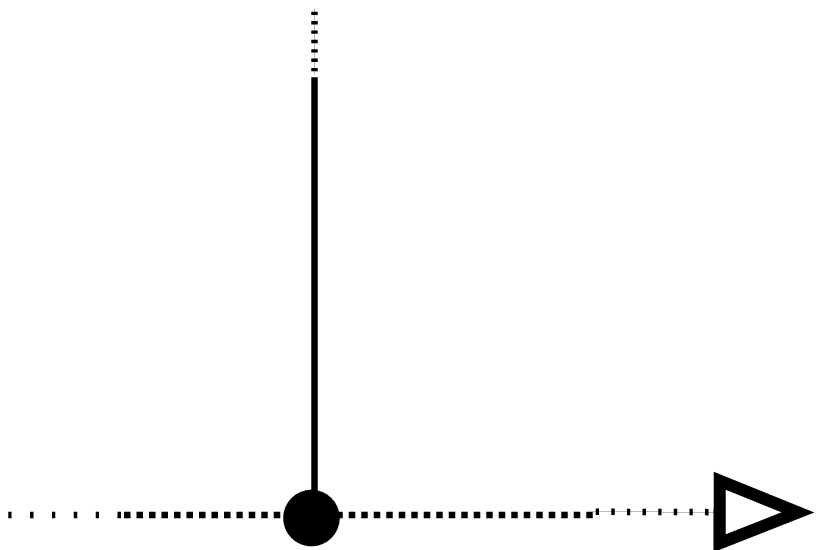}}}
& \mapsto & \ \
\raisebox{-3.5ex}{\scalebox{0.23}{\includegraphics{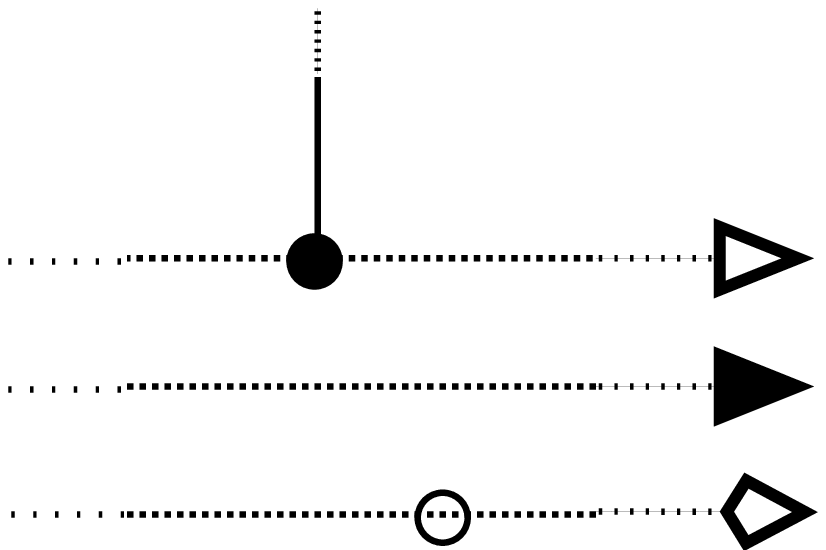}}}\
-\
\raisebox{-3.5ex}{\scalebox{0.23}{\includegraphics{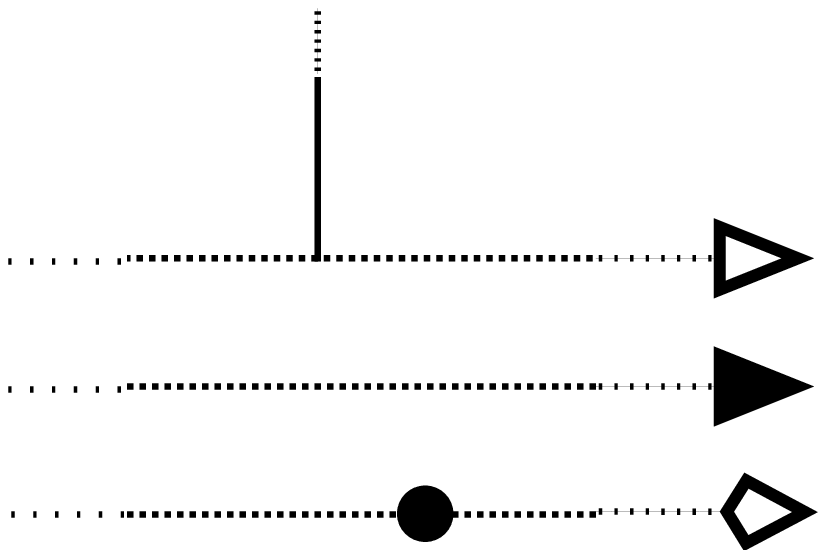}}}\
+\
\raisebox{-3.5ex}{\scalebox{0.23}{\includegraphics{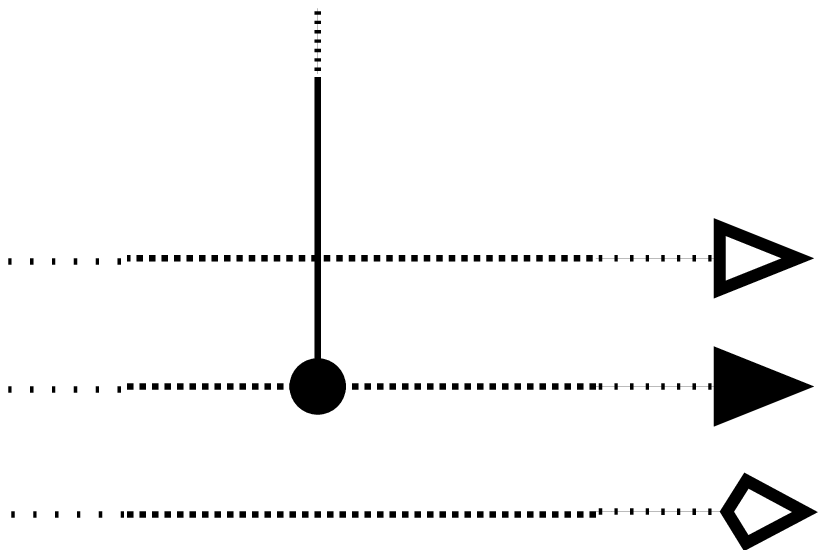}}}
\\[0.25cm]
& & -
\raisebox{-3.5ex}{\scalebox{0.23}{\includegraphics{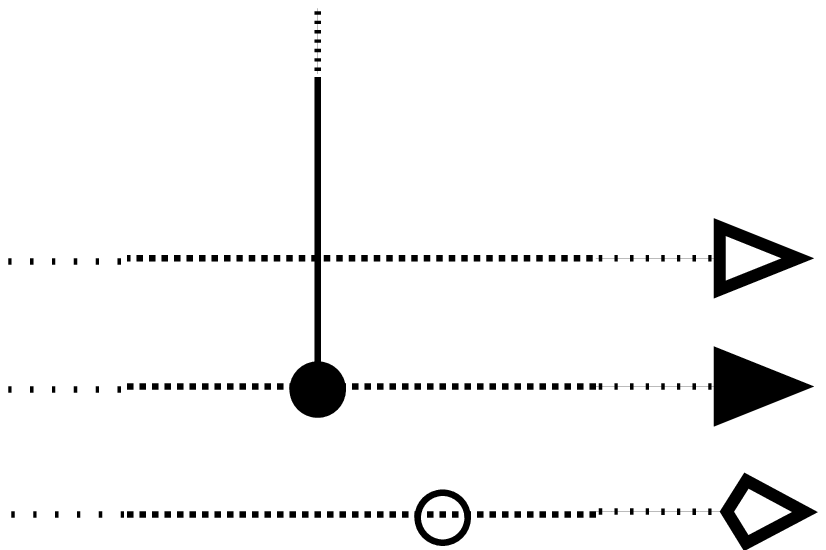}}}\
+\
\raisebox{-3.5ex}{\scalebox{0.23}{\includegraphics{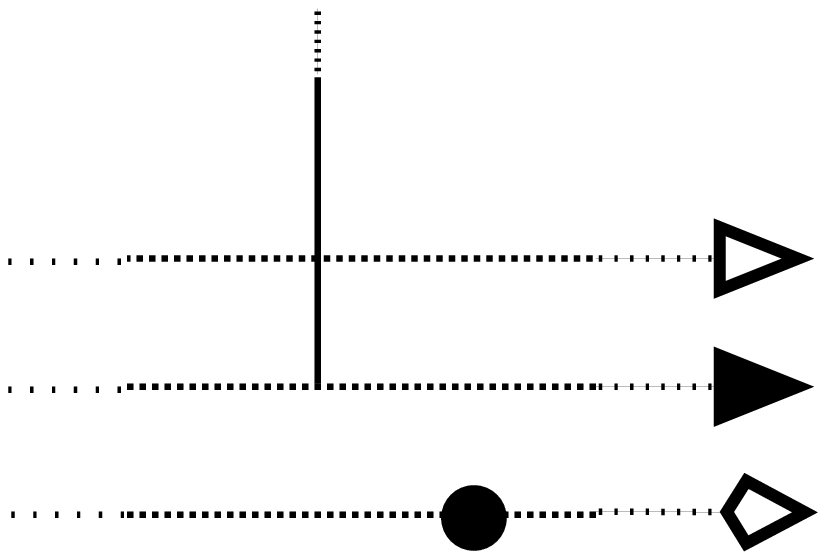}}}\
\ \ . \\
\end{eqnarray*}
\end{figure}

So define a map of $\iota$-complexes $\theta : \widetilde{\mathcal{W}} \rightarrow
\mathcal{T}_{\mathrm{dR}}$
by replacing legs according to the rules shown in Figure \ref{thetarules}.
It is not hard to motivate this seemingly complicated map. The
maps of degree 1 legs are a formal representation of the
combination:
\[
t*\left(\begin{array}{l} \text{Place the leg on the} \\
\text{non-commutative line.} \end{array}\right) +
(1-t)*\left(\begin{array}{l} \text{Place the leg on the} \\
\text{commutative line.} \end{array}\right)
\]
The degree 2 legs then map in the only way possible that will
result in a map of complexes.
\begin{prop}\label{thetaisamapoficomplexes}
This defines a map of $\iota$-complexes.
\end{prop}
Checking that this is true is an instructive exercise. Now let's
meet the combinatorial maps which represent the act of ``setting
the parameter to $0$" and ``setting the parameter to $1$".

\begin{defn}
Define a map of $\iota$-complexes
$\mathrm{Ev}_{\circ\rightarrow
0} : \mathcal{T}_\mathrm{dR} \rightarrow \mathcal{T}$
by declaring its value on some diagram to be zero if the third
line of the diagram is labelled with anything, otherwise declaring
its value to be the diagram that remains when that unlabelled
third line is removed.
\end{defn}

\begin{defn}
Define a map of $\iota$-complexes $
\mathrm{Ev}_{\circ\rightarrow 1} : \mathcal{T}_\mathrm{dR}
\rightarrow \mathcal{T}$
by declaring its value on a diagram to be 0 if the third line of
that diagram has any grade 1 (i.e. filled) discs on it, otherwise
declaring its value to be the diagram that remains when the third
line is removed.
\end{defn}
For example:
\[
\mathrm{Ev}_{\circ\rightarrow 1}\left(
\raisebox{-5ex}{\scalebox{0.24}{\includegraphics{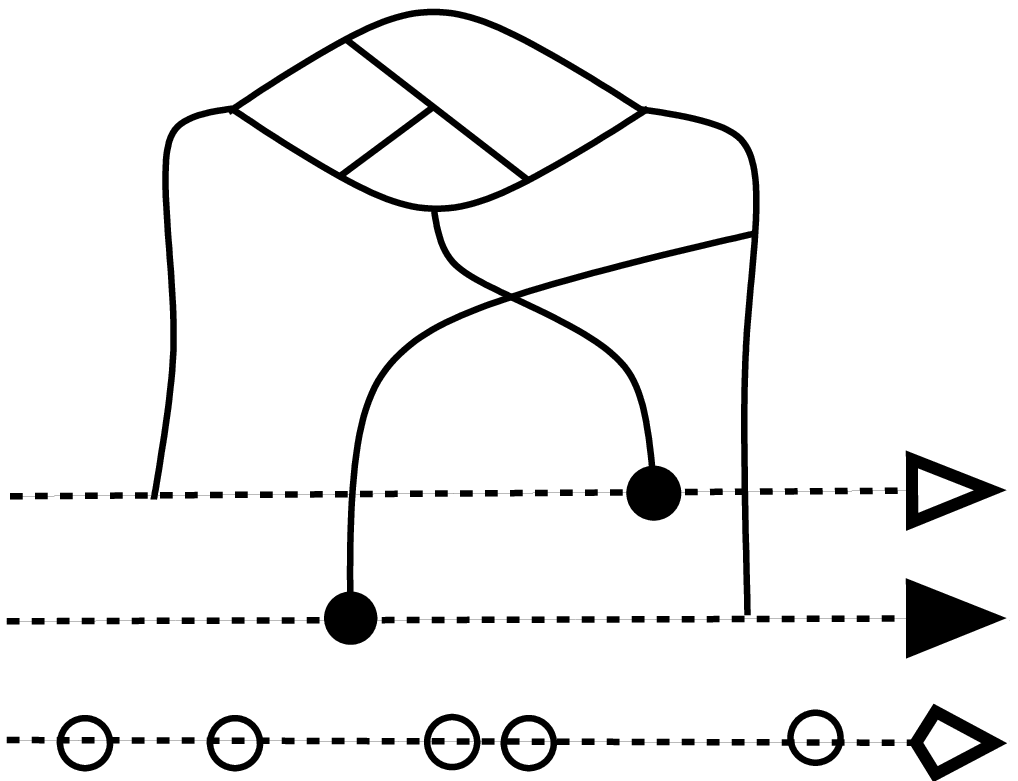}}} \ +\
\raisebox{-5ex}{\scalebox{0.24}{\includegraphics{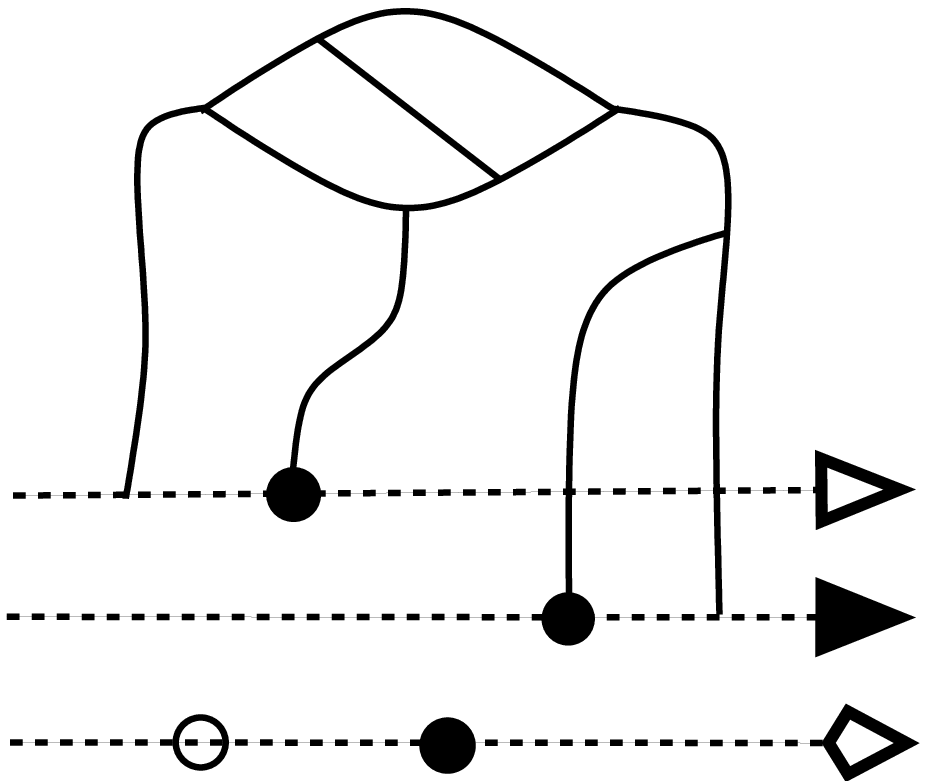}}}
\right)\ =\
\raisebox{-4ex}{\scalebox{0.24}{\includegraphics{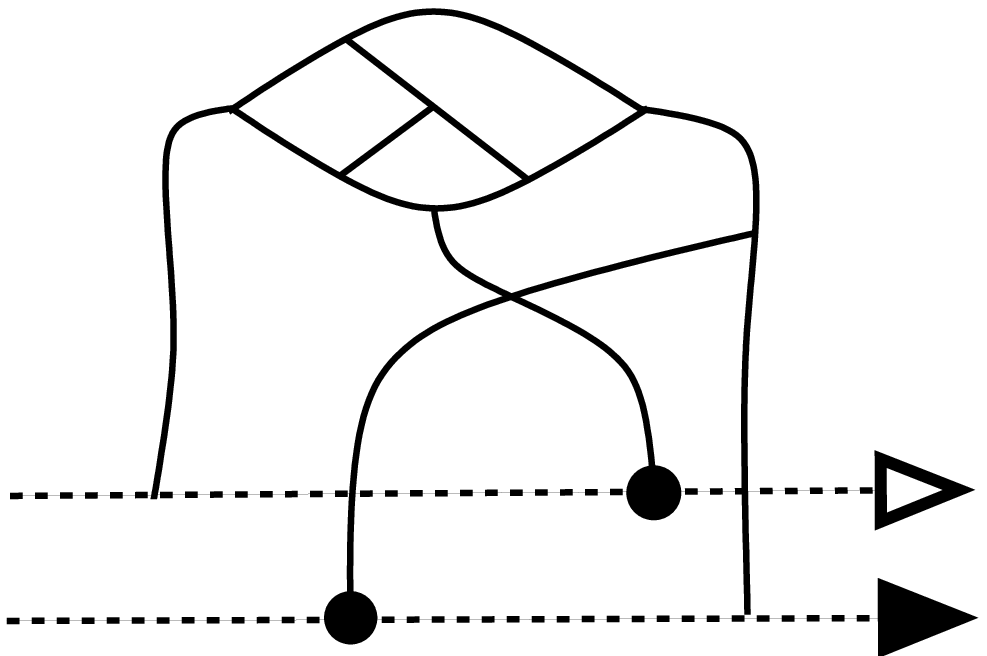}}}\ \ .
\]
The reader can check that, as desired:
\begin{eqnarray}
i_c & = & \mathrm{Ev}_{\circ\rightarrow 0}\circ \theta, \label{pleasingproperty1}\\
i_n & = & \mathrm{Ev}_{\circ\rightarrow 1}\circ \theta.
 \label{pleasingproperty2}
\end{eqnarray}

Now what we want to study is the difference $i_n-i_c$, which we
have just shown is equal to:
\[
i_n-i_c = \left(\mathrm{Ev}_{\circ\rightarrow
1}-\mathrm{Ev}_{\circ\rightarrow 0}\right)\circ\theta.
\]
Intuitively, then, we want to take a difference of what we get
when we evaluate some expression at $1$ with what we get when we
evaluate it at $0$. A formal version of the fundamental theorem of
calculus will will do that for us. First we need to learn to
integrate:
\begin{defn}
Define linear maps $\intoper^i:\mathcal{T}_\mathrm{dR}^i
\rightarrow \mathcal{T}^i$ (and also $\iota$ versions)
 by the following procedure:
\begin{enumerate}
\item{The value of $\intoper$ on some diagram is zero unless the
diagram has precisely one grade $1$ (i.e. filled-in) disc on its
third line.} \item{If a diagram has precisely one such disc, then
to calculate the value of $\intoper$ on the diagram, begin by
permuting (with the appropriate signs) the disc to the far-left
end of the third line.} \item{Then remove the third line and
multiply the resulting diagram by $\frac{1}{n+1}$, where $n$ is
the number of grade $0$ discs on that line.}
\end{enumerate}
\end{defn}
Here is an example of the calculation of $\intoper$:
\begin{multline*}
\intoper\left(
\raisebox{-4ex}{\scalebox{0.23}{\includegraphics{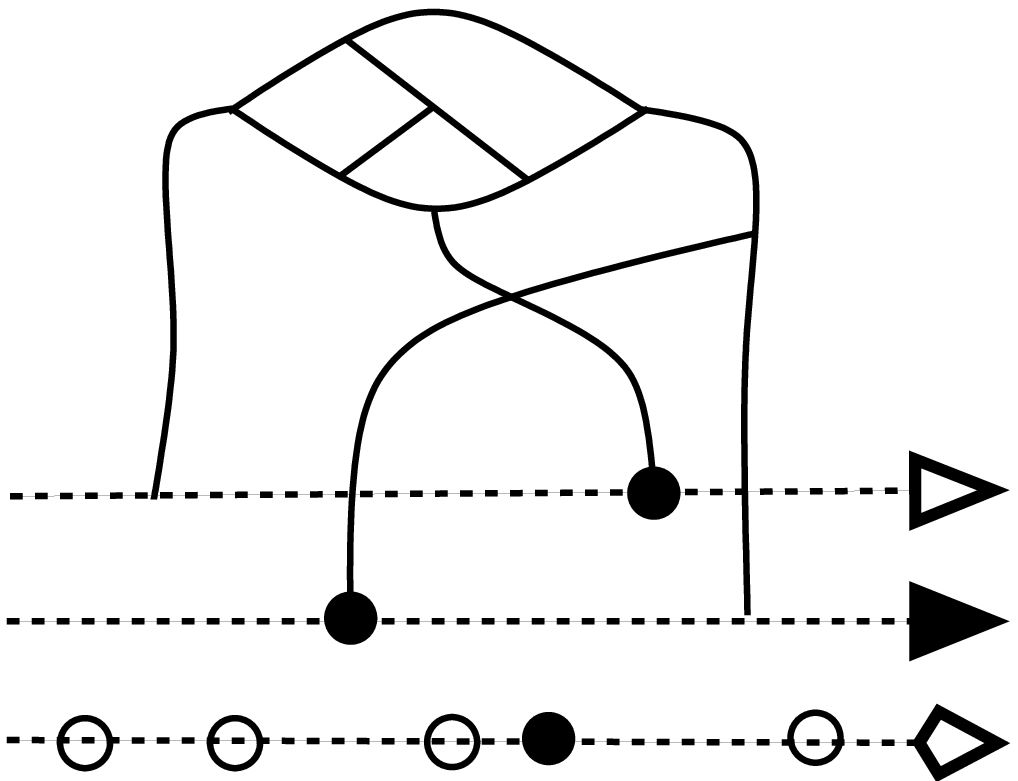}}}\
\right.\left. +\
\raisebox{-4ex}{\scalebox{0.23}{\includegraphics{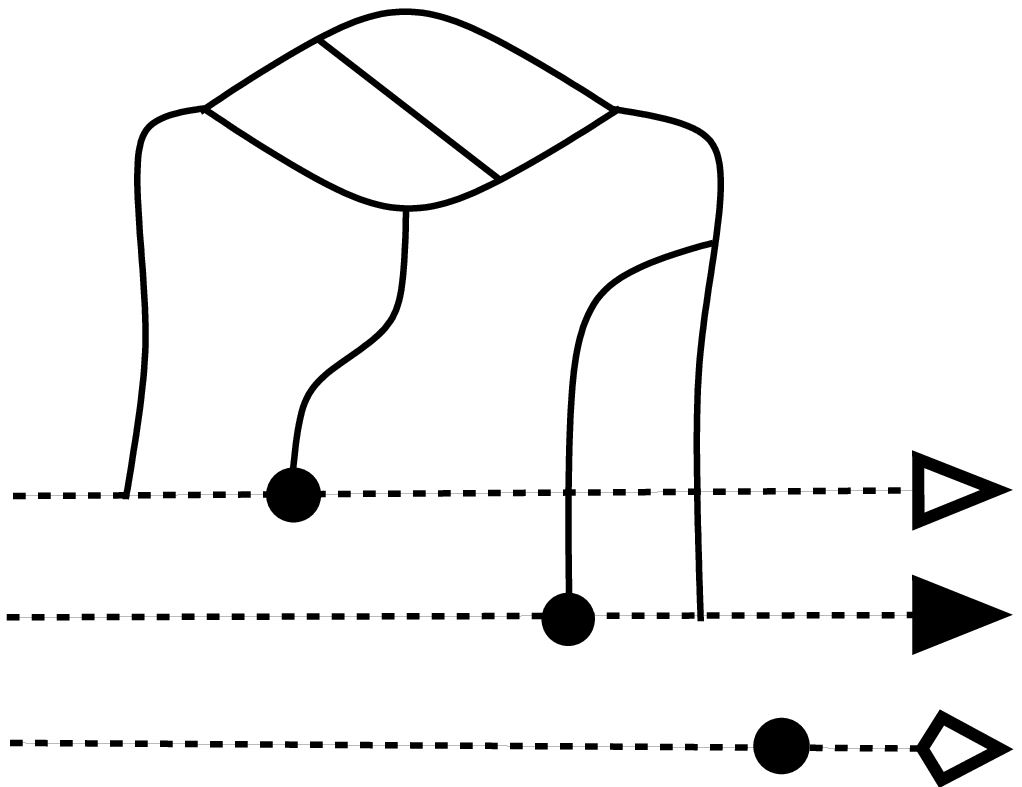}}}\ +\ \
\raisebox{-4ex}{\scalebox{0.23}{\includegraphics{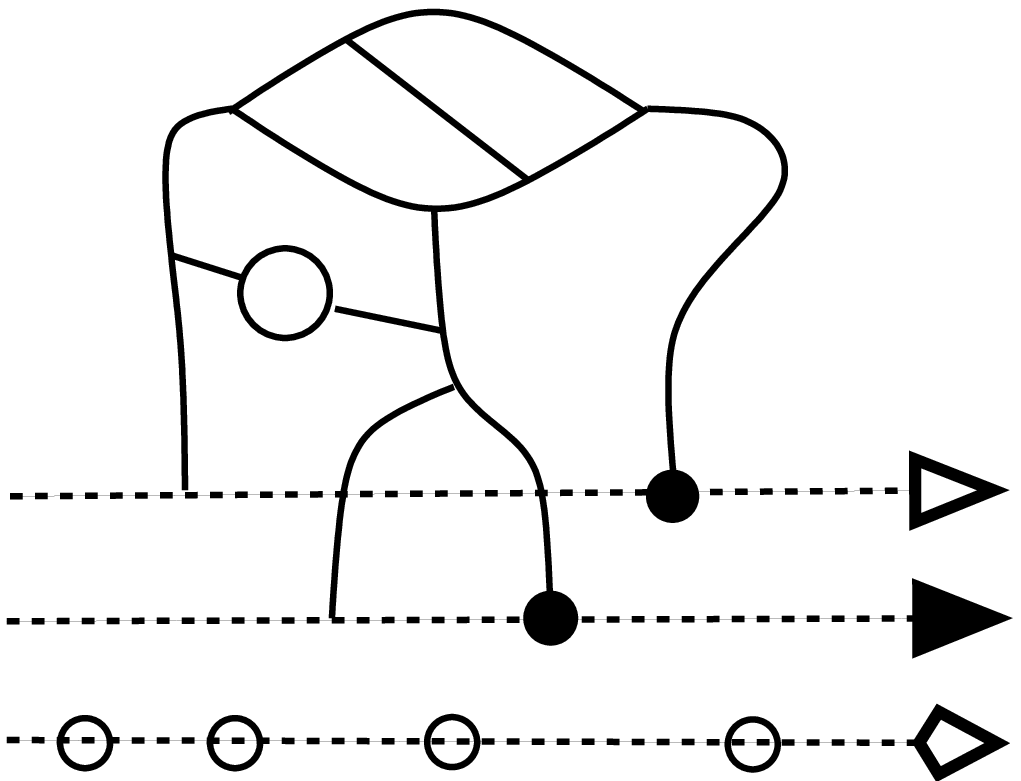}}}
\right)\\
 = -\frac{1}{5}\,
\raisebox{-3ex}{\scalebox{0.23}{\includegraphics{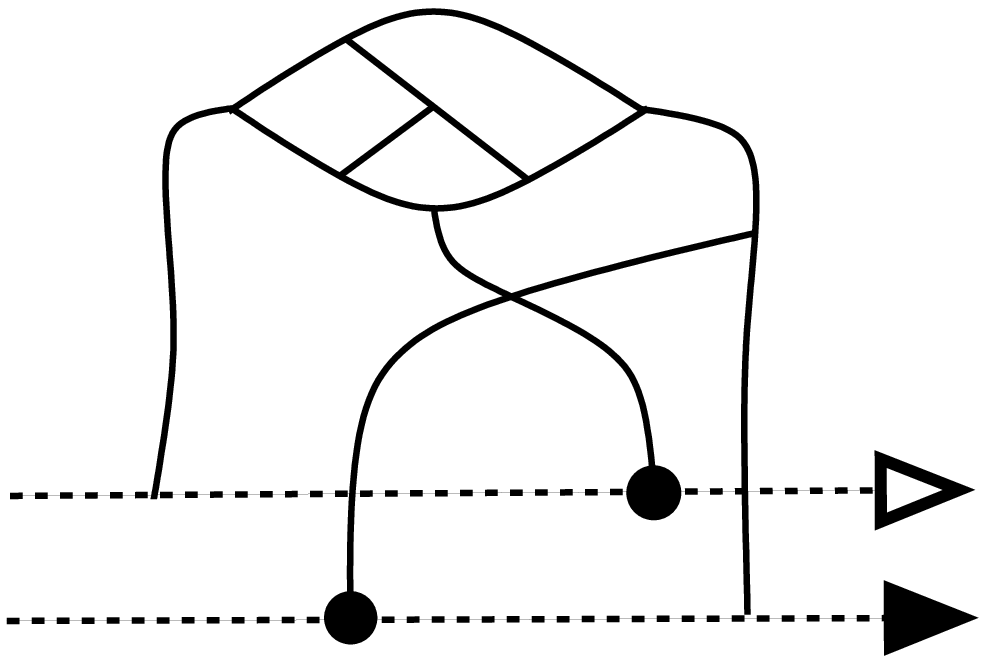}}}\ +
\raisebox{-3ex}{\scalebox{0.23}{\includegraphics{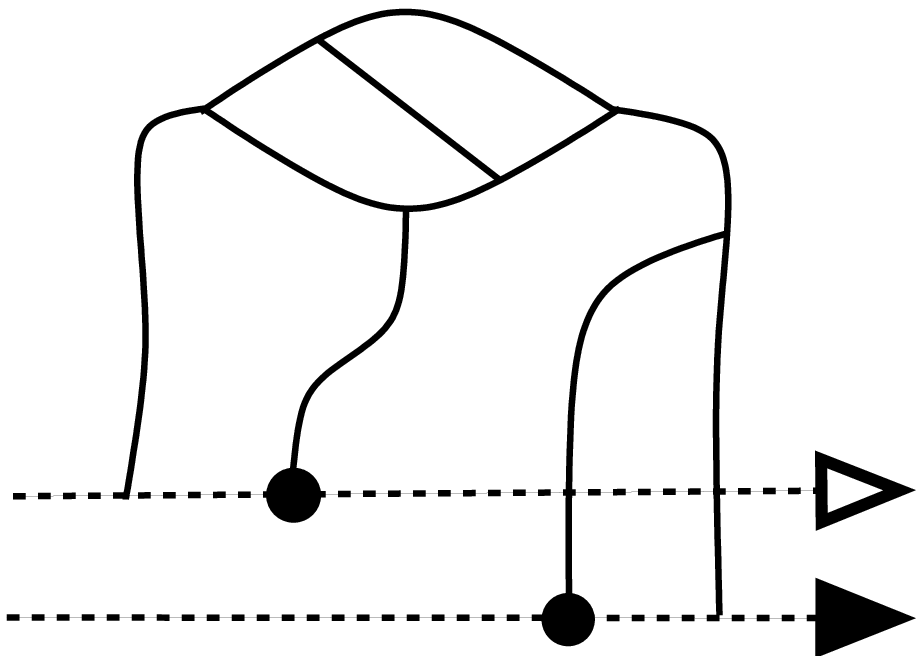}}}\ \ .
\end{multline*}

These grade $-1$ maps do not {\it quite} commute with all the
involved differential operators; observe that the following lemma
is not stated in terms of the full differential on
$\mathcal{T}_{\mathrm{dR}}$. (Equation \ref{derivdefn} defines the
differential $d_{\mathcal{T}_{\mathrm{dR}}}$ as a sum of two grade
$1$ maps: $d_{\mathrm{T}} + d_\bullet$.)

\begin{lem}\label{iotaandscommute}
\begin{eqnarray}
\intoper\circ\,\iota_{\mathcal{T}_{\mathrm{dR}}} +
\iota_{\mathcal{T}}\circ\intoper & = & 0. \label{sandiotaA}
\\
\intoper\circ\, d_{\mathrm{T}} + d_\mathcal{T}\circ\intoper & = &
0. \label{sandiotaB}
\end{eqnarray}
\end{lem}
These equations are straightforward to check - the operation
$\int_{\mathcal{T}_{\mathrm{dR}}}$ is well-defined on
$\mathcal{T}_{\mathrm{dR}}$ so we may assume that the single
filled-in disc lies at the far left hand end of the diagram.

\begin{lem}[``The fundamental theorem of the calculus.'']\label{keylem}
\[
\mathrm{Ev}_{\circ\rightarrow 1} - \mathrm{Ev}_{\circ\rightarrow
0} = \intoper \circ\, d_\bullet\ \ .
\]
\end{lem}
{\it Proof.} We'll simply check that this equation is true on the
different classes of generators.
The only non-trivial case when the third line has no
grade $1$ disc and some positive number $n$ of grade $0$ discs on
its third line. Let's do an example to observe why the equation is
true in this case.

First the left-hand side:
\[
\left(\mathrm{Ev}_{\circ\rightarrow 1} -
\mathrm{Ev}_{\circ\rightarrow 0}\right)\left(
\raisebox{-4ex}{\scalebox{0.23}{\includegraphics{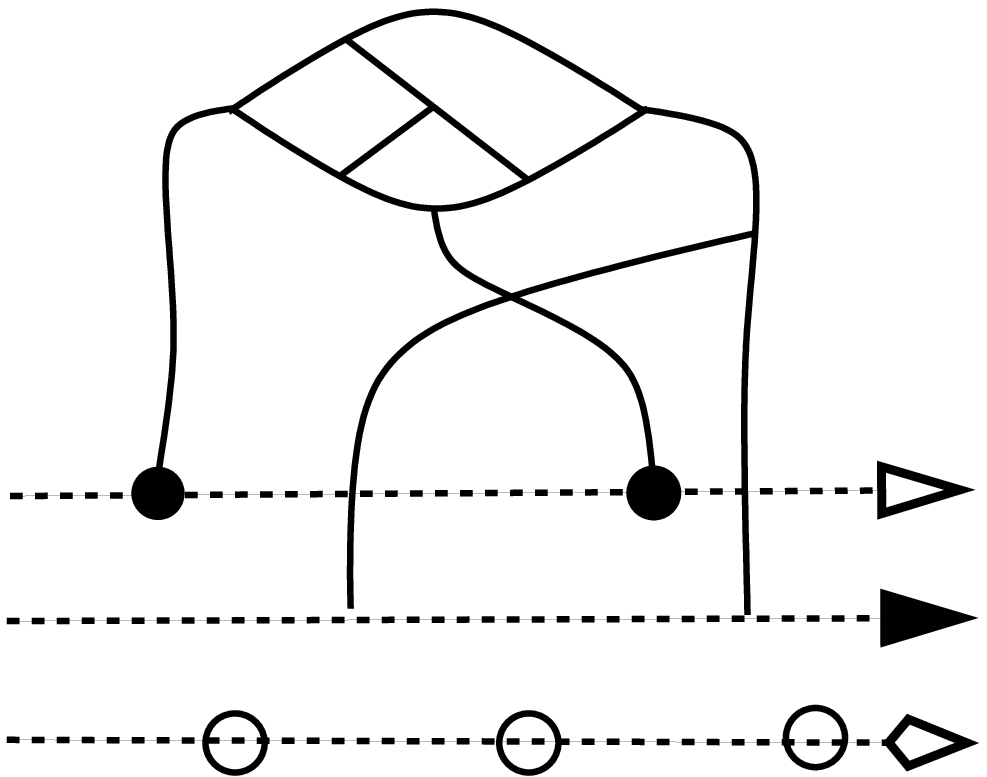}}}
\right)\ \ =\
\raisebox{-4ex}{\scalebox{0.23}{\includegraphics{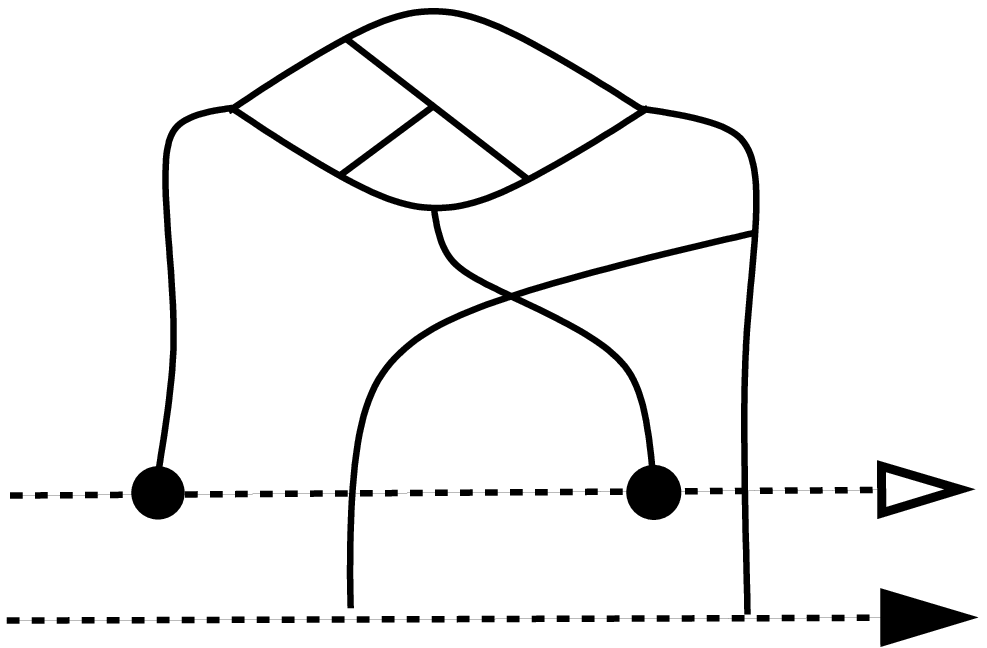}}}\ \
.
\]
Now consider the right-hand side. Applying $d^6_\bullet$ to this diagram gives:
\[
\raisebox{-4ex}{\scalebox{0.23}{\includegraphics{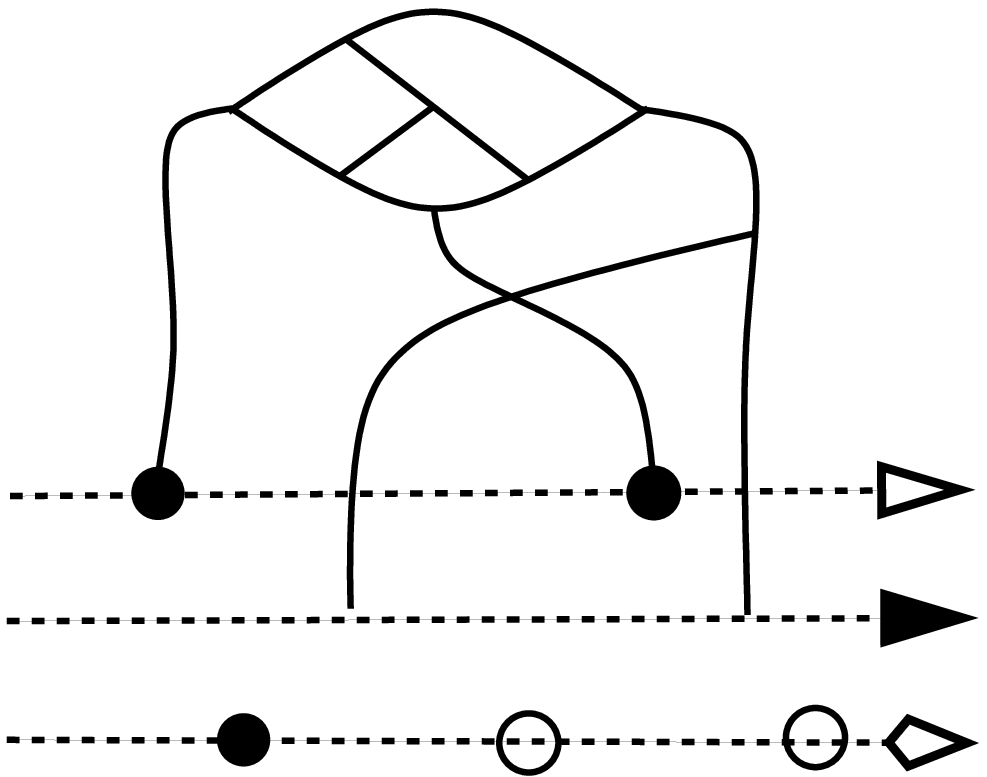}}}\
-\
\raisebox{-4ex}{\scalebox{0.23}{\includegraphics{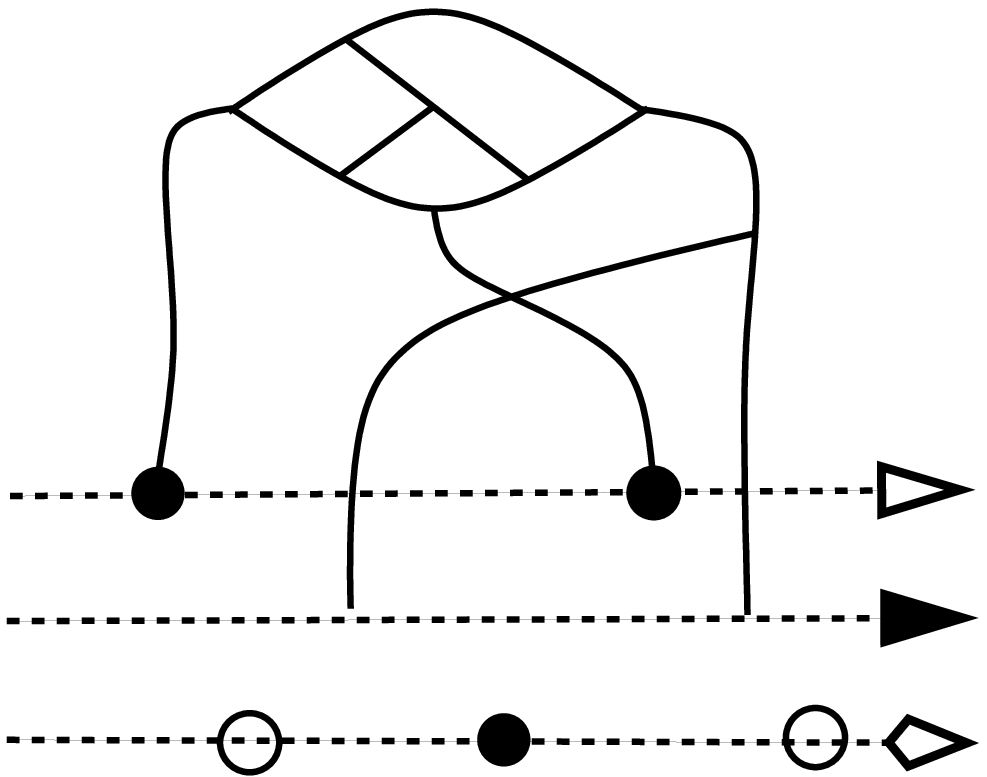}}}\
+\
\raisebox{-4ex}{\scalebox{0.23}{\includegraphics{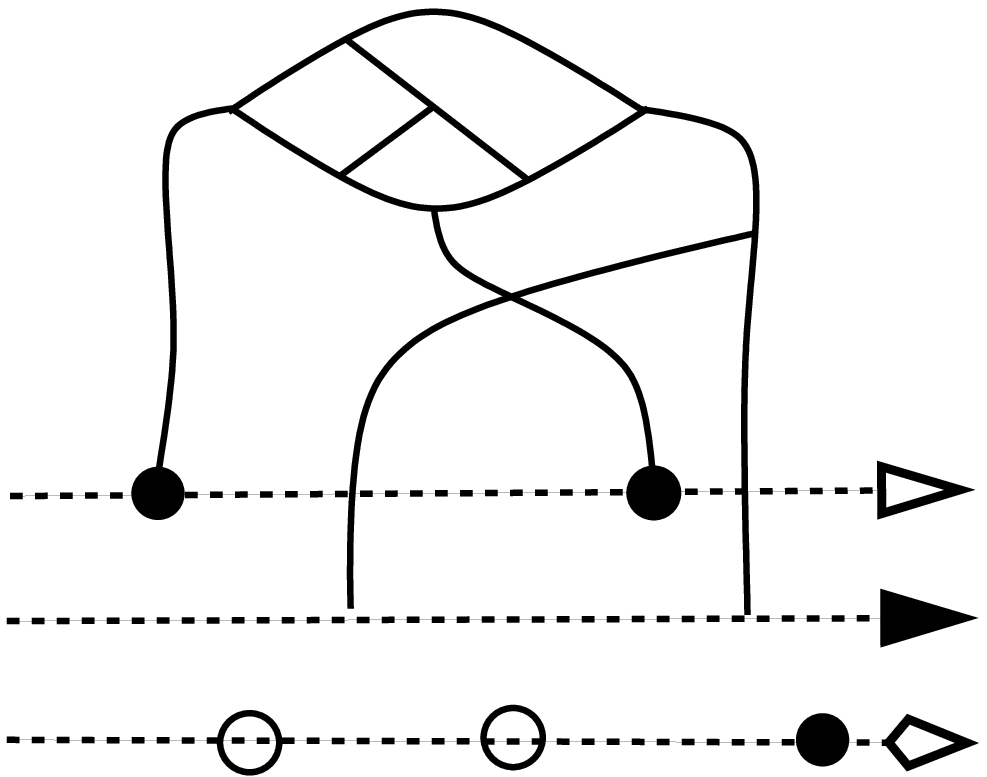}}}\ =\ 3
\raisebox{-4ex}{\scalebox{0.23}{\includegraphics{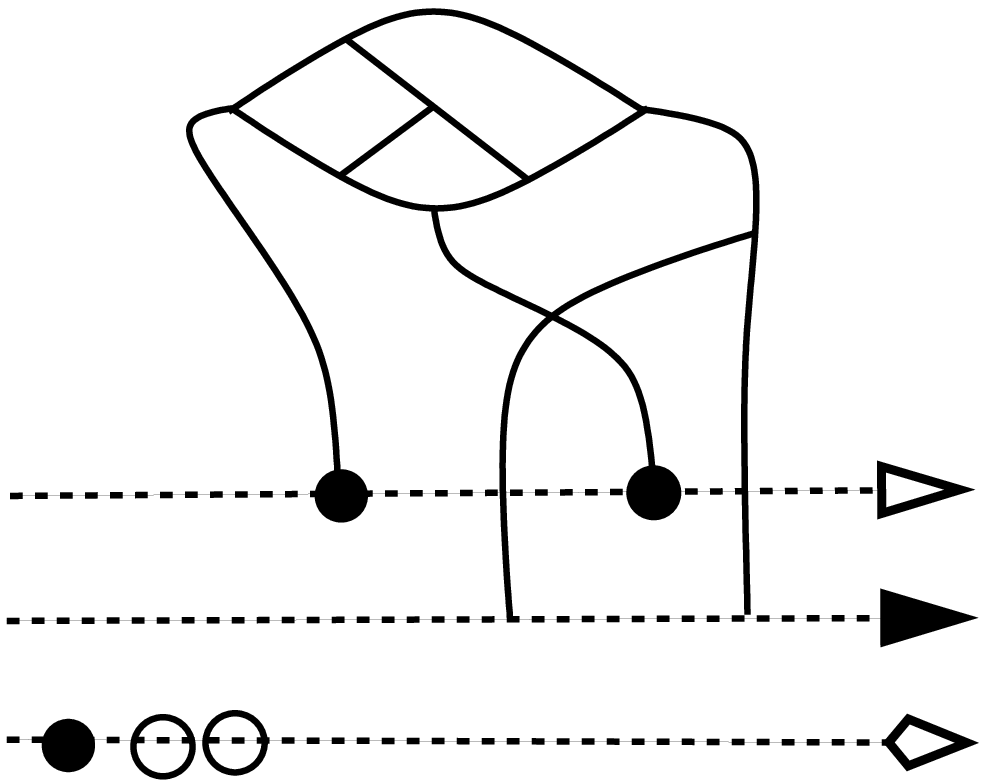}}}\ \ .
\]

Thus:
\begin{eqnarray*}
\left(\intoper\circ\,d_\bullet\right)\left(
\raisebox{-4ex}{\scalebox{0.23}{\includegraphics{lemmaexampA}}}
\right) &  = & 3\,\intoper\left(
\raisebox{-4ex}{\scalebox{0.23}{\includegraphics{lemmaexampY}}}
\right) \\ &  = & (3)\left(\frac{1}{3}\right)
\raisebox{-4ex}{\scalebox{0.23}{\includegraphics{lemmaexampZ}}}\ \
.
\end{eqnarray*}
The reader should not have any difficulties constructing a general
argument from this example.
\begin{flushright}
$\Box$
\end{flushright}

{\it Proof of Theorem \ref{keytheorem}.} Define the required
homotopy $s_\mathcal{T}$ by the formula
\[
s_\mathcal{T} = \intoper\circ\,\theta.
\]
Observe:
\[
\begin{array}{cclp{1cm}l}
i_n - i_c & = & \left(\mathrm{Ev}_{\circ\rightarrow 1} -
\mathrm{Ev}_{\circ\rightarrow 0}\right)\circ\theta, & &
\mbox{(Eqns. \ref{pleasingproperty1} and
\ref{pleasingproperty2})},
\\
& = & \intoper\circ\, d_\bullet \circ \theta, & &
\mbox{(Lemma \ref{keylem})}, \\
& = & \intoper\circ\, (d_{\mathcal{T}_{\mathrm{dR}}} -
d_\mathrm{T}) \circ \theta, &
& \mbox{(Eqn. \ref{derivdefn})}, \\
& = & (\intoper\circ\, d_{\mathcal{T}_\mathrm{dR}}\circ\theta) -
(\intoper\circ d_\mathrm{T}\circ\theta), & & \\
& = & (\intoper\circ\,\theta\circ d_{\ncw}) - (
\intoper\circ\, d_\mathrm{T}\circ\theta), & & \mbox{(Prop. \ref{thetaisamapoficomplexes})}, \\
& = & (\intoper\circ\,\theta\circ d_{\ncw}) + (
d_\mathcal{T}\circ\intoper\circ\,\theta), & & \mbox{(Eqn.
\ref{sandiotaA})},
\\ & = & s_\mathcal{T}\circ
d_{\ncw} + d_\mathcal{T}\circ s_\mathcal{T}.
\end{array}
\]
\begin{flushright}
$\Box$
\end{flushright}

\section{How to obtain Wheeling from Homological
Wheeling.}\label{gettingwheeling}
\subsection{\(\di\)-pairs.}
For the remainder of this work we'll be working with a number of
systems without a natural $\mathbb{Z}$-grading, so it will be
necessary to adjust our language a little. Define a $\di$-pair to
be a pair of vector spaces with maps, \[ \xymatrix{ \mathcal{V}
\ar[r]^{\iota} \ar@(ul,dl)[]_d & \mathcal{V}_\iota
\ar@(ur,dr)[]^{d_\iota} }
\]
such that $d^2 = 0$, $d_\iota^2 = 0$ and $d\circ\iota
= -\iota\circ d$.
An $\iota$-complex, such as $\widetilde{\mathcal{W}}$, may be
viewed as a $(d,\iota)$-pair in the following way:
\[
\xymatrix{*{\parbox{1cm}{\ \ }\oplus_{i=0}^{\infty}\ncw^i\ }
\ar[rr]^(0.5){\oplus\iota} \ar@(ul,dl)[]_{\oplus d^i} & & *{\
\oplus_{i=0}^\infty\ncw_\iota^i\parbox{1cm}{\ }}
 \ar@(ur,dr)[]^{\oplus d^i_\iota}}\ .
\]
Write $\ncw = \oplus_{i=0}^{\infty}\ncw^i$ and $\ncw_\iota=
\oplus_{i=0}^\infty\ncw_\iota^i$.

\subsection{The roadmap}
Recall the statement of Homological Wheeling: \label{gettinghw}
Given elements
$v\in\Bspace^i$ and $w\in\Bspace^j$, there exists an element
$x_{v,w}$ of $\ncw^{i+j-1}$ such that $\iota(x_{v,w})=0$ and such
that
\[
\left(\chi_\Wspace\circ \Upsilon\right)^i(v)\#\left(
\chi_\Wspace\circ \Upsilon\right)^j(w) = \left(\chi_\Wspace\circ
\Upsilon\right)^{i+j}(v\sqcup w) + d(x_{v,w})\ \ \in\ \ncw^{i+j}.
\]
This statement holds in $\ncw$, whereas we want a statement in
$\Aspace$, the space of ordered Jacobi diagrams. Well, the
non-commutative space $\ncw$ has no relations which concern the
legs of diagrams. We can introduce whatever relations we can think
of, and see what consequences may be derived.
For example, in
attempting to derive a statement in $\Aspace$, we could introduce
the STU relations (among degree 2 legs), which is precisely what we will do in the next section:
\begin{equation}\label{introducedSTU}
\raisebox{-2ex}{\scalebox{0.22}{\includegraphics{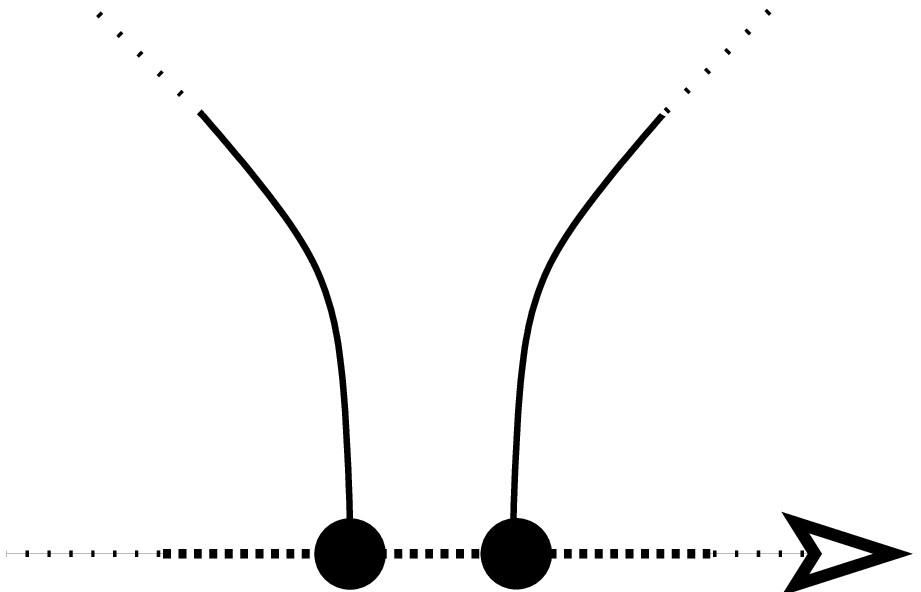}}}  -
 \raisebox{-2ex}{\scalebox{0.22}{\includegraphics{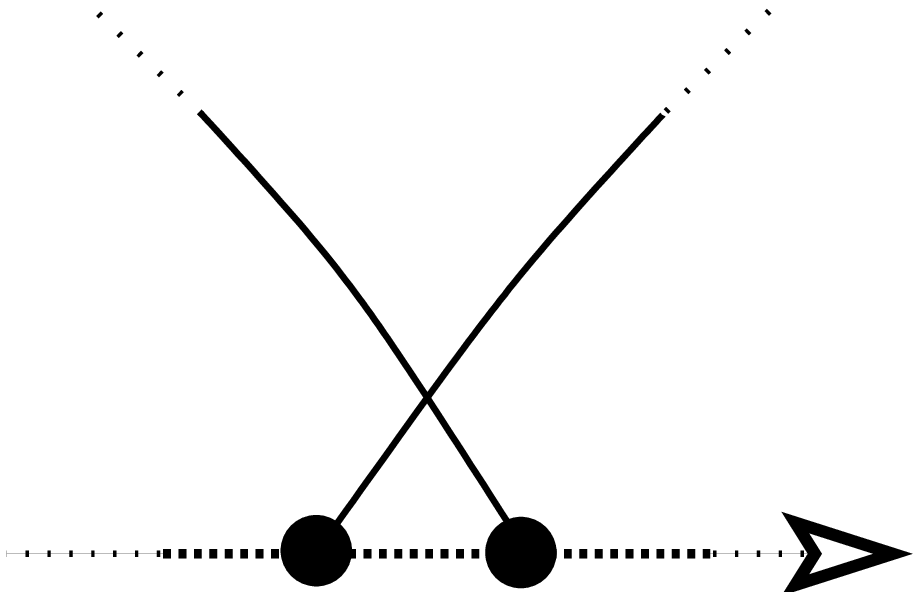}}}
=  \raisebox{-2ex}{\scalebox{0.22}{\includegraphics{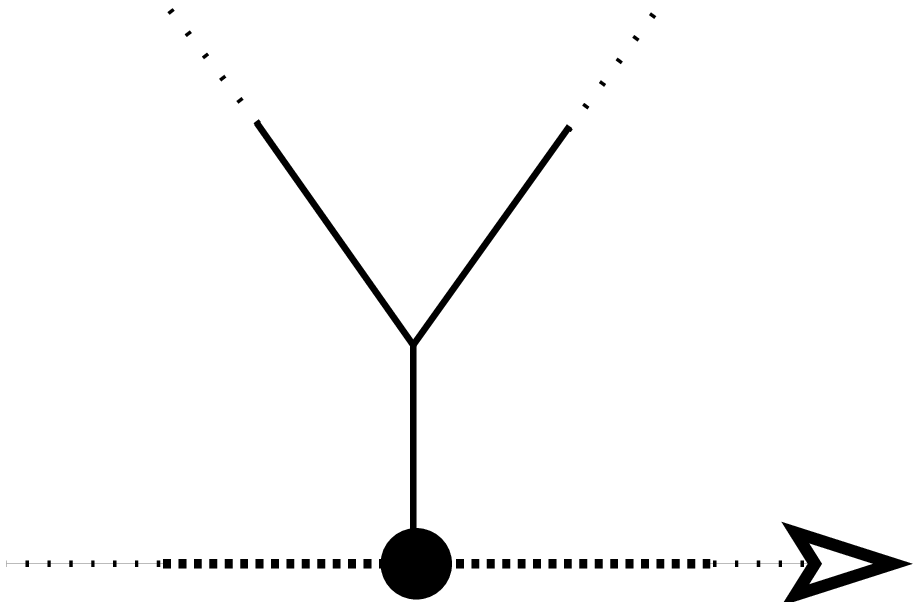}}}\
\ .
\end{equation}

To deduce the usual Wheeling Theorem in $\Aspace$, we
will take the following route:
\[
\ncw \xrightarrow[\pi]{\begin{array}{c}\mbox{Introducing} \\
\mbox{STU and other} \\ \mbox{relations.}\end{array}} \What
\xrightarrow[\basebulltoF]{\begin{array}{c}\text{Changing to} \\ \text{the ``curvatures"} \\
\text{basis.} \end{array}} \WhatF
\xrightarrow[\lambda]{\begin{array}{c} \text{Symmetrising} \\
\text{the grade $1$} \\ \text{legs.} \end{array}}
\Whatwedge=\Aspace\oplus\ldots\ \ .
\]

The $(d,\iota)$-pair $\WhatF$ is the combinatorial analogue of the complex $\mathcal{W}_G=U\mathfrak{g}\otimes \mathrm{Cl}\,g$ that the work \cite{AleksMein} refers to as the ``non-commutative Weil complex". This was the first structure discovered in the theory.

\subsection{The \(\di\)-pair \((\What,\What_\iota,d,\iota)\)}
The maps $d$ and $\iota$ play a crucial role in the statement of
HW. But simply introducing the STU relations leaves these maps
ill-defined.
So we must introduce a
number of other relations at the same time as we introduce STU, to
ensure that the maps $d$ and $\iota$ descend to the quotient
spaces.

Define $\widehat{\mathcal{W}}$ to be the vector space generated by
Weil diagrams modulo \ASreln, \IHXreln, the STU relation (Equation \ref{introducedSTU}), and also the following pair of relations:
\[
\begin{array}{ccccc}
\raisebox{-2ex}{\scalebox{0.22}{\includegraphics{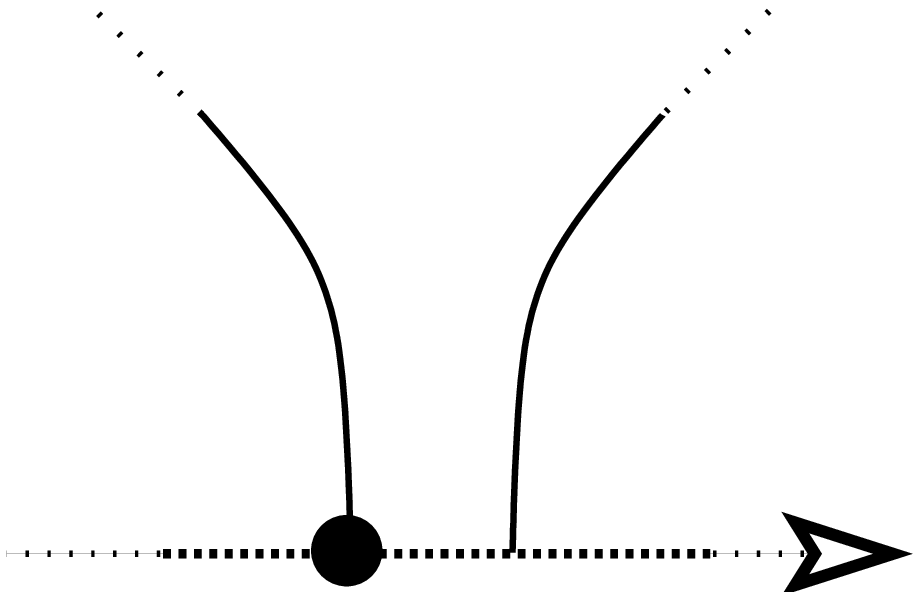}}} & -
& \raisebox{-2ex}{\scalebox{0.22}{\includegraphics{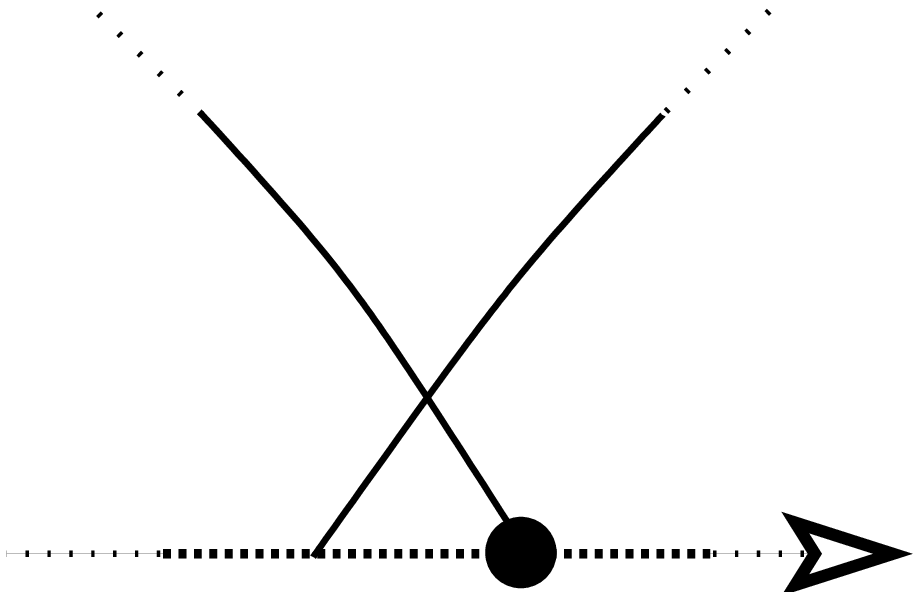}}} &
= &
\raisebox{-2ex}{\scalebox{0.22}{\includegraphics{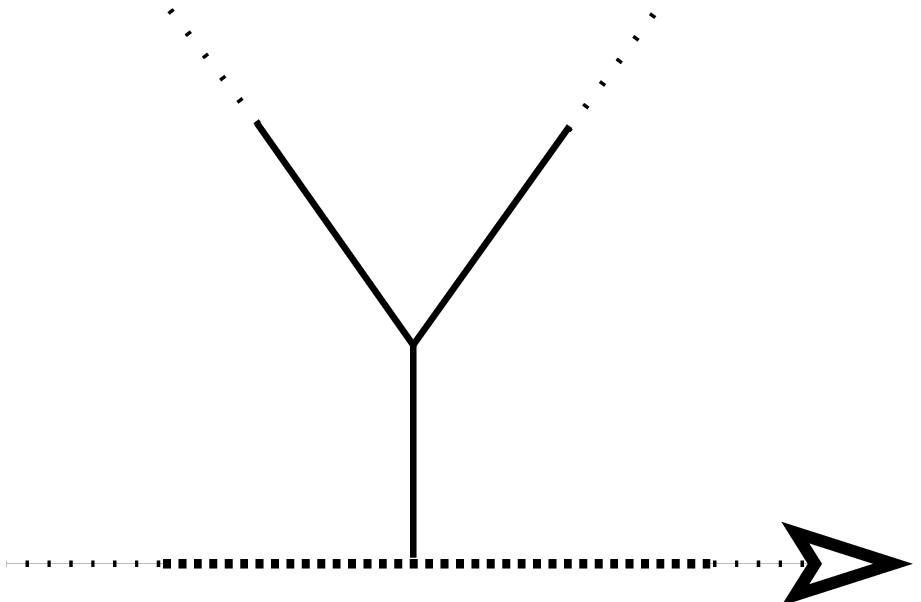}}}\ \
,
\\[0.55cm]
\raisebox{-2ex}{\scalebox{0.22}{\includegraphics{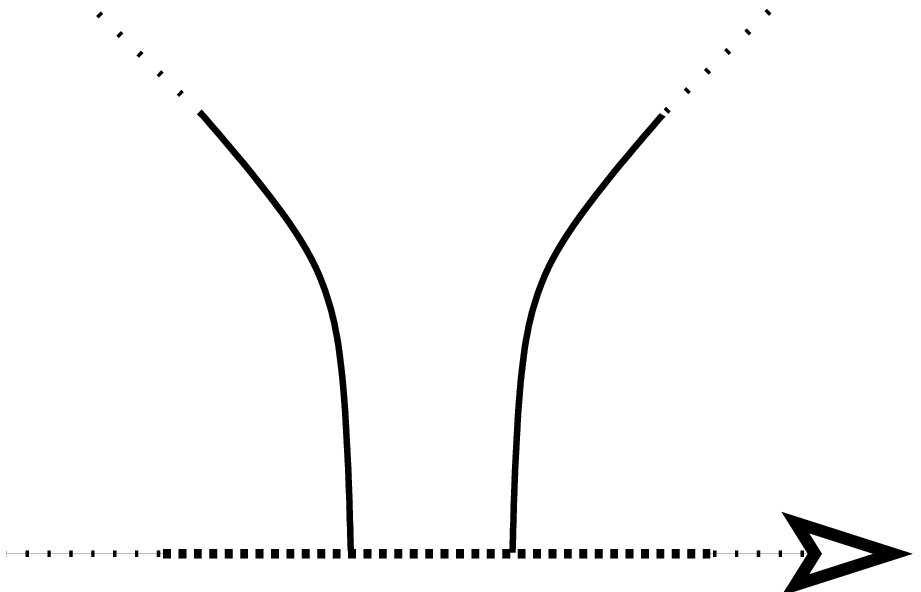}}} &
+ &
\raisebox{-2ex}{\scalebox{0.22}{\includegraphics{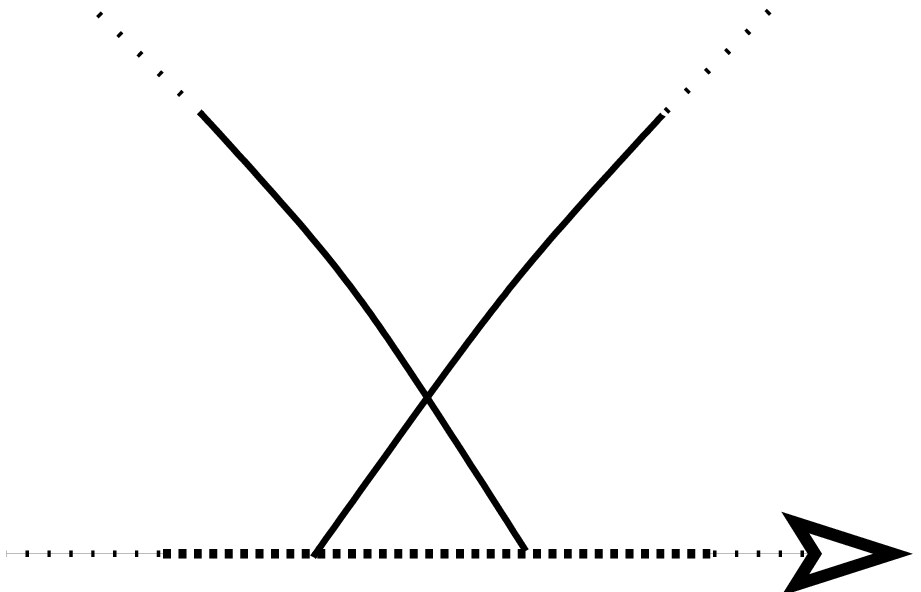}}} &
= &
\raisebox{-2ex}{\scalebox{0.22}{\includegraphics{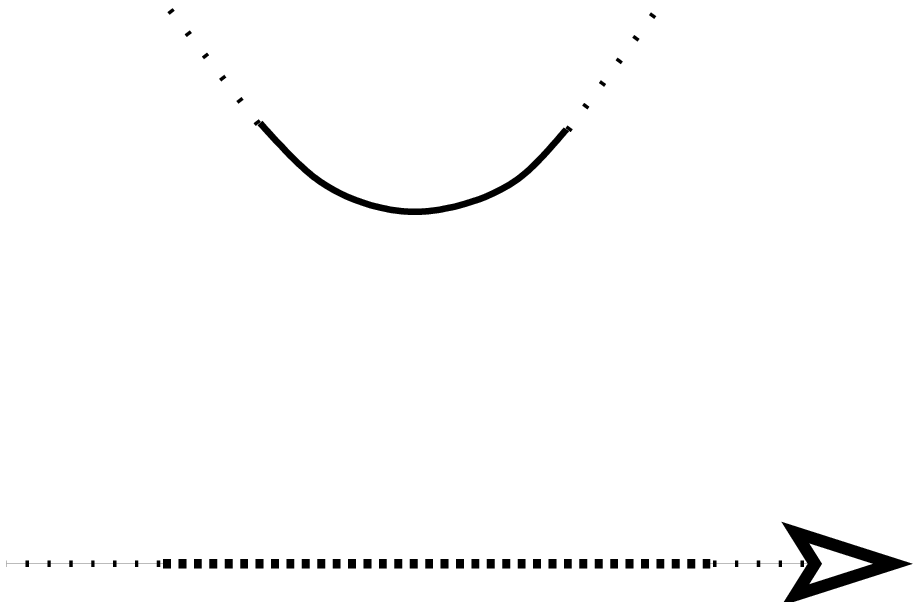}}}\ \
.
\end{array}
\]
Similarly define the vector space $\widehat{\mathcal{W}}_\iota.$

 Recall that the STU relation is a formal analogue of the relation in a universal enveloping algebra which equates a commutator with the corresponding bracket. And note further that the relation just introduced amongst grade 1 legs is a formal analogue of the defining relation of a {\bf Clifford algebra}.

\begin{prop}
The maps $d$ and $\iota$ descend to the spaces
$\widehat{\mathcal{W}}$ and $\widehat{\mathcal{W}}_\iota$. Thus,
quotienting gives a map of $\di$-pairs $
\pi : (\ncw,\ncw_\iota,d,\iota) \rightarrow
(\What,\What_\iota,d,\iota)$.
\end{prop}

\subsection{Some clarifications concerning the diagrams in this section.}\label{loopclarification}
The Weil diagrams generating the spaces introduced in this section will be allowed to have components consisting of closed loops without vertices. (Even though such components don't fit into the usual definition of a graph, it is straightforward to set this generality up formally.) The reason we need such generality is that the `Clifford algebra' relation may produce such components when applied to the two ends of a trivial chord with two grade 1 legs.\footnote{A referee alerted the author to this detail.}

We remark that this is only a formal device we employ to make our definitions logically consistent; the expressions that actually arise when a HW statement is inserted into this composition never use such components.

In addition: let $\Aspace_l$ be constructed in exactly the same way as the space $\Aspace$ but where we allow the ordered Jacobi diagrams to have components consisting of closed loops without vertices. There are no relations between these components and the rest of the graph. Observe that $\Aspace$ is obviously a subspace of $\Aspace_l$.

\subsection{A change of basis: The \(\di\)-pair $\widehat{\mathcal{W}}_{\mathrm{F}}$.}
\label{theamncwc}
We'll now change back to the ``curvatures" basis. Here it has the
salutary effect of making the two kinds of legs invisible to each
other. To be precise: define $\widehat{\mathcal{W}}_\mathrm{F}$ to
be the vector space generated by {\bf F}-Weil diagrams
taken modulo \ASreln, \IHXreln, and the following three classes of
relations:
\[
\begin{array}{ccccc}
\raisebox{-2ex}{\scalebox{0.22}{\includegraphics{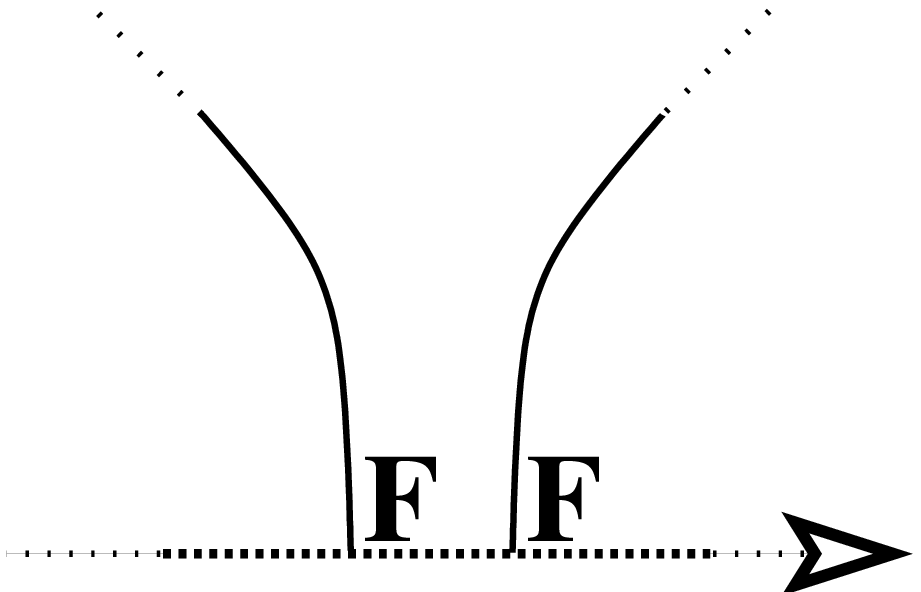}}} & - &
\raisebox{-2ex}{\scalebox{0.22}{\includegraphics{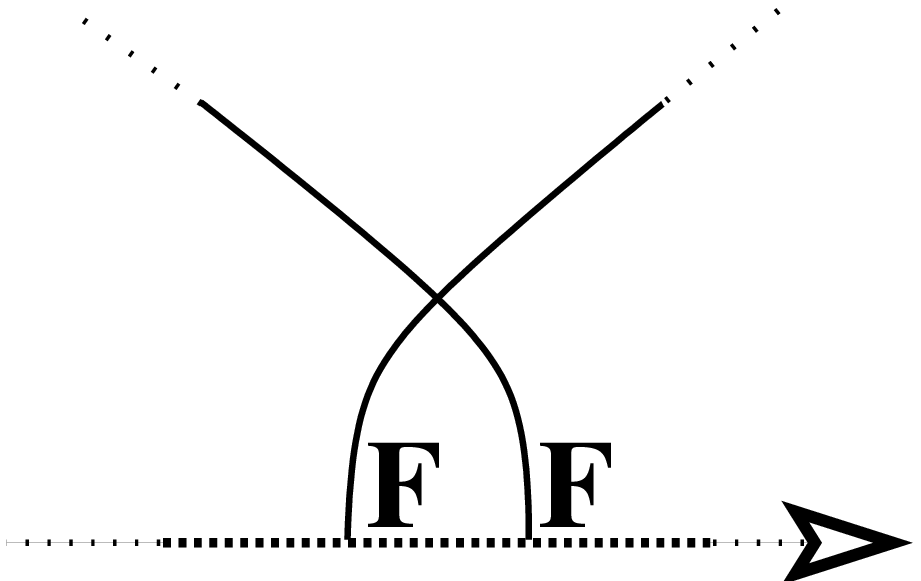}}} & = &
\raisebox{-2ex}{\scalebox{0.22}{\includegraphics{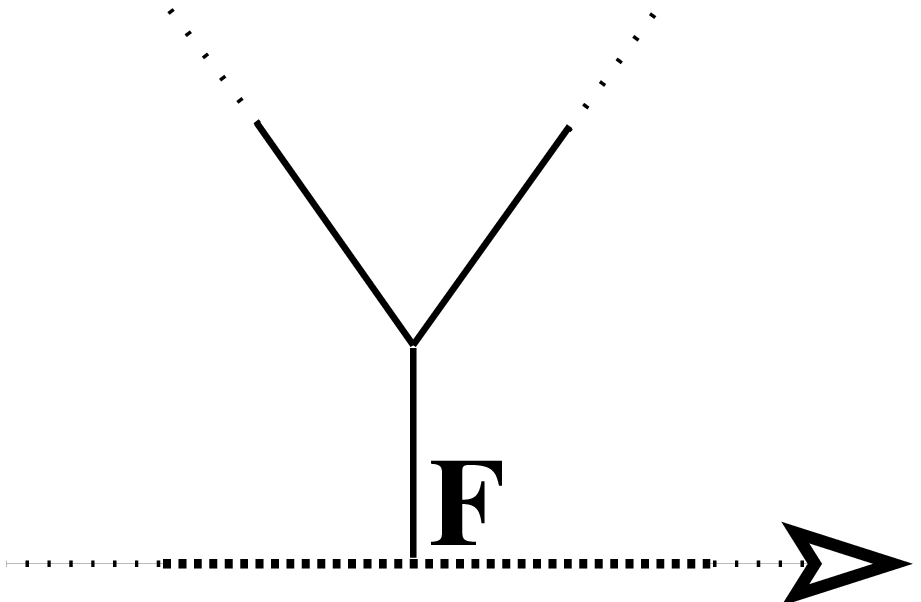}}}\ \ ,
\\[0.55cm]
\raisebox{-2ex}{\scalebox{0.22}{\includegraphics{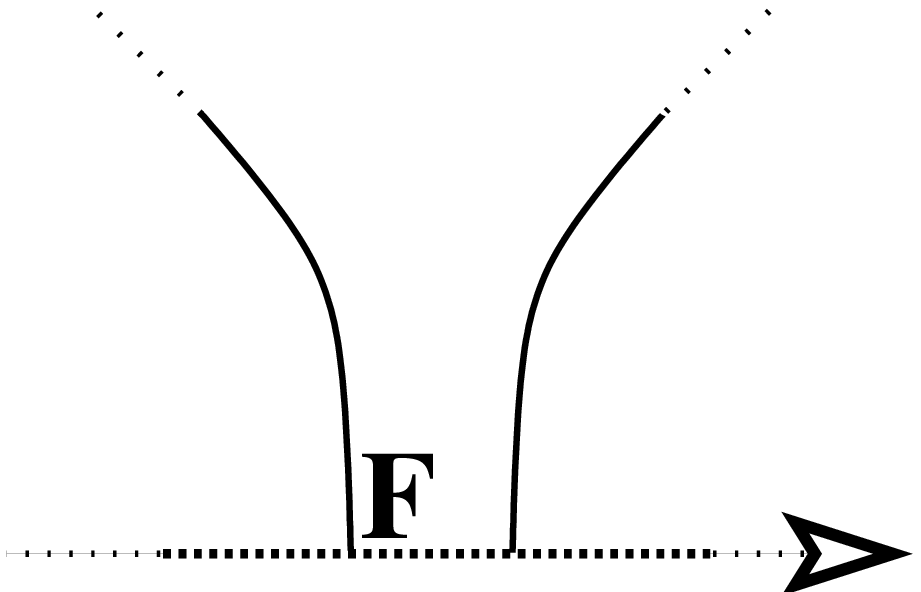}}} &  -
& \raisebox{-2ex}{\scalebox{0.22}{\includegraphics{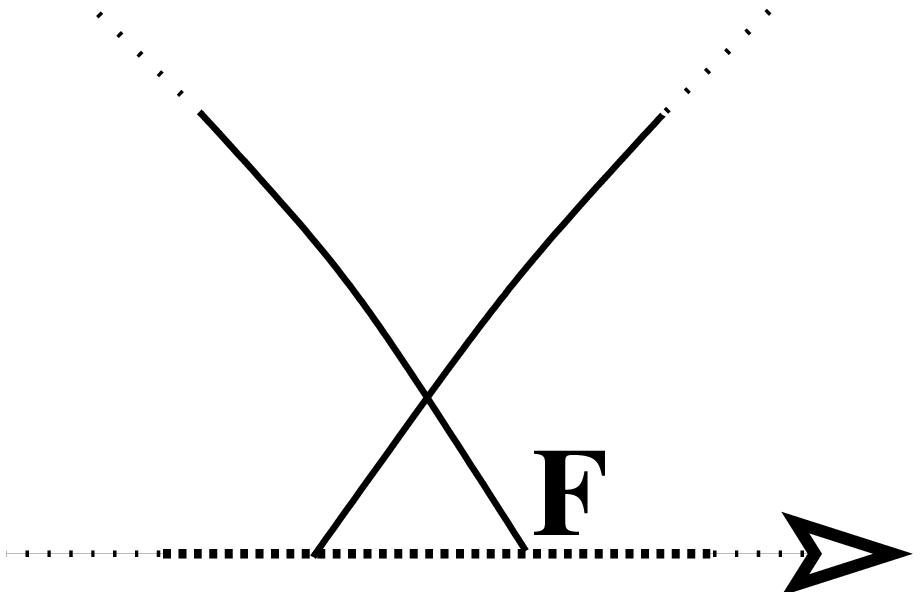}}} & =
& 0\ \ ,
\\[0.55cm]
\raisebox{-2ex}{\scalebox{0.22}{\includegraphics{whatrelnAPP}}} &
+ &
\raisebox{-2ex}{\scalebox{0.22}{\includegraphics{whatrelnBPP}}} &
= &
\raisebox{-2ex}{\scalebox{0.22}{\includegraphics{whatrelnCPP}}}\ \
.
\end{array}
\]
Similarly define $\widehat{\mathcal{W}}_{\mathrm{F}\iota}.$ If we
equip this pair of spaces with maps $d$ and $\iota$ using the
usual substitution rules for {\bf F}-Weil diagrams then we have a
$\di$-pair. (See Section \ref{changebasistoF} for the reason that
$[d,\iota]=0$. And the calculation that $d^2=0$ is a little
surprising in this basis, as we invite the reader to see.) Define
maps $B_{\bullet\rightarrow F}:\What\rightarrow \WhatF$ and
$B_{{\bullet\rightarrow F},\iota}: \Whatiota\rightarrow
\WhatFiota$ by expanding grade $2$ legs in the usual way:
\[
\raisebox{-3ex}{\scalebox{0.25}{\includegraphics{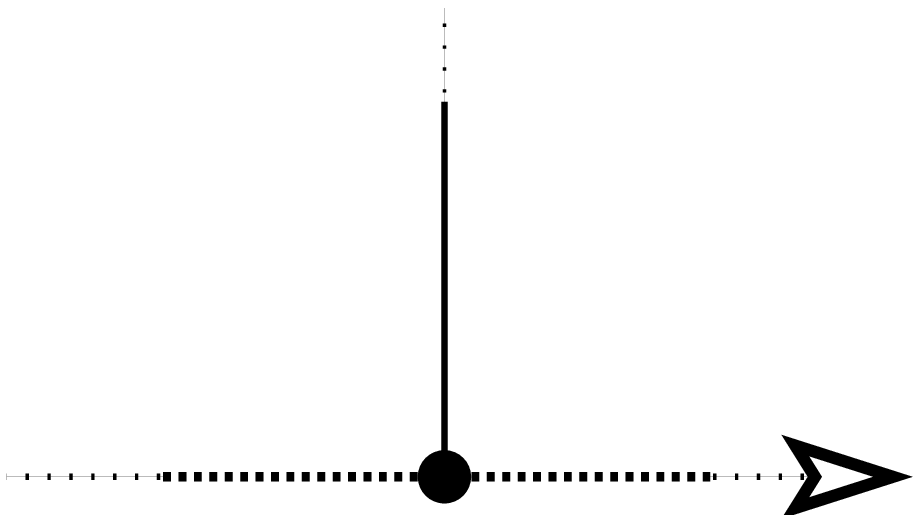}}}\
\mapsto
\raisebox{-3ex}{\scalebox{0.25}{\includegraphics{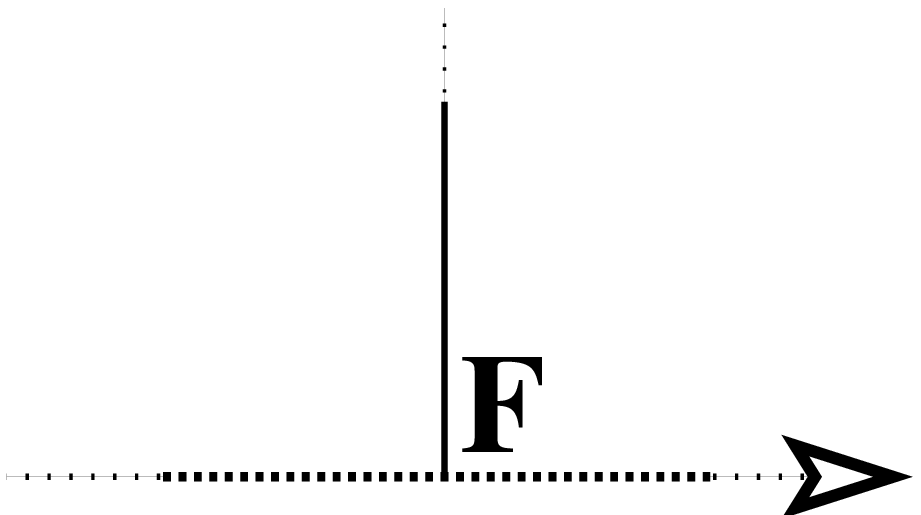}}}\ \
+\frac{1}{2}\,
\raisebox{-3ex}{\scalebox{0.25}{\includegraphics{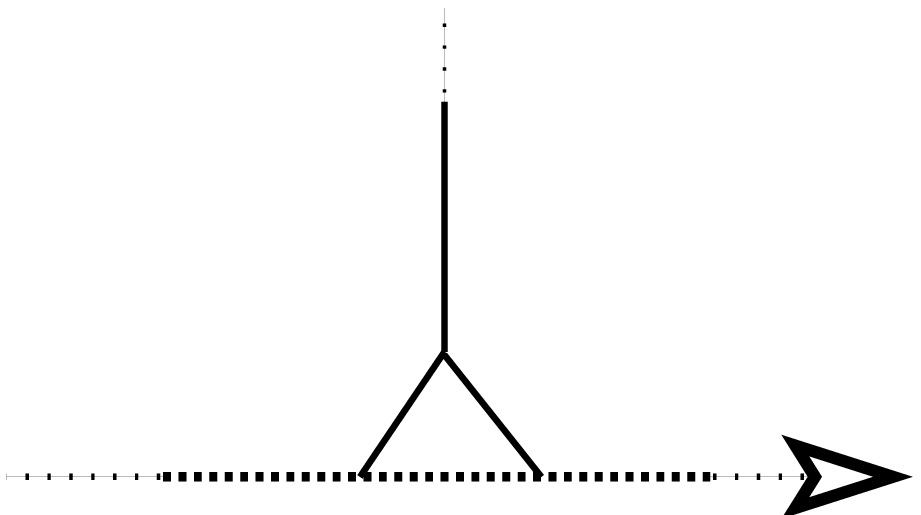}}}\ .
\]
These maps give a map of $(d,\iota)$-pairs $B_{\bullet\rightarrow F}:
(\widehat{\mathcal{W}},\widehat{\mathcal{W}}_\iota,d,\iota)\rightarrow
(\widehat{\mathcal{W}}_\mathrm{F},\widehat{\mathcal{W}}_\mathrm{F\iota},d,\iota)$.
\subsection{Symmetrizing the grade 1 legs.}

The final step on our journey to $\Aspace$ is to (graded)
symmetrize the grade 1 legs. Intuitively speaking, we are just
taking a special choice of basis in $\WhatF$ (and so in $\What$).
Formally: we'll define new vector spaces $\Whatwedge$ and
$\Whatwedgeiota$ and equip them with ``averaging" maps into
$\WhatF$ and $\WhatFiota$.

The reason to symmetrize is that this final space $\Whatwedge$
admits a direct-sum decomposition $\Whatwedge = \bigoplus_{i=0}^{\infty}\Whatwedge^i$
according to the number of grade 1 legs in a diagram.  We'll
locate $\Aspace_l$ in the grade 0 position of this decomposition:
$\Aspace_l \cong \Whatwedge^0$.
Finally, we'll show that $\ker \iota = \Whatwedge^0\cong \Aspace_l$.

\subsubsection{The vector spaces \(\widehat{\mathcal{W}}_{\wedge}\) and
\(\widehat{\mathcal{W}}_{\wedge\iota}\)}
 Define
$\widehat{\mathcal{W}}_{\wedge}$ to be the vector space generated
by {\bf F}-Weil diagrams modulo AS, IHX, STU (amongst F-legs), and the following two
classes of relations:
\[
\begin{array}{ccccc}
\raisebox{-2ex}{\scalebox{0.23}{\includegraphics{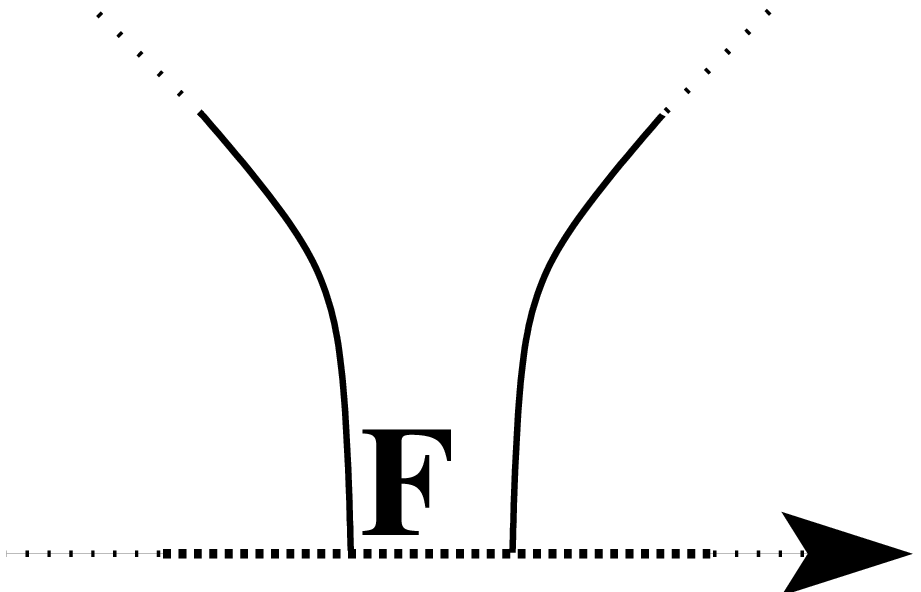}}} &
- &
\raisebox{-2ex}{\scalebox{0.23}{\includegraphics{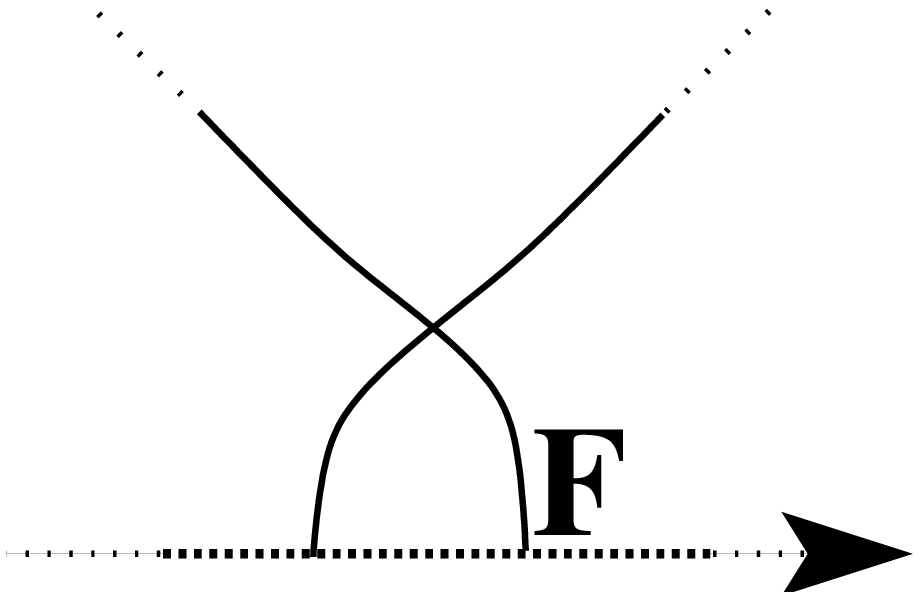}}} &
= & 0\ ,
\\[0.5cm]
\raisebox{-2ex}{\scalebox{0.23}{\includegraphics{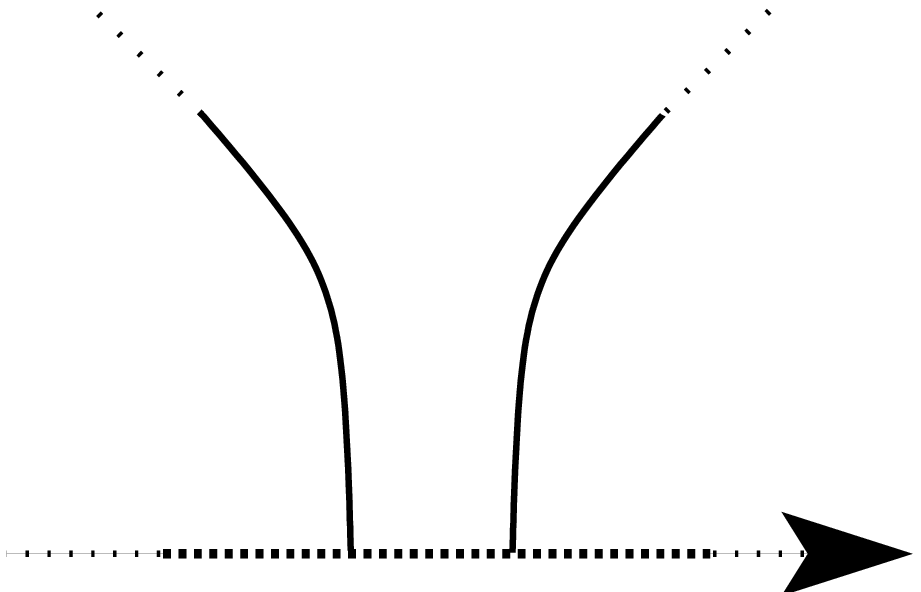}}} &
+ &
\raisebox{-2ex}{\scalebox{0.23}{\includegraphics{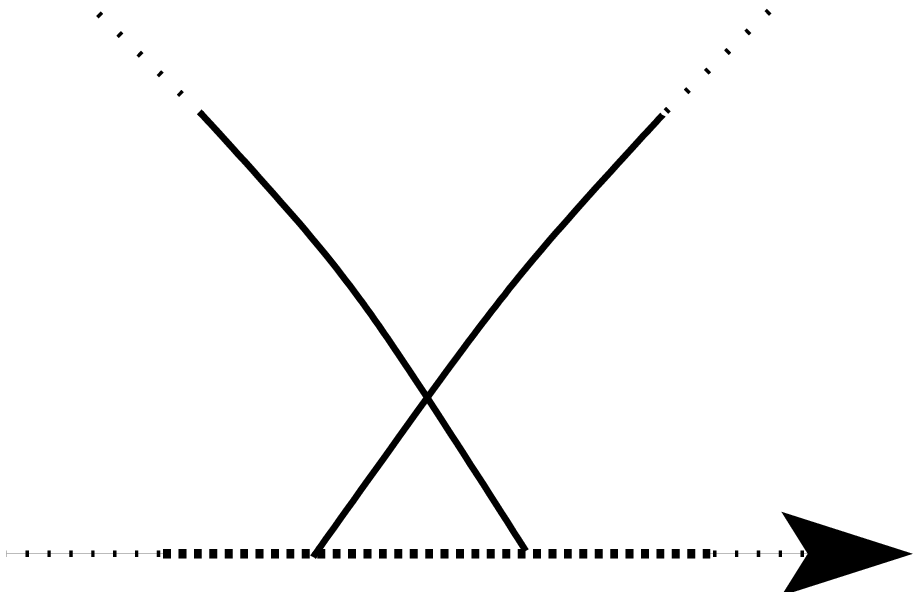}}} &
= & 0\ .
\end{array}
\]
 Define
$\widehat{\mathcal{W}}_{\wedge\iota}$ similarly. In these spaces,
then, the grade 1 legs may be moved about freely (up to sign).

We remark, as discussed in Section \ref{loopclarification}, that we are formally allowing
closed loop components. In $\widehat{\mathcal{W}}_{\wedge}$ there are no relations that involve such components at all.

We will {\bf not} introduce a differential into these spaces. (The
reason is that, in $\WhatF$, when the differential acts on a grade
1 leg, it produces terms with new grade 1 legs. So the action of
the differential in the ``symmetrized" basis will look very
complicated.)

We will only define the map $\iota :
\widehat{\mathcal{W}}_{\wedge}\rightarrow\widehat{\mathcal{W}}_{\wedge\iota}$.
This is defined as a formal linear differential operator, in the
usual way. The following theorem will be established in some
detail in Section \ref{bigwedgetheoremproof}:
\begin{thm}\label{bigwedgetheorem}
There exist vector space isomorphisms $\chi_\wedge :
\widehat{\mathcal{W}}_{\wedge} \rightarrow
\widehat{\mathcal{W}}_{\mathrm{F}}$  and
$\chi_{\wedge\iota} : \widehat{\mathcal{W}}_{\wedge\iota}
\rightarrow \widehat{\mathcal{W}}_{\mathrm{F}\iota}$
commuting with the action of $\iota$.
%
\end{thm}

\begin{figure}
\caption{An example of the map $\chi_\wedge$.
\label{chiwedgeexamp}}
\parbox{12cm}{
\begin{eqnarray*}
\chi_\wedge\left(
\raisebox{-2.75ex}{\scalebox{0.25}{\includegraphics{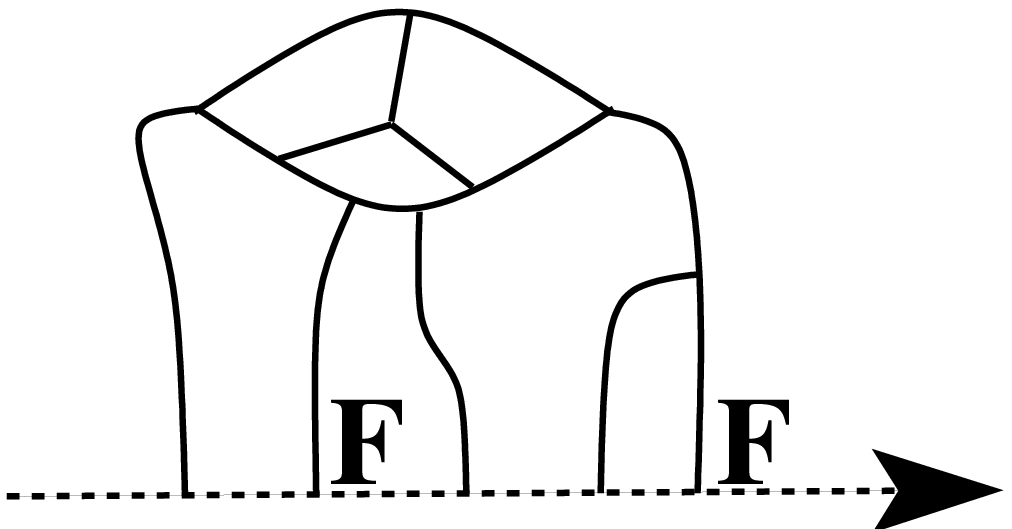}}}
\right) & = &\frac{1}{3!}\left(
\raisebox{-2.75ex}{\scalebox{0.25}{\includegraphics{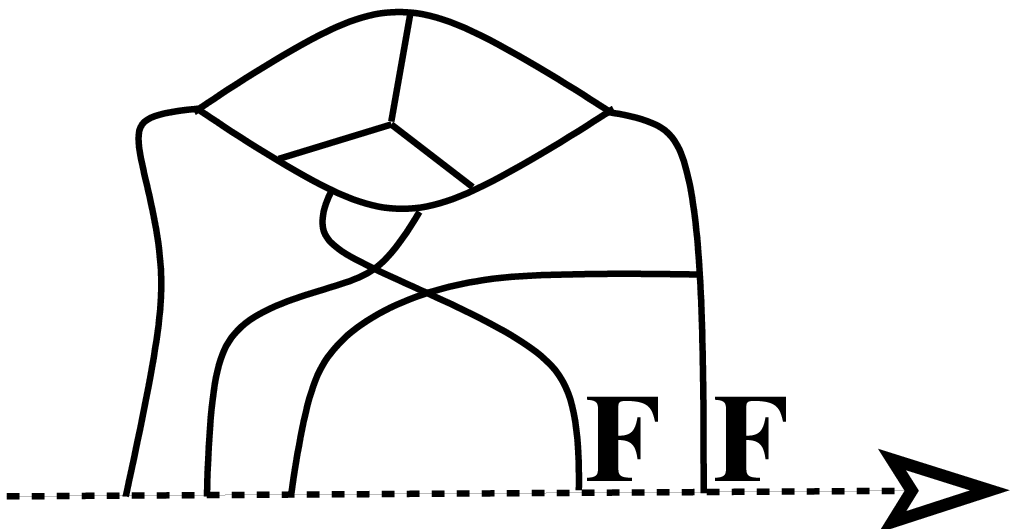}}}
\ -\
\raisebox{-2.75ex}{\scalebox{0.25}{\includegraphics{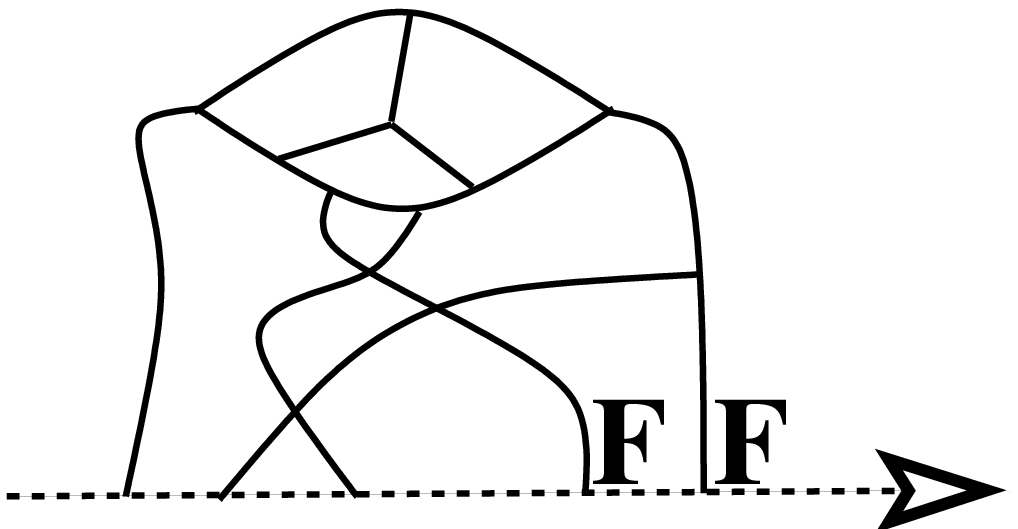}}}\right.
\\
& & \ +\
\raisebox{-2.75ex}{\scalebox{0.25}{\includegraphics{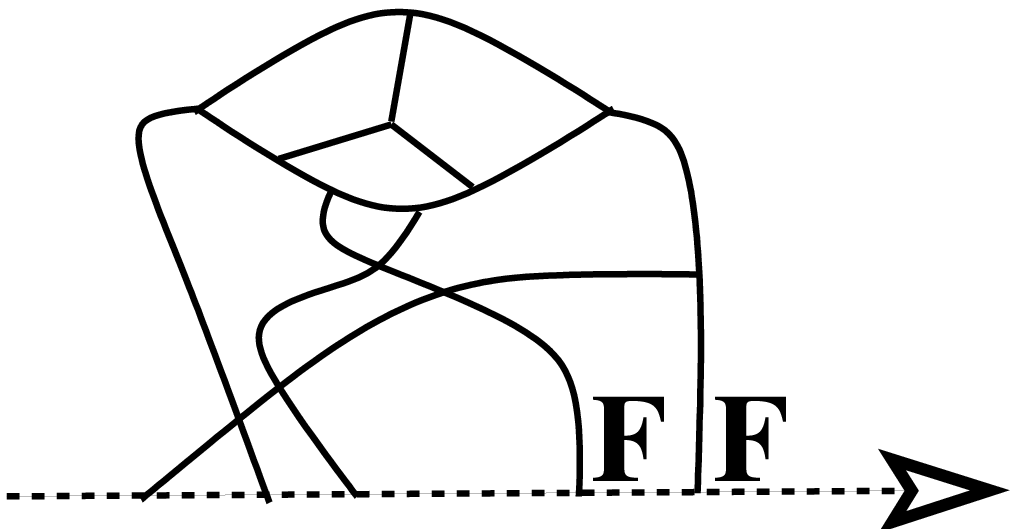}}}
\ -\
\raisebox{-2.75ex}{\scalebox{0.25}{\includegraphics{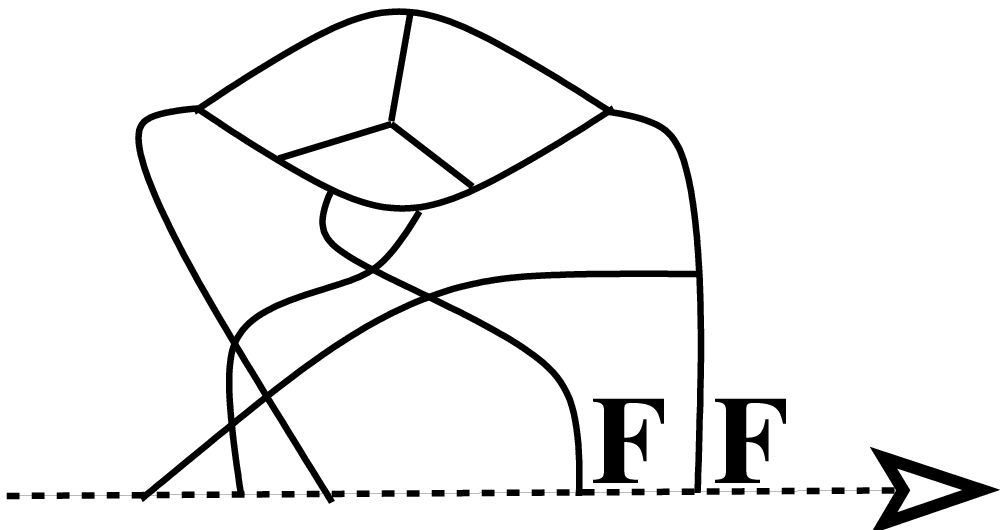}}}
\\[0.4cm]
& &\left. \ +\
\raisebox{-2.75ex}{\scalebox{0.25}{\includegraphics{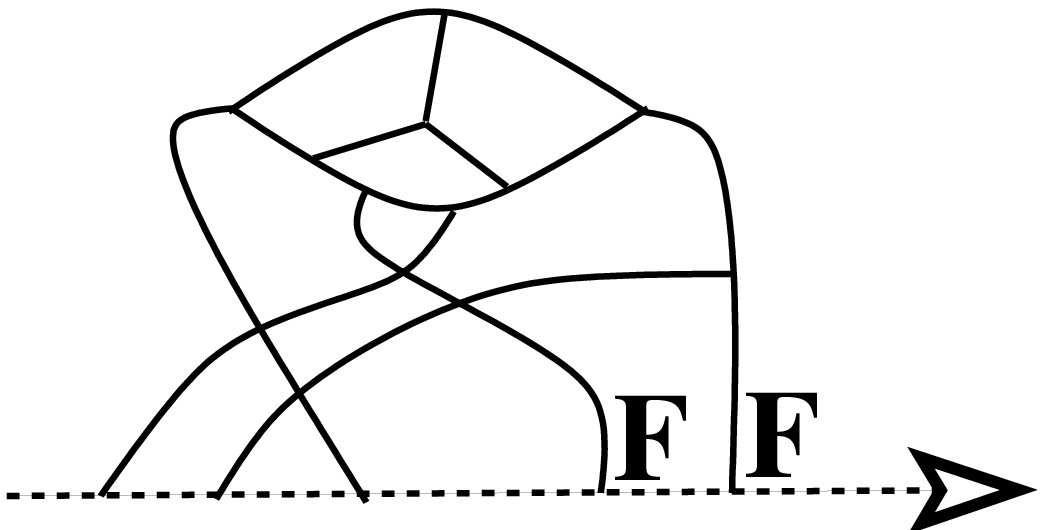}}}
\ -\
\raisebox{-2.75ex}{\scalebox{0.25}{\includegraphics{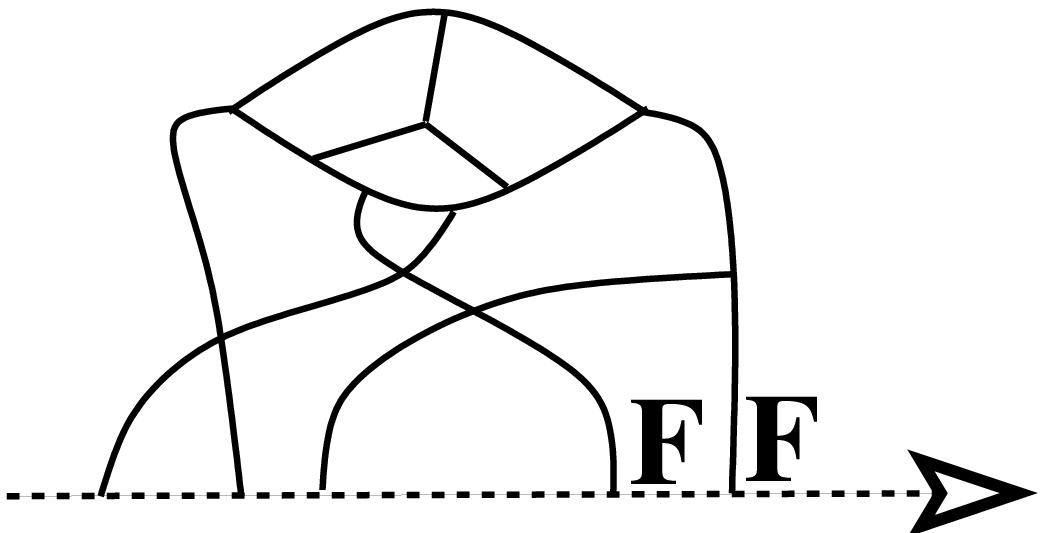}}}\right)\
.
\end{eqnarray*}
} \underline{\hspace{7cm}}
\end{figure}

 The map $\chi_\wedge :
\widehat{\mathcal{W}}_{\wedge} \rightarrow
\widehat{\mathcal{W}}_{\mathrm{F}}$ is defined by declaring the
value of the map on some diagram to be the average of the terms
that can be obtained by doing (signed) permutations of the degree
1 legs of the diagram. $\chi_{\wedge\iota}$ is defined similarly.
An example is shown in Figure \ref{chiwedgeexamp}.

That the maps commute with $\iota$ is straightforward.
The difficult part of the theorem is showing that $\chi_\wedge$ is
a vector space isomorphism. This will be done by explicitly
constructing an inverse.
\begin{defn}
The inverse to $\chi_{\wedge}$ will be denoted $\lambda: \WhatF \to
\Whatwedge.$ \end{defn} Similarly we have $\lambda_\iota$. This
map $\lambda$ is a crucial map which will be constructed
explicitly later, in Section \ref{lambdaconstruct}. After
constructing $\lambda$ there we'll go through the details of
Theorem \ref{bigwedgetheorem}. For the present we'll just assume
that such maps exist.

\subsubsection{Decomposing according to the number of grade $1$ legs.}
Observe that there are no relations in $\Whatwedge$ which involve
diagrams with differing numbers of the degree 1 legs. Thus the
space $\Whatwedge$ is graded by that quantity. Let $\Whatwedge^i$
denote the subspace of $\Whatwedge$ generated by diagrams which
have exactly $i$ grade 1 legs. We have a direct-sum decomposition into these subspaces:
$\Whatwedge = \bigoplus_{i=0}^{\infty} \Whatwedge^i\ $.
Observe now that there is a canonical isomorphism (``label every
leg with an {\bf F}") from  $\Aspace_l$, the familiar space of
ordered Jacobi diagrams (recall the $l$ subscript allows the possibility of closed loop components), to $\Whatwedge^0$. We'll denote it $\phi_\Aspace : \Aspace_l \xrightarrow{\cong} \Whatwedge^0$.
For example:
\[
\phi_\Aspace\left(
\raisebox{-4.5ex}{\scalebox{0.22}{\includegraphics{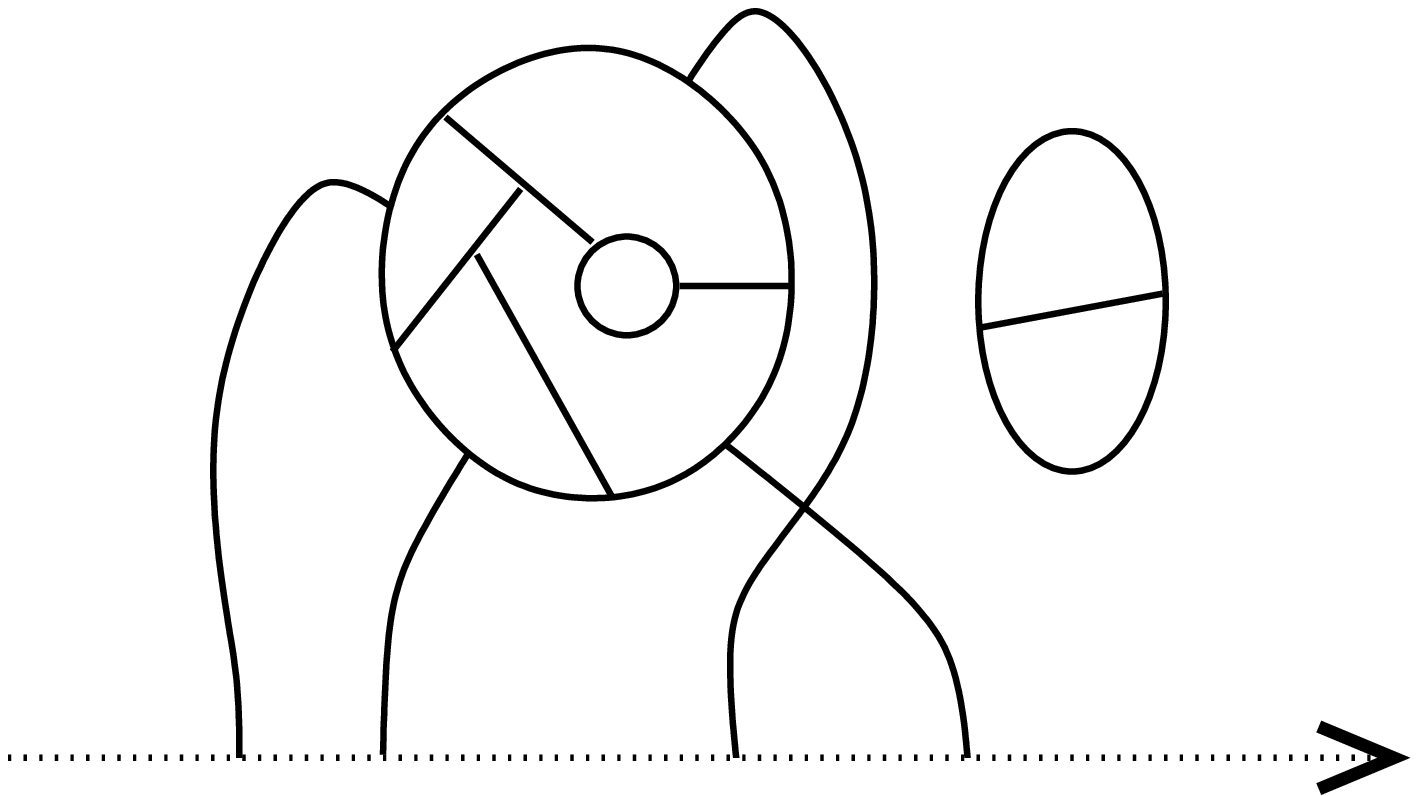}}}
\right) \ =\
\raisebox{-4.5ex}{\scalebox{0.22}{\includegraphics{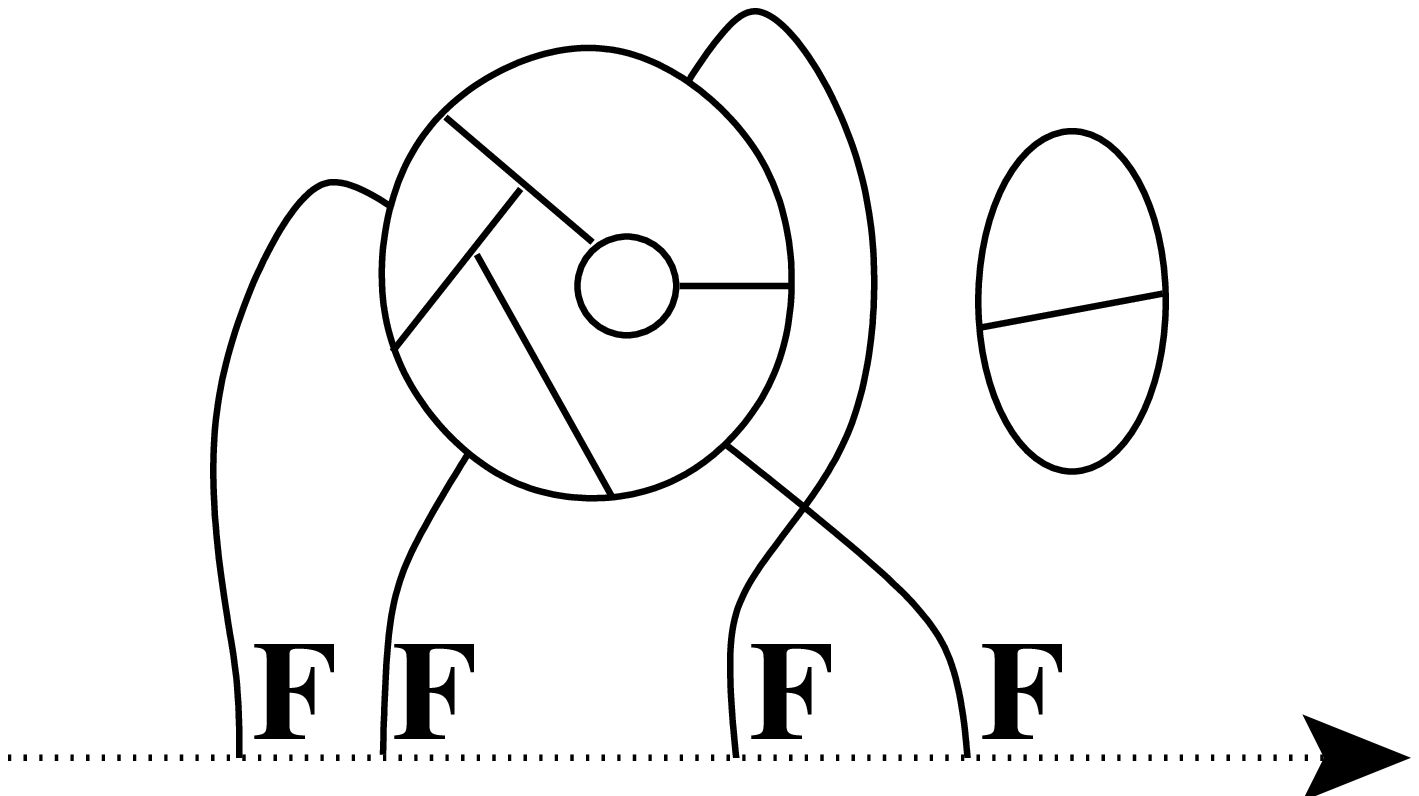}}}\ \ .
\]

\subsubsection{The kernel of $\iota$ on $\Whatwedge$.}
Here we will establish that the kernel of the map $\iota$ is precisely the zeroth summand in this decomposition:
\begin{equation}\label{keriotanolegs}
\ker \iota = \Whatwedge^0 \subset \Whatwedge.
\end{equation}
This fact is a consequence of the map
$\hat{\iota}:\Whatwedgeiota\to\Whatwedge$ that is defined by
placing the special $\iota$-labelled leg on the far left-hand end
of the orienting line.
Because, for some generator $v\in\Whatwedge$,
\[
\left(\hat{\iota}\circ\iota\right)(v) = \left(\begin{array}{c} \text{Number of} \\
\text{grade 1 legs of $v$} \end{array} \right)v\ ,
\]
an element of $\Whatwedge$ lies in the kernel of
$\iota$ if and only if it lies in $\Whatwedge^0$.

\subsection{Wheeling - the key lemmas.}
Now that we have introduced the spaces involved, we can turn to
the explanation of how Wheeling arises from HW. Recall the
statement of HW: Given elements $v\in\Bspace^i$ and
$w\in\Bspace^j$, there exists an element $x_{v,w}$ of
$\ncw^{i+j-1}$ such that $\iota(x_{v,w})=0$ and such that
\[
\left(\chi_\Wspace\circ \Upsilon\right)^i(v)\#\left(
\chi_\Wspace\circ \Upsilon\right)^j(w) = \left(\chi_\Wspace\circ
\Upsilon\right)^{i+j}(v\sqcup w) + d(x_{v,w})\ \ \in\ \ncw^{i+j}.
\]
We'll now take the three summands of this equation and
insert them into the composition of linear maps
\[
\ncw \xrightarrow{\pi} \What \xrightarrow{\basebulltoF} \WhatF
\xrightarrow{\lambda} \Whatwedge\ ,
\]
so as to produce a statement in $\Whatwedge$.

To begin, note that each of the three summands
lies in the
kernel of $\iota_{\ncw}$. Thus, each of their images under the map
$\left(\lambda\circ\basebulltoF\circ\pi\right)$ lies in the kernel
of $\iota_{\Whatwedge}$, which we know is isomorphic to $\Aspace_l$, via the linear isomorphism $\phi_\Aspace$.

So we move on to asking: if we pull these elements back to
$\Aspace_l$, what do we get? The following, crucial, lemma says that
we just get {\bf The Wheeling Map} composed with the averaging map
from $\Bspace$ to $\Aspace\subset \Aspace_l$ (just as appears in the Wheeling
Theorem). The proof is a quite involved (but conceptually straightforward) combinatorial calculation which will
appear in a separate publication \cite{K}.
\begin{lem}\label{howwheelsappearlem}
For all $x\in \Bspace$,
\[
\left(\phi_{\Aspace}^{-1}\circ\lambda\circ\basebulltoF\circ\pi\right)\left(\left(\chi_{\Wspace}\circ\Upsilon\right)(x)\right)
=\left(\chi_{\Bspace}\circ\partial_{\Omega}\right)(x)\ .
\]
\end{lem}
The next lemma is a quick corollary, as will be explained at the
end of this subsection:
\begin{lem}\label{smalllem}
The element
$\left(\phi_{\Aspace}^{-1}\circ\lambda\circ\basebulltoF\circ\pi\right)
\left(\left(\chi_{\Wspace}\circ\Upsilon\right)(v)\#
\left(\chi_{\Wspace}\circ\Upsilon\right)(w)\right)$ is equal to
$\left(\chi_{\Bspace}\circ\partial_{\Omega}\right)(v) \#
\left(\chi_{\Bspace}\circ\partial_{\Omega}\right)(w)$.
\end{lem}

Assuming these lemmas, the only thing which stands in
the way of the Wheeling Theorem is the error term:
$d(x_{v,w})$, where $\iota(x_{v,w})=0$.
This error term is dispatched by the following lemma.
\begin{lem}\label{dispatchlem}
If $z\in \ncw$ is such that $\iota_{\ncw}(z)=0$, then
\[
\left(\basebulltoF\circ\pi\right)\left(\,d(z)\,\right) = 0\ \in
\WhatF .
\]
\end{lem}
\begin{proof}
We'll start by establishing a small claim. We claim that if an
element $w\in \WhatF$ is in the kernel of $\iota_{\WhatF}$, then
it can be expressed as a linear combination of diagrams each of
whose legs is an {\bf F}-labelled leg. The deduction is that:
\begin{itemize}
\item{$\lambda(w)$ is in the kernel of $\iota_{\Whatwedge}$
(because $\lambda$ commutes with $\iota$, by Theorem
\ref{bigwedgetheorem}).} \item{Thus $\lambda(w)$
lies in $\Whatwedge^0$ (by Equation \ref{keriotanolegs}).}
\item{Thus $\lambda(w)$ is a linear combination of diagrams all of
whose legs are {\bf F}-labelled.} \item{Thus
$w=\chi_{\Whatwedge}\left(\lambda(w)\right)$ must also consist of
entirely-{\bf F}-labelled diagrams. }
\end{itemize}
This establishes the claim. Then:
\begin{eqnarray*}
\left(\basebulltoF\circ\pi\right)\left(\,d(z)\,\right) & = &
d\left(\left(\basebulltoF\circ\pi\right)(z)\right) \\
& = & d\left(
\begin{array}{c}
\text{A combination of diagrams in $\WhatF$} \\
\text{every leg of which is labelled by an {\bf F}.}
\end{array} \right)
\\
& = & 0.
\end{eqnarray*}
The first equality follows because the maps $\basebulltoF$ and
$\pi$ are both maps of $(d,\iota)$-pairs. The second equality
follows from the claim that started this proof (noting that
$\iota\left(\left(\basebulltoF\circ\pi\right)(z)\right)=0$).
The
third equality we leave as an exercize. (Compare with the proof of
Proposition \ref{phiisaniotamap}.)
\end{proof}

{\it Proof of Lemma \ref{smalllem} using Lemma
\ref{howwheelsappearlem}.} Consider the expression
that is the subject of this lemma.
Briefly considering the constructions of the maps $\basebulltoF$
and $\pi$, we realize that we can write it:
\[
\left(\phi_{\Aspace}^{-1}\circ\lambda\right)
\left(\left(\basebulltoF\circ\pi\circ\chi_{\Wspace}\circ\Upsilon\right)(v)\#
\left(\basebulltoF\circ\pi\circ\chi_{\Wspace}\circ\Upsilon\right)(w)\right).
\]
Now, just as in the claim used in the proof of Lemma
\ref{dispatchlem}, we observe that because the factors
$\left(\basebulltoF\circ\pi\circ\chi_{\Wspace}\circ\Upsilon\right)(v)$
and
$\left(\basebulltoF\circ\pi\circ\chi_{\Wspace}\circ\Upsilon\right)(w)$
are both in the kernel of $\iota_{\WhatF}$, they can be expressed
using diagrams whose legs are entirely ${\bf F}$-legs. The map
$\lambda$ does not do anything to such diagrams. Thus, the
expression can be written:
$\left(\phi_{\Aspace}^{-1}\circ\lambda\circ\basebulltoF\circ\pi\circ\chi_{\Wspace}\circ\Upsilon\right)(v)\#
\left(\phi_{\Aspace}^{-1}\circ\lambda\circ\basebulltoF\circ\pi\circ\chi_{\Wspace}\circ\Upsilon\right)(w).
$
\begin{flushright}
$\Box$
\end{flushright}

\subsection{The map \(\lambda\).}\label{lambdadefnsect}
\label{lambdaconstruct} The remaining business of this paper is to
prove Theorem \ref{bigwedgetheorem}.
This subsection will define maps
$\lambda$ and $\lambda_\iota$,
prove that they are well-defined,
and prove that they left-invert the averaging maps: $\lambda\circ\chi_\wedge=\mathrm{id}_{\widehat{\mathcal{W}}_\wedge}$
and
$\lambda_\iota\circ\chi_{\wedge\iota}=\mathrm{id}_{\widehat{\mathcal{W}}_{\wedge\iota}}$.
The next and final subsection will combine these facts with the
easy facts that $\chi_{\wedge}$ and $\chi_{\wedge\iota}$ are
surjective, and that $\iota$ commutes with the $\chi_\wedge$, to
obtain the theorem.

The constructions of $\lambda$ and $\lambda_\iota$ are quite
natural and easily remembered: ``glue chords (with certain easily
determined signs and coefficients) into the grade 1 legs in all
possible ways".

\subsubsection{The formal definition of $\lambda$.}
\begin{figure}
\caption{An example of the action of a pairing.
\label{pairingexamp}}
\parbox{12cm}{
\begin{multline*}
\mathcal{D}_{\{\{2,4\},\{1,3\}\}}\left(
\raisebox{-6ex}{\scalebox{0.23}{\includegraphics{lambdaopexampA}}}
\right) =  \mathcal{D}_{\{2,4\}}\left(\mathcal{D}_{\{1,3\}}
\left(
\raisebox{-6ex}{\scalebox{0.23}{\includegraphics{lambdaopexampA}}}
\right) \right) \\[0.25cm]
  =  -\frac{1}{2}\ \mathcal{D}_{\{2,4\}}\left(
\raisebox{-6ex}{\scalebox{0.23}{\includegraphics{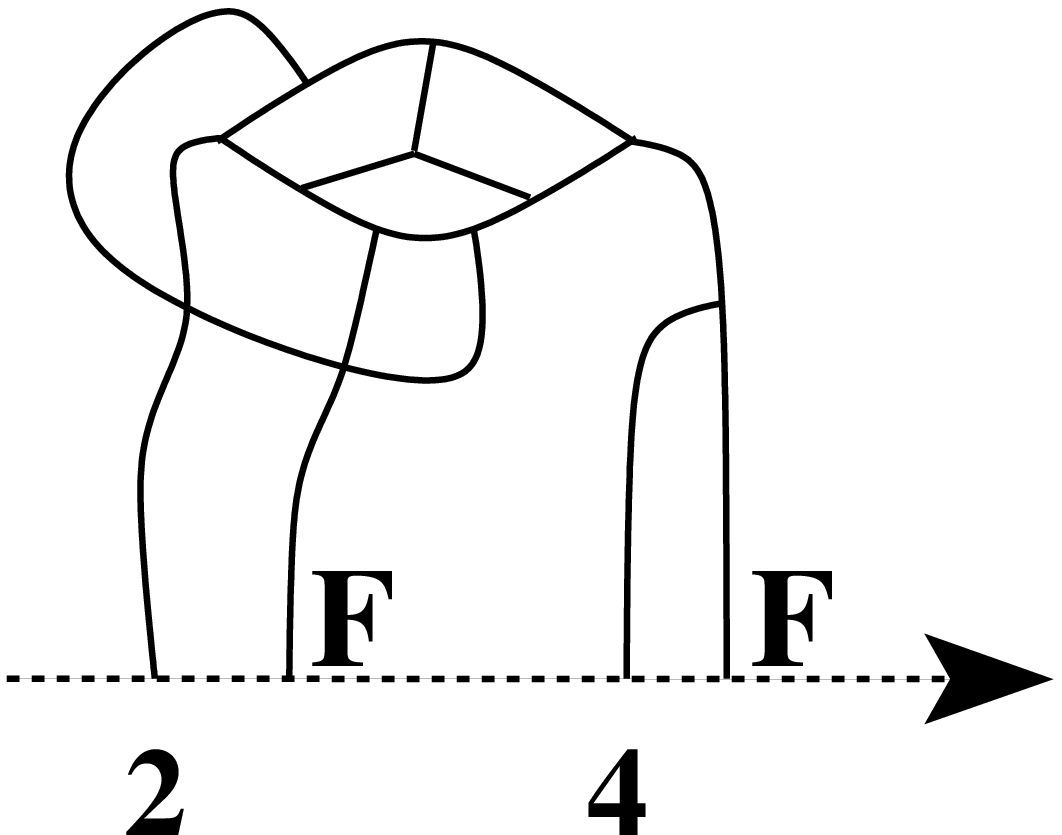}}}\
\right)\ =\  -\left(\frac{1}{2}\right)^2\
\raisebox{-4ex}{\scalebox{0.23}{\includegraphics{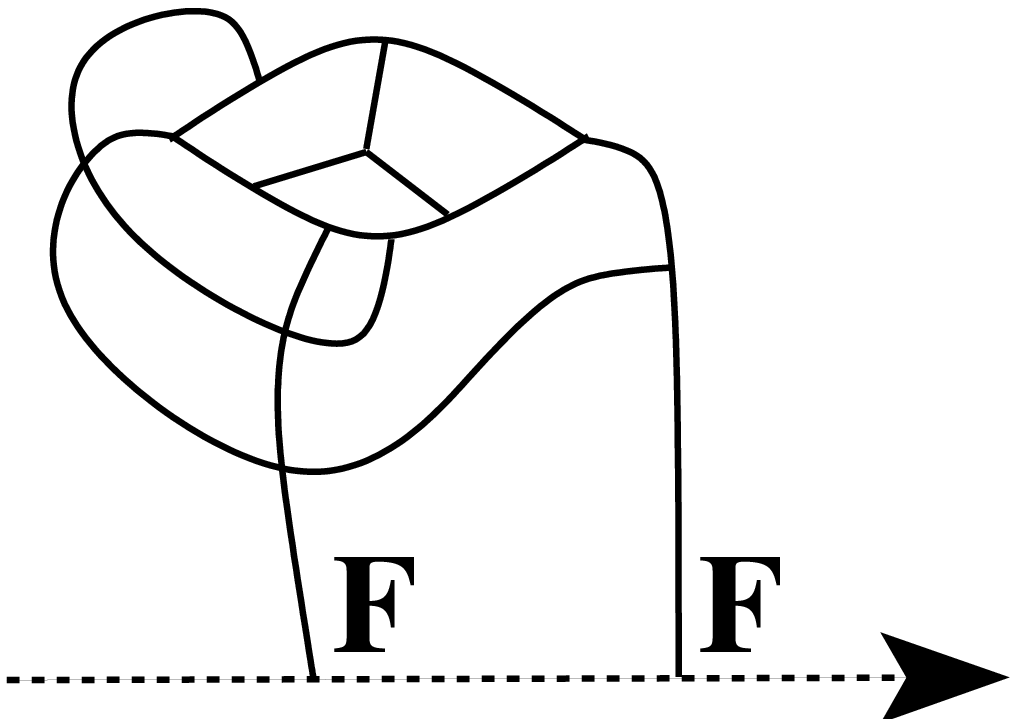}}}\
.
\end{multline*}
} \underline{\hspace{7cm}}
\end{figure}

Consider a diagram $w$ from $\WhatF$. Let $\mathcal{L}_{\bot}(w)$
denote the set of grade 1 legs of $w$. To begin, we'll explain how
to operate with a two-element subset $\mathcal{S}$ of
$\mathcal{L}_\bot(w)$ on $w$. The resulting term, which will be
denoted $\mathcal{D}_{\mathcal{S}}(w)$, is obtained from $w$ by:
\begin{enumerate}
\item{Doing signed permutations to move one of the legs in
$\mathcal{S}$ until it is adjacent to the other leg.} \item{Then
gluing those two legs together and multiplying the result by
$\frac{1}{2}$\,.}
\end{enumerate}
For example, consider the following diagram $w$, with the elements
of its set $\mathcal{L}_\bot(w)$ enumerated.
To illustrate this operation we will calculate
$\mathcal{D}_{\{2,4\}}(w).$ First we do signed permutations to
make leg $2$ adjacent to leg $4$:
\[
\raisebox{-6ex}{\scalebox{0.25}{\includegraphics{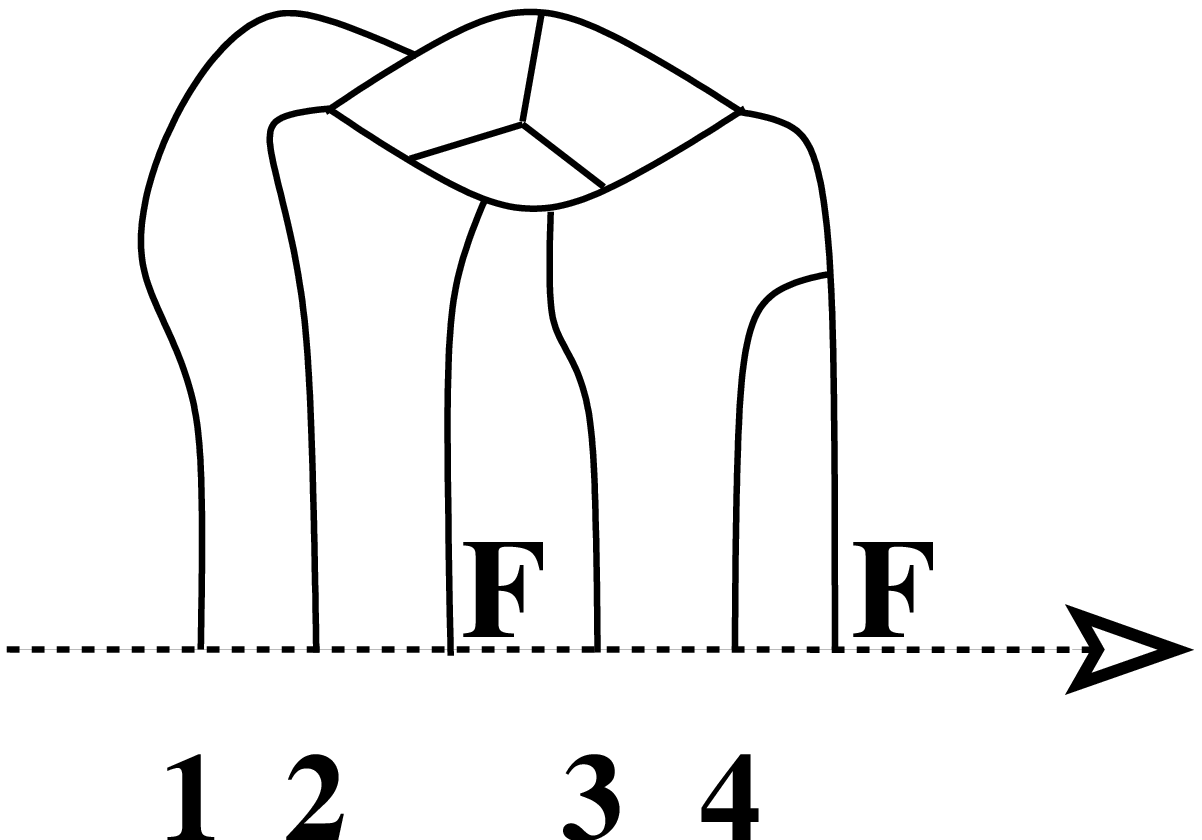}}}\
\leadsto\
\raisebox{-6ex}{\scalebox{0.25}{\includegraphics{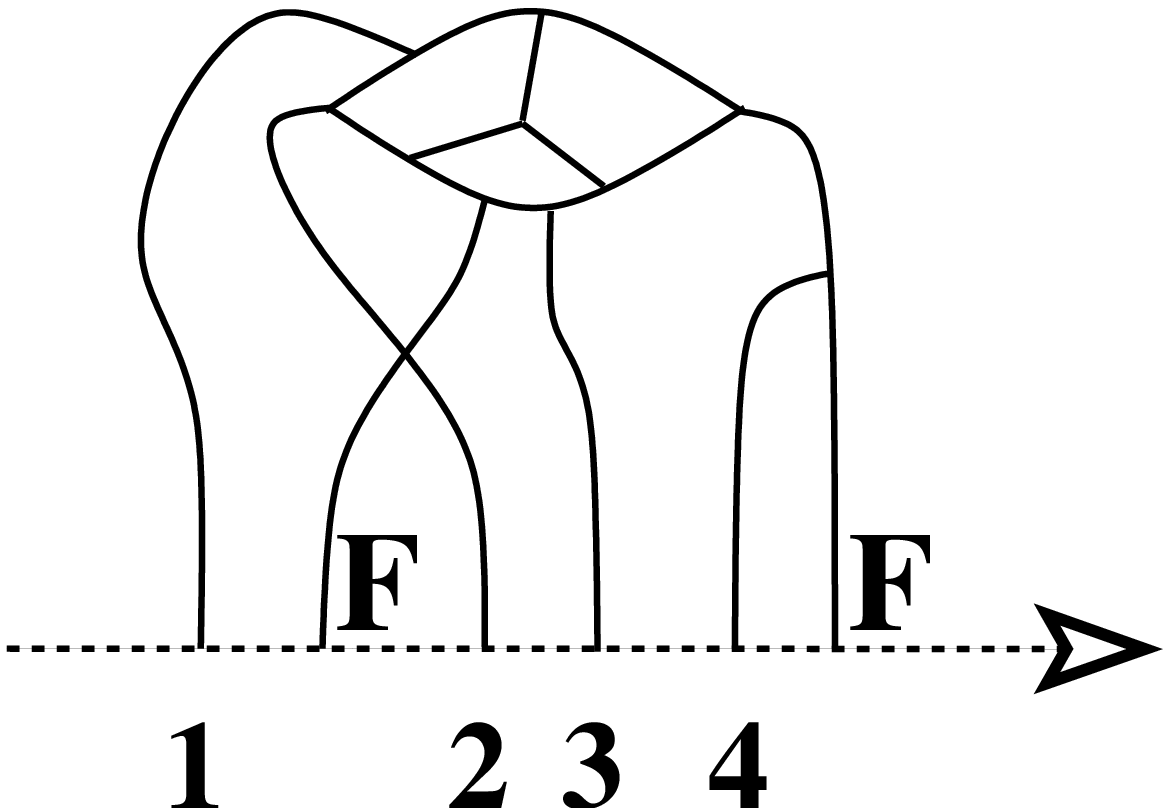}}}\
\leadsto\
-\raisebox{-6ex}{\scalebox{0.25}{\includegraphics{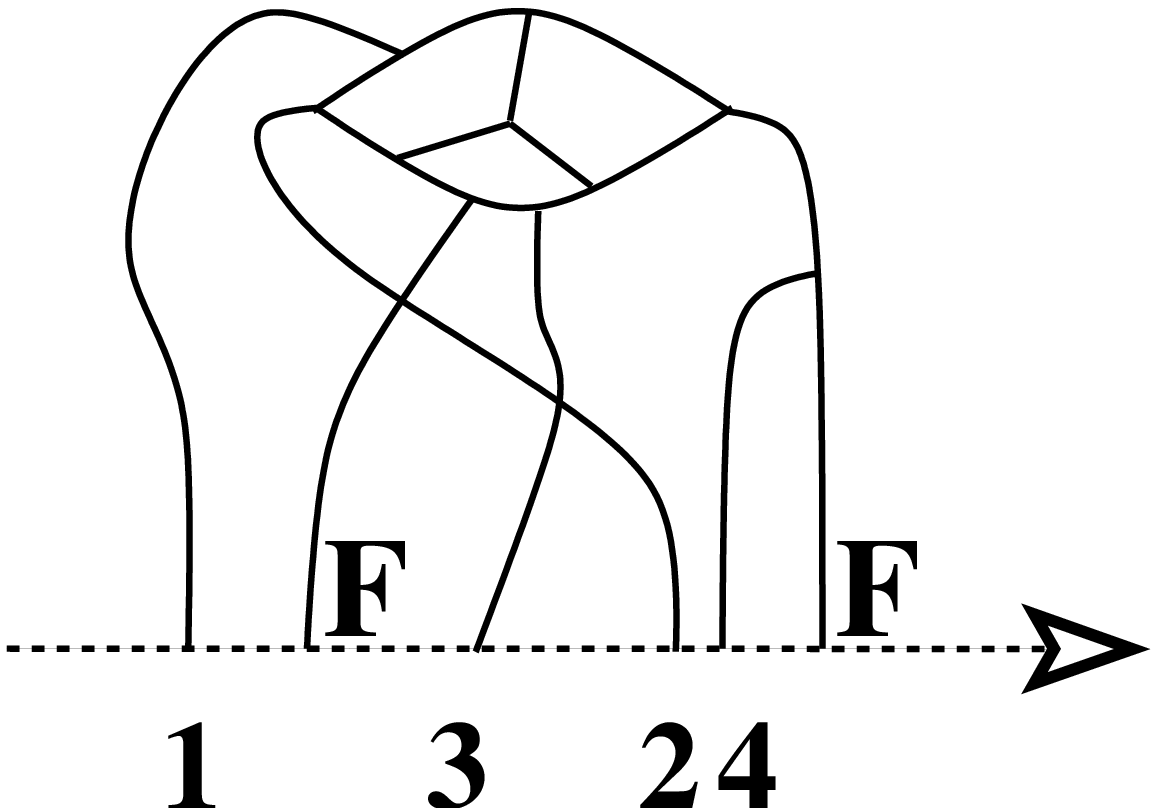}}}\
.
\]
Then we glue those two legs together, then multiply by
$\frac{1}{2}$. Thus:
\[
\mathcal{D}_{\{2,4\}}\left(
\raisebox{-6ex}{\scalebox{0.25}{\includegraphics{lambdaopexampA}}}\right)
=\ -\frac{1}{2}\
\raisebox{-6ex}{\scalebox{0.25}{\includegraphics{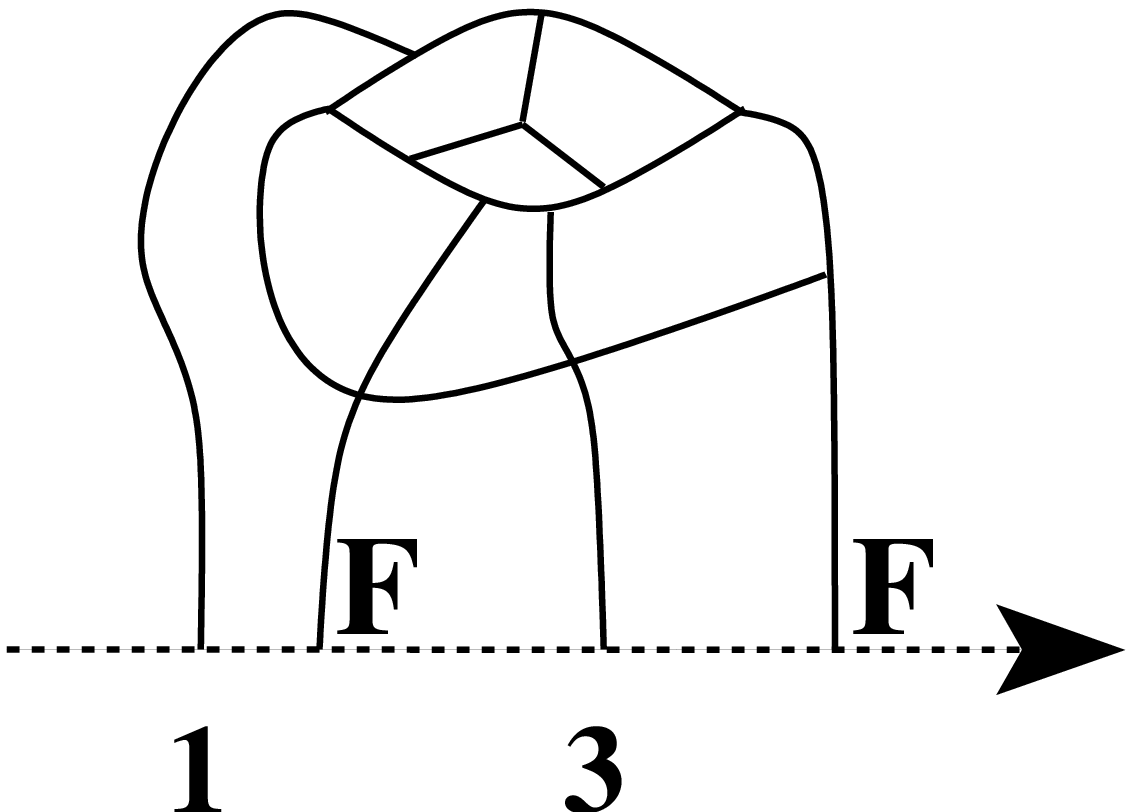}}}\
.
\]
Now we'll extend this to an action on $w$ of a {\it pairing} of
the legs of $w$. The desired map $\lambda$ will then be defined to
be a sum over one term for every pairing.
\begin{defn}
A {\bf pairing} of the legs of a diagram $w$ is defined to be a
set of mutually disjoint 2-element subsets of
$\mathcal{L}_\bot(w)$. The empty set will be regarded as a
pairing. Let $\mathcal{P}(w)$ denote the set of pairings of the
diagram $w$. Define the action of a pairing
$\wp=\{\mathcal{S}_1,\ldots,\mathcal{S}_n\}$ on the diagram $w$,
denoted $\mathcal{D}_\wp(w)$, to be:
\[
\mathcal{D}_\wp(w) = \left(\mathcal{D}_{\mathcal{S}_1}\circ \ldots
\circ \mathcal{D}_{\mathcal{S}_n}\right)(w).
\]
\end{defn}
Figure \ref{pairingexamp} gives an example of the action of a pairing on a diagram.

For this to give a well-defined operation on diagrams it is
necessary that a different choice of ordering of the elements of
$\wp$ will yield the same result. It is sufficient to show that
the equation $(\mathcal{D}_{\{a,b\}}\circ \mathcal{D}_{\{c,d\}})(w) =
(\mathcal{D}_{\{c,d\}}\circ \mathcal{D}_{\{a,b\}})(w)$ always holds,
where $a,b,c$ and $d$ are different legs of some diagram $w$. To
show this equation we just have to observe that it is true for all
possible arrangements of the positions of these legs. Let's just
check one possibility, where the two pairs are `linked':
\[
\raisebox{-5ex}{\scalebox{0.28}{\includegraphics{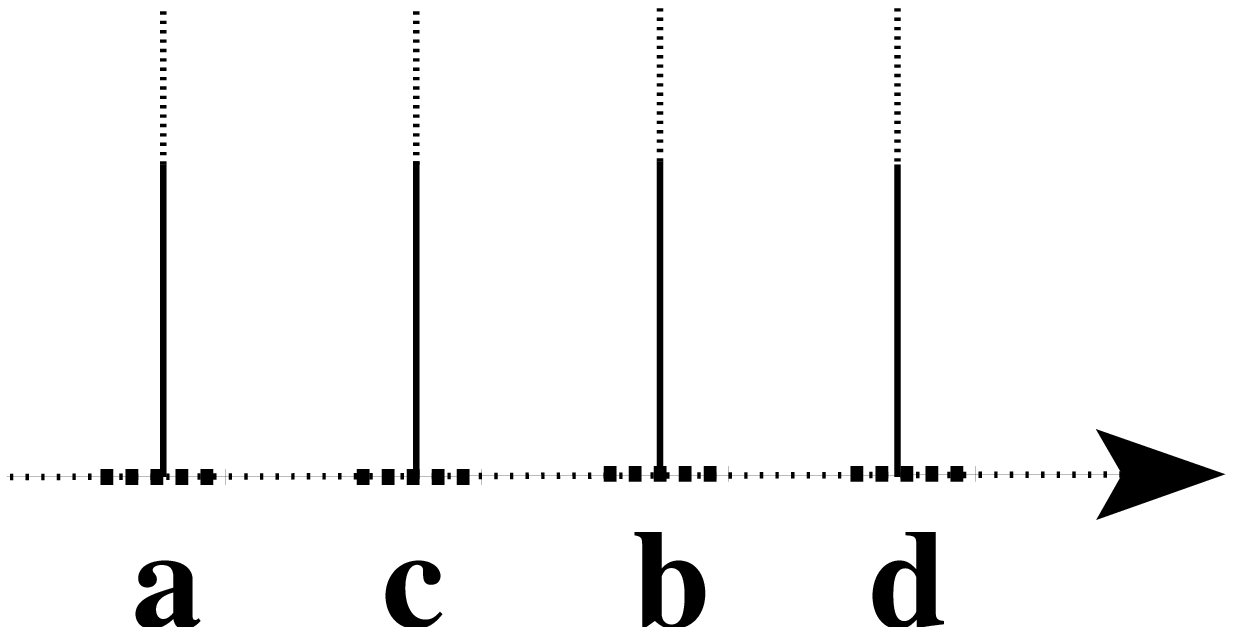}}}\ \
.
\]
Assume that the total grade of the legs between $a$ and $c$ in the
above diagram is $p$, between $c$ and $b$ is $q$, and between $b$
and $d$ is $r$. Then the result of applying
$\mathcal{D}_{\{a,b\}}\circ \mathcal{D}_{\{c,d\}}$ to the above
diagram is readily calculated to be
\[
(-1)^{q+1+r}(-1)^{p+q}\frac{1}{4}\
\raisebox{-1ex}{\scalebox{0.28}{\includegraphics{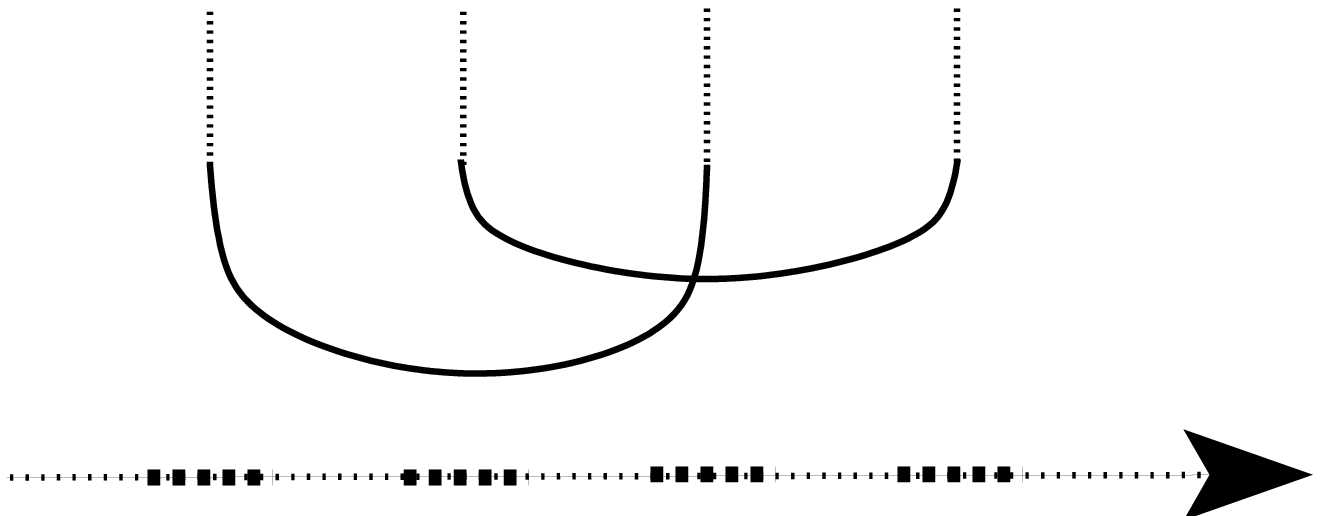}}}\ \
\ .
\]
On the other hand, the result of applying
$(\mathcal{D}_{\{c,d\}}\circ \mathcal{D}_{\{a,b\}})$ is
\[
(-1)^{p+q+1}(-1)^{q+r}\frac{1}{4}\
\raisebox{-1ex}{\scalebox{0.28}{\includegraphics{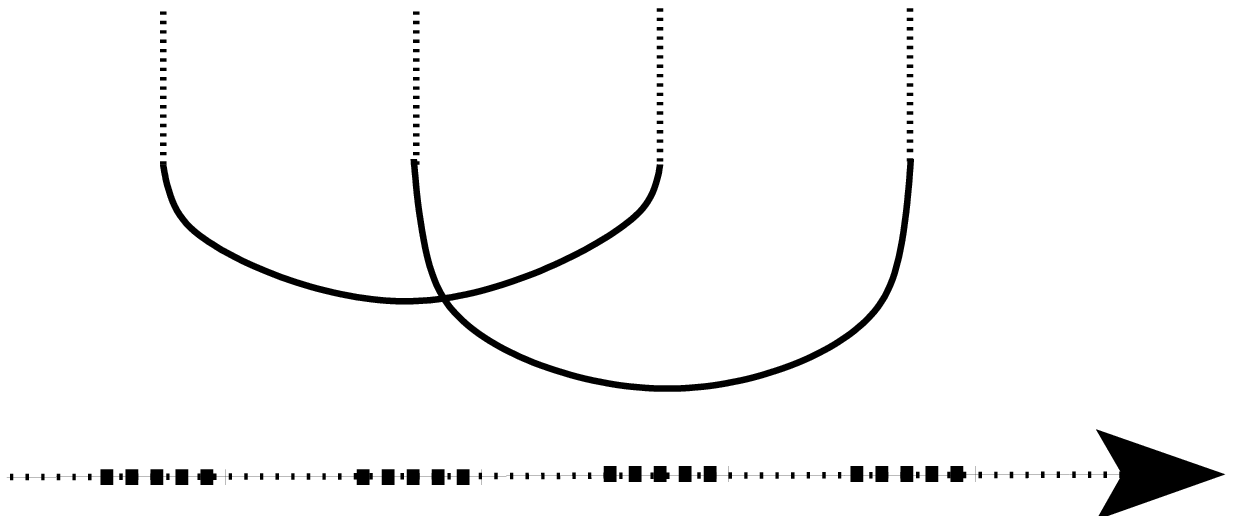}}}\ \
\ ,
\]
which is clearly the same thing. Nothing surprising happens in the
other cases.

\begin{figure}
\caption{An example of the full calculation of $\lambda$.
\label{lambexamp}}
\parbox{12cm}{
\begin{eqnarray*}
& & \\[-0.3cm]
 \lefteqn{\lambda\left(
\raisebox{-4ex}{\scalebox{0.25}{\includegraphics{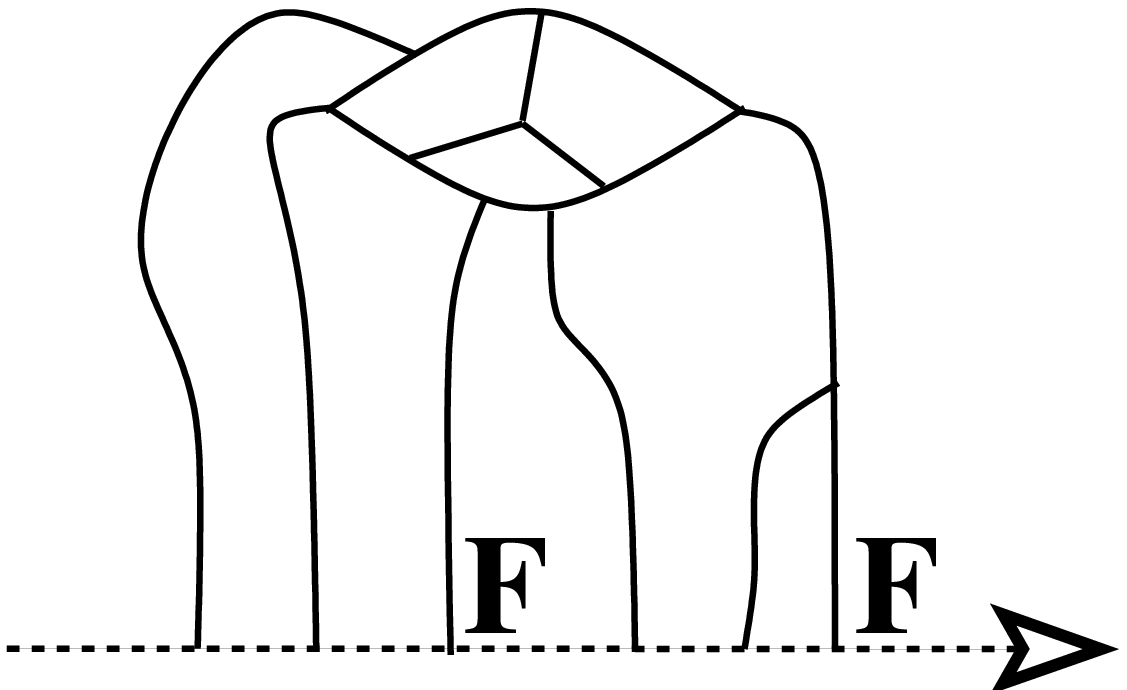}}}\right)\
\ =\ \
\raisebox{-4ex}{\scalebox{0.25}{\includegraphics{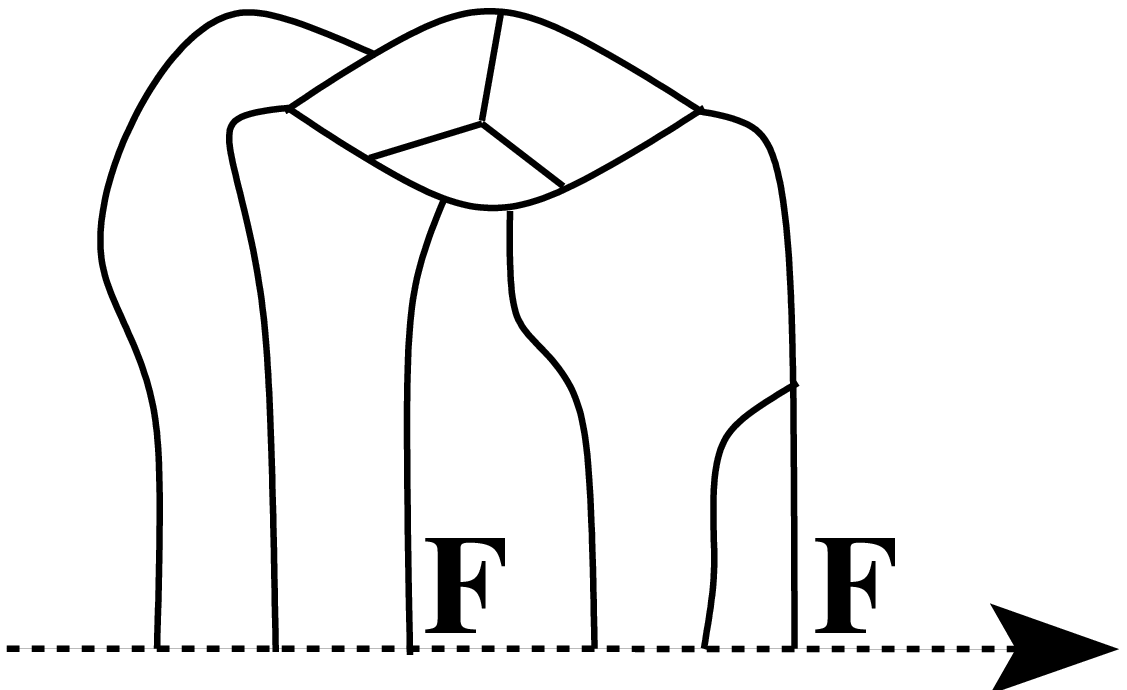}}}}\\[-0.1cm]
& & + \frac{1}{2}\
\raisebox{-3ex}{\scalebox{0.25}{\includegraphics{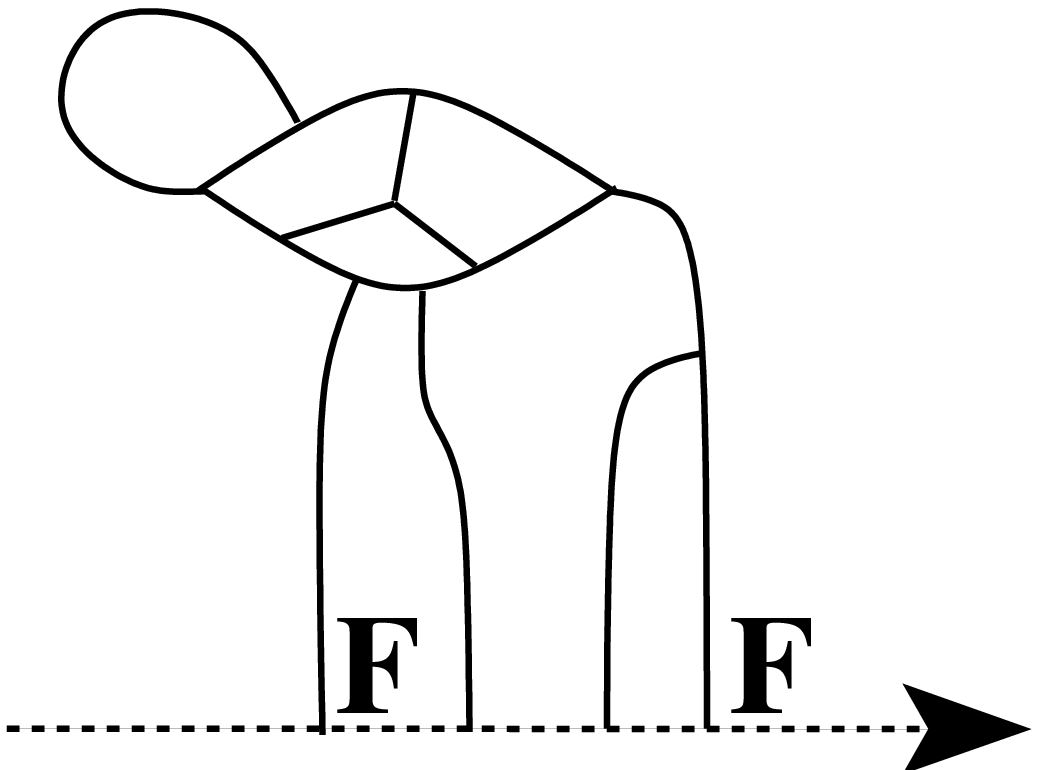}}}\
-\frac{1}{2}\
\raisebox{-3ex}{\scalebox{0.25}{\includegraphics{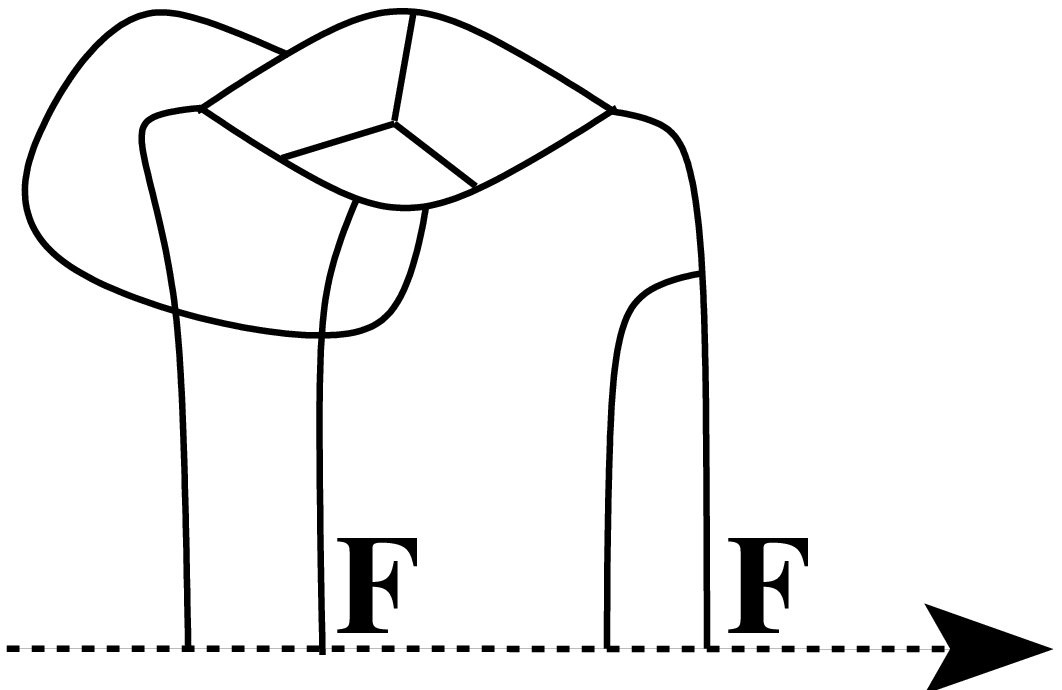}}}\
+\frac{1}{2}\
\raisebox{-3ex}{\scalebox{0.25}{\includegraphics{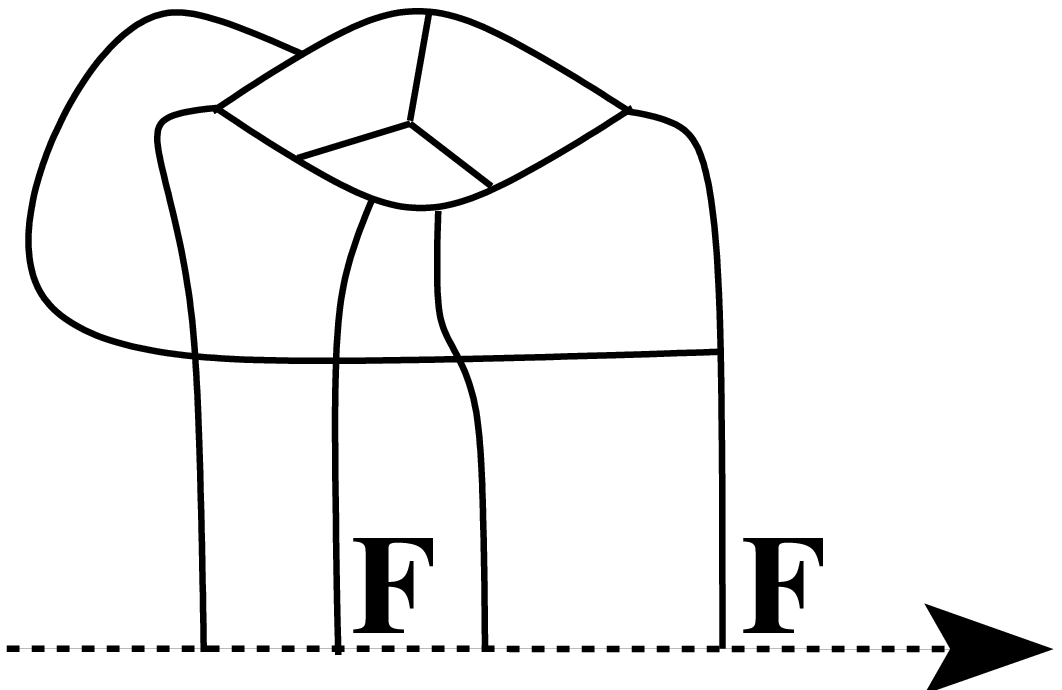}}}\
\\[0.23cm]
& & +\frac{1}{2}\
\raisebox{-3ex}{\scalebox{0.25}{\includegraphics{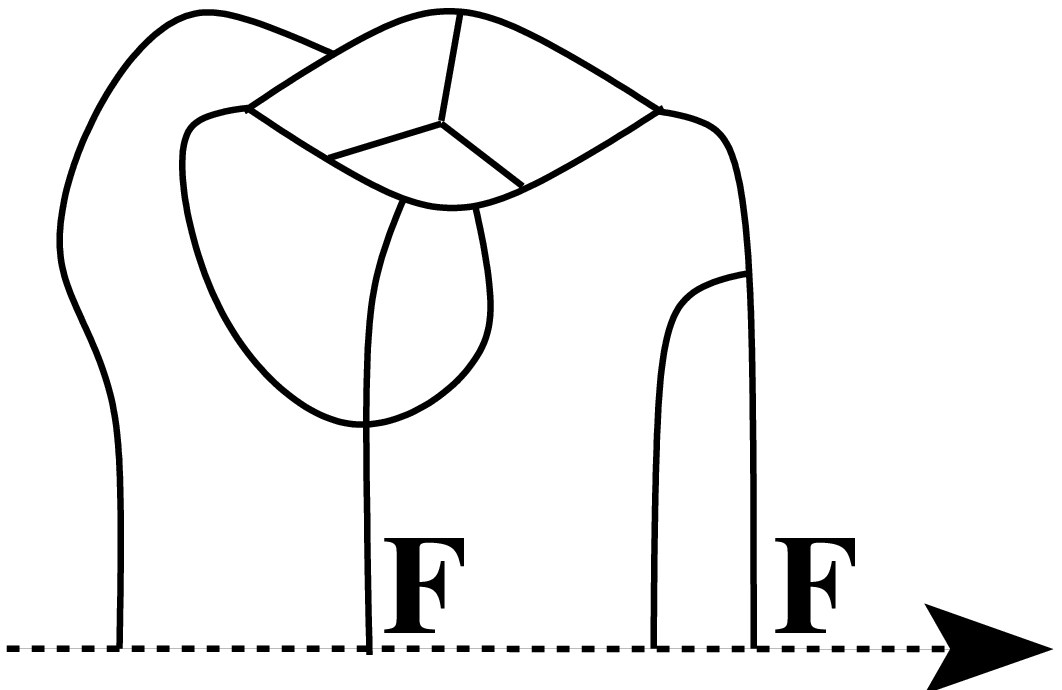}}}\
-\frac{1}{2}\
\raisebox{-3ex}{\scalebox{0.25}{\includegraphics{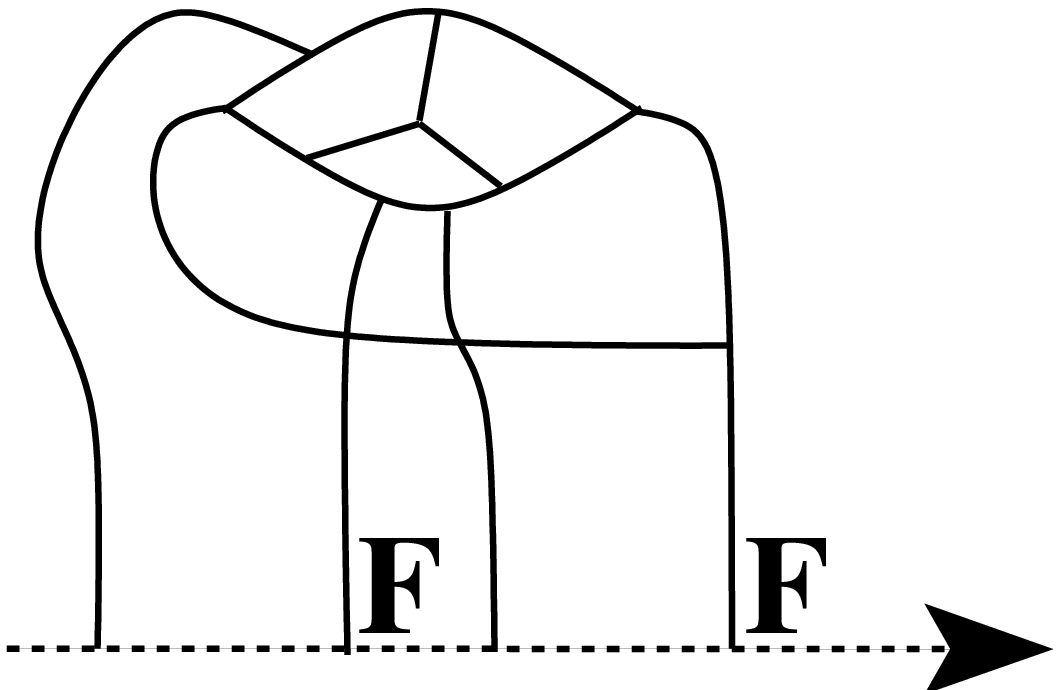}}}\
+\frac{1}{2}\
\raisebox{-3ex}{\scalebox{0.25}{\includegraphics{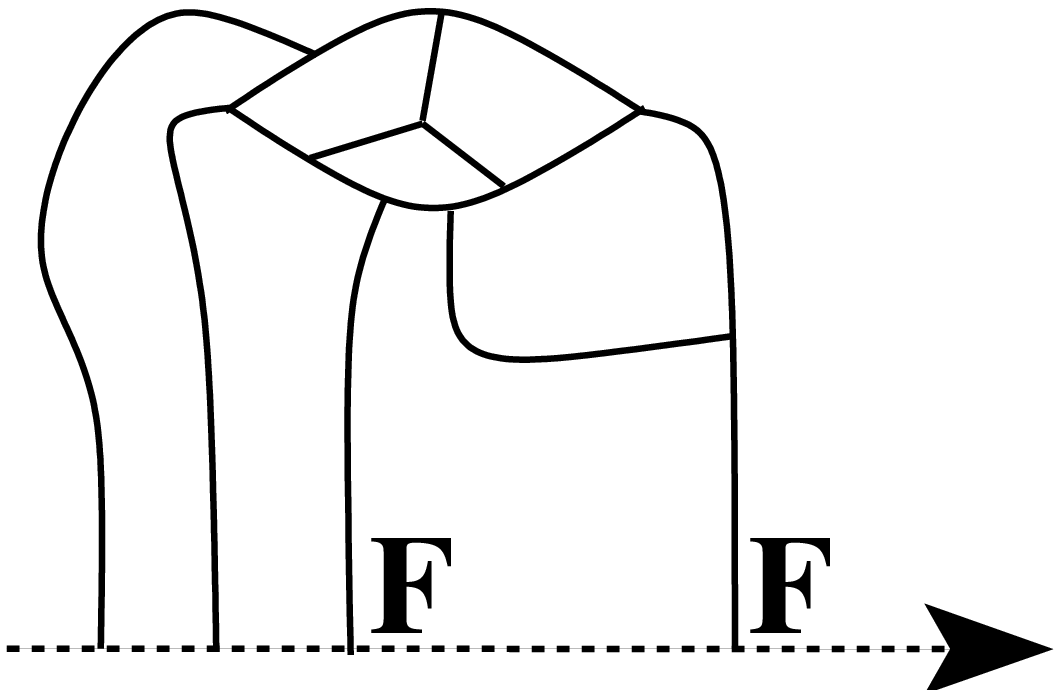}}}
\\[0.23cm]
& & +\frac{1}{4}\
\raisebox{-3ex}{\scalebox{0.25}{\includegraphics{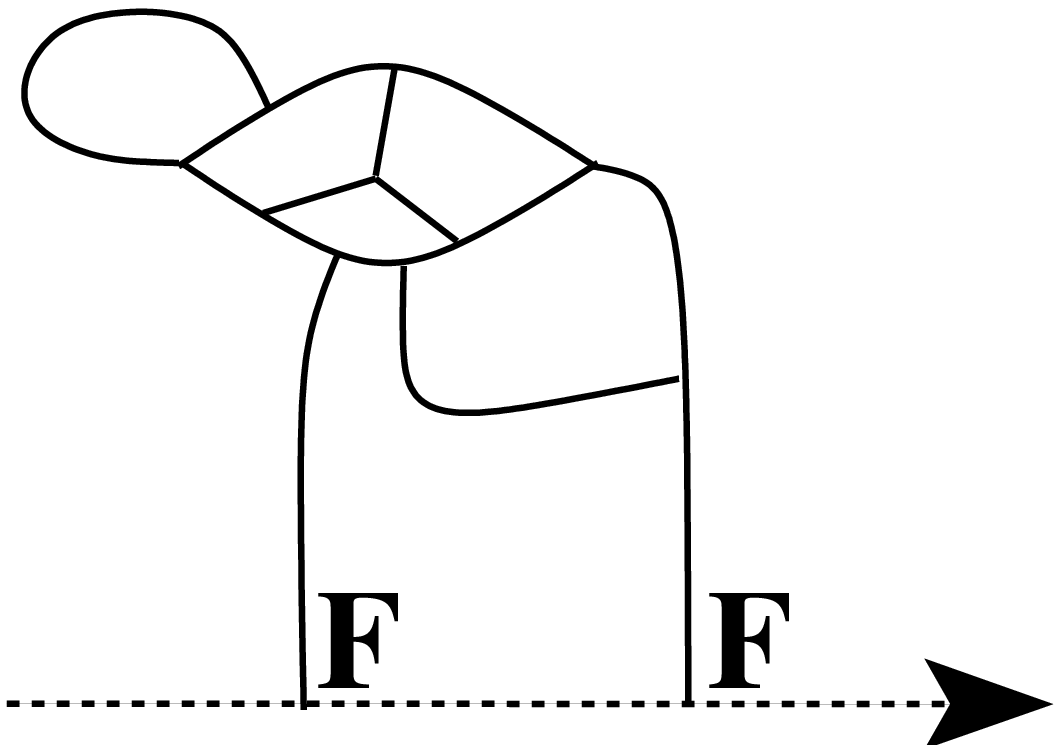}}}\
-\frac{1}{4}\
\raisebox{-3ex}{\scalebox{0.25}{\includegraphics{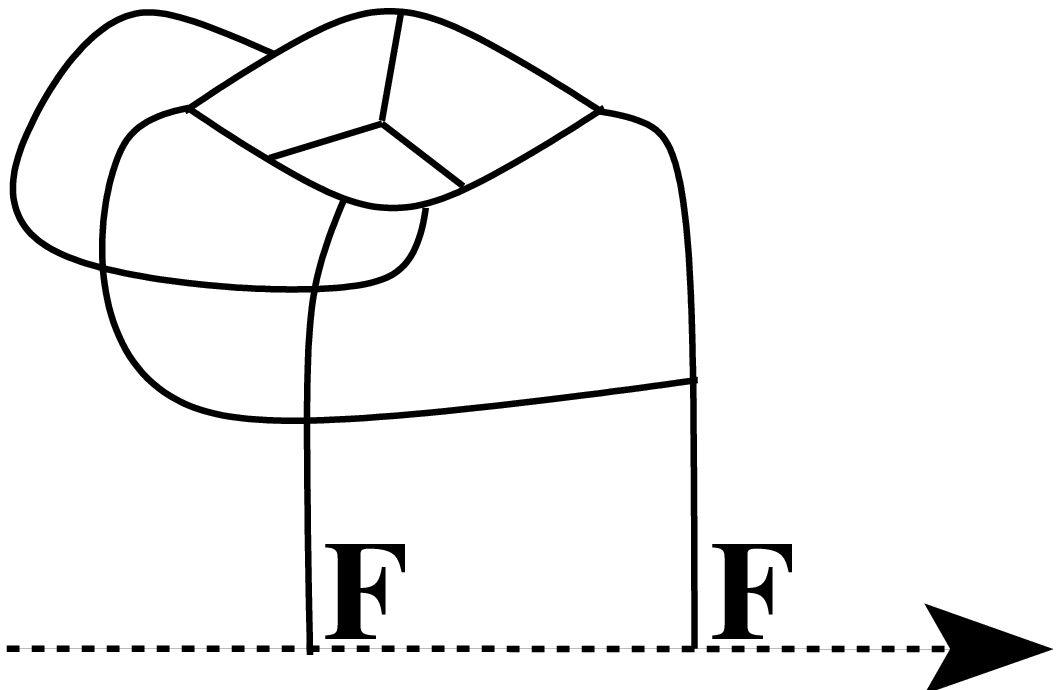}}}\
+\frac{1}{4}\
\raisebox{-3ex}{\scalebox{0.25}{\includegraphics{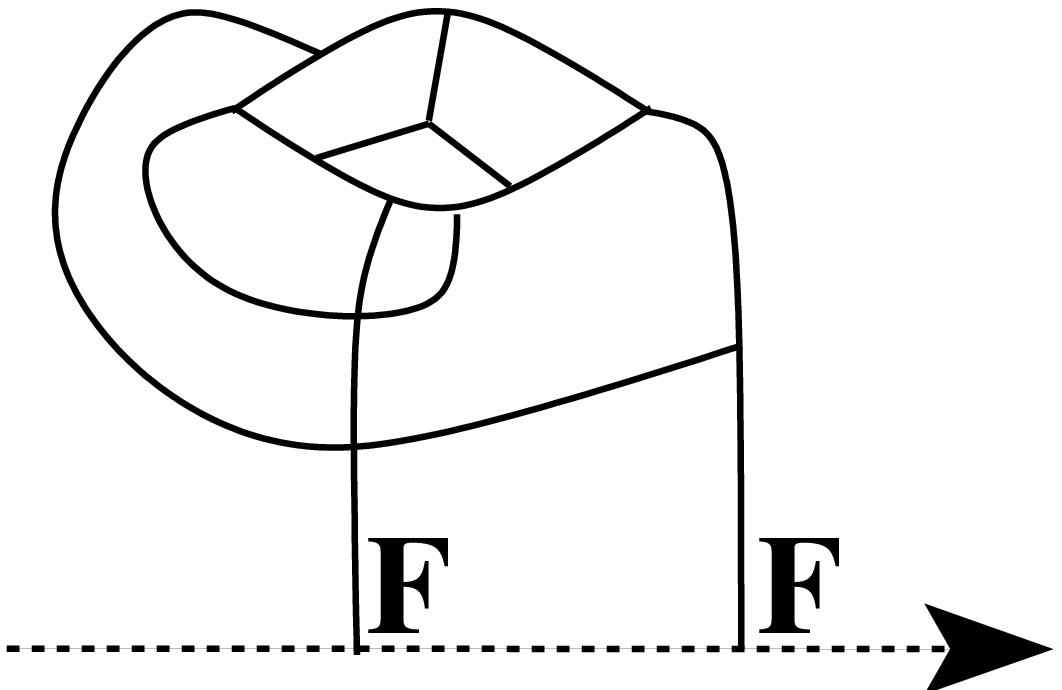}}}\
.\\[0.3cm]
\end{eqnarray*}
} \underline{\hspace{7cm}}
\end{figure}

\begin{defn}
Define a linear map $\lambda : \widehat{\mathcal{W}}_{\mathrm{F}}
\rightarrow \widehat{\mathcal{W}}_{\wedge}$ by declaring the value
of $\lambda$ on some diagram $w$ to be
$\lambda(w) = \sum_{\wp\in\mathcal{P}(w)}\mathcal{D}_\wp(w)$.
Define a map $\lambda_\iota$ similarly.
\end{defn}
An example of the result of a full computation is shown in Figure
\ref{lambexamp}.

\subsubsection{The map $\lambda$ respects
relations.}\label{lambdarespect}
\begin{prop}
The map $\lambda$ is well-defined.
\end{prop}
\begin{proof}
We must check that $\lambda$ maps relations to relations. The map
$\lambda$ only `sees' the grade 1 legs, so our problem is to show
that $\lambda$ maps the `Clifford algebra' relations to zero.
So consider a particular case of these relations, and let $D_1$,
$D_2$ and $D_3$ denote the three involved diagrams. For example:
\[
\begin{array}{c}
\raisebox{-3ex}{\scalebox{0.24}{\includegraphics{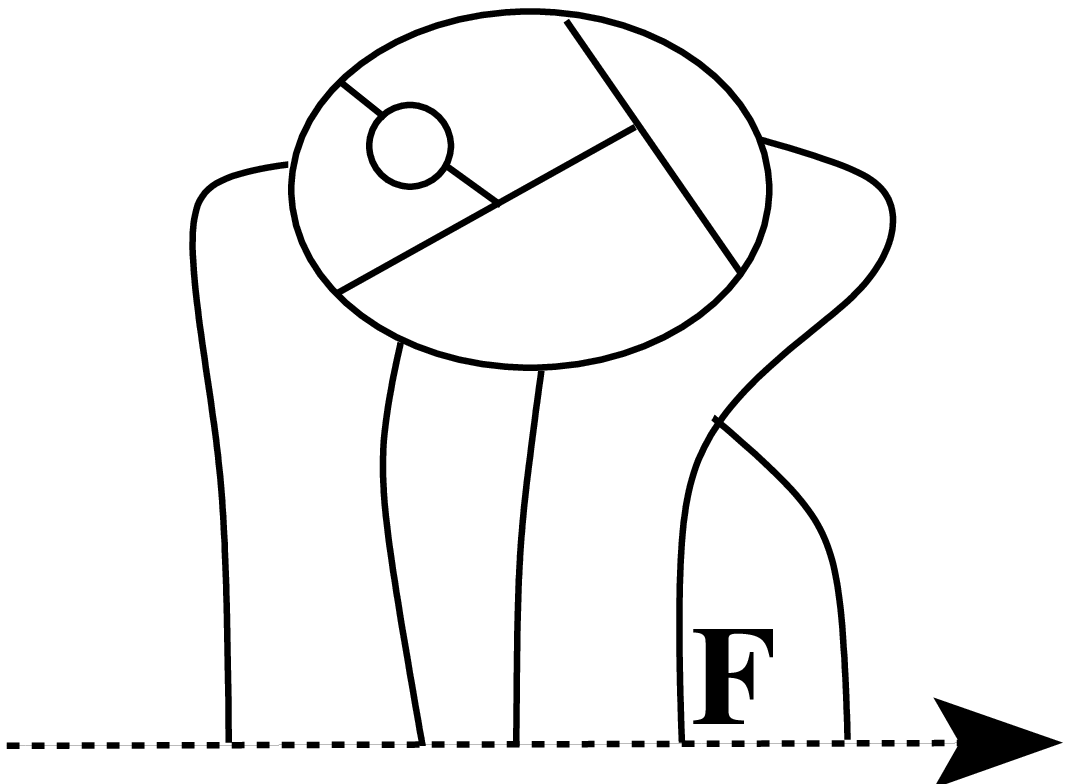}}}
\\[0.7cm]
D_1
\end{array}
 +
\begin{array}{c}
\raisebox{-3ex}{\scalebox{0.24}{\includegraphics{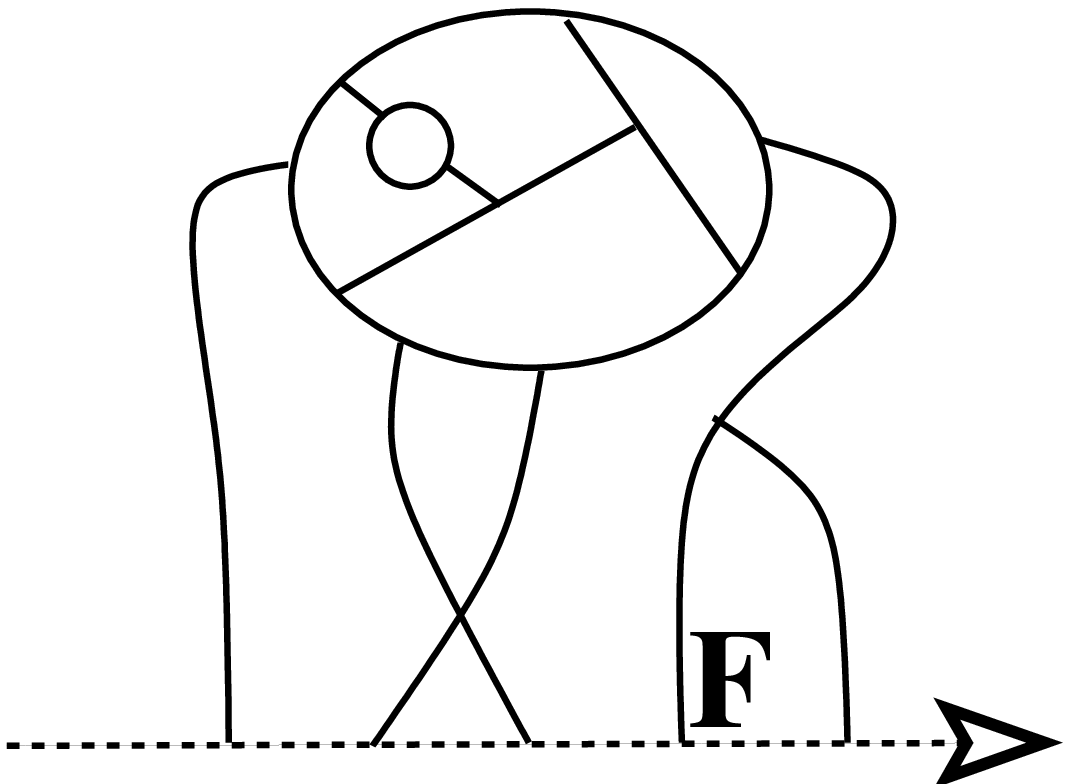}}}
\\[0.7cm]
D_2
\end{array}-
\begin{array}{c}
\raisebox{-3ex}{\scalebox{0.24}{\includegraphics{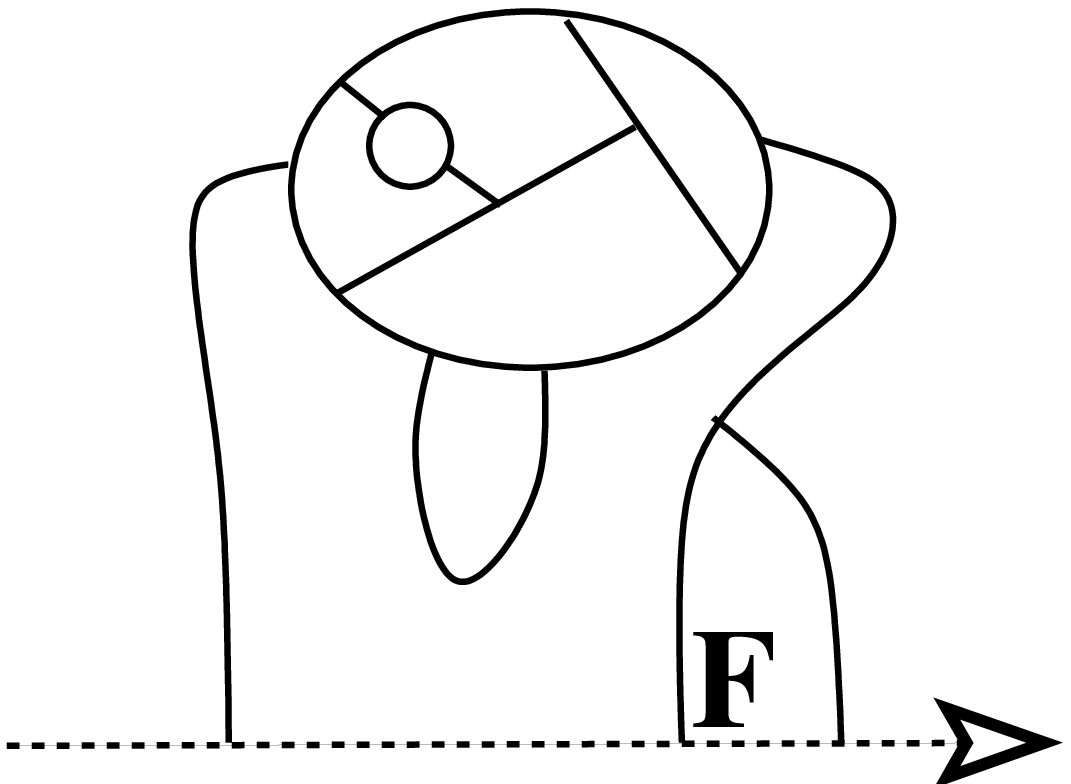}}}
\\[0.7cm]
D_3
\end{array}\ .
\]
Identify the sets of grade 1 legs $\mathcal{L}_\bot(D_1)$ and
$\mathcal{L}_\bot(D_2)$ in the natural way determined by the rest
of the graph. (So in the above example the 2nd leg of $D_1$ is
identified with the third leg of $D_2$.) Thus we can write:
\[
\lambda(D_1+D_2) =
\sum_{\wp\in\mathcal{P}(D_1)}\left(\mathcal{D}_\wp(D_1)+\mathcal{D}_\wp(D_2)\right).
\]
We'll compute the terms of this sum according to the following
three possibilities for the pairing $\wp$. Call the two grade 1
legs that change position in the expression the {\it active} legs.
The other grade 1 legs are called the {\it inactive} legs. The
grade 2 legs can be ignored completely.

{\it Possibility \#1: Neither of the active legs appear in $\wp$.}
In this case it is clear that
$\mathcal{D}_\wp(D_1)+\mathcal{D}_\wp(D_2) = 0$,
using a leg permutation relation in $\Whatwedge$.

 {\it Possibility \#2: At least one of the active legs are
paired by $\wp$ with an inactive leg.} In this case, we again have
that $\mathcal{D}_\wp(D_1)+\mathcal{D}_\wp(D_2) = 0$.
To see this, write
\[
\mathcal{D}_\wp =
\mathcal{D}_{\mathcal{S}_1}\circ\ldots\circ\mathcal{D}_{\{a,i\}}
\]
where $a$ denotes the active leg referred to, and $i$ denotes the
inactive leg it is paired to. Now calculate: \[
\mathcal{D}_{\{a,i\}}(D_1)=\mathcal{D}_{\{a,i\}}\left(
\raisebox{-5ex}{\scalebox{0.27}{\includegraphics{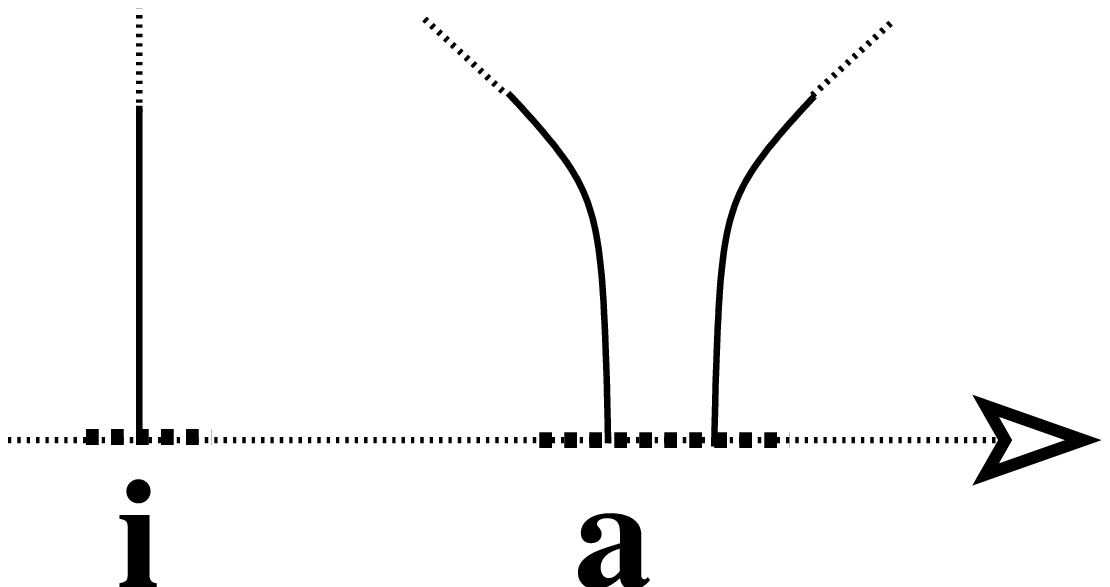}}}
\right)\ \ =\ \ \frac{1}{2}(-1)^x
\raisebox{-3ex}{\scalebox{0.27}{\includegraphics{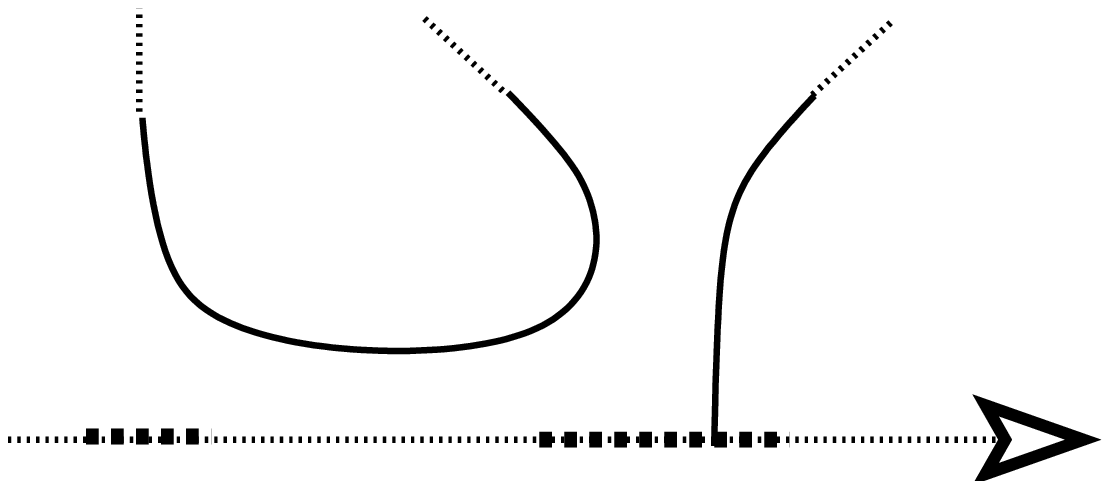}}}\ \
,
\]
where $x$ denotes the total grade of the legs between $i$ and $a$.
Also:
\[
\mathcal{D}_{\{a,i\}}(D_2)= \mathcal{D}_{\{a,i\}}\left(
\raisebox{-5ex}{\scalebox{0.27}{\includegraphics{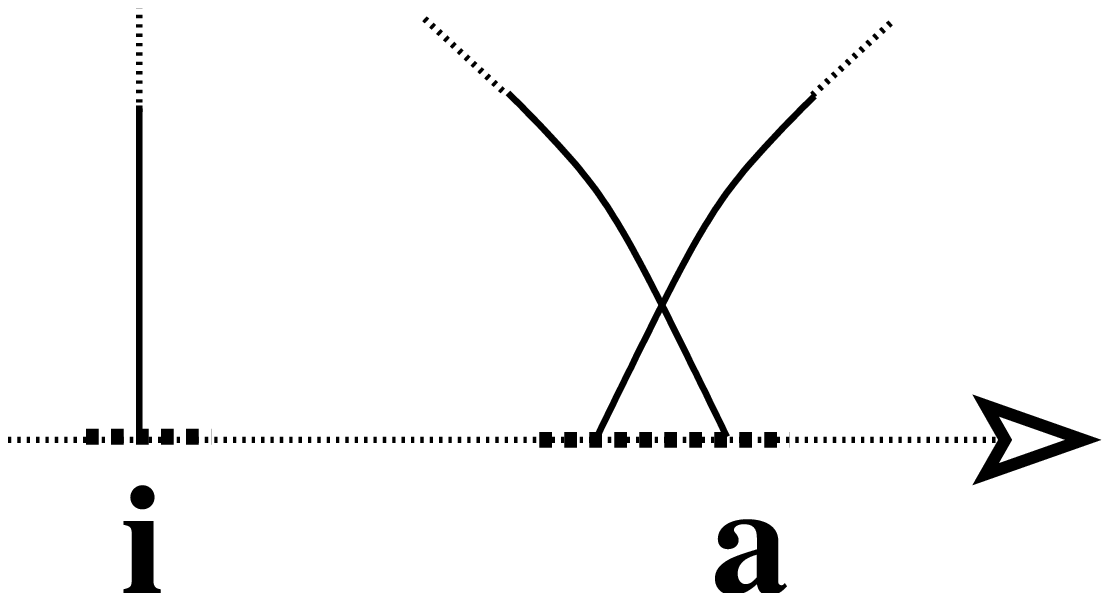}}}
\right)\ \ =\ \ \frac{1}{2}(-1)^{x+1}
\raisebox{-3ex}{\scalebox{0.27}{\includegraphics{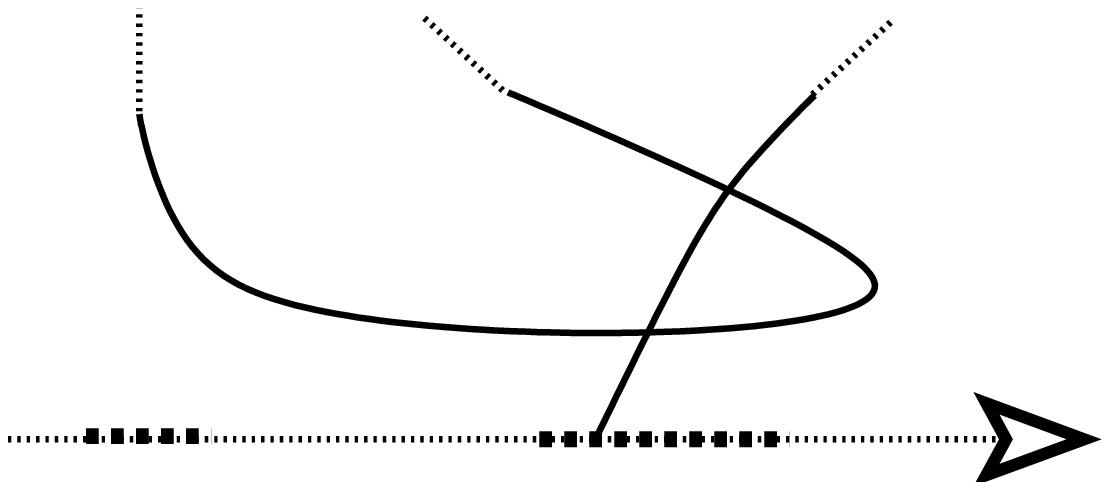}}}.
\]

{\it Possibility \#3: The two active legs are paired with each
other.} Let $\mathcal{P}^a(D_1)\subset \mathcal{P}(D_1)$ denote
the set of pairings of $D_1$ which pair the two active legs
together. There is a bijection $r: \mathcal{P}^a(D_1) \rightarrow \mathcal{P}(D_3)$.
The map is just to remove the pair $\{a_1,a_2\}$ from the pairing.
A quick calculation gives the equation
\[
\mathcal{D}_{\wp}(D_1)+\mathcal{D}_\wp(D_2) =
\mathcal{D}_{r(\wp)}(D_3).
\]
Note: this equation is why we need a ``$\frac{1}{2}$ for every
chord" in the definition of $\lambda$.
\end{proof}

\subsubsection{The map $\lambda$ does indeed left-invert $\chi_\wedge$.}

\begin{proof} Consider some diagram $w$ with $n$ degree 1 legs. Write
the signed average of the grade $1$ legs
$\chi_\wedge(w)=\sum_{\sigma\in \text{Perm}_n} w_\sigma$, where
$\text{Perm}_n$ denotes the set of permutations of the set
$\{1,\ldots,n\}$. For example:
\[
\mbox{If}\ \ w=
\raisebox{-3ex}{\scalebox{0.25}{\includegraphics{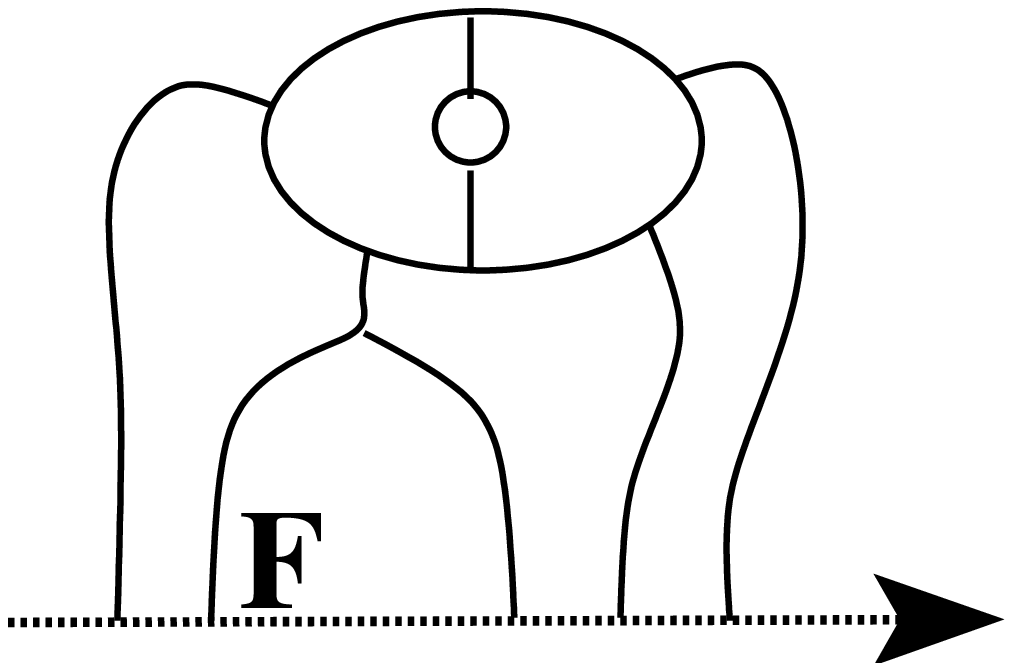}}} \ \
\mbox{then}\ \ w_{\left({ 1 2 3  4 \atop 1  4  3  2 } \right)} \
=\ \ -\frac{1}{4!}\
\raisebox{-3ex}{\scalebox{0.25}{\includegraphics{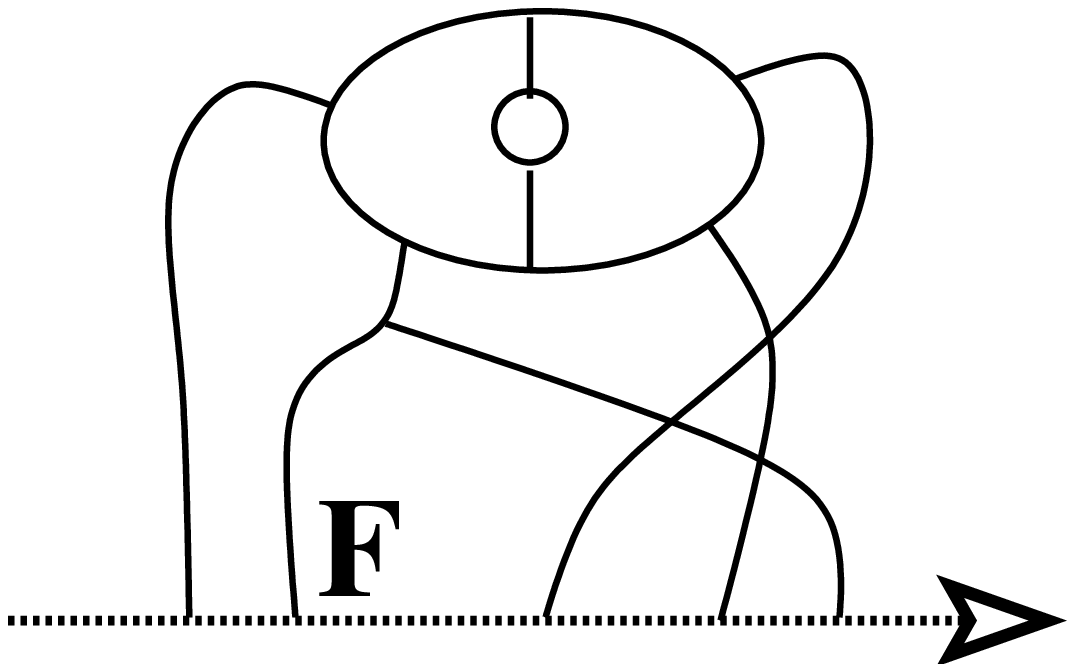}}} \ \
.
\]
For each $\sigma$, identify the set of grade 1 legs
$\mathcal{L}_\bot(w_\sigma)$ with $\mathcal{L}_\bot(w)$ using the
order of the legs as they appear along the orienting line from
left-to-right. Given a pairing $\wp\in\mathcal{P}(w),$ we'll show
that
$\mathcal{D}_\wp\left(\sum_{\sigma\in \text{Perm}_n}
w_\sigma\right) = 0
$
unless $\wp=\phi$, from which the proposition follows.

To begin, choose some
decomposition $\mathcal{D}_{\wp} =
\mathcal{D}_{\mathcal{S}_1}\circ\ldots\mathcal\circ{D}_{\{l_1,l_2\}}$.
Assuming w.l.o.g. that $l_1<l_2$, let
$\text{Perm}_n^{||}(l_1,l_2)\subset \text{Perm}_n$ denote the
subset of permutations $\sigma$ with the property that if
$\sigma(x)=l_1$ and $\sigma(y)=l_2$ then $x<y$. Graphically
speaking, $\text{Perm}_n^{||}(l_1,l_2)$ consists of those
permutations where the straight lines ending at positions $l_1$
and $l_2$ do not cross.
Similarly, define ${\text{Perm}}_n^{X}(l_1,l_2)\subset
\text{Perm}_n$ to be the subset of permutations where they do
cross. There is clearly a bijection
$\vartheta:\text{Perm}_n^{||}(l_1,l_2)\stackrel{\cong}{\longrightarrow}
\text{Perm}_n^{X}(l_1,l_2)$:
\[
\vartheta\left(
\raisebox{-5.5ex}{\scalebox{0.22}{\includegraphics{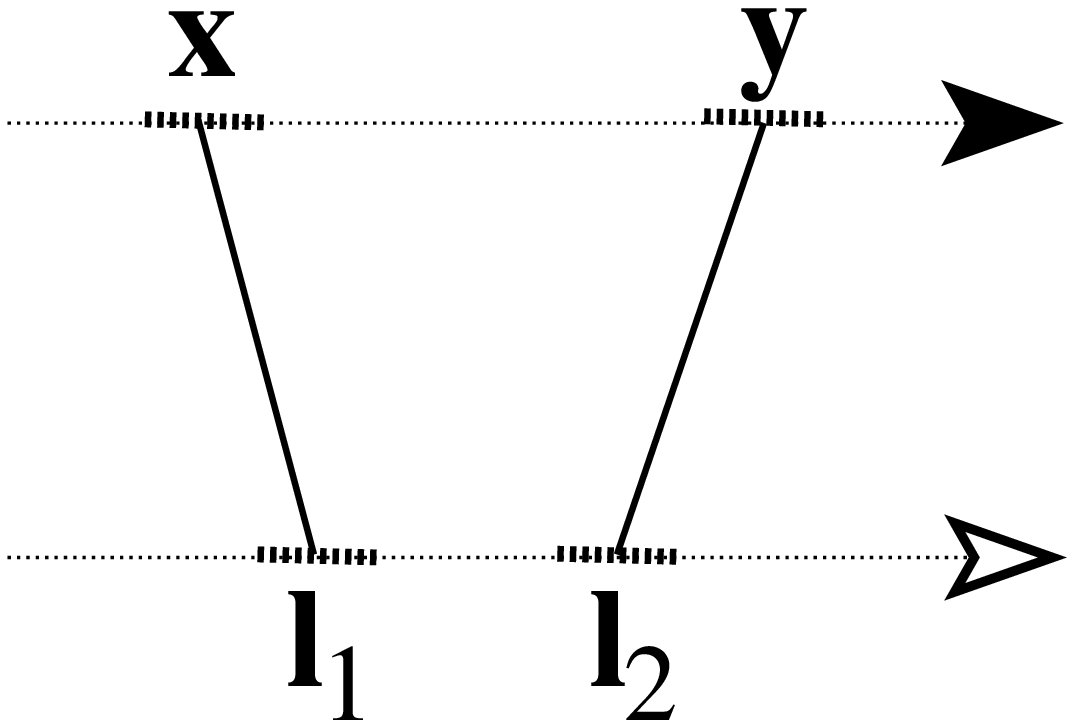}}}
\right)\ \ =\ \
\raisebox{-5.5ex}{\scalebox{0.22}{\includegraphics{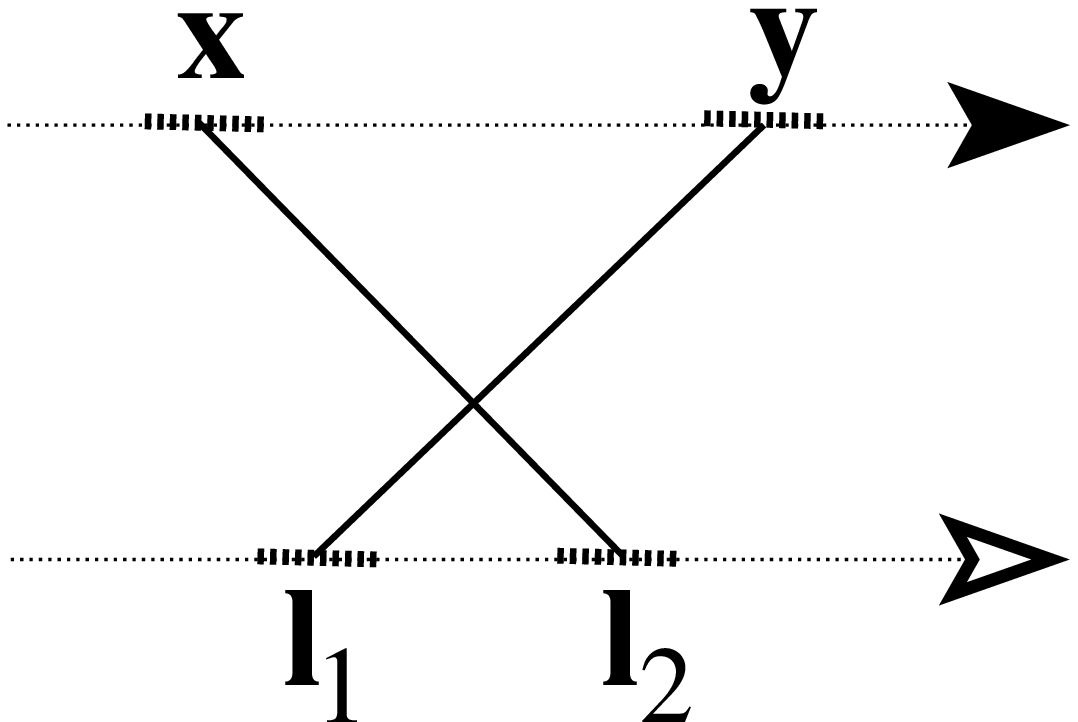}}}\
.
\]
Now calculate:
\[
\mathcal{D}_{\wp}\left(\sum_{\sigma\in \text{Perm}_n} w_\sigma\right) =
\left(\mathcal{D}_{\mathcal{S}_1}\circ\ldots\mathcal\circ{D}_{\{l_1,l_2\}}\right)\left(\sum_{\sigma\in \text{Perm}_n^{||}(l_1,l_2)} (w_\sigma+ w_{\vartheta(\sigma)})\right).
\]
The claim now follows from the observation that $
\mathcal{D}_{\{l_1,l_2\}}(w_\sigma) = -
\mathcal{D}_{\{l_1,l_2\}}(w_{\vartheta(\sigma)})$.
\end{proof}

\subsection{Proof of Theorem \ref{bigwedgetheorem}}
\label{bigwedgetheoremproof} To prove this theorem we must:
\begin{enumerate}
\item{Prove that $\iota\circ \chi_\wedge=\chi_{\wedge\iota}\circ\iota$.} \item{Prove
that $\chi_\wedge$ and $\chi_{\wedge\iota}$ are surjective.}
\item{Observe that because $\lambda\circ\chi_\wedge =
\mathrm{id}_{\widehat{\mathcal{W}}_{\wedge}}$ and
$\lambda_\iota\circ\chi_{\wedge\iota} =
\mathrm{id}_{\widehat{\mathcal{W}}_{\wedge\iota}}$, it follows
that $\chi_\wedge$ and $\chi_{\wedge\iota}$ are also injective.}
\end{enumerate}

The first fact is a (straightforward if a little bit technical) combinatorial exercise.

\begin{prop}
The maps $\chi_\wedge$ and $\chi_{\wedge\iota}$ are surjective.
\end{prop}
{\it Proof.} This argument will be familiar. Take $w$, an
arbitrary {\bf F}-Weil diagram. Let $w_\wedge$ be the same
diagram, but regarded as a generator of
$\widehat{\mathcal{W}}_\wedge$. Then:
\[
w - \chi_\wedge(w_\wedge) = \left(\parbox{4cm}{A sum of diagrams
with less grade 1 legs than $w$.}\right).
\]
By induction, $w$ lies in the image of $\chi_\wedge$.
\begin{flushright}
$\Box$
\end{flushright}

\end{document}